\documentclass[11pt]{article}

\usepackage{amssymb,latexsym}

\input pictex.sty

\makeatletter
\@addtoreset{equation}{section}
\makeatother
\renewcommand{\theequation}{\arabic{section}.\arabic{equation}}
\def\proof{\par{\it Proof}. \ignorespaces}
\def\endproof{{\ \vbox{\hrule\hbox{%
   \vrule height1.3ex\hskip0.8ex\vrule}\hrule }}\par}

\makeatletter
\def\@thmcounterend{}

\def\newtheorem{\@ifstar{\@sthm}{\@Sthm}}

\def\@Sthm#1{\@ifnextchar[{\@othm{#1}}{\@nthm{#1}}}

\def\@xnthm#1#2[#3]#4#5{\expandafter\@ifdefinable\csname
#1\endcsname
   {\@definecounter{#1}\@addtoreset{#1}{#3}%
   \expandafter\xdef\csname the#1\endcsname{\expandafter\noexpand
     \csname the#3\endcsname \@thmcountersep \@thmcounter{#1}}%
   \expandafter\xdef\csname #1name\endcsname{#2}%
   \global\@namedef{#1}{\@thm{#1}{\csname
#1name\endcsname}{#4}{#5}}%
                             
\global\@namedef{end#1}{\@endtheorem}}}

\def\@ynthm#1#2#3#4{\expandafter\@ifdefinable\csname #1\endcsname
   {\@definecounter{#1}%
   \expandafter\xdef\csname the#1\endcsname{\@thmcounter{#1}}%
   \expandafter\xdef\csname #1name\endcsname{#2}%
   \global\@namedef{#1}{\@thm{#1}{\csname
#1name\endcsname}{#3}{#4}}%
                              
\global\@namedef{end#1}{\@endtheorem}}}

\def\@othm#1[#2]#3#4#5{%
  \@ifundefined{c@#2}{\@latexerr{No theorem environment `#2'
defined}\@eha}%
  {\expandafter\@ifdefinable\csname #1\endcsname
  {\global\@namedef{the#1}{\@nameuse{the#2}}%
  \expandafter\xdef\csname #1name\endcsname{#3}%
  \global\@namedef{#1}{\@thm{#2}{\csname
#1name\endcsname}{#4}{#5}}%
  \global\@namedef{end#1}{\@endtheorem}}}}

\def\@thm#1#2#3#4{\refstepcounter
{#1}\@ifnextchar[{\@ythm{#1}{#2}{#3}{#4}}{\@xthm{#1}{#2}{#3}{#4}}}

\def\@xthm#1#2#3#4{\@begintheorem{#2}{\csname
the#1\endcsname}{#3}{#4}%
                    \ignorespaces}

\def\@ythm#1#2#3#4[#5]{\@opargbegintheorem{#2}{\csname
       the#1\endcsname}{#5}{#3}{#4}\ignorespaces}

\def\@begintheorem#1#2#3#4{\trivlist
                 \item[\hskip\labelsep{#3#1\ #2\@thmcounterend}]#4}

\def\@opargbegintheorem#1#2#3#4#5{\trivlist
      \item[\hskip\labelsep{#4#1\ #2\ (#3)\@thmcounterend}]#5}

\def\@sthm#1#2{\@Ynthm{#1}{#2}}

\def\@Ynthm#1#2#3#4{\expandafter\@ifdefinable\csname #1\endcsname
   {\global\@namedef{#1}{\@Thm{\csname #1name\endcsname}{#3}{#4}}%
    \expandafter\xdef\csname #1name\endcsname{#2}%
    \global\@namedef{end#1}{\@endtheorem}}}

\def\@Thm#1#2#3{\@ifnextchar[{\@Ythm{#1}{#2}{#3}}{\@Xthm{#1}{#2}{
#3}}}

\def\@Xthm#1#2#3{\@Begintheorem{#1}{#2}{#3}\ignorespaces}

\def\@Ythm#1#2#3[#4]{\@Opargbegintheorem{#1}
       {#4}{#2}{#3}\ignorespaces}

\def\@Begintheorem#1#2#3{#3\trivlist
                          
\item[\hskip\labelsep{#2#1\@thmcounterend}]}

\def\@Opargbegintheorem#1#2#3#4{#4\trivlist
      \item[\hskip\labelsep{#3#1\ (#2)\@thmcounterend}]}

\makeatother
\makeatletter\newif\if@ceqnsw\newskip\baselineskipp\newskip\cabovedisplayskip
\newdimen\cpamarraycolsep
\newcounter{abceqn}

\def\veqnarray{\baselineskipp=\baselineskip%\baselineskip=-1000\p@%
\cabovedisplayskip=\abovedisplayskip%
\advance\abovedisplayskip by -.3\cabovedisplayskip%
\advance\belowdisplayskip by -.3\cabovedisplayskip%
$$\global\@ceqnswtrue\refstepcounter{equation}%
\edef\@currentlabel{\theequation}%
\m@th\global\@eqcnt\z@\tabskip\@centering\let\\\@eqncr%
\vcenter\bgroup\everycr{}\vskip .7\cabovedisplayskip%
\baselineskip=\baselineskipp%
\halign\bgroup\hskip\@centering\global\@eqnswfalse
  $\displaystyle\tabskip\z@skip{##}$\@eqnsel&\global\@eqcnt\@ne
  \hskip \tw@\cpamarraycolsep \hfil${##}$\hfil
  &\global\@eqcnt\tw@ \hskip \tw@\cpamarraycolsep
  $\displaystyle{##}$\hfil%
   \tabskip\@centering&\global\@eqcnt\thr@@
       \hbox to\z@\bgroup\hss##\egroup\tabskip\z@skip\cr}

\def\endveqnarray{\@@eeqncr\egroup\vskip .7\cabovedisplayskip\egroup%
\eqno\if@ceqnsw\hbox{\@eqnnum}\else \global\advance\c@equation\m@ne\fi%
\global\@ignoretrue$$}

\def\nonumber{\global\@eqnswfalse\global\@ceqnswfalse}
\def\@@eeqncr{\let\@tempa\relax\ifcase\@eqcnt\def\@tempa{& &}\or\def\@tempa{&}%
   \or \else
       \let\@tempa\@empty
       \@latex@error{Too many columns in eqnarray environment}\@ehc\fi
    \global\@eqcnt\z@\cr}

\makeatother

\newtheorem{Theorem}[equation]{Theorem}{\bf}{\it}
\newtheorem{Lemma}[equation]{Lemma}{\bf}{\it}
{\bf}{\it}
\newtheorem{Corollary}[equation]{Corollary}{\bf}{\it}
\newtheorem{Proposition}[equation]{Proposition}{\bf}{\it}
{\bf}{\it}
\newtheorem{Remark}[equation]{Remark}{\bf}{\it}
\newtheorem{Definition}[equation]{Definition}{\bf}{\it}
\newtheorem{Notation}[equation]{Notation}{\bf}{\it}

\begin{document}
\title{{\bf{ On the existence of quasi-periodic lattice oscillations
\footnote{This paper was submitted as Habilitaionsschrift to the 
Ludwig-Maximilians-Universit\"at M\"unchen in 1999.}}}}

\author{{\bf{T. Kriecherbauer }} \\  Ruhr-Universit\"at Bochum, Germany}
\maketitle

%\pagenumbering{roman}
%\evensidemargin 0.0in
%\oddsidemargin 0.0in
%\textwidth 6.5in
%\topmargin 0.0in
%\textheight 8.8in
%\thispagestyle{empty}

%\begin{center}
%{\LARGE \bf On the existence of quasi-periodic lattice oscillations} 
%\vspace{0.5cm} \vspace{6cm}
%\end{center}

%\begin{center} 
%Habilitationsschrift\\
%zur Feststellung der Lehrbef\"ahigung\\
%f\"ur das Fachgebiet Mathematik\\
%5in der Fakult\"at f\"ur Mathematik und Informatik\\
%der Ludwig-Maximilians-Universit\"at M\"unchen \vspace{5cm}
%\end{center}

%\begin{center}
%vorgelegt von \\
%Thomas Kriecherbauer\\
%im November 1999
%\end{center}

%\clearpage

\evensidemargin 0.0in
\oddsidemargin 0.0in
\textwidth 6.5in
\topmargin -0.5in
\textheight 8.8in

%\setcounter{page}{3}
%\begin{center}
%{\Large \bf Summary} \vspace{0.5cm}\\
%\end{center}
\begin{abstract}
We construct multi-phase travelling waves which are periodic in time
and quasi-periodic in the spatial variable for infinite nonlinear chains with
nearest neighbor interactions. Such solutions have been observed numerically
in \cite{DKV}. Their existence has so far been established only
for linear systems and in the integrable case of the Toda lattice
\cite{Krc}, \cite{DKV}. For a general class of nonlinear chains
one can show the existence of single-phase travelling waves using either
perturbation theory \cite{DKV} or topological methods \cite{FV}.
For travelling waves with more than one phase, however, a
small divisor problem occurs and the above mentioned constructions fail.

In this paper we prove a KAM-type theorem which yields the existence
of a large family of multi-phase waves for a general class of 
chains. In fact, we show that  
most small amplitude travelling wave solutions
of a linear system persist if a generic (analytic) nonlinearity is 
added to the 
force law (see theorems \ref{Tsmr3.1} and \ref{Tass.1}). Furthermore,
we describe the exceptional set of nonlinear force laws which are not
covered by our theorem. They form a countable collection of four-dimensional
sub-manifolds in the space of real analytic functions.

During the past ten years a number of KAM theorems have been developed
for infinite-dimensional systems
(see e.g. \cite{Kuk}, \cite{Po1}, \cite{CW1}, \cite{B1}). 
Although none of them implies our result, our proof uses many
ideas introduced by
Craig and Wayne \cite{CW1} and further developed by Bourgain 
\cite{B1}, \cite{B4}. In particular, we follow the analysis of 
\cite{B1} and \cite{B4}. In order to show our result
the techniques of \cite{B1}, \cite{B4} had to be modified and refined.
For example, we extend the analysis of \cite{B1}
to obtain enough control on the geometry of
certain resonant sets such that the amplitude -- frequency modulation,
provided by the nonlinear dispersion relation, can be used to avoid strong
resonances. Therefore we are not forced to introduce external
parameters in our equation which would have no physical justification. Furthermore, we present proofs
for a number of statements which were formulated in \cite{B1}, \cite{B4}
without proof.

From a technical point of view the essential part of the proof is to
control the inverse of families of matrices $T = D + R$ where
$D$ denotes a diagonal matrix and $R$ is a Toeplitz matrix. Although the
matrix $R$ is small, it is 
not dominated by $D$ due to the small divisor problem.
To obtain estimates on $T^{-1}$ one employs a multi-scale analysis which 
is similar to the one 
introduced by Fr\"ohlich and Spencer \cite{FS} in the context
of Anderson localization. 
We use a new version of the Fr\"ohlich -- Spencer technique containing 
ideas of Bourgain \cite{B1} and of the author \cite{Kri}.  
\end{abstract}

%\clearpage
%\setcounter{page}{5}
\tableofcontents

\clearpage
%\pagenumbering{arabic}

\begin{center}
{\huge \bf Chapter I} \vspace{1.2cm}\\
{\huge \bf Introduction} \vspace{2cm}
\end{center}

\section{The dynamics of nonlinear chains}
\label{int}

In this paper we study the dynamics of nonlinear, infinite 
chains with nearest neighbor interactions which are governed by a
system of differential
equations of the form
\begin{eqnarray}
\label{int.10}
\frac{d^2}{dt^2} x_n(t) = F(x_{n-1}(t) - x_{n}(t)) -
F(x_{n}(t) - x_{n+1}(t)), \;\;\; n \in {\Bbb Z}.
\end{eqnarray}
The simplest physical model for such a system is a chain of particles
each having unit mass, where neighboring particles are connected by 
identical springs.
The displacement of the $n$-th particle at time $t$ is described by $x_n(t)$.
The force transmitted by each spring depends only on the distance of 
the corresponding particles and the force law is prescribed by
the function $F$. In the case that the springs obey Hooke' s law the force 
$F$ is a linear function with $F' > 0$. In this paper we will consider
nonlinear force laws $F$ which are strictly
increasing functions. As one can easily see from (\ref{int.10}) the 
monotonicity
implies that $F$ represents a restoring force, i.e. at all times 
the force on the $n$-th particle drives $x_n$ toward the 
midpoint $(x_{n-1} + x_{n+1})/2$
of its neighboring particles. 

Mechanical systems are not the only applications of
(\ref{int.10}). For example, there is a corresponding electric model
of chains of coupled LC-circuits. In fact, lattice models are quite
frequently used in solid state physics, nonlinear optics, material
sciences, biology and chemistry (see \cite{Mal}, \cite{Tod} and 
references therein).
Also, from a mathematical point of view, lattice systems are of considerable
interest. For example, they arise naturally in the numerical analysis
of partial differential equations either 
by expanding the equation in a suitable basis in function space or
by discretisation of the spatial
variable.

Another attractive feature of lattice models is that they are usually
easy to analyze numerically. They are therefore well suited to discover and 
study nonlinear effects for differential equations. 
In fact, this
was the motivation for the well-known Fermi-Pasta-Ulam experiments \cite{FPU}.
Their rather unexpected findings became the starting point for many 
interesting investigations on the dynamics of nonlinear lattices.
The idea of Fermi, Pasta and Ulam was to apply the emerging
powers of the electronic computing machines to simple nonlinear 
problems which could not be treated analytically. As a first problem
they chose to investigate the dynamics of chains
of 65 particles for three different types of force laws.
They assumed that the first particle $x_0$ and the last particle $x_{64}$
are pinned, i.e. $x_0(t) = 0$, $x_{64}(t) = 0$ for all times.
Recall that in the case of a linear force law
$F(x) = \alpha x$, $\alpha > 0$, 
any initial value problem can be solved explicitly by expanding the solution
in a basis of eigenmodes
$s^{(k)}$, $c^{(k)}$, $1 \leq k \leq 63$, 
where
\begin{eqnarray}
\nonumber
s^{(k)}_n (t) &=& \sin (\lambda_k t) \sin (n \frac{k \pi}{64} ) , \;\;\;
c^{(k)}_n (t) = \cos (\lambda_k t) \sin (n \frac{k \pi}{64} ), \;\;\;
\mbox{ and } \\
\nonumber
\lambda_k &:=& 2 \sqrt{\alpha} \sin (\frac{k \pi}{128}).
\end{eqnarray}
For linear systems energy is not exchanged between modes corresponding
to different indices $k$. 
The general belief among physicists at the time of the FPU experiments 
was that the solutions of generic nonlinear systems exhibit 
ergodic behavior (ergodicity hypothesis). 
To study this effect
Fermi, Pasta and Ulam decided to investigate 
thermalization, i.e. the exchange of energy between
the (linear) eigenmodes which according to the ergodicity hypothesis
should eventually lead to equidistribution of energy
between all the modes. For three classes of nonlinear force laws
$F$ (linear + quadratic, linear + cubic, linear + piecewise linear)
they studied the distribution of energy between the different modes,
where at time $t=0$ all energy was placed in the lowest eigenmode. 
They were surprised to find that thermalization did not occur.
In fact, energy was only exchanged between a few modes and in a
rather regular fashion. Moreover, after some intermediate time (at the order
of about 100 cycles, i.e. $t \approx 100 (2 \pi/\lambda_1)$) the system
was very close to its initial state. They concluded that
the ``prevalent beliefs in the universality of mixing and thermalization
in nonlinear systems may not be always justified'' \cite{FPU}. 
Indeed, at almost the same time as the Fermi-Pasta-Ulam experiments
were conducted the ergodicity hypothesis received another blow.
Kolmogorov announced at the ICM 1954 his result which 
implies, roughly speaking, that there are open sets of Hamiltonian systems 
(in the neighborhood of completely integrable systems) for which
the dynamics cannot be ergodic since there exist invariant sets of small but
positive measure.
Therefore ergodicity is not a generic property of nonlinear Hamiltonian
systems.
Since Fermi, Pasta and Ulam worked in a regime where the amplitudes
and hence the effects of the nonlinearity were moderate it is conceivable
\cite{Way}
that they were observing traces of the quasi-periodic motion whose existence
Kolmogorov had announced in 1954 (the details of the proof of Kolmogorov's 
result were subsequently
given by Arnold in the analytic case and by Moser in the $C^{\infty}$ -- case
leading to what is now called KAM theory). It was later observed numerically
(see \cite{Tod}) 
that energy sharing between modes is enhanced by increasing the 
total energy of the system. However, until today no general
result on the long time behavior of nonlinear Hamiltonian systems has been 
proved and
despite a number of interesting examples and results the original questions
of Fermi, Pasta and Ulam on ergodicity, mixing and thermalization of nonlinear
conservative systems are not resolved.

The findings of
Fermi, Pasta and Ulam led to a series of numerical studies by
a number of people. One of the striking features which were observed was 
that nonlinear lattices allow not only the propagation of periodic waves 
but also
the propagation of localized pulses (solitons), a genuinely nonlinear 
phenomenon. Although some of the numerical results 
could be explained by considering the 
Korteweg -de Vries equation as a continuum limit of the lattice equation 
there was hardly any rigorous analysis.
In order to overcome the lack of analysis, Toda
set out to find a particular nonlinear force law for which one could write
down explicit travelling wave solutions. By some ingenious 
considerations \cite{Tod2} (see also \cite{Tod}) he found that for 
$F(x) = e^x$ (Toda lattice) 
there exist periodic wave solutions which can be expressed
in terms of elliptic functions. Toda's success was more thorough than one 
would have expected. Flaschka \cite{Fla}, \cite{Fla2} and 
independently Manakov \cite{Man}
found that the Toda lattice belongs to the select 
class of completely integrable systems. Not only the dynamics of
finite lattices can be analyzed completely (see e.g. \cite{Mos})
but also the dynamics of infinite chains can be investigated via the
inverse spectral transformation. As Toda had hoped, many rigorous
results on the dynamics of the Toda lattice have been established,
shedding light on the behavior of general nonlinear
lattices. During the last ten years, motivated in part by the 
success of the Toda lattice, analytic results were also obtained for
non-integrable lattices. 
Most of these results concern the propagation of waves. In the 
non-integrable case the main focus so far has been in establishing
the existence of travelling wave solutions (as an exception see \cite{FP}
for a remarkable stability result for solitons). In the integrable case,
in contrast, one can analyze the long time behavior of 
solutions for a large class of
initial value problems.

We will now describe some of the results
obtained for nonlinear lattices. Hereby we restrict ourselves
to infinite lattices with nearest neighbor interactions. 
In particular, our discussion will not include another famous lattice model,
the Calogero -- Moser lattice \cite{Cal}, \cite{Mos2}
where all the particles interact with each other
and which again constitutes an integrable system.
Also, we will not discuss the dynamics of higher dimensional lattices
(see \cite{Mal} and references therein).

We start by briefly explaining the integrability of the Toda lattice
as discovered by Flaschka \cite{Fla} and Manakov \cite{Man}.
We introduce Flaschka' s variables
\begin{eqnarray}
\nonumber
a_n(t) := - \frac{\dot{x}_n(t)}{2}, \quad 
b_n(t) := \frac{1}{2}e^{\frac{1}{2}(x_n(t) - x_{n+1}(t))}.
\end{eqnarray}
Set $L(t)$ to be the symmetric tridiagonal matrix with entries $a_n(t)$ on the 
diagonal and off-diagonal entries $b_n(t)$. Furthermore, set
$B(t) := L_+(t) - L_+^T(t)$ where $L_+$ denotes the upper triangular part of 
$L$. The equation
for the Toda lattice can be written in Lax-pair form, i.e.
\begin{eqnarray}
\nonumber
(x_n, \dot{x}_n) \mbox{ solves Toda } 
\Leftrightarrow
\dot{L} = [B, L] \; (:= BL - LB).
\end{eqnarray}
The consequences of this observation are striking. It is not difficult to
see that a family of operators $L(t)$ which satisfies 
$\dot{L}(t) = [B(t), L(t)]$
for some family of skew-symmetric operators $B(t)$
is unitarily equivalent, i.e. 
for every $t$ there exists an unitary operator $Q(t)$ such that
$L(t) = Q(t) L(0) Q^{*}(t)$. This implies in particular that the 
spectrum of $L(t)$ is independent of $t$.
Moreover, recall that tridiagonal symmetric operators are 
determined by the spectrum
together with a set of additional data (the first components of the 
eigenvectors in the matrix case, the auxiliary spectrum in the 
periodic case and scattering data in the case of infinite operators 
with entries $a_n$, $b_n$ converging as $|n| \to \infty$).
In all three cases one can show that {\em all} spectral data 
needed to reconstruct the operator $L(t)$ evolve simply in time
and can be deduced from the corresponding data at time $t=0$
by explicit formulae. This implies, for example, that we can use the
spectral transformation
to solve initial value problems. Starting with initial
conditions $x_n(0)$, $\dot{x}_n(0)$ one defines the operator $L(0)$
via Flaschka's variables. Then one has to determine the spectral 
data for $L(0)$. The corresponding spectral 
data for $L(t)$ are obtained by explicit formulae. Applying 
the inverse spectral transformation yields the values 
for $x_n(t)$. Instead of solving a dynamics problem we now have to study
the spectral transformation and its inverse. 
Such a remarkable method for solving nonlinear differential systems was
first discovered through the pioneering work of Kruskal et al. \cite{GGKM}
and Lax \cite{Lax} for the Korteweg - de Vries equation.
Nonlinear partial differential equations which can be written in 
Lax-pair form are often called integrable partial differential equations.
There is quite a number of such equations 
(see e.g. \cite{AC}).
In addition to the Toda lattice and the Korteweg -- de Vries equation,
the nonlinear Schr\"odinger equation and the Boussinesq 
equation belong to this class.
It turned out that in the scattering case one can analyze the 
inverse spectral transformation very efficiently by formulating
it as a matrix -- valued Riemann -- Hilbert problem. Deift and
Zhou \cite{DZ1}, \cite{DZ2} have introduced a nonlinear steepest 
descent -- type
method to analyze parameter dependent Riemann -- Hilbert problems.
Their technique was used and further developed in \cite{DIZ} -- \cite{DKMVZ3},
\cite{DVZ}, \cite{BDJ}. 
Matrix -- valued
Riemann -- Hilbert problems not only appear in the theory of
integrable systems but have recently been used to establish
new asymptotic results in such diverse fields as statistical mechanics 
\cite{DIZ}, combinatorics \cite{BDJ}, 
orthogonal polynomials \cite{Dei}, \cite{DKMVZ2}, \cite{KM}
and random matrix theory
\cite{DKMVZ1}, \cite{DKMVZ3}.

Let us return to the results for the Toda lattice.
In one direction, Deift and T-R McLaughlin \cite{DMcL} 
study the continuum limit of the Toda lattice. Their work yields
an important case study
to understand the relation between solutions of hyperbolic systems
which develop shocks and the solutions of the corresponding spatially
discretized systems.
In another direction, towards an understanding of wave propagation
in nonlinear lattices, the dynamics of driven semi-infinite lattices
was analyzed. Since this particular class of initial
boundary value problems is also the starting point for the present
work, we explain it now in more detail.

Consider the following class of initial boundary value problems
\begin{eqnarray}
\nonumber
\ddot{x}_n(t) &=& F(x_{n-1}(t) - x_{n}(t)) - F(x_{n}(t) - x_{n+1}(t)),
\;\;\; n \geq 1,\\
\label{int.100}
x_{0}(t) &=& h(t), \;\;\; t \geq 0\\
\nonumber
x_n(0) - nd &=& {\cal O}(e^{- \alpha n}), \;\;\; \dot{x}_n(0) = 
{\cal O}(e^{- \alpha n}) \quad (d \in {\Bbb R}, \alpha > 0),
\end{eqnarray}
which describes the motion of a
semi-infinite lattice which is driven from one end
by particle $x_0$ of prescribed motion $h(t)$.
In numerical experiments Holian and Straub \cite{HS}
made a remarkable discovery.
They investigated the shock problem where the driver moves
toward the lattice at constant speed $a$ (i.e. $h(t) = a t$, $a > 0$). 
They observed that for a large class of nonlinear force laws $F$,
there was a critical speed $a^{*} \in (0, \infty)$.
For sub-critical driving speeds $0 < a < a^{*}$, the $n$-th particle 
$x_n(t) \to a t + n d'$ as time $t \to \infty$. This means, in the
frame of the driving particle the lattice eventually tends to a quiescent
lattice with some new lattice spacing $d'$.
On the other hand, for supercritical driving speeds ($a > a^{*}$) the 
particles, viewed again in the frame of the driver, keep oscillating for all
times. The energy generated behind the shock front is never
entirely dispersed. 
This effect is genuinely nonlinear;
linear lattices exhibit sub-critical behavior for all driving speeds $a > 0$.
In the case of the Toda lattice the shock problem could be analyzed in
detail using integrability
\cite{HFMcL}, \cite{VDO}, \cite{Kam1}, \cite{Kam2}. 
In fact, not only the existence of 
a critical shock speed could be explained, but also a detailed description
of the propagation of the shock, and of the region of modulated waves 
behind the shock was derived. 
More recently,  
the rarefaction problem was investigated as well. Here
one withdraws the driver from the lattice at constant speed $a < 0$.
Again there exists a critical speed $a_{*} \in (-\infty, 0)$. For
sub-critical speeds $a_{*} < a < 0$ the lattice tends to a rarefied but
quiescent lattice moving at the speed of the driver. For  
supercritical speeds $a < a_{*}$, however, the driver moves away from 
the bulk of the lattice which then behaves like a free semi-infinite lattice.
The distances between neighboring particles tend to infinity at a
logarithmic rate. For the Toda lattice this genuinely nonlinear behavior
was established in \cite{DKKZ} for the supercritical case. 

Let us now assume that the driver in system (\ref{int.100}) 
undergoes a periodic motion
\begin{eqnarray}
\label{int.150}
h(t) = p(\gamma t),
\end{eqnarray}
where $p$ denotes a $2 \pi$-periodic function and the parameter
$\gamma > 0$ represents the frequency of the oscillations of the
driver. Although the corresponding initial boundary value problem
seems very similar to the shock and rarefaction problem described above, 
the crucial difference for the analysis is that the periodically driven
semi-infinite lattice is
not integrable for any nonlinear force law, including the Toda lattice. 
In \cite{DKV} problem (\ref{int.100}), (\ref{int.150})
was studied numerically. 
We observed the following time asymptotic behavior. 
There exists a sequence of thresholds
$\infty =: \gamma_0 > \gamma_1 > \gamma_2 > \ldots > 0$, such that for 
driving frequencies $\gamma_{\nu} > \gamma > \gamma_{\nu+1}$ the solutions 
of (\ref{int.100}), (\ref{int.150}) are well described by travelling
$\nu$-phase waves, as time becomes large, except for a boundary layer close 
to the driver. More precisely, there exists a function $\chi$ defined on
the $\nu$-dimensional torus ${\Bbb T}^{\nu} \equiv 
({\Bbb R}/ 2 \pi {\Bbb Z})^{\nu}$ and real numbers $\tilde{d}$,
$\omega_1, \ldots, \omega_{\nu}$ such that for $t/n$ large 
\begin{eqnarray}
\nonumber
x_n(t) \approx x_n^{as}(t) \equiv
\tilde{d} n +
\chi(\omega_1 n - \gamma t, \omega_2 n - 2 \gamma t, 
\ldots, \omega_{\nu} n -
\nu \gamma t) + {\cal O}(e^{-\beta n}), \;\;\; n \geq 1,
\end{eqnarray}
for some $\beta > 0$. One of the goals of \cite{DKV} was to construct
the numerically observed limit cycles $x^{as}$ for small amplitudes
of the driving function $h$ by perturbative methods. It turned out that
the main difficulty in the construction were resonances 
($l \gamma$ is contained in the continuous spectrum
of some associated linearized operator for $l = 0, \pm 1, \ldots, \pm \nu$)
leading to an over-determined equation. In order to overcome this 
problem it was necessary to obtain sufficiently ample families
of travelling wave solutions of the {\em doubly infinite lattice}, i.e.
solutions of the form
\begin{eqnarray}
\label{int.200}
x_n(t) = \chi(\omega_1 n - \gamma t, \omega_2 n - 2 \gamma t, 
\ldots, \omega_{\nu} n - \nu \gamma t), \;\;\; n \in {\Bbb Z},
\end{eqnarray}
where again $\chi$ is a function on the $\nu$ -- dimensional torus 
${\Bbb T}^{\nu}$. In the case $\nu=1$ such families of travelling
wave solutions (with small amplitudes) can be constructed
by a standard Lyapunov -- Schmidt reduction argument (see
\cite{DKV}). However, for $\nu \geq 2$ a small divisor problem
occurs. It is the goal of the present paper to overcome this 
small divisor problem and to show that for generic nonlinear real analytic
force laws $F$ there exist uncountably many solutions of (\ref{int.10})
of type (\ref{int.200}). Note that for linear lattices such solutions
can be obtained by a simple superposition of travelling one-phase
waves (see section \ref{smr1} below) and in the case of the Toda lattice
($F(x) = e^{x}$) such solutions can be expressed explicitly in 
terms of theta-functions (see \cite{Krc}, \cite{DKV}).

Small divisor problems have a long history in the theory of 
nonlinear dynamical systems. In fact, many investigations in the field
of celestial mechanics were aimed at understanding the effects
of resonances or near resonances which are the physical reason behind the  
occurrence of small divisors (see \cite{Gre} for the history of the
three body problem in celestial mechanics).
A general method to overcome such problems was presented only 
in the middle of this century, when the above mentioned
KAM theory was developed for finite dimensional
Hamiltonian systems. The basic idea is to use a rapidly
converging iteration scheme to balance the effects of the
small divisors. In recent years this idea has also been extended
to infinite dimensional systems, establishing the 
existence of finite dimensional invariant tori.
In infinite dimensions, two somewhat
different methods have been most successful.
The first method, developed by Kuksin \cite{Kuk} and P\"oschel 
\cite{Po1} (see also \cite{KP},
\cite{Po2})
extends the 
classical approach of KAM-theory to iteratively construct a normal form
using a sequence of symplectic transformations 
(cf. \cite{Eli}, \cite{Po3}, \cite{W3}). 
The second method, introduced by 
Craig and Wayne in \cite{CW1}, uses Fourier analysis to obtain 
a nonlinear equation on a lattice. They apply a Lyapunov -- Schmidt 
reduction to this equation. The infinite dimensional part of the
reduction contains the small divisor problem which is overcome
by a modified Newton scheme in the spirit of the Nash -- Moser 
technique. 
The main difficulty in this procedure lies in inverting
matrices of the form $D+R$, where $D$ is a diagonal matrix dominating 
$R$, {\it except} for some diagonal entries $D(n, n)$ of small absolute
values which represent the 
small divisors. The lattice sites $n$ where $|D(n, n)|$
is of order ${\cal O}(|R|)$ are called singular sites. To control 
$(D+R)^{-1}$ the separation of singular sites plays an important role.
In the work of Craig and Wayne the existence of 
families of time-periodic solutions
was established for $(1 + 1)$ -- dimensional nonlinear wave equations
\cite{CW1}
and for $(1 + 1)$-dimensional nonlinear Schr\"odinger equations \cite{CW2}.
For the construction of quasi-periodic solutions (in time), however,
the singular sites are less separated and a new idea was needed.
Such an idea was provided by Bourgain \cite{B1} who introduced a multi-scale 
analysis to overcome this difficulty. He  
established the existence of quasi-periodic solutions for 
$(1+1)$ -- dimensional wave equations \cite{B1} and for $(1+1)$ -- and 
$(1+2)$ -- dimensional nonlinear Schr\"odinger equations \cite{B1}, \cite{B4}
(cf. \cite{B3} for a finite dimensional application of his technique).
The estimates on the inverse matrices $(D+R)^{-1}$ in the 
work of Craig -- Wayne and in the work of Bourgain are based on a technique
which was first introduced by Fr\"ohlich and Spencer \cite{FS} in the 
theory of Anderson localization, where similar inversion problems occur.
Note furthermore that the idea of using a Lyapunov-Schmidt reduction
to prove the existence of quasi-periodic solutions was already used
by Scheurle \cite{Sch1}, \cite{Sch2} in the context of finite dimensional
systems.

In our proof of the existence of travelling multi-phase waves
we follow the construction of Bourgain \cite{B1}. In order to obtain
our result we will modify and refine his method at various points.
Hereby an important role is played by a new version of the 
coupling lemma in the Fr\"ohlich -- Spencer technique which the
author first presented in \cite{Kri}. A more detailed description 
of our proof will be given in section \ref{ova} after the statement
of the main result in section \ref{smr}.

We conclude this section by briefly describing related results on
the existence of travelling waves in nonlinear (non-integrable) 
infinite chains with nearest neighbor interactions. 
In contrast to our theorem, all results
obtained previously concern the construction of one-phase travelling waves
\begin{eqnarray}
\label{int.300}
x_n(t) = \chi(\omega n - \gamma t). 
\end{eqnarray}
These results are either of perturbative nature
(i.e. small amplitude solutions) or they are derived by variational
methods. The existence of solitary travelling wave solutions 
with prescribed potential energy
was shown by Friesecke and Wattis \cite{FW} for a large class of
force laws. Later Smets and Willem \cite{SW} used a different 
variational formulation to construct solitary waves with given
wave speeds (see \cite{Sme} for extensions).
Periodic travelling waves were constructed 
by Filip and Venakides \cite{FV}. All of the just mentioned results
do not require smallness of the solutions.

Regarding small amplitude solutions 
the above mentioned paper of Friesecke and Pego \cite{FP} yields a 
remarkably detailed
analysis of solitary waves at near sonic speeds, including their
stability properties. Small amplitude periodic travelling waves 
were constructed by Deift, Kriecherbauer, Venakides 
in \cite{DKV}. Georgieva, Kriecherbauer and Venakides \cite{GKV1}, \cite{GKV2}
have extended the analysis of \cite{DKV} to diatomic chains with 
periodically varying masses (of period 2) discovering interesting
resonance phenomena which might provide a mechanism for frequency doubling.

Finally, we mention that the existence of travelling single-phase waves 
has also been established for a related system of coupled oscillators
which is described by 
\begin{eqnarray}
\label{int.310}
\ddot{x}_n = x_{n-1} + f(x_n) + x_{n+1}.
\end{eqnarray}
%c%The existence of solitary waves was proved by Zinner \cite{Zin}
%c%for a special class of functions $f$ via topological methods.
A general and powerful method to determine all (small amplitude) solutions
of (\ref{int.310}) of the form (\ref{int.300}) was recently presented
by Iooss and Kirchg\"assner \cite{IK}. They use their technique to construct
solutions with interesting and unexpected profile functions $\chi$.

\section{Statement of the main result}
\label{smr}
In this section we formulate our main result. We will construct
solutions of 
\begin{eqnarray}
\label{smr.10}
\ddot{x}_n(t) = F(x_{n-1}(t) - x_n(t)) - F(x_n(t) - x_{n+1}(t)), \;\; 
n \in {\Bbb Z}, \; t \in {\Bbb R},
\end{eqnarray}
which are of the form
\begin{eqnarray}
\label{smr.20}
x_n(t) = n b + \chi (n \omega_1 - \gamma t, \ldots, n \omega_{\nu} - 
\nu \gamma t)
\end{eqnarray}
for some $\nu \in {\Bbb N} \setminus \{ 1 \}$, 
$b$, $\gamma$, $\omega_i \in {\Bbb R}$ and
$\chi: ({\Bbb R} / {2 \pi \Bbb Z})^{\nu} \equiv {\Bbb T}^{\nu} 
\to {\Bbb R}$.
Solutions of this form describe time-periodic travelling multi-phase waves 
where
$\nu$ denotes the number of phases. 
As explained in section \ref{int} there is numerical evidence \cite{DKV}
that such solutions exist for a large class of force laws $F$.
Note that we exclude the case $\nu = 1$ from
our considerations because in this case there is no small divisor problem  
and the existence of periodic single-phase waves 
has already been proved in \cite{DKV}. 
Roughly speaking, our
main result shows that for generic analytic 
force laws $F$, sufficiently
low frequencies $\gamma > 0$ and most values of averaged
lattice spacings $b$ there exist uncountably many solutions of type
(\ref{smr.20}) with small amplitudes. These solutions can be  
parameterized smoothly by $z \in {\cal Z} \subset {\Bbb C}^{\nu}$
where ${\cal Z}$ is a set of relative large measure
in a neighborhood of the origin. In the following
we will make this statement precise. In order to understand the origin of
this class of travelling waves we first discuss the linear case.   

\subsection{The linear lattice}
\label{smr1}

Suppose $F(x) = \alpha x$ with $\alpha > 0$. A function
\begin{eqnarray}
\label{smr1.5}
x_n(t) = e^{i(\omega n - \gamma t)}
\end{eqnarray} 
solves (\ref{smr.10}) if and 
only if the frequencies $\gamma$ and $\omega$ satisfy the 
following dispersion relation
\begin{eqnarray}
\label{smr1.10}
\gamma^2 = 4 \alpha \sin^2 \frac{\omega}{2}.
\end{eqnarray}
For given $\gamma > 0$ there exists an unique 
$\nu \in {\Bbb N}_0$ with 
\begin{eqnarray}
\label{smr1.20}
(\nu \gamma)^2 \leq 4 \alpha < ((\nu+1) \gamma)^2.
\end{eqnarray}
Note that in the case $\nu = 0$ (i.e. $\gamma > 2 \sqrt{\alpha}$) 
no bounded solution of (\ref{smr.10}) of the form (\ref{smr1.5}) exists.
Assume $\nu \geq 1$ and define
\begin{eqnarray}
\label{smr1.30}
\omega_j^{(0)} :=  2 \arcsin \left( \frac{j \gamma}{2 \sqrt{\alpha}} \right)
\in (0, \pi]
\;\;\; \mbox{ for } j \in \{1, \ldots, \nu \}.
\end{eqnarray}
For $z \in {\Bbb C}^{\nu}$ , $b \in {\Bbb R}$ set
\begin{eqnarray}
\label{smr1.40}
x_n(t;b,z) &:=& n b + \chi_z^{(lin)} (\omega^{(0)} n - g \gamma t), \;\;\;
\mbox{ where }
\\
\label{smr1.50}
\chi_z^{(lin)} &:& {\Bbb T}^{\nu} \to {\Bbb R}, \;\;
\chi_z^{(lin)}(\xi) 
= \sum_{j=1}^{\nu} (z_j e^{i \xi_j} + \overline{z_j} 
e^{-i \xi_j}),   
\\
\label{smr1.60}
\omega^{(0)} &:=& (\omega_1^{(0)}, \ldots, \omega_{\nu}^{(0)}),
\\
\label{smr1.70}
g &:=& (1, 2, \ldots, \nu).
\end{eqnarray} 
The functions $(x_n(\cdot \; ;  b, z))_{n \in {\Bbb Z}}$ provide 
a family of solutions
of (\ref{smr.10}) of type (\ref{smr.20}).

\subsection{Assumptions}
\label{smr2}

In this section we state four assumptions 
which will constitute the hypothesis of our main theorem
\ref{Tsmr3.1}. They concern the force law $F$, the averaged lattice
spacing $b$ and the (time-) frequency $\gamma$.

{\em
{\bf Assumption A1:}
$F$ is real analytic in a neighborhood of $-b$ with $F'(-b) > 0$.
}

%\noindent
%Assumption 1 implies that there exist coefficients $\alpha_j \in {\Bbb R}$
%with
%\begin{eqnarray}
%\label{smr2.10}
%F(y+b) = \sum_{j > 0} \alpha_j y^j \;\; \mbox{ for $y$ close to $-b$ }
%\end{eqnarray}
%and
%\begin{eqnarray}
%\label{smr2.15}
%\alpha_1 > 0.
%\end{eqnarray} 

{\em
{\bf Assumption A2:}
$0 < \gamma < \sqrt{F'(-b)}$ and $2 \sqrt{F'(-b)} \notin \gamma {\Bbb Z}$.
}

\noindent
Assumption A2 implies the existence of $\nu \in {\Bbb N} \setminus \{ 1 \}$
satisfying
\begin{eqnarray}
\label{smr2.20} 
\nu \gamma < 2 \sqrt{F'(-b)} < (\nu+1) \gamma.
\end{eqnarray}
We set 
\begin{eqnarray}
\label{smr2.30}
\omega_j^{(0)} &:=&  2 \arcsin \left( \frac{j \gamma}{2 \sqrt{F'(-b)}} 
\right) \in (0, \pi)
\;\;\; \mbox{ for } j \in \{1, \ldots, \nu \}. \\
\label{smr2.40}
\omega^{(0)} &:=& (\omega_1^{(0)}, \ldots, \omega_{\nu}^{(0)}),
\end{eqnarray}
Note that these definitions are consistent with the corresponding 
definitions in section \ref{smr1}.

{\em {\bf Assumption A3} (diophantine condition):
There exist positive constants $s$ and $\tau$ such that 
\begin{eqnarray}
\label{smr2.50}
\mbox{ dist}(<\omega^{(0)}, m>, 2 \pi {\Bbb Z}) > \frac{s}{|m|^{\tau}}
\quad \mbox{ for all } m \in 
{\Bbb Z}^{\nu} \setminus \{ 0 \}.
\end{eqnarray}
}
We have used the notation $<a, b> = \sum_{i=1}^{\nu} a_i b_i$ for $a$, 
$b \in {\Bbb R}^{\nu}$.
Before we can state assumption A4, we need a few more definitions.
\begin{eqnarray}
\label{smr2.60}
\Lambda_j := \frac{(j \gamma)^2 \cos (\frac{\omega_j^{(0)}}{2})}
{4 \sin^3(\frac{\omega_j^{(0)}}{2})}, \;\;\mbox{ for } 1 \leq j \leq \nu.
\end{eqnarray}
Recall the definition of the integer vector $g$ in (\ref{smr1.70}).
For $m \in {\Bbb Z}^{\nu}$ and $\omega \in {\Bbb R}^{\nu}$ we set
\begin{eqnarray}
\label{smr2.70}
V(\omega)(m) := 
\left\{
\begin{array}{ll}
F'(-b) & \mbox{ if } <g, m> = 0 \\
F'(-b) - \frac{<g, m>^2 \gamma^2}
{4 \sin^2 \frac{<\omega, m>}{2}} & \mbox{ else }.
\end{array}
\right. 
\end{eqnarray}
Here we understand $V(\omega)(m) = \infty$ in the case $<g, m> \neq 0$ and
$<\omega, m> \in 2 \pi {\Bbb Z}$.
%Note that it follows from assumption A3 that $V(\omega^{(0)})(m)$ is finite
%for all $m \in {\Bbb Z}^{\nu}$.
Denote by
\begin{eqnarray}
\label{smr2.80} 
e_j := (0, \ldots, 0, 1, 0, \ldots, 0) \in {\Bbb Z}^{\nu}.
\end{eqnarray}
Set $\Omega^{(\nu)} \equiv \Omega$ 
to be the $\nu \times \nu$ matrix with entries
\begin{eqnarray}
\label{smr2.90}
\Omega^{(\nu)}_{j, l} \equiv
\Omega_{j, l} :=
\left\{
\begin{array}{ll}
3 \Lambda_l^{-1} \left(
F''(-b)^2 \frac{1}{V(\omega^{(0)})(2 e_l)} - F'''(-b)
\right) & \mbox{ if } j=l, \\
2 \Lambda_l^{-1} \left(
F''(-b)^2 \left(
\frac{1}{V(\omega^{(0)})(e_j+e_l)} 
+ \frac{1}{V(\omega^{(0)})(e_j-e_l)}
\right)
- F'''(-b)
\right) & \mbox{ if } j \neq l.
\end{array}
\right.
\end{eqnarray}
It is shown in lemma \ref{Lass.1} that
the entries of $\Omega$ are well defined real numbers.  

{\em 
{\bf Assumption A4} (non-degeneracy condition):
No entry of $\Omega$ vanishes and $\det \Omega \neq 0$.
}
 
In theorem \ref{Tass.1} 
we formulate and prove our claim that 
assumptions A1 -- A4 are generically
satisfied.
 
\subsection{The Theorem}
\label{smr3}

\begin{Theorem}
\label{Tsmr3.1}
Suppose $F: {\Bbb R} \to {\Bbb R}$, $b \in {\Bbb R}$ and 
$\gamma \in {\Bbb R}_+$ satisfy assumptions A1 -- A4.
Let $g$, $\nu$, $\omega^{(0)}$ be defined as in (\ref{smr1.70}) --
(\ref{smr2.40}) above and
let $0 < \kappa < 1$. Then there exist $\rho > 0$, a measurable set
${\cal Z} \subset U_{\rho}(0) \equiv 
\{ z \in {\Bbb C}^{\nu}: |z_j| < \rho\; , \; 1 \leq j \leq \nu \}$ 
and functions
$\omega: U_{\rho}(0) \to {\Bbb R}^{\nu}$,
$\chi: {\Bbb T}^{\nu} \times U_{\rho}(0) \to {\Bbb R}$
such that for every $z \in {\cal Z}$ the
functions
\begin{eqnarray}
\label{smr3.10}
x_n(t) = n b + \chi(\omega(z)n - g \gamma t, z) \;\;\; \quad 
n \in {\Bbb Z}, \;
t \in {\Bbb R},
\end{eqnarray}
define solutions of system (\ref{smr.10}) of type (\ref{smr.20}) and
statements (\ref{smr3.20}) -- (\ref{smr3.50}) below hold. 
\begin{eqnarray}
\label{smr3.20}
&\mbox{ vol}({\cal Z}) \geq \kappa \mbox{ vol}(U_{\rho}(0)) \;\;\; 
\mbox{ (volume with respect to Lebesgue measure in } {\Bbb R}^{2\nu}),& \\
\label{smr3.30}
&\mbox{ the map } z \mapsto \chi(\cdot, z) \mbox{ is one to one },& \\
\label{smr3.40}
&|\omega(z) - \omega^{(0)}| = {\cal O}(|z|^2),& \\
\label{smr3.50}
&\sup_{\xi \in {\Bbb T}^{\nu}} 
\left|\chi(\xi, z) - \chi_{z'}^{(lin)}(\xi) \right| = {\cal O}(|z|^2), \;\;\;
\mbox{ (see (\ref{smr1.50})) and where for } z \in U_{\rho}(0)& \\
\label{smr3.15}
&z' := \frac{i \sqrt{F'(-b)}}{\gamma} \left(
\frac{z_1}{1}, \frac{z_2}{2}, \ldots, \frac{z_{\nu}}{\nu} \right).&
\end{eqnarray}
\end{Theorem}

\begin{Remark}
\label{Rsmr3.1}
{\em
Theorem \ref{Tsmr3.1} together with theorem \ref{Tass.1} implies that
most small amplitude travelling wave solutions (\ref{smr1.40})
of the linear system persist
if generic higher order terms are added to the force law. 
The profile function $\chi(\cdot, z)$ as well as the frequency
vector $\omega(z)$ differ from the corresponding linear quantities
only to order $|z|^2$, where $|z|$ is a measure of the amplitude of
the oscillations (see (\ref{smr3.40}) and (\ref{smr3.50})). 
Statements (\ref{smr3.20}) and (\ref{smr3.30}) imply that we have constructed 
uncountably many solutions
of (\ref{smr.10}) of type (\ref{smr.20}).
}
\end{Remark}

\begin{Remark}
\label{Rsmr3.2}
{\em 
We have required in assumption $A1$ that the force law is analytic. This
will ensure rapid decay on the sequence of Fourier coefficients. 
As in classical KAM theory one might expect
that it suffices to assume some finite regularity. A first step to support
this view for infinite dimensional systems 
has been taken in \cite{CS}. In this paper we do not investigate 
the question of minimal regularity and treat only the analytic case.
The second assumption in A1, $F'(-b)  > 0$ is necessary for the existence
of (bounded) travelling waves in the linear case
as the analysis of section \ref{smr1}
shows. Recall from the first paragraph of section \ref{int} that 
this assumption also has the physical interpretation that $F$ represents a
restoring forces on the lattice. 

Assumption $A2$ can again be understood from the linear case. It implies
that the linearized equation carries multi-phase waves with two or more
phases. Furthermore, assumption A2 excludes those frequencies where
the number $\nu$ of phases changes.

The diophantine assumption $A3$ implies that 
small divisors do not appear for low Fourier modes. 
Such a condition (possibly in a weaker form) is needed to 
start the Newton iteration scheme.

Assumption $A4$ will be used to show that the nonlinear dispersion relation
$\omega(z)$ is not degenerate to second order. This implies,
roughly speaking, that the dependence of $\omega(z)$ on $z$ is sufficiently
large to ensure that we can avoid resonant
parameters $(z, \omega(z))$ by small changes in $z$.
}
\end{Remark}

\section{Overview of the analysis}
\label{ova}

Our proof of theorem \ref{Tsmr3.1} is based on a method which
was first introduced by Craig and Wayne to construct
Cantor families of time-periodic, small amplitude
solutions for nonlinear wave equations \cite{CW1} 
and for  nonlinear Schr\"odinger
equations \cite{CW2}, where the spatial variable lives on a finite 
interval and
either Dirichlet or periodic boundary conditions are prescribed.
The basic strategy of Craig and Wayne is to expand the solution
with respect to the eigenfunction basis of the corresponding
linearized operator. This leads to a nonlinear equation
on a two-dimensional lattice. Then they apply a Lyapunov -- Schmidt 
reduction. The small divisor problem is contained in
the infinite dimensional part of the reduction. In order to solve this 
equation Craig and Wayne devise a modified Newton scheme  
in the spirit of Nash -- Moser implicit function theorems. 
At each step of the iteration scheme the crucial difficulty lies in 
deriving estimates on the inverse of matrices $T(\lambda)$
where $\lambda$ denotes a parameter. More precisely, 
at the $k-th$ step of the iteration $T_k(\lambda)$ represents a family of
$L_k \times L_k$ matrices ($L_k$ grows exponentially with $k$) and
one has to show that, except for a small set of parameters $\lambda$, 
the inverse matrices $T_{k}^{-1}(\lambda)$ can be bounded in a suitable 
norm by $\Phi(L_k)$ for some function $\Phi$ of sub-exponential growth
at $\infty$ (e.g. $\Phi$ polynomial).
The matrices $T_k(\lambda)$ are of the form $D + R$, where $D$ is
a diagonal matrix and $R$ is small with respect to some norm.
Clearly, if the entries of the diagonal matrix
$D$ were bounded away from zero and $R$ were sufficiently small such that 
$\| D^{-1} R \| < 1$, i.e. if we were in the case of diagonal dominance, the
inverse of $D+R$ would exist and could be expanded in a Neumann series. 
However, the existence of small divisors is reflected in the fact that 
for all parameter values $\lambda$ the matrices $T_k(\lambda)$  
have diagonal entries $T_k(\lambda)(n, n)$ of small absolute value
and we are therefore {\em not} in the case of diagonal dominance. 
The corresponding 
lattice points $n$ with $|T_k(\lambda)(n, n)| = {\cal O}(|R|)$ 
are called singular sites.
The location of these singular sites, in particular their mutual
distances, determines how difficult it is to obtain the required estimates
on $T_k^{-1}(\lambda)$, 
if it is possible at all.
In \cite{CW1}, \cite{CW2} Craig and Wayne proved that for most
parameter values $\lambda$, singular sites $n$ (respectively
clusters of two singular sites) are well separated. As a consequence
the effects of different singular sites decouple, i.e. 
one can restrict oneself to investigating submatrices of $T_k(\lambda)$
which contain at most one singular site (respectively one
cluster of singular sites). The inverse of such submatrices can be controlled
by studying the dependence of their small eigenvalues on $\lambda$.

For the construction of quasi-periodic solutions, however, there is less 
separation of singular sites than in the periodic case. This has the effect
that clusters of singular sites cannot be decoupled. The analysis
of Craig and Wayne had to be developed further. This was done by
Bourgain who added a multi-scale analysis to the procedure of Craig and
Wayne. Bourgain was able to prove the existence of quasi-periodic solutions of
nonlinear Schr\"odinger equations and nonlinear wave equations
\cite{B1}, \cite{B4}. 

The results of Craig -- Wayne and Bourgain are based in some form on 
proving exponential localization of the eigenvectors corresponding
to small eigenvalues.
A powerful technique for proving such localization results
for a similar class
of operators was introduced by Fr\"ohlich and Spencer \cite{FS} in the 
theory of Anderson localization. Their method has been used and further
developed by many authors (see e.g. 
\cite{DK}, \cite{FSW}, \cite{P2}, \cite{Sur}, \cite{CW1}, \cite{B1},
\cite{Gol}, \cite{B4}, \cite{Kri}) including the above mentioned
work of Craig -- Wayne and Bourgain.

In the remainder of this section we explain in more detail how the method of
Craig -- Wayne and Bourgain can be used to prove theorem \ref{Tsmr3.1}.
In order to keep the overview of
the analysis as non-technical as possible we will introduce 
some simplifications.
In remarks \ref{Rova.1} \ref{Rova.2}, \ref{Rova.3}, \ref{Rova.5} 
below we describe the 
modifications which are necessary to treat the full problem. At the
end of this section we outline the plan of this paper.

The first step of our analysis is to perform a Fourier transformation
of the equation. To this end we expand the function $\chi(\cdot, z)$
(see (\ref{smr3.10})) 
\begin{eqnarray}
\nonumber
\chi(\xi, z) = \sum_{m \in {\Bbb Z}^{\nu}} \tilde{u}_m(z) e^{i<m, \xi>}
= \sum_{m \in {\Bbb Z}^{\nu}} \frac{u_m(z)}{-2i \sin(<\omega(z), m>/2)} 
e^{i<m, \xi>}.
\end{eqnarray}
Note that we will
work with the modified Fourier coefficients $u_m(z)$ rather than
with the sequence $\tilde{u}_m(z)$ itself, mainly for the reason that the 
nonlinear part of the equation takes a simpler form when regarded as a 
function of $(u_m)$. The letter $z$ denotes a parameter of the construction.
The interpretation of this parameter is given by the relation
\begin{eqnarray}
\nonumber
u_{e_j}(z) = z_j, \quad 1 \leq j \leq \nu,
\end{eqnarray}
which will hold by definition.
It is our task to determine the nonlinear dispersion relation $\omega(z)$ and
the sequence of modified Fourier coefficients 
$(u_m(z))_{m \in {\Bbb Z}^{\nu}}$, for a large family of 
parameters $z \in {\Bbb C}^{\nu}$, in such a way that the
corresponding multi-phase waves (defined through (\ref{smr3.10})) solve
equation (\ref{smr.10}). In section \ref{fa} below we derive
the equations for $\omega$ and $(u_m)$ 
(see proposition \ref{Pfeq}) which can be written in the form
\begin{eqnarray}
\label{ova.80}
(D(\omega) u + W(u))_m = 0 \;\;\; \mbox{ for } m \in {\Bbb Z}^{\nu}
\setminus \{ 0 \}.
\end{eqnarray}
Here $D(\omega)$ denotes a diagonal matrix,
\begin{eqnarray}
\label{ova.85}
D(\omega)(m, n) = \delta_{m, n} 
\left(
F'(-b) - \frac{\gamma^2 <m, g>^2}{4 \sin^2 \frac{<\omega, m>}{2}}
\right),
\end{eqnarray}
and represents the linearized equation, whereas $W(u)$ contains the 
nonlinear part of the equation with $W(u)$ = ${\cal O}(|u|^2)$.
Observe that $D(\omega)(m, m)$ is not defined for 
$<\omega, m> \in 2 \pi {\Bbb Z}$. As it turns out this causes only minor
difficulties in our analysis (see remark \ref{Rova.5} below). In case
$<\omega, m>$ $\in 2 \pi {\Bbb Z}$ we shall
define $D(\omega)(m, m) = \infty$ if $<m, g> \neq 0$, and
$D(\omega)(m, m) = F'(-b)$ if $<m, g> = 0$. 
Using the symmetries of the system (see remark \ref{Rova.1}) we may
restrict ourselves to search for solutions
$(u_m)_{m \in {\Bbb Z}^{\nu}} \equiv
(u(m))_{m \in {\Bbb Z}^{\nu}}$ of (\ref{ova.80}) which lie in the set
\begin{eqnarray}
\label{ova.87}
\{ u: {\Bbb Z}^{\nu} \to {\Bbb C} \; : \; u(0)=0, \;u(m)=u(-m) \in {\Bbb R}
\; \mbox{ for } m \in {\Bbb Z}^{\nu} \}.
\end{eqnarray}

In order to apply a Lyapunov -- Schmidt reduction to (\ref{ova.80})
we determine the kernel of the diagonal operator $D(\omega^{(0)})$
where $\omega^{(0)}$ (see (\ref{smr2.30}), (\ref{smr2.40}))
satisfies the dispersion relations of the linearized equation. 
It follows from the definition of $\omega^{(0)}$ that $D(\omega^{(0)})(m, m)
 = 0$ for all $m \in {\cal S}$ where
\begin{eqnarray}
\label{ova.90}
{\cal S} := \{e_j, -e_j: 1 \leq j \leq \nu \}.
\end{eqnarray}
In order to show that $D(\omega^{(0)})(m, m) \neq 0$ for all $m \in
{\Bbb Z}^{\nu} \setminus {\cal S}$ one needs to employ the 
diophantine condition of assumption A3. Hence the dimension of the kernel of 
$D(\omega^{(0)})$ is $2 \nu$ and
\begin{eqnarray}
\nonumber
\ker D(\omega^{(0)}) = \{u: u(m) = 0 \mbox{ for } m \in
{\Bbb Z}^{\nu} \setminus {\cal S} \}.
\end{eqnarray}
Define $Q$ to be the projection on the kernel
\begin{eqnarray}
\label{ova.100}
(Q u)(m) :=
\left\{
\begin{array}{ll}
u(m) & \mbox{ if } m \in {\cal S}, \\
0 & \mbox{ else }.
\end{array}
\right.
\end{eqnarray}
We denote by $P$ the complementary projection
\begin{eqnarray}
\label{ova.110}
(P u)(m) :=
\left\{
\begin{array}{ll}
u(m) & \mbox{ if } m \in {\Bbb Z}^{\nu} \setminus ({\cal S} \cup \{ 0 \}), \\
0 & \mbox{ if } m \in {\cal S} \cup \{ 0 \}.
\end{array}
\right.
\end{eqnarray}
Note that $P = Id - Q$ holds for all sequence spaces which require
$u(0) = 0$. Since we look for solutions in such a class (see (\ref{ova.87}))
we may decompose  
\begin{eqnarray}
\nonumber
u = Qu + Pu \equiv x + v.
\end{eqnarray}
With this notation one can write equation (\ref{ova.80}) in the following way
\begin{eqnarray}
\label{ova.112}
D(\omega) x + Q W (x + v) = 0&&\mbox{(Q -- equation)},\\ 
\label{ova.115}
D(\omega) v + P W (x + v) = 0&&\mbox{(P -- equation)}.
\end{eqnarray}
In a standard Lyapunov -- Schmidt reduction the operator $D(\omega^{(0)})$
has a bounded inverse on the $P$ -- projection of the space.
Therefore the second equation can be solved by the implicit
function theorem which determines $v = v(x, \omega)$ as a function of 
$x$ and $\omega$. Inserting $v(x, \omega)$ into the first equation
one obtains a degenerate, finite dimensional equation depending only
on $x$ and $\omega$. This equation is often called the bifurcation equation.
Typically, the bifurcation equation can be solved by determining
the nonlinear dispersion relation
$\omega = \omega(x)$. In this way one has constructed
a family of solutions, parameterized by (small) $x$ which 
correspond to solutions of the linearized equation.
This procedure works without any difficulties in the case $\nu = 1$ 
leading to small amplitude single-phase travelling waves \cite{DKV}.
However, in the case $\nu \geq 2$ considered in this paper it
is clear from (\ref{ova.85}) that $D(\omega^{(0)})$ does not have a
bounded inverse on the $P$ -- projected space because the diagonal entries 
of $D(\omega^{(0)})$ accumulate at $0$, reflecting the small divisor problem.
As a result we will not be able to solve the $P$ -- equation
for all values $(x, \omega)$ in a neighborhood of $(0, \omega^{(0)})$
but for a Cantor -- type subset ${\cal N}^{\infty}$ of relative large measure. 
Since we also need to satisfy the bifurcation equation (\ref{ova.112}), 
not only the measure
but also the geometry
of the set ${\cal N}^{\infty}$ must be controlled. Indeed, we have to show
that the subset $\{ (x, \omega(x)): x \mbox{ small } \}$ 
of codimension $\nu$ lies in 
${\cal N}^{\infty}$ for most values of $x$.

We now describe the analysis of the $P$ -- equation in more detail. 
In the reduced space (\ref{ova.87}) we can parameterize the 
kernel of $D(\omega^{(0)})$ by $\varphi(a)$, $a \in {\Bbb R}^{\nu}$, where
\begin{eqnarray}
\label{ova.120}
\varphi(a)(m) := 
\left\{
\begin{array}{ll}
a_j & \mbox{ if } m \in \{e_j, -e_j\} \\
0   & \mbox{ else } 
\end{array}
\right.
\end{eqnarray}
Furthermore we denote
\begin{eqnarray}
\label{ova.130}
\lambda &:=& (a, \omega), \\
\label{ova.140}
\lambda^{(0)} &:=& (0, \omega^{(0)}), \\
\label{ova.150}
{\cal Q}(v, \lambda) &:=& D(\omega) \varphi(a) + Q W( \varphi(a) + v ),\\
\label{ova.160}
{\cal P}(v, \lambda) &:=& D(\omega) v + P W( \varphi(a) + v ).
\end{eqnarray}

Chapter II is devoted to solving the $P$ -- equation. More precisely, we
show that there exist a positive number $\rho_1$, a smooth function
$v(\lambda)$, defined for $\lambda \in {\cal N}^{(1)} \equiv 
B_{\rho_1}(\lambda^{(0)})$, 
and a subset ${\cal N}^{\infty}$ of ${\cal N}^{(1)}$
such that $v(\lambda) = P v(\lambda)$ and 
${\cal P}(v(\lambda), \lambda) = 0$ for all $\lambda \in {\cal N}^{\infty}$
(cf. theorem \ref{TP}).
Following Craig and Wayne we construct $v$ by
a modified Newton scheme, $ v :=  \lim_{j \to \infty} v_j$ where
\begin{eqnarray}
\nonumber
v_1 (\lambda) &=& 0, \\
\label{ova.200}
v_{j+1}(\lambda) &=& v_j(\lambda) - T_j^{-1}(\lambda)
{\cal P}(v_j(\lambda), \lambda), \;\;\; \mbox{ with } \\
\nonumber
T_j (\lambda) &=& [D(\omega) + DW(\varphi(a) + v_j(\lambda)]|_{\Lambda_j}
\end{eqnarray}
and $\Lambda_j$ is a nested sequence of finite lattices
exhausting the set 
${\Bbb Z}^{\nu} \setminus (\{ 0 \} \cup {\cal S})$. 
Roughly speaking, if one can show that 
$T_j^{-1}(\lambda)$ is bounded by $\Phi(\# \Lambda_j)$ in a suitable norm
for some function $\Phi$ of sub-exponential growth at infinity 
(e.g. $\Phi$ polynomial) and
$\# \Lambda_j$ grows exponentially with $j$, then the sequence 
$v_j(\lambda)$ will converge to some function $v(\lambda)$ which solves
${\cal P}(v(\lambda), \lambda) = 0$. Note that the modification of the
Newton scheme consists in restricting the linearized operator $T_j$
to some sub-lattice. The truncation of the high Fourier modes corresponds
to the regularization operations used in the Nash -- Moser technique.
The error introduced by this modification of the Newton scheme 
can be controlled employing the decay of the sequences $v_j(m)$ (in $m$). 
The decay of the $v_j$'s stems from the 
regularity of the equation. To this end we have required 
in assumption A1 the force law to be analytic. In a recent paper
Craig and Su \cite{CS} have shown in the case of 
time-periodic solutions of (1+1) -- dimensional
nonlinear wave equations that analyticity can be replaced by some
finite regularity assumption. This suggests that the method of 
Craig and Wayne can generally be extended to the non-analytic case.
However, we will not pursue the question of minimal regularity in this paper.

The main difficulty of the iteration scheme (\ref{ova.200}) 
is to prove the existence of the inverse matrix 
$T_j^{-1}(\lambda)$ and to derive the corresponding bounds.
It is clear that such bounds cannot hold uniformly for all values 
of the parameter 
$\lambda \in {\cal N}^{(1)}$. Even for the diagonal operator
$D(\omega)$ estimates on the inverse matrices 
$|(D|_{\Lambda_j})^{-1}(\omega)| \leq \Phi(\# \Lambda_j)$ require
some diophantine conditions on $\omega$. Such conditions
cannot hold in any open neighborhood of $\omega^{(0)}$. 
Therefore we may only expect
to bound $T_j^{-1}(\lambda)$ for $\lambda$ belonging to some subset
${\cal N}^{(j)} \subset {\cal N}^{(1)}$. Definition (\ref{ova.200})
only applies for $\lambda \in {\cal N}^{(j)}$. For technical reasons 
we will extend $v_{j+1}$ to a smooth function on all of ${\cal N}^{(1)}$
but for values $\lambda \in {\cal N}^{(1)} \setminus {\cal N}^{(j)}$ we cannot
control $|{\cal P}(v_{j+1}(\lambda), \lambda)|$. The sets of good
parameters ${\cal N}^{(j)}$ form a nested sequence and 
${\cal N}^{\infty}$ denotes the intersection of all these sets. 

We determine the sets ${\cal N}^{(j)}$ of good parameters and 
the estimates on $T_j^{-1}(\lambda)$ by a multi-scale analysis similar
to the one used by Bourgain \cite{B1}. In order to explain the multi-scale
analysis we consider a class of matrices which have a similar
structure as the matrices $T_j(\lambda)$ but which are easier to 
analyze. We show in remark \ref{Rova.2} below how 
to modify the argument to treat the case at hand. Set
\begin{eqnarray}
\label{ova.210}
T(\omega)(m, n) = (<\omega, m>) \delta_{m, n} + r(m - n), \;\;\;
m, n \in {\Bbb Z}^{\nu}.
\end{eqnarray}
Here $r$ denotes a sequence on ${\Bbb Z}^{\nu}$ which has small norm
and $r(m) \to 0$ for $|m| \to \infty$ at a sufficiently fast rate.
Observe that the matrices $T(\omega)$ are of the form $D + R$ 
where $D$ is diagonal and $R$ is a Toeplitz matrix with entries
decaying rapidly with the distance from the diagonal. 
These are exactly the properties
which are shared by the matrices $T_j(\lambda)$ (see (\ref{ova.200}))
and which are basic requirements for a multi-scale analysis. 

The main idea behind the multi-scale analysis, which was first introduced
by Fr\"ohlich and Spencer \cite{FS}, can be described in the 
following way. Denote
by $N_1 < N_2 < N_3 < \ldots$ a suitably chosen sequence of length scales
and fix a suitable function $\Phi(x)$ of sub-exponential growth at 
infinity. For fixed parameter
$\omega$ we say that
a lattice site $n \in {\Bbb Z}^{\nu}$ is $\omega$-moderate at scale $j$, 
if
\begin{eqnarray}
\nonumber
\| T(\omega)|_{B_{N_j}(n)}^{-1} \| \leq M_j \equiv \Phi(N_j)
\end{eqnarray}
where $B_{N_j}(n) = \{m \in {\Bbb Z}^{\nu}: |m-n| < N_j\}$ and
$\| \cdot \|$ denotes a suitable weighted operator norm 
(which will depend on $j$).
To see the significance of
this definition assume that $\Lambda$ is some sub-lattice of ${\Bbb Z}^{\nu}$
and suppose that there exists an integer $j_0$ 
such that each $n \in \Lambda$ is 
$\omega$-moderate at some scale $1 \leq j \leq j_0$. Then one can construct
an inverse of the restricted matrix $T(\omega)|_{\Lambda}$ and the bound on
the inverse is roughly given by
$M_{j_0}$. In other words, we can paste together the local inverse
matrices $T(\omega)|_{B_{N_j}(n)}^{-1}$ to a global inverse 
$T(\omega)|_{\Lambda}^{-1}$. In section \ref{cl} we state a recent
version of the coupling lemma (lemma \ref{Lcl.1}) which was introduced
by the author in \cite{Kri}. The hypothesis of the coupling
lemma essentially requires that $r(|m|)$ decays faster than 
$(\Phi(|m|))^{-1}$. This implies that for a site 
$n$ which is $\omega$-moderate 
at scale $j$ the product of $M_j$ 
(the bound on $\| T(\omega)|_{B_{N_j}(n)}^{-1} \|$) 
with the interaction term $|R(m, n)| = |r(m-n)|$ 
is small for any lattice site $m$
in $\Lambda \setminus B_{N_j}(n)$.
Observe that lemma \ref{Lcl.1} also requires that the neighborhoods
$B_{N_j}(n)$ are contained in the set $\Lambda$. We address this 
issue in remark \ref{Rova.3} below.

Of course, the question remains how to establish that a
lattice point is $\omega$-moderate at some scale. 
Here we again use the special 
structure of the matrix $T(\omega)$ which implies
\begin{eqnarray}
\label{ova.220}
T(\omega)(m+k, n+k) = \delta_{m, n} (<\omega, k>) + T(\omega)(m, n).
\end{eqnarray}
This motivates the introduction of a new parameter. 
For $\theta \in {\Bbb R}$ define
\begin{eqnarray}
\nonumber 
T(\theta, \omega)(m, n) := (\theta + <\omega, m>) \delta_{m, n} + 
r(m - n), \;\;\; m, n \in {\Bbb Z}^{\nu}.
\end{eqnarray}
Then 
\begin{eqnarray}
\label{ova.230}
T(0, \omega) &=& T(\omega),\\
\label{ova.240}
T(\theta, \omega)(m+k, n+k) &=&  T(\theta + <\omega, k>, \omega)(m, n).
\end{eqnarray} 
Relation (\ref{ova.240}) implies that at each scale $j$ it suffices to 
investigate the invertibility of 
$T(\theta, \omega)|_{B_{N_j}(0)}$ for all 
$\theta \in {\Bbb R}$, since 
\begin{eqnarray}
\label{ova.250}
T(\omega)|_{B_{N_j}(n)} = T(<\omega, n>, \omega)|_{B_{N_j}(0)}.
\end{eqnarray}
Denote $\Gamma_j := B_{N_j}(0)$.
Our goal to determine whether a
site $n$ is $\omega$-moderate can be expressed in the following way.
Construct a nested sequence of sets
\begin{eqnarray}
\nonumber
{\Bbb R} \supset I^{(1)}_{\omega} \supset I^{(2)}_{\omega} \supset
\ldots \; ,
\end{eqnarray}
such that
\begin{eqnarray}
\label{ova.300}
\| T(\theta, \omega)|_{\Gamma_1}^{-1} \| \leq M_1 &&
\mbox{ if } \theta \in {\Bbb R} \setminus I^{(1)}_{\omega}, \\
\label{ova.310}
\| T(\theta, \omega)|_{\Gamma_j}^{-1} \| \leq M_j &&
\mbox{ if } \theta \in  I^{(j-1)}_{\omega} \setminus I^{(j)}_{\omega},
\;\; j \geq 2.
\end{eqnarray}
By (\ref{ova.250}) a site $n$ is $\omega$-moderate for some scale
$1 \leq j \leq j_0$ if $<n, \omega> \in {\Bbb R} \setminus I^{(j_0)}_{\omega}$.
Now the objective is to construct the sets $I^{(j)}_{\omega}$ and to control
the location and measure of these sets.
The sets $I^{(j)}_{\omega}$ will be constructed inductively (in $j$).
The following separation property is essential.
\begin{eqnarray}
\label{ova.320}
\mbox{ If } \theta + <\omega, m> \in I^{(j)}_{\omega}
\mbox{ and } \theta + <\omega, n> \in I^{(j)}_{\omega}, \mbox{ then either }
m = n \mbox{ or } |m - n| > 2 N_{j+1}.
\end{eqnarray}
Next we present an outline of the construction of the sets
$I^{(j)}_{\omega}$. Let us assume that the smallest length scale $N_1 = 1$.
Then $\Gamma_1 = \{ 0 \}$ and the question of invertibility of the 
$1 \times 1$ matrix $T(\theta, \omega)|_{\Gamma_1}$ reduces to the 
question whether $\theta + r(0)$ can be inverted. Recall that we have
assumed that $r$ is small. We denote $\epsilon := \|r\| \geq |r(0)|$ 
for some suitable norm $\| \cdot \|$. For $|\theta| \geq 2 \epsilon$, say, 
we have $|T(\theta, \omega)|_{\Gamma_0}^{-1}| \leq 1/ \epsilon =: M_0$.
We can therefore satisfy (\ref{ova.300}) by setting
\begin{eqnarray}
\nonumber
I^{(1)}_{\omega} := (-2 \epsilon, 2 \epsilon)
\end{eqnarray}
Observe that the separation property (\ref{ova.320}) for $I^{(1)}_{\omega}$ is
equivalent to showing that for $k \in {\Bbb Z}^{\nu}$ satisfying
$0 < |k| \leq 2 N_2$ the following diophantine condition holds
\begin{eqnarray}
\label{ova.420}
|<\omega, k>| > 4 \epsilon.
\end{eqnarray}
Such a condition is satisfied for most values of $\omega$, if $\epsilon$ is 
sufficiently small (depending on $N_2$).

Suppose now that the sets $I^{(l)}_{\omega}$ have been defined for $1 \leq
l \leq j$ satisfying (\ref{ova.300}) -- (\ref{ova.320}). In order to construct
$I^{(j+1)}_{\omega}$ we need to determine those values $\theta \in 
I^{(j)}_{\omega}$ for which we can establish 
$\| T(\theta, \omega)|_{\Gamma_{j+1}}^{-1} \| \leq M_{j+1}$.  
Fix $\theta \in I^{(j)}_{\omega}$. Since $\theta = \theta + <\omega, 0>
\in I^{(j)}_{\omega}$ we conclude from the separation property
(\ref{ova.320}) that $\theta + <\omega, n> \in {\Bbb R} \setminus
I^{(j)}_{\omega}$ for all lattice sites $n$ satisfying
$0 < |n| \leq 2 N_{j+1}$. This implies in particular, that for each
$n \in \Lambda := \Gamma_{j+1} \setminus \{ 0 \}$ there exists
some scale $1 \leq l \leq j$ such that 
$\| T(\theta, \omega)|_{B_l(n)}^{-1} \| \leq M_l$. 
As we have described above, the coupling
lemma (see lemma \ref{Lcl.1}, cf. remark \ref{Rova.3}) 
then shows that $T(\omega)|_{\Lambda}$ is invertible and 
$\| T(\omega)|_{\Lambda}^{-1} \| \lesssim M_j$.
Consider the block decomposition
\begin{eqnarray}
\label{ova.350}
T|_{\Gamma_{j+1}} =
\left(
\begin{array}{cc}
T|_{\Lambda} & P_1 \\ P_2 & \theta + r(0)
\end{array}
\right).
\end{eqnarray}
Define
\begin{eqnarray}
\label{ova.400}
b(\theta, \omega) := \theta + r(0) 
- P_2 (T|_{\Lambda})^{-1}(\theta, \omega) P_1.
\end{eqnarray}
Suppose that $b(\theta, \omega) \neq 0$. Then 
we can invert $T|_{\Gamma_{j+1}}$,
\begin{eqnarray}
\nonumber
T|_{\Gamma_{j+1}}^{-1} =
\left(
\begin{array}{cc}
T|_{\Lambda}^{-1} +  (T|_{\Lambda}^{-1}) P_1 b^{-1} P_2 (T|_{\Lambda}^{-1})& 
- (T|_{\Lambda}^{-1}) P_1 b^{-1} \\ - b^{-1} P_2 (T|_{\Lambda}^{-1}) & b^{-1}
\end{array}
\right).
\end{eqnarray}
Thus we define 
\begin{eqnarray}
\label{ova.410}
I^{(j+1)}_{\omega} := \{ \theta \in I^{(j)}_{\omega}: 
|b(\theta, \omega)| < \delta_{j+1} \}, \;\;\; \mbox{ with }
\delta_{j+1} \sim M_{j+1}^{-1}.
\end{eqnarray} 
Since the off-diagonal parts $P_1$, $P_2$ are assumed to be small with
entries decaying rapidly with the distance from the diagonal we 
expect that $b(\theta, \omega) \sim \theta$.
This would imply that the set 
$I^{(j+1)}_{\omega}$ has length of order $\delta_{j+1} \sim 
\Phi(N_{j+1})^{-1}$ which tends to zero as $j \to \infty$.
We still need to investigate the separation condition (\ref{ova.320}) for
the set $I^{(j+1)}_{\omega}$. Assuming again that $b(\theta, \omega) \sim
\theta$ condition (\ref{ova.320}) essentially leads to a 
diophantine condition
\begin{eqnarray}
\label{ova.412}
|<\omega, k>| \gtrsim \delta_{j+1}
\end{eqnarray}
for all $k \in {\Bbb Z}^{\nu}$ with $0 < |k| \leq 2N_{j+2}$. 
This condition
can be satisfied for most values of $\omega$ if $\delta_{j+1}$ is 
sufficiently small (depending on $N_{j+2}$). 

Suppose one can choose the parameters $N_j$, $\Phi$, $\delta_j$ of the
construction such that the above described procedure works. Then 
for most parameter values $\omega$
(which satisfy some diophantine type conditions to ensure that
separation properties (\ref{ova.320}) are satisfied) we have obtained
a nested sequence of open sets
$I^{(j)}_{\omega}$ (of length $\to 0$ as $j \to \infty$) such that for 
each $n \in {\Bbb Z}^{\nu}$ with $<\omega, n> \in {\Bbb R} \setminus
I_{\omega}^{(j)}$ the site
$n$ is $\omega$-moderate at some scale $1 \leq l \leq j$.
In view of the coupling lemma we can therefore prove
bounds $\| T(\omega)|_{\Lambda}^{-1} \| \lesssim M_j$ for subsets $\Lambda
\subset {\Bbb Z}^{\nu}$,
if we can show that $<\omega, n> \in {\Bbb R} \setminus
I_{\omega}^{(j)}$ for all $n \in \Lambda$. Since the inductive construction
of the sets $I_{\omega}^{(j)}$ was rather explicit (involving the function $b$)
the multi-scale analysis as described above yields a powerful tool to
investigate for a given sub-lattice $\Lambda \subset {\Bbb Z}^{\nu}$ the set
of parameters $\omega$ for which the inverse of $T(\omega)|_{\Lambda}$ 
is bounded by $\Phi(\# \Lambda)$. 
Recall that this is the kind of control
required for the modified Newton scheme (\ref{ova.200}) to solve the 
$P$ -- equation (\ref{ova.115}).

We briefly turn to the solution of the bifurcation equation 
($Q$ -- equation).
This equation is degenerate. However, we transform the equation 
in such a way that it can be solved by a standard
implicit function theorem,
determining the frequency vector $\omega$ as a function of the parameter $a$. 
What is more delicate to prove is the fact that
the set $\{ (a, \omega(a) ): a \mbox{ small } \}$ is contained
in ${\cal N}^{\infty}$ for most values of $a$. To do this
we must find a good description of the set of 
resonant parameters $\lambda \in {\cal N}^{(1)} \setminus {\cal N}^{\infty}$, 
for which either the separation property
(cf. (\ref{ova.320})) is violated or the inverse of $T_j(\lambda)$ cannot be
bounded by $\Phi (\# \Lambda_j)$. Following Bourgain \cite{B1} we
construct sets of polynomials 
with coefficients depending
on $\lambda$ which can be used to characterize the set of 
resonant parameters
(see the second part of remark \ref{Rova.2} for a motivation of these 
polynomials).
To obtain an lower estimate on the set of parameters $a$ 
for which $(a, \omega(a)) \in {\cal N}^{\infty}$ we employ assumption
$A4$ which ensures that to second order $\omega$ depends in a non-degenerate
way on $a$.

\begin{Remark}
\label{Rova.1} Symmetries of the equation.
{\em
For given $F$, $b$, $\gamma$ we say that the pair
$(\chi, \omega)$ defines a travelling wave solution 
($\chi: {\Bbb T}^{\nu} \to {\Bbb R}$, $\omega \in {\Bbb R}^{\nu}$) 
if 
\begin{eqnarray}
\label{ova.500}
\gamma^2 <g, (D^2 \chi)(\xi) g> = F(-b +\chi(\xi - \omega) - \chi(\xi))
- F(-b + \chi(\xi) - \chi(\xi + \omega)) \; \mbox{ for all }
\xi \in {\Bbb T}^{\nu}.
\end{eqnarray}
Observe that (\ref{ova.500}) implies that 
the system $(x_n)$ defined through (\ref{smr.20}) solves (\ref{smr.10}).
Furthermore, it will become clear in section \ref{fa} that
equation (\ref{ova.80}) is equivalent to (\ref{ova.500}).

Suppose now that $(\chi, \omega)$ defines a travelling wave solution. Then
\begin{itemize}
\item[(S1)] {\em Additive constants} \newline
$(\chi_x, \omega)$ defines a travelling wave solution for $x \in {\Bbb R}$
where $\chi_x(\xi) = \chi(\xi) + x$.

\item[(S2)] {\em Phase shifts} \newline
$(\chi_{\zeta}, \omega)$ defines a travelling wave solution for 
$\zeta \in {\Bbb T}^{\nu}$ where $\chi_{\zeta} (\xi) = \chi(\xi + \zeta)$.

\item[(S3)] {\em Odd reflection} \newline
$(\tilde{\chi}, \omega)$ defines a travelling wave solution 
where $\tilde{\chi}(\xi) = - \chi(-\xi)$.
\end{itemize}
These symmetries allow us to restrict our attention to the reduced space 
presented in (\ref{ova.87}). The restriction $u(0) = 0$ is stated in 
proposition
\ref{Pfeq} and is related to (S1). 
Due to symmetry (S3) we may search for
solutions ($\chi$, $\omega$) with $\chi(\xi) = - \chi(-\xi)$. For the 
corresponding sequence $u$ this implies that $u(m) \in {\Bbb R}$ for all
$m \in {\Bbb Z}^{\nu}$. Since we are only interested in real-valued functions
$\chi$ (i.e. $u(m) = \overline{u(-m)}$) we obtain $u(m) = u(-m) \in {\Bbb R}$.

Suppose that we have constructed a solution ($\chi$, $\omega$) where $\chi$ is
an odd function. Denote $a_j = u(e_j) \in {\Bbb R}$. For $\zeta \in 
{\Bbb T}^{\nu}$ define $\chi_{\zeta}$ as in (S2) and denote by $u_{\zeta}$
the corresponding sequence. Then $u_{\zeta}$ again solves (\ref{ova.80})
by (S2) and $u_{\zeta}(e_j) = a_j e^{i \zeta_j}$. We therefore obtain solutions
with respect to complex parameters $z_j = u(e_j)$ (as claimed in theorem
\ref{Tsmr3.1}) through a phase shift from
the corresponding odd solutions.

\noindent
{\em Warning}: In our analysis we will also consider 
complex parameters $\lambda=(a, \omega)$ because this simplifies some of the 
estimates, e.g. by using the Cauchy integral formula. However, for 
values of $a$ (and $\omega$) with non-vanishing imaginary parts these
solutions do {\em not} correspond to physical solutions.   
In particular they do not correspond to some solution $u_{\zeta}$ defined
above. 
}
\end{Remark}

\begin{Remark}
\label{Rova.2}
{\em
The matrices $T_j(\lambda)$ which appear in the modified Newton scheme 
(\ref{ova.200}) are of the form 
\begin{eqnarray}
\nonumber
T_j(\lambda)(m, n) = V(\omega)(m) \delta_{m, n} + r_j(\lambda)(m - n)
\end{eqnarray}
where $V(\omega)$ was defined in (\ref{smr2.70}). Regarding the 
multi-scale analysis described above for operators (\ref{ova.210}) 
a number of modifications are needed. In this remark we describe 
three of the more significant modifications. 

Firstly, the diagonal entries
have a more complicated behavior with respect to translations
(cf. (\ref{ova.220})). In fact, we may only consider translations with
respect to vectors $k$ which satisfy $<k, g> = 0$. Define 
\begin{eqnarray}
\label{ova.520}
V(\theta, \omega)(m) &:=&
\left\{
\begin{array}{ll}
F'(-b) & \mbox{ if } <g, m> = 0, \\
F'(-b) - \frac{<g, m>^2 \gamma^2}
{4 \sin^2 \frac{\theta + <\omega, m>}{2}} & \mbox{ else },
\end{array}
\right. \\
\nonumber
T_j(\theta, \lambda)(m, n) &:=& V(\theta, \omega)(m) 
\delta_{m, n} + r_j(\lambda)(m - n).
\end{eqnarray}
Then 
\begin{eqnarray}
\nonumber
T_j(0, \lambda) &=& T_j(\lambda), \\
\nonumber
T_j(\theta, \lambda)(m + k, n + k) &=& T_j(\theta + <\omega, k>, \lambda)
(m, n) \;\;\; \mbox{ if } <k , g> = 0, 
\end{eqnarray}
replacing relations (\ref{ova.230}), (\ref{ova.240}) of the simplified model. 
In principle, the restriction
$<k, g> = 0$ on the directions $k$ of translation could create great difficulties.
However, in our case all singular sites are located within the strip
$\{ m \in {\Bbb Z}^{\nu}: |<m, g>| \leq  \nu \}$ and the restricted directions
are sufficient to translate between different singular sites.

Secondly, singular sites appear in clusters of up to $2 \nu$ lattice points.
Such clusters always appear in the construction of quasi-periodic solutions
and we may essentially think of them as translates of the set ${\cal S}$ 
defined
in (\ref{ova.90}). This effects the inductive construction of the sets
$I_{\omega}^{(j+1)}$ in the following way. The decomposition of the
set $\Gamma_{j+1} = \Lambda \cup \{0\}$ in (\ref{ova.350}) has to be
replaced by  $\Gamma_{j+1} = \Lambda \cup S$, where $S$ is some subset
of ${\cal S}$.
Bourgain \cite{B1} introduced a powerful and general technique to deal
with such a situation: Instead of 
the scalar function $b$ (see (\ref{ova.400})) he defines a
matrix-valued function 
\begin{eqnarray}
\nonumber
b(\theta, \lambda) = (T_j|_{S})(\theta, \lambda) - P_2  
(T_j|_{\Lambda})^{-1}(\theta, \lambda) P_1,
\end{eqnarray}
and the definition of the set $I_{\omega}^{(j+1)}$ (now 
$I_{\lambda}^{(j+1)}$, since $T_j$ depends on $\lambda$) is replaced by
\begin{eqnarray}
\nonumber
I_{\lambda}^{(j+1)} := \{ \theta \in I_{\lambda}^{(j)}: 
|\det b(\theta, \lambda)| < \delta_{j+1} \}.
\end{eqnarray} 
The function $\det b(\theta, \lambda)$ is more complicated to analyze 
than the function $b(\theta, \omega)$ defined for our simplified version
in (\ref{ova.400}). Recall that in
the discussion following (\ref{ova.410}) we have argued that 
$b(\theta, \omega) \sim \theta$ in order to obtain information on 
the sets $I_{\omega}^{(j)}$ and to prove the separation property (\ref{ova.320}).
Bourgain proposes to apply the Weierstrass preparation theorem
to locally replace the function $\det b (\theta, \lambda)$ by polynomials. 
In our situation we will construct a finite set of pairs $(p, \vartheta)$ 
where each $p(z, \lambda)$ is a polynomial in $z$ 
of degree $\leq 2 \nu$ with leading
coefficient 1 (the other coefficients are small and depend on $\lambda$), 
such that for $\theta$ close to $\vartheta$ we can estimate
\begin{eqnarray}
\nonumber
|\det b(\theta, \lambda)| \geq c |p(\theta - \vartheta, \lambda)|
\end{eqnarray}
for some positive constant $c$. 
This implies, for example, that a lattice site $n$ is $\lambda$-moderate,
at some scale $\leq j+1$, if
\begin{eqnarray}
\nonumber
|p(<\omega, n> - \vartheta, \lambda)| > \frac{\delta_{j+1}}{c}
\end{eqnarray}
for some suitable pair $(p, \vartheta)$. Using estimates on the 
coefficients of the polynomials one can analyze the set of parameters
$\lambda$ for which such a condition holds.

The power of Bourgain's idea is best demonstrated by investigating the 
separation property (\ref{ova.320}).
In fact, suppose that there exist lattice sites $m \neq n$ and $\theta \in 
{\Bbb R}$, 
such that $\theta + <\omega, n> \in I_{\lambda}^{(j+1)}$ and 
$\theta + <\omega, m> \in I_{\lambda}^{(j+1)}$. By definition there
exist pairs
$(p_1, \vartheta_1)$, $(p_2, \vartheta_2)$ such that
\begin{eqnarray}
\nonumber
|p_1(\theta + <\omega, m> - \vartheta_1, \lambda)| &\leq& 
\frac{\delta_{j+1}}{c},
\\
\nonumber
|p_2(\theta + <\omega, n> - \vartheta_2, \lambda)| &\leq& 
\frac{\delta_{j+1}}{c}
\end{eqnarray}
A resultant type construction (see e.g. section \ref{res}) then implies
that there exists a polynomial $p \equiv p_1 \ominus p_2$ of degree
$\leq 4 \nu^2$ such that
\begin{eqnarray}
\nonumber
|p(<\omega, m - n> - (\vartheta_1 - \vartheta_2), \lambda)| =
{\cal O}(\delta_{j+1}).
\end{eqnarray}
Observe that the dependence on $\theta$ has dropped out. We can replace
the diophantine condition (\ref{ova.412}) which ensured the separation property
in the simplified model above by a condition
\begin{eqnarray}
\nonumber
|p(<\omega, m - n> - (\vartheta_1 - \vartheta_2), \lambda)| >
C \delta_{j+1}
\end{eqnarray}
for some positive constant $C$. With the help of suitable estimates
on the coefficients of $p$ one may show that such a 
condition is satisfied for most values of the parameter $\lambda$.

We turn to the third modification. Observe that the matrices $T_j(\theta,
\lambda)$ depend on the scale $j$. One might expect that one has to 
perform a multi-scale analysis for each operator $T_j$ separately,
leading to sequences $I_{\lambda}^{(l, j)}$ for $l \leq j$.
However, we will be able to show that
\begin{eqnarray}
\label{ova.543}
|T_j(\theta, \lambda) - T_l(\theta, \lambda)| << M_l^{-1} \;\;\;
\mbox{ for } j \geq l.
\end{eqnarray}
Suppose we have constructed sets $I^{(l)}_{\lambda}$ such that
\begin{eqnarray}
\label{ova.547}
\| T_1(\theta, \lambda)|_{\Gamma_1}^{-1} \| \leq M_1 &&
\mbox{ if } \theta \in {\Bbb R} \setminus I^{(1)}_{\lambda}, \\
\label{ova.549}
\| T_l(\theta, \lambda)|_{\Gamma_l}^{-1} \| \leq M_l &&
\mbox{ if } \theta \in  I^{(l-1)}_{\lambda} \setminus I^{(l)}_{\lambda},
\;\; l \geq 2
\end{eqnarray}
(cf. (\ref{ova.300}), (\ref{ova.310})). Then it follows from (\ref{ova.543})
by a standard application of the Neumann series that
\begin{eqnarray}
\nonumber
\| T_j(\theta, \lambda)|_{\Gamma_1}^{-1} \| \leq 2 M_1 
&\mbox{ for } \theta \in {\Bbb R} \setminus I^{(1)}_{\lambda}, 
& j \geq 1\\
\nonumber
\| T_j(\theta, \lambda)|_{\Gamma_l}^{-1} \| \leq 2 M_l 
&\mbox{ for } \theta \in  I^{(l-1)}_{\lambda} \setminus I^{(l)}_{\lambda},
&j \geq l \geq 2.
\end{eqnarray}
This is the reason why it suffices to construct one sequence of nested
sets $I_{\lambda}^{(l)}$ which covers all matrices $T_j(\theta, \lambda)$
simultaneously.
}
\end{Remark}

\begin{Remark}
\label{Rova.3}
{\em
In our discussion of the multi-scale analysis above we have used
the coupling lemma in order to construct inductively the sets 
$I_{\omega}^{(j)}$. The coupling lemma \ref{Lcl.1} allows to paste together
local inverse matrices $T|_{B_{N_l}(n)}^{-1}$, $n \in \Lambda$ 
to an inverse matrix
of $T|_{\Lambda}$. However, one requirement of this procedure is that
the local neighborhoods $B_{N_l}(n)$ are contained in the set $\Lambda$.
Clearly, this condition is violated for lattice sites $n$ which lie close
to the boundary of the set $\Lambda$. We address this problem in the
following way. At scale $j$ we not only consider the invertibility
of $T_j|_{\Gamma_j}(\theta, \lambda)$ (recall $\Gamma_j = B_{N_j}(0)$)
but also the invertibility of $T_j|_{C}(\theta, \lambda)$ for $C \in
{\cal C}^{(j)}$ where ${\cal C}^{(j)}$
denotes a certain class of subsets of $\Gamma_j$.
Correspondingly, we need to construct sets $I_{C, \lambda}^{(j)}$ for
all sets $C \in {\cal C}^{(j)}$.    
}
\end{Remark}

\begin{Remark}
\label{Rova.4}
{\em
Observe that the function $V(\theta, \omega)(m)$ defined in (\ref{ova.520})
is $2 \pi$ -- periodic in the variable $\theta$. Hence we can restrict
ourselves to define the sets $I^{(j)}_{C, \lambda}$ as subsets of
$(-\pi, \pi]$. The periodic extension of this set to the real line
will be denoted by
\begin{eqnarray}
\label{ova.550}
\tilde{I}^{(j)}_{C, \lambda} := I^{(j)}_{C, \lambda} + 2 \pi {\Bbb Z}.
\end{eqnarray}
In this context the following map will be useful.
\begin{eqnarray}
\label{ova.560}
[ \cdot ]&:& {\Bbb C} \to 
\{ w \in {\Bbb C}: - \pi < \mbox{ Re}(w) \leq \pi \},\\
\nonumber
&& \mbox{ where } [z] \mbox{ is uniquely defined by the condition }
z - [z] \in 2 \pi {\Bbb Z}.
\end{eqnarray}
}
\end{Remark}

\begin{Remark}
\label{Rova.5}
{\em
In this remark we discuss the effects of the singularity in the definitions
of $T_j(\lambda)$ and $T_j(\theta, \lambda)$. For any
subset $C$ of ${\Bbb Z}^{\nu}$ we define the sets
\begin{eqnarray}
\label{ova.600}
Z_C &:=& \{\lambda \in {\Bbb C}^{2 \nu}: < \omega, m> \in 2 \pi {\Bbb Z}\;
\mbox{ for some } m \in C \mbox{ with } <m, g> \neq 0 \}, \\
\label{ova.610}
\tilde{Z}_C &:=& \{(\theta, \lambda) \in {\Bbb C}^{2 \nu+1}:
\theta + < \omega, m> \in 2 \pi {\Bbb Z}\;
\mbox{ for some } m \in C \mbox{ with } <m, g> \neq 0 \}.
\end{eqnarray}
Observe that the restricted matrices $T_j|_{C}(\lambda)$ are well defined
if $\lambda \notin Z_C$ and the matrices $T_j|_{C}(\theta, \lambda)$ 
are well defined if $(\theta, \lambda) \notin \tilde{Z}_C$. These
sets of singularities, however, create no difficulties in our analysis. Recall
that the Newton scheme only uses the inverse of the matrices $T_j|_C$.
We will show that the inverse matrices $(T_j|_C)^{-1}(\lambda)$, 
$(T_j|_C)^{-1}(\theta, \lambda)$ can be continued analytically across
$Z_C$, respectively across $\tilde{Z}_C$.
}
\end{Remark}

We conclude this section by describing the plan of the paper.

The proof of theorem \ref{Tsmr3.1} contains many parameters, like the 
size of the length scales $N_j$, the bounds $M_j$ on the inverse matrices,
the choice of $\delta_j$ in the definition of the sets $I^{(j)}_{C, \lambda}$
and many more. In section \ref{npc} we show that we can choose
all these parameters of the construction in such a way that the above 
described multi-scale analysis works (see lemma \ref{hyp}). 
In section \ref{npc} we also collect
and explain all the notation used in chapters II and III. 
For the convenience of the reader 
we also include a table of notation at the end of the appendix 
(section \ref{tnot}).

In section \ref{fa} we derive the equations for the frequency vector 
$\omega$ and 
for the modified Fourier coefficients $u(m)$.
Our analysis is greatly facilitated by using appropriate norms for
the sequences of Fourier coefficients (cf. \cite{Kri}). We choose weighted 
$\ell_1$-norms. The weights are either of exponential growth
or lie in the Gevrey class. Similar weights have already been used
in the context of small divisor problems by DeLatte \cite{DeL},
Craig and Wayne \cite{CW2} and Bourgain \cite{B1}. 
In the first part of section \ref{fa} we introduce
the weight functions and norms used in this paper and we prove 
their basic properties.

In sections \ref{fis} -- \ref{sp} we deal with the small divisor problem
and solve the $P$ -- equation (see (\ref{ova.115}) and theorem \ref{TP}).
As explained above we proceed inductively. At every step $j$ of the 
induction we construct the approximate solutions $v_j$, the set of 
good (i.e. non-resonant parameters ${\cal N}^{(j)})$ and the sets
$I^{(j)}_{C, \lambda}$ of the multi-scale analysis. The induction
statements $({\cal IS})_j$, which take a slightly different form for $j=1$ and 
$j \geq 2$,
are stated in sections \ref{fis} and \ref{is}. Essentially they contain the 
following claims.  

\newpage
\noindent
{\bf Induction Statement $({\cal IS})_j$}
\begin{itemize}
\item[(1)]
\begin{itemize}
\item[] Regularity properties of $v_j(\lambda)$.
\item[] Estimates on $v_j(\lambda)$, $\partial^{\beta}v_j(\lambda)$ and
on ${\cal P}(v_j(\lambda), \lambda)$.
\item[] Symmetries and support of $v_j(\lambda)$.
\end{itemize}
\item[(2)] 
\begin{itemize}
\item[] Description of the set of ``good'' parameters ${\cal N}^{(j)}$ 
at scale $j$
as a disjoint union of cubes.
\item[] Description of the set of parameters 
${\cal N}^{(j-1)} \setminus {\cal N}^{(j)}$ 
which is removed at the $j$-th step 
in terms of a finite set of 
polynomials (cf. remark \ref{Rova.2}
above). Estimates on the number of polynomials and on the coefficients of
the polynomials.  
\end{itemize}
\item[(3)]
\begin{itemize}
\item[(a)] Location of the set $I^{(j)}_{C, \lambda}$ (cf remarks \ref{Rova.3} and
\ref{Rova.7}).
\item[(b)] Diophantine estimates for lattice sites $m$ with 
$2 N_{j-1} < |m| \leq 2 N_j$ and frequencies $\omega$ 
satisfying $(a, \omega) \in {\cal N}^{(j)}$.
\item[(c)]  For all $\lambda \in {\cal N}^{(j)}$
the lattice sites $m$ satisfying $2 N_{j-1} < |m| \leq 2 N_j$
are $\lambda$-moderate at some scale $1 \leq l \leq j$ 
with respect to $T_j(\lambda)$. 
\item[(d)] For all $\lambda \in {\cal N}^{(j)}$
the separation property (cf. (\ref{ova.320})) 
holds for the sets $I^{(j)}_{C, \lambda}$.
\end{itemize}
\item[(4)]
\begin{itemize}
\item[]
The analog of statements (\ref{ova.547}), respectively (\ref{ova.549}) hold
for the matrices 
\newline
$T^{(j)}(\theta, \lambda) = \mbox{ diag} ( V(\theta, \omega) ) +
(DW)(\varphi(a) + v_j(\lambda))$.
\end{itemize}
\end{itemize}

\begin{Remark}
\label{Rova.10}
{\em
When we actually perform the induction in chapter II, there are a 
few differences from the just described induction statements which 
might be somewhat confusing (e.g. at step $j$ we will consider length
scales $N_{j-1}$ rather than length scales $N_j$). It should be
emphasized that the discussion of section \ref{ova} is formal
and although it is very close to the actual analysis of chapters II and
III, the accurate definitions and statements are found in those chapters.
}
\end{Remark}

\begin{Remark}
\label{Rova.6}
{\em
The induction statements of sections \ref{fis} and \ref{is} are slightly
stronger than just described. 
They are required to hold also for parameters which lie
in small complex neighborhoods of the sets ${\cal N}^{(j)}$ and
$I^{(j)}_{C, \lambda}$.
}
\end{Remark}

\begin{Remark}
\label{Rova.7}
{\em
We find it convenient to use cube decompositions to control the
geometry of the sets ${\cal N}^{(j)}$. The cubes are numbered
by some $k \in K^{(j)}$. We choose the sets $I^{(j)}_{C, \lambda}$ 
in such a way that they agree for all values of $\lambda$
belonging to the same sub-cube. Hence we can index these
sets by $I^{(j)}_{C, k}$ with $k \in K^{(j)}$.
}
\end{Remark}

\noindent
The induction statement simplifies in the case $j=1$ since we choose
$v_1(\lambda) \equiv 0$ and ${\cal N}^{(1)}$ to be some cubic neighborhood
of $\lambda^{(0)}$. Statements (3) and (4) follow essentially from the 
diophantine assumption A3. Observe that statement (4) is formulated in terms
of a lower bound on the diagonal entries of $T^{(1)}$.
The details of the proof of $({\cal IS})_{j=1}$ are presented in section
\ref{fis}.

The proof of the induction step $j \to j+1$ 
stretches over sections \ref{cv2} -- 
\ref{comp}. First we construct the new iterate $v_{j+1}$. Recall that
we denote the scales corresponding to the multi-scale analysis by $N_j$.
These scales grow super-exponentially in $j$ and hence they increase 
too fast for
a Newton scheme. For this reason we need to perform many Newton iterations 
to obtain $v_{j+1}(\lambda)$ from $v_j(\lambda)$. For each
iteration the  bounds on the inverse matrices follow via the coupling lemma
from the induction hypothesis (3c) and from estimates 
$\|T^{(j)} - T^{(l)} \| << \|(T^{(l)})^{-1}\|$, $j \geq l$, which follow from
the fast convergence of the sequence $v_j$ (see induction statement (4)
above for a definition of $T^{(j)}$).

For the construction of $v_{j+1}$ we distinguish the cases $j = 1$ 
(section \ref{cv2}) and $j \geq 2$ (section \ref{cvj}). The reason
for this distinction is that we have chosen the parameters 
of the construction in such a way that
in the step $1 \to 2$ we 
do not encounter small divisors. This enables us to use norms
with exponentially growing weights (rather than the slower growing
Gevrey type weights which we use in the context of the multi-scale analysis) 
leading to better estimates on $v_2$ which are much needed in chapter III.
In section \ref{cpol} we construct the polynomials which are used to 
define the sets $I^{(j+1)}_{C, \lambda}$ and ${\cal N}^{(j+1)}$
(cf. the second part of remark \ref{Rova.2}). We complete the
proof of the induction step in section \ref{comp}. In section \ref{sp}
the solution of the $P$ -- equation and certain properties of the
solution are stated in theorem \ref{TP} which
follows immediately from the induction statements. This
ends chapter II.

In chapter III we solve the bifurcation equation ($Q$ -- equation) by
determining $\omega$ as a function of $a$ (section \ref{solq}). Section
\ref{est} is devoted to deriving the estimates on the derivatives of
the function $\omega = \omega(a)$ which are necessary to prove
lower bounds on the measure of 
the set $\{ a: (a, \omega(a)) \in {\cal N}^{\infty} \}$.
In section \ref{pmt} we are finally ready to prove our main result, theorem
\ref{Tsmr3.1}.

We have collected those lemmata in the appendix which can be 
formulated somewhat independently from the constructions of chapters II
and III. Some of these tools might be interesting by themselves.

In section \ref{ass} we study the genericity of assumptions
A1 -- A4 of our main result. We find a fairly explicit description
of the exceptional set of force laws $F$ for which we cannot verify assumptions
A1 -- A4 (see lemma \ref{Lass.3} and theorem \ref{Tass.1}).
Section \ref{wf} contains the essential part of the proof of
property (\ref{wf.120}) for the weight functions we consider in this paper.
Property (\ref{wf.120}) together with the definition of the corresponding
norms on sequence and matrix spaces facilitate a number of estimates
needed in the multi-scale analysis. In sections \ref{enp} and \ref{lop}
we prove basic facts about the nonlinear, respectively linear part of 
our lattice equation. In section \ref{cl} we formulate and prove
a version of the coupling lemma which has been introduced in \cite{Kri}.

We mentioned in remark \ref{Rova.2} above that we have to estimate
the coefficients of polynomials which are
obtained by an application of the Weierstrass preparation theorem.
Hence one needs to carefully investigate the construction of these polynomials
in the proof of the Weierstrass preparation theorem. Such an analysis
was provided by Bourgain \cite{B4} in the even more complicated situation
that the degree of the polynomials may be arbitrarily large. 
For completeness sake
we present Bourgain's nice proof in section \ref{wpt}.

The estimates on the sets 
${\cal N}^{(j-1)} \setminus {\cal N}^{(j)}$
of resonant parameters are all based on an
elementary lemma for real valued functions defined on an interval. The lemma
concerns functions with the property that there exists a $k \in {\Bbb N}$
such that the $k$-th derivative of this function is continuous and 
uniformly bounded away from $0$. For such functions $g$ one can derive
good upper bounds on the sets $\{x: |g(x)| \leq \delta \}$. This is 
made precise and proved in section \ref{el}. 

In section \ref{res}
we perform a resultant type construction. For given polynomials 
$p$, $q \in {\Bbb C}[z]$ it provides three polynomials
$p \ominus q \in {\Bbb C}[z]$ and $R_1$, $R_2 \in {\Bbb C}[x, y]$
such that
\begin{eqnarray}
\nonumber
(p \ominus q) (x-y) = p(x) R_1(x,y) + q(y) R_2(x, y).
\end{eqnarray}
As mentioned in remark \ref{Rova.2} above it was the idea of Bourgain
to use such an algebraic construction to prove the separation property
(\ref{ova.320}). In our analysis we not only use the existence of such
a representation but we also need estimates on the 
coefficients of $p \ominus q$, $R_1$ and $R_2$. Therefore we investigate
how the coefficients of these polynomials can be obtained from
the coefficients of $p$ and $q$. 

Throughout our analysis we need to estimate higher order derivatives
of composed functions or of functions which are defined implicitly.
In many cases we cannot use the Cauchy integral formula to obtain
such bounds. In such situations we employ somewhat explicit
formulae for the higher order chain rule which are adapted to
the desired estimates. We have collected these different versions
of the chain rule in section \ref{cr}.

\noindent
{\bf Acknowledgments.}
It is my great pleasure to thank Professor P. Deift and Professor E.
Wienholtz for their constant support, advice and encouragement. 
I thank J. Bourgain, L. Chierchia, 
W. Craig, P. Deift, R. de la Llave, K. T-R McLaughlin,
J. P\"oschel, J. Shatah, T. Spencer and C. E. Wayne for useful
conversations on small divisor problems. Part of the work was
completed while I was visiting Courant Institute, New York,
and the Max-Planck-Institut f\"ur Mathematik in den Naturwissenschaften,
Leipzig. I am grateful to both institutes for their  hospitality and
for providing a stimulating research environment. The author was 
supported in part by the Deutsche Forschungsgemeinschaft grants
\# Kr 1673/1-1 and \# Kr 1673/2-1.

\newpage
\noindent
\begin{center}
{\huge \bf Chapter II} \vspace{1.2cm}\\
{\huge \bf The Small Divisor Problem} \vspace{2cm}
\end{center}

The goal of chapter II is to solve the $P$ -- equation (see theorem
\ref{TP}). Throughout chapter II and chapter III
we assume that the hypothesis of theorem \ref{Tsmr3.1} holds,
i.e. $F:{\Bbb R} \to {\Bbb R}$, $b \in {\Bbb R}$ and $\gamma \in 
{\Bbb R}_+$, $0 < \kappa < 1$ are given and assumptions A1 -- A4 
are satisfied.

\section{Notation and definitions}
\label{npc}

In this section we explain the notation and definitions 
used in chapters II and  III to prove theorem \ref{Tsmr3.1}.
For quick reference we have also included a table of notation
in section \ref{tnot} at the end of this paper. There is a reason why
we introduce all the notation before we begin the proof. As we will 
see there are many parameters in the construction which must satisfy
a large number of conditions. If we would introduce these parameters only
when they first appear in the proof it would be a difficult task to
verify that no circular or contradictory definitions have been made.
This issue is resolved in the following way. We present the definitions
in such an order that by starting from the quantities given in theorem 
\ref{Tsmr3.1} and in assumptions A1 -- A4 one may verify one by one
that the definitions make sense. However, from  part {\bf H)} on 
the definitions will depend on a parameter $n_0 \in {\Bbb N}$. We
show at the end of the present section (lemma \ref{hyp}) that 
$n_0$ can be chosen in such a way that all conditions on the 
parameters which appear in the proof of theorem \ref{Tsmr3.1} 
are satisfied. To facilitate the proof of lemma \ref{hyp} we collect at the 
beginning of each of the sections \ref{fis} - \ref{pmt}
the conditions which are used in that particular section.  

\noindent
{\bf A) General Notation} \label{npc.10}

For any $d \in {\Bbb R}$ we denote by $| \cdot |$ the maximum norm
on ${\Bbb Z}^d$, ${\Bbb R}^d$ and ${\Bbb C}^d$. Sometimes we also use
$\ell_1$ -- or $\ell_2$ -- norms on these spaces which are denoted by
$| \cdot |_1$, respectively $| \cdot|_2$. 
In our notation we distinguish real and complex neighborhoods.
For sets $X \subset {\Bbb R}^d$, $Z \subset {\Bbb C}^d$ and 
positive numbers $\rho$ we define
\begin{eqnarray}
\label{npc.20}
B_{\rho}(X) &:=& B(X, \rho) := \{y \in {\Bbb R}^d: |y-x| < \rho \;
\mbox{ for some } \; x \in X \}, \\
\label{npc.30}
U_{\rho}(Z) &:=& U(Z, \rho) := \{w \in {\Bbb C}^d: |w - z| < \rho \;
\mbox{ for some } \; z \in Z \}.
\end{eqnarray}
In case the sets $X = \{x\}$, $Z = \{z\}$ contain only a single point
we usually replace in (\ref{npc.20}), (\ref{npc.30}) 
the sets $X$, $Z$ by the elements $x$, $z$.
Note that in this case $B_{\rho}(x)$ defines a cube and $U_{\rho}(z)$
defines a polydisc since we use the maximum norm for the definition of
the neighborhoods. By a slight abuse of notation we sometimes
restrict the neighborhoods to the lattice ${\Bbb Z}^d$, i.e.
$B_N(n) = \{m \in {\Bbb Z}^d: |m - n| < N \}$ 
for lattice sites $n \in {\Bbb Z}^d$ and $N > 0$. 

For $a$, $b \in {\Bbb C}^d$ we define 
\begin{eqnarray}
\label{npc.40}
<a, b> := \sum_{i=1}^d a_i b_i.
\end{eqnarray}
Note that there are no complex conjugates in this definition.

The function $[ \cdot ]$ was defined in (\ref{ova.560}). We also use
the floor function
\begin{eqnarray}
\label{npc.50}
\lfloor \cdot \rfloor \; :\; {\Bbb R} \to {\Bbb Z}\;\; ; \;\;\; 
\lfloor x \rfloor = \max \{ n \in {\Bbb Z}: n \leq x \}.
\end{eqnarray}

\noindent
{\bf B) Notation related to the statement of theorem \ref{Tsmr3.1}}

We are given the function 
$F$ and the numbers $b$, $\gamma$, $\kappa$ through the hypothesis
of theorem \ref{Tsmr3.1}. The following quantities are defined
due to assumptions A1 -- A4: the number of phases $\nu$ (\ref{smr2.20}),
the frequency vector of the linearized system $\omega^{(0)}$ (\ref{smr2.30}),
(\ref{smr2.40}), the vector $g$ (\ref{smr1.70}), positive constants
$s$ and $\tau$ related to assumption A3 (\ref{smr2.50}), and the 
matrix $\Omega$ (\ref{smr2.60}) -- (\ref{smr2.90}). 

\noindent
{\bf C) Notation related to the weight functions}

In section \ref{fa} we define families of weight functions on ${\Bbb Z}^{\nu}$
and corresponding weighted norms for sequence and matrix spaces. We define
the following parameters
\begin{eqnarray}
\label{npc.52}
c &:=& 0.01 \\
\label{npc.54}
D_N &:=& \max (D_{\nu, c}, D_{\nu, 1}) \geq 1,
\end{eqnarray}
where $D_{\nu, c}$ and $D_{\nu, 1}$ are defined in (\ref{wf.10}) and
shown to be finite in proposition \ref{Pwf.1}.
Weight functions $w_{\sigma, c}$, $w_{\sigma, 1}$ are defined for
$\sigma \geq 1/4$ in \ref{DFa1}. The corresponding sequence
spaces $(X_{\sigma, c}, \| \cdot \|_{\sigma, c})$, 
$(X_{\sigma, 1}, \| \cdot \|_{\sigma, 1})$ and matrix spaces
$({\cal L}_{\sigma, c}, \| \cdot \|_{\sigma, c})$, 
$({\cal L}_{\sigma, 1}, \| \cdot \|_{\sigma, 1})$ are defined 
in \ref{DFa2} (see also  (\ref{fa1.20}), (\ref{fa1.60})). 

With each application of the coupling lemma \ref{Lcl.1} 
we loose some of the decay properties of the matrices. Therefore 
we use a scale of spaces ${\cal L}_{\sigma_j, c}$ with
\begin{eqnarray}
\label{npc.56}
\sigma_j := \frac{1}{2}(1 + \frac{1}{j})
\;\;\; \mbox{ for } j \geq 1.
\end{eqnarray}

\noindent
{\bf D) Notation related to the nonlinear part of the equation}

The function $F$  is assumed to be real analytic in a neighborhood of
$-b$. We denote by $\tilde{r}_{F, b}$ the radius of convergence of
the corresponding power series at $-b$ and set 
\begin{eqnarray}
\label{npc.60}
r_{F, b} &:=& \min(1, \tilde{r}_{F, b}) , \\
\label{npc.70}
F(y) &=& \sum_{k=0}^{\infty} \alpha_k (y+b)^k \;\;\; \mbox{ for }
|y + b| < r_{F, b}
\end{eqnarray}
which defines the sequence $(\alpha_k)_{k \geq 0}$ 
($\alpha_k = F^{(k)}(-b)/(k!)$).
The nonlinear part of the Fourier equation is given by
\begin{eqnarray}
\nonumber
W(u) = \sum_{k=2}^{\infty} \alpha_k u^{*k},
\end{eqnarray}
where $u^{*k}$ denotes the $k$-th convolution power of the 
sequence $u$ (see notation \ref{Nfa1.1}). The definition of $W(u)$ is purely
formal at this point. In lemma \ref{Lenp.1} the definition is made 
rigorous for sequences $u$ in $X_{\sigma, c}$, resp. $X_{\sigma, 1}$. Furthermore,
the constant $D_W$ is defined through lemma \ref{Lenp.1}.

\noindent
{\bf E) Notation related to the Lyapunov -- Schmidt decomposition}

In order to describe the Lyapunov -- Schmidt reduction we defined the 
diagonal operators $D(\omega)$ (\ref{ova.85}),
the set ${\cal S}$ (\ref{ova.90}), the projections $Q$ (\ref{ova.100})
and $P$ (\ref{ova.110}), the parameterization 
$\varphi(a)$ (\ref{ova.120}), the parameter $\lambda^{(0)}$ 
(\ref{ova.140}) and the functions ${\cal P}$, ${\cal Q}$ (\ref{ova.160}),
(\ref{ova.150}) which describe the $P$ -- equation, 
respectively $Q$ -- equation of
the Lyapunov - Schmidt reduction.

Recall that we first solve the $P$ -- equation for $v$. The parameters
of this equation $(a, \omega)$ are denoted by $\lambda$ 
(see (\ref{ova.130})).

\noindent
{\bf F) Notation related to the linearized equations}

For $j \in {\Bbb N}$ we denote by $v_j(\lambda)$ the $j$-th
approximation of the solution $v(\lambda)$ of the $P$ -- equation.
As explained in section \ref{ova} we need to investigate the invertibility
of the following matrices
\begin{eqnarray}
\label{npc.63}
T^{(j)}(\lambda)(m, n) &=& V(\omega)(m) \delta_{m, n} +
DW(\varphi(a) + v_j(\lambda)) (m, n), \\
\label{npc.64}
T^{(j)}(\theta, \lambda)(m, n) &=& V(\theta, \omega)(m) \delta_{m, n} +
DW(\varphi(a) + v_j(\lambda)) (m, n),
\end{eqnarray}
where $V(\omega)$, $V(\theta, \omega)$ were defined in (\ref{smr2.70}),
(\ref{ova.520}). For $l \in {\Bbb N}$ set
\begin{eqnarray}
\label{lop.50}
V_l: {\Bbb C} \to {\Bbb C} \cup \{ \infty \} \; ; \;\;\;
V_l(\vartheta) := \alpha_1 - 
\frac{l^2 \gamma^2}{4 \sin^2 \frac{\vartheta}{2}}.
\end{eqnarray}
Observe that $V(\theta, \omega)(m) = V_{|<m, g>|}(\theta + <\omega, m>)$.
In our analysis we use a few estimates on the functions $V_l$, mainly in
a neighborhood of their zeros. It is easy to see from assumption
A2 that $V_l$ has zeros only for $1 \leq l \leq \nu$. In this case the
set of zeros of $V_l$ is given by $\{ \pm \omega_l^{(0)} \} + 2 \pi {\Bbb Z}$.
The properties of the functions $V_l$, $1 \leq l \leq \nu$, which are needed
in the proof of theorem \ref{Tsmr3.1} are collected in proposition 
\ref{Plop.1}. In this proposition one also finds the definition of the
constants $d_V$, $\delta_V$ and $D_V$.

The quantities $V(\omega)(m)$, $V(\theta, \omega)(m)$ are not well
defined if $<\omega, m> \in 2 \pi {\Bbb Z}$, respectively 
$\theta + <\omega, m> \in 2 \pi {\Bbb Z}$. Hence for $C \subset 
{\Bbb Z}^{\nu}$ the restricted matrices $T^{(j)}_C(\lambda)$,
$T^{(j)}_C(\theta, \lambda)$ are not well defined for $\lambda \in 
Z_C$, respectively $(\theta, \lambda) \in \tilde{Z}_C$. The sets
$Z_C$, $\tilde{Z}_C$ were defined in (\ref{ova.600}), (\ref{ova.610}).

\noindent
{\bf G) Special constants and functions}

In the course of proving theorem \ref{Tsmr3.1} a number
of constants and functions are used. We list them here. Of particular interest
are definitions (\ref{npc.72}) -- (\ref{npc.77}) as they are closely
related to the choice of parameters in the multi-scale analysis (see
(\ref{npc.410}) -- (\ref{npc.440})).

\begin{eqnarray}
\label{npc.72}
A &:=& 25600 \nu^8 \left(2 \lfloor \tau \rfloor + 5 \right) \\
\label{npc.73}
A_1 &:=& 4 \nu \\
\label{npc.75}
E_{\rho} &:=& 32 \nu^3 (\lfloor \tau \rfloor + 3) \\
\label{npc.76}
E_{\delta} &:=& 32 \nu^3 (\lfloor \tau \rfloor + 1) \\
\label{npc.77}
E_{M} &:=& 32 \nu^3 (2\lfloor \tau \rfloor + 5) \\
\label{npc.78}
q &:=& 1 + A(E_M + E_{\rho}) \\
\label{npc.80}
B_0 &:=& 22 \nu^3 \\
\label{npc.90}
B_1 &:=& 8 \nu^2 \\
\label{npc.100}
B_2 &:=& 18 \nu^3 \\
\label{npc.110}
D_{\tau} &:=& \sum_{k=1}^{\infty}
5^{k(\tau+1)} e^{-1.5^k} \\
\label{npc.118}
D_{E}(x) &:=& \sum_{k=1}^{\infty} 5^{k E_{\rho} (x-1)} e^{-1.5^k} \;\;
\mbox{ for } x \geq 1 \\
\label{npc.120}
d_{\tau, c} &:=& \min_{x \geq 1} \frac{e^{\frac{1}{2} x^c}}{2 x^{\tau}} \\
\label{npc.130}
d_{min} &:=& \min_{\xi \in {\Bbb R}^{\nu}, |\xi|_2 = 1} |\Omega \xi| \\
\label{npc.140}
\tilde{d}_{min} &:=& \min_{1 \leq j, l \leq \nu} |\Omega_{j, l}| \\
\label{npc.150}
D_{1,1} &:=& 4 D_N \left( \sum_{p=1}^{B_0+1} p^{B_0} \right) 
\left( \frac{16 \alpha_1}{d_V \gamma^2} \right)^{B_0+1} 2^{E_{\rho} B_0}
\\
\label{npc.160}
D_{1,2} &:=& 2 \cdot 8^{B_0} (B_0!)  
\left( \sum_{p=1}^{B_0+1} p^{B_0} D_N^p \right)
\\
\label{npc.170}
D_{1,3} &:=& (8 D_N \max(D_{1,1}, D_{1,2}) )^{B_0+1} 
\left( \sum_{p=1}^{B_0+1} p^{B_0} \right)
\left( \max_{y \geq 0} (1 +y)^{q(2 B_0 +1)} e^{-\frac{1}{4} y^c} \right)
\\
\label{npc.180}
D_{1,4} &:=& 2^{B_0} D_V + 4^{B_0} D_V D_{1, 3} 
\left( \max_{y \geq 0} (1 +y)^{q(2 B_0 +1)} e^{-\frac{1}{4} |y-1|^c} \right)
\\
\label{npc.190}
D_{1,5} &:=& (2 \nu)! 2^{2 \nu} (2 \nu)^{B_0}
\\
\label{npc.200}
D_{1} &:=& (2 \nu)! (2 D_{1, 4})^{2 \nu} (2 \nu)^{B_0}
\\
\label{npc.210}
D_{2} &:=& (2 \nu + B_1)^{B_1^2} (75 D_1 (B_1!) 2^{B_1})^{2 B_1 -1}
\\
\label{npc.220}
D_{3} &:=& (B_1 + 1)! 37 e^6 \nu^2 D_W D_N^2 D_{\tau} D_E(B_1 + 1)
\\
\label{npc.230}
D_{4,2} &:=& D_W D_3^{B_1 + 1}
\left( \sum_{p=1}^{B_1+1} \frac{1}{p!} p^{B_1 + 1} \right) 
\\
\label{npc.250}
D_{4,1} &:=& D_{V} + D_{4, 2} 
\\
\label{npc.260}
D_{4} &:=& D_{4, 1} (\nu + B_1)^{B_1} \frac{2}{d_V} \\
\label{npc.480}
D_P &:=& (4 \nu^2)! 2^{4 \nu} (4 \nu)^{B_1} D_2^{4 \nu} \\
\label{npc.115}
D_K &:=& \nu^2 D_P + 10 \nu^2 4^{(2 \nu+1)^2} [(2\nu + 1)^2 !] 
\end{eqnarray}

\noindent
{\em Definition of the mollifier $\psi$}.
We choose a function $\psi: {\Bbb R}^{4 \nu} \to {\Bbb R}$ with the 
following four properties.
\begin{itemize}
\item[(1)] $\psi \in C^{\infty}({\Bbb R}^{4 \nu})$
\item[(2)] supp$(\psi) \subset \{ x \in {\Bbb R}^{4 \nu}: |x| < 1 \}$.
\item[(3)] $\psi(x) \geq 0$ for all $x \in {\Bbb R}^{4 \nu}$.
\item[(4)] $\int_{{\Bbb R}^{4 \nu}} \psi(x) dx = 1$.
\end{itemize}
For $k \in {\Bbb N}_0$ define
\begin{eqnarray}
\label{npc.265}
D_{\psi}(k) := \sup_{|\beta|_1 \leq k} \int_{{\Bbb R}^{4 \nu}} 
|\partial^{\beta} \psi (x)| dx.
\end{eqnarray}

\noindent
{\em Definition of $\tilde{\kappa}$}.
The set ${\cal Z}$ in the statement of theorem \ref{Tsmr3.1} turns out to be 
of the form
\begin{eqnarray}
\nonumber
{\cal Z} = \{(z_1, \ldots, z_{\nu}) \in {\Bbb C}^{\nu}: 
(|z_1|, \ldots, |z_{\nu}|) \in {\cal M}^{\infty} \}
\end{eqnarray}
for some suitable set ${\cal M}^{\infty} \subset {\Bbb R}^{\nu}$.
We need to understand how lower bounds on the $\nu$-dimensional 
Lebesgue measure of the set ${\cal M}^{\infty}$ translate to lower bounds
on the $2\nu$-dimensional 
Lebesgue measure of the set ${\cal Z}$ regarded as a subset of 
${\Bbb R}^{2 \nu}$. For a set $M \subset {\Bbb R}^{\nu}$ denote
\begin{eqnarray}
\nonumber
M^{{\Bbb C}} := \{(z_1, \ldots, z_{\nu}) \in {\Bbb C}^{\nu}: 
(|z_1|, \ldots, |z_{\nu}|) \in M \} \subset {\Bbb C}^{\nu} \equiv
{\Bbb R}^{2 \nu}.
\end{eqnarray}
Recall that $\kappa \in (0, 1)$ was given in theorem \ref{Tsmr3.1}. 
We define
\begin{eqnarray}
\label{npc.270}
\tilde{\kappa} := \inf \{
\mbox{ vol}^{\nu}(B_{1}(0) \setminus M) : M \subset B_{1}(0) 
\mbox{ and vol}^{2\nu}(M^{\Bbb C}) \leq \kappa \mbox{ vol}^{2 \nu}
(U_{1}(0)) \}.
\end{eqnarray}
It is elementary to show the following proposition.
\begin{Proposition}
\label{Pnpc.1}
Given $\nu \in {\Bbb N}$, $0 < \kappa < 1$ and let $\tilde{\kappa}$
be defined as in (\ref{npc.270}). Then
\begin{itemize}
\item[(a)] $\tilde{\kappa} > 0$.
\item[(b)] For any $\rho > 0$ and any measurable set $M \subset B_{\rho}(0)$
with vol$^{\nu}(B_{\rho}(0) \setminus M) \leq \tilde{\kappa} \rho^{\nu}$
it follows that
vol$^{2 \nu}(M^{{\Bbb C}}) \geq \kappa \mbox{ vol}^{2 \nu}(U_{\rho}(0))$.
\end{itemize}
\end{Proposition}

\noindent 
{\bf H) Parameters of the multi-scale analysis depending on $n_0$}
\begin{eqnarray}
\label{npc.400}
N_0 &:=& 5 ^{n_0}\\
\label{npc.410}
N_j &:=& N_0^{A^j}\\
\label{npc.420}
\rho_j &:=& N_j^{-E_{\rho}}\\
\label{npc.425}
\tilde{\rho}_j &:=& N_j^{-A_1 E_{\rho}} \\
\label{npc.430}
\delta_j &:=& N_j^{-E_{\delta}}\\
\label{npc.440}
M_j &:=& N_j^{E_M}
\end{eqnarray}

\newpage
\noindent 
{\bf I) The sets ${\cal C}^{(j)}$ and ${\cal POL}$}

We define the following classes of subsets in ${\Bbb Z}^{\nu}$ (cf. remark
\ref{Rova.3}).
\begin{eqnarray}
\label{npc.450}
{\cal C}^{(1)} &:=& 
\{ C \subset B_{1.5}(0): C \cap {\cal S} \neq \emptyset \}, \\
\label{npc.460}
{\cal C}^{(j)} &:=&
\{ C \subset B_{N_{j-1}}(0): C \cap {\cal S} \neq \emptyset \mbox{ and }
C \mbox{ interval } \}, \;\; j \geq 2,
\end{eqnarray}
where we say that a set $C \subset {\Bbb Z}^{\nu}$ is an interval, if
there exist real numbers $a_i$, $b_i$, $1 \leq i \leq \nu$ such that
$C = \{ n \in {\Bbb Z}^{\nu}: a_i < n_i < b_i \mbox{ for all }
1 \leq i \leq \nu \}$.

Recall the definition of the constant $D_P$ in (\ref{npc.480}).
We call a map $p: {\Bbb C} \times U \to {\Bbb C}$ an {\em admissible
polynomial}, if $\emptyset \neq U \subset {\Bbb C}^{2 \nu}$ and
there exist $d \in {\Bbb N}$, and analytic functions
$c_i : U \to {\Bbb C}$, $0 \leq i < d$, such that 
\begin{eqnarray}
\nonumber
p(\theta, \lambda) = \theta^d + \sum_{0 \leq i< d} c_i(\lambda) \theta^i
\end{eqnarray}
and the following three conditions are satisfied:
\begin{itemize}
\item[(1)] $1 \leq d \leq 4 \nu^2$,
\item[(2)] $c_i(\lambda) \in {\Bbb R}$ for all
$\lambda \in U \cap {\Bbb R}^{2 \nu}$, $0 \leq i < d$,
\item[(3)] $|\partial^{\beta} c_i (\lambda)| \leq D_P M_0^{4(B_0 + \nu)
|\beta|_1}$ for all $\lambda \in U$, $0 \leq i < d$, 
$0 \leq |\beta|_1 \leq B_1$.
\end{itemize}
We define
\begin{eqnarray}
\label{npc.490}
{\cal POL} := \{ (p, \vartheta): p \mbox{ is an admissible polynomial and }
\vartheta \in {\Bbb R} \}.
\end{eqnarray}

\noindent 
{\bf J) Cube decompositions and the projections $\pi^{(j)}_l$}

Let $\rho$, $\tilde{\rho} > 0$ with $\rho/\tilde{\rho} \in {\Bbb N}$ and 
let $x \in {\Bbb R}^{2 \nu}$. Then there exists a set 
$\{ y_k \in {\Bbb R}^{2 \nu}: 1 \leq k \leq (\rho/\tilde{\rho})^{2 \nu} \}$
which is uniquely defined
by the following two conditions.
\begin{eqnarray}
\nonumber
\overline{B_{\rho}(x)} &=& \bigcup_{k} \overline{B_{\tilde{\rho}}(y_k)} \\
\nonumber
B_{\tilde{\rho}}(y_k) \cap B_{\tilde{\rho}}(y_l) &=& \emptyset \; \mbox{ for }
l \neq k. 
\end{eqnarray}
We use such cube decompositions for the sets of parameters ${\cal N}^{(j)}$
which appear in the induction statements in sections \ref{fis} and \ref{is}.
\begin{eqnarray}
\label{npc.495}
\overline{{\cal N}^{(j)}} = \bigcup_{k \in K^{(j)}} 
\overline{B_{\rho_j}(\lambda_k^{(j)})} =
\bigcup_{k \in K^{(j + 0.5)}} 
\overline{B_{\tilde{\rho}_j}(\lambda_k^{(j+0.5)})}
\end{eqnarray}
The sets $K^{(l)}$, $\{ \lambda_k^{(l)}: k \in K^{(l)} \}$,
$K^{(l+0.5)}$, $\{ \lambda_k^{(l+0.5)}: k \in K^{(l+0.5)} \}$ are defined 
inductively in $l$.
We will define $K^{(1)} = \{ 1 \}$ and $\lambda_1^{(1)} = \lambda^{(0)}$
(see (\ref{fis.111}) -- (\ref{fis.113}) below). Each induction step 
$j \to j+1$ consists of two parts. First we perform a cube
decomposition of each cube $B_{\rho_j}(\lambda_k^{(j)})$, $k \in K^{(j)}$
into sub-cubes of radius $\tilde{\rho}_j$, yielding sets 
$K^{(j+0.5)}$, $\{ \lambda_k^{(j+0.5)}: k \in K^{(j+0.5)} \}$ (see
also (\ref{npc.495})).
Observe that 
$\rho_j/\tilde{\rho}_j \in {\Bbb N}$, since $A$, $A_1$, $E_{\rho}$, $N_j$
are integers (see (\ref{npc.72}) -- (\ref{npc.75}), (\ref{npc.400}) --
(\ref{npc.425}), $n_0 \in {\Bbb N}$ by lemma \ref{hyp}) and the cube
decomposition is well defined. Then one performs a cube decomposition
on each cube $B_{\tilde{\rho}_j}(\lambda_k^{(j+0.5)})$, $k \in K^{(j+0.5)}$
into sub-cubes of radius $\rho_{j+1}$ ($\tilde{\rho}_j / \rho_{j+1} \in 
{\Bbb N}$ since $A / A_1 \in {\Bbb N}$) obtaining sets
$\hat{K}^{(j+1)}$, $\{ \lambda_k^{(j+1)}: k \in \hat{K}^{(j+1)} \}$.
The set $K^{(j+1)}$ is defined as a subset of $\hat{K}^{(j+1)}$ by removing
those cubes which contain resonant parameters. 
We define the following projection map between 
sets $K^{(l)}$. Let $l$, $j \in \frac{1}{2}({\Bbb N} + 1)$ with $l \leq j$.
\begin{eqnarray}
\label{npc.500}
\pi^{(j)}_l : K^{(j)} \to K^{(l)}, \;\; 
\pi^{(j)}_l(k) \mbox{ is defined by the condition } \;
\lambda_k^{(j)} \in B_{\rho_l} \left(\lambda_{\pi^{(j)}_l(k)}^{(l)} \right).
\end{eqnarray}
In the case that $l$ is not an integer we understand
$\rho_{l} \equiv \tilde{\rho}_{l - 0.5}$.

We conclude this section by determining the parameter $n_0$.
\begin{Lemma}
\label{hyp}
Suppose that $F:{\Bbb R} \to {\Bbb R}$, $b \in {\Bbb R}$, $\gamma \in 
{\Bbb R}_+$ and $0 < \kappa < 1$ 
are given and satisfy assumptions A1 -- A4. 
There exists an integer $n_0 \in {\Bbb N}$  such
that conditions
(\ref{fis5}) -- (\ref{fis25}), (\ref{cv2.5}) -- (\ref{cv2.65}), 
(\ref{sp.50}), (\ref{est3.10}) -- (\ref{est3.20}), 
(\ref{pmt.10}) -- (\ref{pmt.30}) are satisfied;
conditions (\ref{conf.2}) -- (\ref{dpoh.2}), 
(\ref{comp.30}) -- (\ref{comp.70}), (\ref{est1.200}) -- (\ref{est4.9})
are satisfied for $j \geq 1$ and conditions (\ref{cvj.5}) -- (\ref{cvj.80})
are satisfied for $j \geq 2$.
\end{Lemma}
For the remaining part of chapter II and for chapter III we will let 
$n_0$ be a fixed integer satisfying the conditions stated in lemma \ref{hyp}.
This completes the definitions given in parts H) and I) of this section.

\begin{proof}
One proves lemma \ref{hyp} by showing that each condition is satisfied
for sufficiently large values of $n_0$. Of course, in case the condition 
depends on the induction step $j \in {\Bbb N}$ we must ensure that
the lower bound on $n_0$ is uniform in $j$. For the proof it is important
to observe that the only quantities appearing in the conditions which
depend on $n_0$ are $N_j$, $\rho_j$, $\tilde{\rho}_j$, $\delta_j$ and 
$M_j$ ($j \geq 0$). We choose a few examples to demonstrate the method
of proof.

$\bullet$ (\ref{cv2.35}):
It follows from assumption A2 and from 
the definition of $\nu$ in (\ref{smr2.20}) that 
\begin{eqnarray}
\label{npc.600}
\frac{(\nu+1)^2 \gamma^2}{4} - \alpha_1 > 0.
\end{eqnarray}
The only quantity in (\ref{cv2.35}) which depends on $n_0$ is $N_0$.
Since $E_{\rho} > 1$ (\ref{npc.75}) and $\tau > 0$ (assumption A3)
(\ref{npc.600}) implies that (\ref{cv2.35}) holds for sufficiently
large $n_0$.

$\bullet$ (\ref{cv2.40}): The only quantities, depending on $n_0$ are
$N_0$ and $\rho_0 = N_0^{-E_{\rho}}$. Since $E_{\rho} > \tau + 1$
(see (\ref{npc.75})) it is clear that $N_0^{\tau + 1} \rho_0 \to 0$
as $n_0 \to \infty$, proving (\ref{cv2.40})

$\bullet$ (\ref{cvj.45}):
Using (\ref{npc.410}), (\ref{npc.440}) it suffices to show that for all 
$j \geq 1$
\begin{eqnarray}
\label{npc.610}
\log (4 D_N) + A E_M \log (N_j) \leq 
\frac{1}{4(j+1)(j+2)} N_j^c. 
\end{eqnarray}
We obtain the condition on $n_0$ by investigating (\ref{npc.610}) 
inductively. For $j = 1$ we need
\begin{eqnarray}
\nonumber
\log (4 D_N) + n_0 A^2 E_M \log 5 \leq \frac{1}{24} 5^{Ac n_0}
\end{eqnarray}
which is clearly satisfied for sufficiently large $n_0$.
For the induction step $j \to j+1$, $j \geq 1$, we estimate the ratios
\begin{eqnarray}
\nonumber
0 \leq \frac{\log (4 D_N) + A E_M \log (N_{j+1})}
{\log (4 D_N) + A E_M \log (N_j)} &\leq& A, \\
\nonumber
\frac{4(j+1)(j+2) N_{j+1}^c}{4(j+2)(j+3)N_j^c} &\geq& \frac{1}{2} N_j^{(A-1)c}
\geq \frac{1}{2} 5^{A(A-1)c n_0}.
\end{eqnarray}
For $n_0$ sufficiently large we have 
\begin{eqnarray}
\nonumber
A \leq \frac{1}{2} 5^{A(A-1)c n_0}
\end{eqnarray}
completing the proof of (\ref{npc.610})

$\bullet$ (\ref{esf.80}):
We may replace (\ref{esf.80}) by the following three conditions
\begin{eqnarray}
\label{npc.620}
(8/\rho_{j+1})^{B_0} &\leq& e^{\frac{1}{32} N_j^c} \;\;\; \mbox{ for all }
j \geq 1, \\
\label{npc.630}
B_0! D_W e^{-\frac{1}{32} N_1^c} &\leq& \frac{1}{2}, \\
\label{npc.640}
e^{-\frac{1}{32} N_{j+1}^c} &\leq& \frac{1}{2} e^{-\frac{1}{32} N_j^c}
\;\;\; \mbox{ for all }
j \geq 1.
\end{eqnarray}
Observe that there exists a number $\tilde{N} > 0$ such that for all 
$N \geq \tilde{N}$
\begin{eqnarray}
\nonumber
(8 N^{A E_{\rho}})^{B_0} \leq e^{\frac{1}{32} N^c}.
\end{eqnarray}
If $n_0$ is chosen large enough such that $5^{A n_0} \geq \tilde{N}$ then
(\ref{npc.620}) is satisfied. Condition (\ref{npc.630}) can clearly be 
satisfied by choosing $n_0$ large and condition (\ref{npc.640}) is
satisfied if $N_1^c(N_1^{c(A-1)} - 1) \geq 32 \log 2$ which again is achieved
by making $n_0$ large.

$\bullet$ (\ref{awp.8}):
Recall that $k! \geq (k/e)^k$. Since $\nu \geq 2$ we have $B_2/4 \geq 12 \nu$
and $B_2 \geq 2\nu(B_1+1)$ (see (\ref{npc.90}), (\ref{npc.100})). 
Hence $B_2! \geq (12 \nu)^{B_2} \geq  (12 \nu)^{2\nu(B_1+1)}$.

$\bullet$ (\ref{est4.9}): 
Since $N_j = 5^{n_0 A^j}$ and $A \geq 2$ it is clear that $N_j \geq 2^j$
independent of the choice of $n_0 \in {\Bbb N}$. We can therefore
replace (\ref{est4.9}) by
\begin{eqnarray}
\label{npc.650}
N_j^{-A E_{\delta}/B_1 + A_1 E_{\rho} + (\nu + 2)A +(4 \nu +1) + 2 E_{\delta}}
\leq \tilde{\kappa} \;\;\; \mbox{ for all } j \geq 1.
\end{eqnarray}
One verifies from the definitions of $A$, $A_1$, $E_{\delta}$, $E_{\rho}$ and
$B_1$ that the exponent in (\ref{npc.650}) is negative. Since
$N_j \geq N_1$ for all $j \geq 1$ one can satisfy (\ref{est4.9})
simultaneously for all values of $j$ by choosing $n_0$ large enough.
\end{proof}
 
\begin{Remark}
\label{Rnpc.1}
{\em
The conditions which appear in lemma \ref{hyp} not only determine $n_0$
but they are also the reason behind the choices which were made for 
$c$, $B_0$, $B_1$, $B_2$, $A$, $A_1$, $E_{\rho}$, $E_{\delta}$, $E_M$.
To determine these constants on can proceed for example in the following
way. First define $c$ satisfying (\ref{cvj.65}), then choose
$B_0$, $B_1$ and $B_2$ via (\ref{est3.475}), (\ref{awp.4}) and (\ref{awp.1}).
We are then left with conditions on $A$, $A_1$, $E_{\rho}$, 
$E_{\delta}$ and $E_M$. 
The conditions for these quantities are somewhat involved
as one can already see from the condition 
that the exponent in (\ref{npc.650})
is negative appearing in the proof of (\ref{est4.9}) above. 
The following list shows all conditions which contain at least
two quantities of $A$, $A_1$, $E_{\rho}$, $E_{\delta}$, $E_M$.
\\
$\bullet$ $E_{\delta}+1 < E_{\rho}$:
(\ref{fis15}), (\ref{cvj.40}), (\ref{conf.25}), (\ref{conf.37}),
(\ref{conf.45}), (\ref{esf.5}), (\ref{esf.77}), (\ref{esf6.310}),
(\ref{comp.35}), (\ref{comp.40})
\\
$\bullet$ $E_M > E_{\delta} + E_{\rho}$:
(\ref{esf.75}), (\ref{esf.90}), (\ref{cvj.50}), (\ref{conf.55})
\\
$\bullet$ $A E_{\delta} > E_{\rho} - 1$: (\ref{conf.15})
\\
$\bullet$ $A E_{\rho} > E_M(2 B_0 + 1)$: (\ref{esf.95})
\\
$\bullet$ $A E_{\delta} > E_{M}$: (\ref{esf6.320})
\\
$\bullet$ $2 E_M (B_0 + \nu) + 2 \nu A E_{\delta} < E_{\rho} A A_1$:
(\ref{awp.6})
\\
$\bullet$ $A^2( E_{\rho} - E_{\delta}) > 4 E_M (B_0 + \nu)$: (\ref{comp.40})
\\
$\bullet$ $2 E_M + A E_{\delta} + E_{\delta} < A E_M$: (\ref{comp.60})
\\
$\bullet$ $16 \nu^2 B_1 E_{\rho} + 32 \nu^2 (B_0 + \nu) E_M < A$: 
(\ref{est3.490})
\\
$\bullet$ 
$A_1 E_{\rho} + (\nu + 2)A +(4 \nu +1) + 2 E_{\delta} < A E_{\delta}/B_1$:
(\ref{est4.9})
}
\end{Remark}

\section{Fourier analysis}
\label{fa}

In this section we perform a Fourier transform of our original equation.
We begin by defining weighted sequence spaces 
(for the Fourier coefficients)
and corresponding operator spaces. Similar classes of weights
have been used
in the context of small divisor problems by DeLatte \cite{DeL},
Craig and Wayne \cite{CW2} and Bourgain \cite{B1}. The weights and
norms we introduce in section \ref{fa1} are well suited for a
multi-scale analysis and have been presented by the author in 
\cite{Kri}.

\subsection{Sequence and matrix spaces}
\label{fa1}

The Fourier coefficients of the profile functions $\chi( \cdot;z)$
form sequences on the lattice ${\Bbb Z}^{\nu}$. For our analysis it
is convenient to consider certain families of weighted $\ell_1$ spaces.
The corresponding weights belong to the following class which was already
used in \cite{DKV}.
\begin{eqnarray}
\label{fa1.10}
{\cal W} := 
\{
w: {\Bbb Z}^{\nu} \to [1, \infty)\;:\;
w(m) \leq w(n) w(m-n) \; \mbox{ for all } m, n \in {\Bbb Z}^{\nu}
\}
\end{eqnarray}
For weights $w \in {\cal W}$ we define the sequence space
\begin{eqnarray}
\label{fa1.20}
\ell_{1, w} := \left\{
u: {\Bbb Z}^{\nu} \to {\Bbb C} \;:\;
\sum_{m \in {\Bbb Z}^{\nu}} w(m) |u(m)| < \infty
\right\}.
\end{eqnarray}
It is easy to see that the space $\ell_{1, w}$ together with the 
norm
\begin{eqnarray}
\label{fa1.30}
\|u\|_{\ell_{1, w}} := \sum_{m \in {\Bbb Z}^{\nu}} w(m) |u(m)|
\end{eqnarray}
is a Banach space which is contained in $\ell_1$. What will be important 
in our analysis is that the norms defined in (\ref{fa1.30}) are
sub-multiplicative with respect to convolution.
For $u$, $v \in {\ell_1}$ set 
\begin{eqnarray}
\label{fa1.40}
(u * v)(m) := \sum_{n \in {\Bbb Z}^{\nu}} u(m-n) v(n).
\end{eqnarray}
The following proposition is a straightforward consequence of
definition (\ref{fa1.10}) and the well known sub-multiplicativity
of the $\ell_1$ -- norm under convolution.
\begin{Proposition}
\label{Pfa1.1}
Let $w \in {\cal W}$ and $u$, $v \in \ell_{1, w}$. Then $u*v \in
\ell_{1, w}$ and
\begin{eqnarray}
\label{fa1.50}
\| u*v \|_{\ell_{1, w}} \leq \| u \|_{\ell_{1, w}} \| v \|_{\ell_{1, w}}.
\end{eqnarray} 
For $\xi \in {\Bbb R}^{\nu}$ we have
\begin{eqnarray}
\label{fa1.55}
\left(\sum_{m \in {\Bbb Z}^{\nu}} u(m) e^{i<m, \xi>} \right) \cdot
\left(\sum_{m \in {\Bbb Z}^{\nu}} v(m) e^{i<m, \xi>} \right) =
\sum_{m \in {\Bbb Z}^{\nu}} (u*v)(m) e^{i<m, \xi>}
\end{eqnarray}
and
\begin{eqnarray}
\label{fa1.57}
(u e^{i<\cdot, \xi>}) * (v e^{i<\cdot, \xi>}) = (u*v) e^{i<\cdot, \xi>}.
\end{eqnarray} 
\end{Proposition}   

\begin{Notation}
\label{Nfa1.1}
For $u \in \ell_1$ and $k \in {\Bbb N}$ we denote by
\begin{eqnarray}
\nonumber
u^{*k} := u * u * \ldots * u
\end{eqnarray}
the $k$-th convolution power of $u$. 
\end{Notation}

Next we turn to spaces of linear operators on $\ell_{1, w}$. For
$w \in {\cal W}$ we define
\begin{eqnarray}
\label{fa1.60}
{\cal L}_{w} :=
\left\{
R: {\Bbb Z}^{\nu} \times {\Bbb Z}^{\nu} \to {\Bbb C} \; : \;
\| R \|_{w} \equiv \sup_{n \in {\Bbb Z}^{\nu}}
\sum_{m \in {\Bbb Z}^{\nu}} 
w(m-n) |R(m, n)| < \infty \right\}.
\end{eqnarray}
It is not difficult to see that $({\cal L}_{w}, \| \cdot \|_{w})$ is a 
Banach space. Moreover, one can introduce matrix multiplication
in ${\cal L}_{w}$ by setting
\begin{eqnarray}
\label{fa1.70}
(R S)(m, n) := \sum_{p \in {\Bbb Z}^{\nu}} R(m, p) S(p, n).
\end{eqnarray}
Indeed, for $R$, $S \in {\cal L}_{w}$, $w \in {\cal W}$ and 
$n \in {\Bbb Z}^{\nu}$ we conclude from (\ref{fa1.10})
\begin{eqnarray}
\nonumber
\sum_{m \in {\Bbb Z}^{\nu}} w(m-n) |(RS)(m, n)| 
&\leq& \sum_{m, p \in {\Bbb Z}^{\nu}}
w(m-p) |R(m,p)| w(p-n) |S(p, n)| \\
\nonumber
&\leq&
\sum_{p \in {\Bbb Z}^{\nu}} \|R\|_{w} w(p-n) |S(p,n)|  \\
\nonumber
&\leq&
\| R \|_{w} \| S \|_{w}.
\end{eqnarray}
Thus the matrix $RS$ lies in ${\cal L}_{w}$, satisfying
\begin{eqnarray}
\label{fa1.80}
\| R S \|_{w} \leq \| R \|_{w} \| S \|_{w}.
\end{eqnarray}
One readily verifies that ${\cal L}_{w}$ forms an algebra
(see, e.g. \cite{Lar} for a
definition) with unity $I$, $I(m, n) = \delta_{m, n}$. 
Furthermore, for $R \in {\cal L}_w$ and $u \in \ell_{1, w}$
we define $(R u)(m) := \sum_{n \in {\Bbb Z}^{\nu}} R(m, n) u(n)$ and obtain
\begin{eqnarray}
\nonumber
\| R u\|_{\ell_{1, w}} &=& 
\sum_{m} w(m) \left|
\sum_{n} R(m, n) u(n)
\right| \leq \sum_{m, n} w(m-n) |R(m, n)| w(n) |u(n)| \\
\label{fa1.90}
&\leq&
\| R \|_w \| u \|_{\ell_{1, w}}.
\end{eqnarray}
We summarize
\begin{Proposition}
\label{Pfa1.2}
Let $w \in {\cal W}$. Then ${\cal L}_w$ is a Banach algebra 
with respect to the
norm $\| \cdot \|_{w}$. The space ${\cal L}_w$ is contained in the space 
${\cal B}(\ell_{1, w})$ of bounded linear operators on 
$\ell_{1, w}$ and the corresponding operator norm in 
${\cal B}(\ell_{1, w})$ is bounded above by 
$\| \cdot \|_w$.
\end{Proposition}
 
Finally, we introduce two particular families of weight 
functions in ${\cal W}$.
\begin{Definition}
\label{DFa1}
Let $\nu$, $c$ and $D_N$ be given by (\ref{smr2.20}), (\ref{npc.52}) and
(\ref{npc.54}). For
$\sigma \geq 1/4$ and $n \in {\Bbb Z}^{\nu}$ we set
\begin{eqnarray}
\nonumber
w_{\sigma, c}(n) &:=&
D_N e^{\sigma |n|^c},\\
\nonumber
w_{\sigma, 1}(n) &:=&
D_N (1 + |n|)^{\nu+1} e^{\sigma |n|}.
\end{eqnarray}
\end{Definition}

The following properties of the just defined 
weight functions will be used throughout our analysis.

\begin{Lemma}
\label{Lwf.1}
For $\sigma \geq 1/4$ and $x \in \{c, 1 \}$ the following holds.
\begin{eqnarray}
\label{wf.110}
w_{\sigma, x} &\in& {\cal W}
\\
\label{wf.120}
\frac{1}{w_{\sigma, x}(n)} &\geq&
\sum_{m \in {\Bbb Z}^{\nu}}
\frac{1}{w_{\sigma, x}(m) 
w_{\sigma, x}(n-m)} \;\; \mbox{ for all }
n \in {\Bbb Z}^{\nu}
\\
\label{wf.130}
w_{\sigma, x}(n) &=& w_{\sigma + \mu, x}(n) 
e^{-\mu |n|^x}
\;\; \mbox{ for all }
n \in {\Bbb Z}^{\nu}, \; \mu > 0 \\
\label{wf.135}
\| u \|_{\sigma_1, x_1} &\leq& \| u \|_{\sigma_2, x_2} \;\; 
\mbox{ for all } \sigma_1 \leq \sigma_2, \; x_1 \leq x_2, \;
u \in X_{\sigma_2, x_2}.
\end{eqnarray}
\end{Lemma}

\begin{proof}
Properties (\ref{wf.120}), (\ref{wf.130}) and (\ref{wf.135})
follow immediately from the
definition of $w_{\sigma, x}$ and from the definition of $D_N$ 
(see (\ref{npc.54}), proposition \ref{Pwf.1}). It is also clear that 
$w_{\sigma, x}(n) \geq 1$ for all $n \in {\Bbb Z}^{\nu}$. Finally, the 
sub-multiplicative property $w(n) \leq w(m) w(n-m)$ which appears in the
definition of ${\cal W}$ is weaker than (\ref{wf.120}) and hence 
(\ref{wf.110}) follows.
\end{proof}

Property (\ref{wf.110}) shows that the sequence and matrix spaces corresponding
to the weights $w_{\sigma, x}$ satisfy propositions \ref{Pfa1.1} and 
\ref{Pfa1.2}. We introduce the following notation.
\begin{Definition}
\label{DFa2}
\begin{eqnarray}
\nonumber
X_{\sigma, c} 
&:=& \ell_{1, w_{\sigma, c}}, \quad X_{\sigma, 1} 
:= \ell_{1, w_{\sigma, 1}}, \\
\nonumber
{\cal L}_{\sigma, c} 
&:=& {\cal L}_{w_{\sigma, c}}, \quad {\cal L}_{\sigma, 1} 
:= {\cal L}_{w_{\sigma, 1}}.
\end{eqnarray}
By a slight abuse of notation the corresponding norms of the 
sequence spaces $X_{\sigma, c}$, $X_{\sigma, 1}$ and of the matrix spaces 
${\cal L}_{\sigma, c}$, ${\cal L}_{\sigma, 1}$ are both denoted by
$\| \cdot \|_{\sigma, c}$, $\| \cdot \|_{\sigma, 1}$.
\end{Definition}

\subsection{The equation for the Fourier coefficients}
\label{feq}

In this subsection we derive the equation for the Fourier coefficients.

\begin{Proposition}
\label{Pfeq}
Let $F$, $b$, $\gamma$ be given as in theorem \ref{Tsmr3.1}, satisfying
assumptions A1 -- A4. Recall the definitions of $\nu$, $g$, $c$, $r_{F, b}$,
$\alpha_k$ (see section \ref{npc} B) -- D) ).
Suppose that for some
$s_0 > 0$, $\omega \in {\Bbb R}^{\nu}$, and 
$u \in X_{1, c}$ the 
following five conditions hold:
\begin{eqnarray}
\label{feq10}
\| u \|_{1, c} &<& r_{F, b},
\\
\label{feq15}
u(0)&=&0,
\\
\label{feq17}
u(-m)&=&\overline{u(m)} \;\; \mbox{ for all } m \in {\Bbb Z}^{\nu},
\\
\label{feq25}
\mbox{ dist }( <\omega, m>, 2 \pi {\Bbb Z} ) &\geq& s_0 
e^{-\frac{1}{2} |m|^c} \quad \mbox{ for all } 
m \in {\Bbb Z}^{\nu}  \setminus \{ 0 \},
\\
\label{feq20}
\left(
\alpha_1 - \frac{\gamma^2 <m, g>^2}{4 \sin^2 \frac{<\omega, m>}{2}}
\right) u(m) &=& - \sum_{k=2}^{\infty} \alpha_k u^{*k}(m)
\quad \mbox{ for all } 
m \in {\Bbb Z}^{\nu}  \setminus \{ 0 \}.
\end{eqnarray}
Set
\begin{eqnarray}
\label{feq30}
\tilde{u}(m) := 
\left\{
\begin{array}{ll}
\frac{u(m)}{-2i \sin \frac{<\omega, m>}{2}}& ,\mbox{ if } 
m \in {\Bbb Z}^{\nu}  \setminus \{ 0 \}, \\
0 &,\mbox{ if } m = 0,
\end{array}
\right.
\end{eqnarray}
and define
\begin{eqnarray}
\label{feq35}
\chi (\xi) &:=& \sum_{m \in {\Bbb Z}^{\nu}} 
\tilde{u}(m) e^{i <m, \xi>} \;\;
\mbox{ for } \xi \in {\Bbb R}^{\nu},\\
\label{feq40}
x_n(t) &:=& nb + \chi (n \omega - \gamma t g) \;\;
\mbox{ for } t \in {\Bbb R}, n \in {\Bbb Z}.
\end{eqnarray}
Then for every $n \in {\Bbb Z}$ the function $x_n$ is real valued,
smooth and 
$\frac{2 \pi}{\gamma}$-periodic. The chain $(x_n)_{n \in {\Bbb Z}}$ 
satisfies
\begin{eqnarray}
\label{feq45}
\ddot{x}_n (t) = 
F(x_{n-1}(t) - x_n(t)) - F(x_n(t)-x_{n+1}(t)) \;\; 
\mbox{ for all } n \in {\Bbb Z}, t \in {\Bbb R}.
\end{eqnarray}
\end{Proposition}

\begin{proof}
From condition (\ref{feq25}) we deduce \begin{eqnarray}
\label{feq50}
\left| \sin \frac{<\omega, m>}{2} \right| \geq
\frac{1}{\pi} s_0 e^{-\frac{1}{2} |m|^c} \;\;
\mbox{ for } m \in {\Bbb Z}^{\nu} \setminus
\{ 0 \}.
\end{eqnarray}
Using (\ref{fa1.30}), definitions \ref{DFa1}, \ref{DFa2} and
$D_N \geq 1$ (see (\ref{npc.54})) it follows that
\begin{eqnarray}
\label{feq55}
| \tilde{u} (m)| 
\leq \frac{\pi \| u \|_{1, c}}{ 2 s_0 w_{1, c}(m)} 
e^{\frac{1}{2} |m|^c} 
\leq \frac{\pi \| u \|_{1, c}}{ 2 s_0} e^{-\frac{1}{2} |m|^c} 
\;\; \mbox{ for } m \in {\Bbb Z}^{\nu}.
\end{eqnarray}
This estimate together with condition (\ref{feq17}) implies that
$\chi$ as defined in (\ref{feq35}) is a real valued $C^{\infty}$-function
on the torus ${\Bbb T}^{\nu}$.
The time periodicity
of the functions $x_n$ is obvious 
since all components of the vector $g$ are integers. 

It only remains to verify (\ref{feq45}).
From (\ref{feq40}), (\ref{feq35}), (\ref{feq30}), and (\ref{feq15}) 
we conclude
\begin{eqnarray}
\nonumber
x_{n-1}(t) - x_n(t) &=& - b + \sum_{m \in {\Bbb Z}^{\nu}}
\tilde{u}(m) (e^{-i<m, \omega>} - 1) e^{i<m, \omega n - g \gamma t>} \\
\nonumber
&=& -b + \sum_{m \in {\Bbb Z}^{\nu}}
u(m) e^{-\frac{i}{2}<m, \omega>}
e^{i<m, \omega n - g \gamma t>}.
\end{eqnarray}
Using $ w_{1, c}(m) \geq 1$ for all $m \in {\Bbb Z}^{\nu}$
and (\ref{feq10}) we see that the Fourier series of
$x_{n-1} - x_n + b$ converges absolutely with
\begin{eqnarray}
\label{feq60}
|x_{n-1}(t) - x_n(t) + b| < \sum_{m \in {\Bbb Z}^{\nu}} |u(m)| 
\leq \| u \|_{1, c}  < r_{F, b}.
\end{eqnarray}
Thus the following reorderings of summation are justified (see also 
proposition \ref{Pfa1.1}).
\begin{eqnarray}
\nonumber
F(x_{n-1}(t)-x_n(t)) &=&
\sum_{k=0}^{\infty} \alpha_k
\left(
x_{n-1}(t)-x_n(t) +b
\right)^k
\\
\nonumber
&=&
\sum_{k=0}^{\infty} \alpha_k
\left(
\sum_{m \in {\Bbb Z}^{\nu}} u(m) 
e^{-\frac{i}{2}<m, \omega>}
e^{i<m, \omega n - g \gamma t>}
\right)^k
\\
\nonumber
&=& \alpha_0 +
\sum_{k=1}^{\infty} \alpha_k
\left(
\sum_{m \in {\Bbb Z}^{\nu}} u^{*k}(m) 
e^{-\frac{i}{2}<m, \omega>}
e^{i<m, \omega n - g \gamma t>}
\right)
\\
\nonumber
&=& \alpha_0 +
\sum_{m \in {\Bbb Z}^{\nu}}
\left(
\sum_{k=1}^{\infty} \alpha_k
u^{*k}(m) 
\right)
e^{-\frac{i}{2}<m, \omega>}
e^{i<m, \omega n - g \gamma t>}.
\end{eqnarray}
Consequently,
\begin{eqnarray}
\nonumber
F(x_{n}(t)-x_{n+1}(t)) =
\alpha_0 +
\sum_{m \in {\Bbb Z}^{\nu}}
\left(
\sum_{k=1}^{\infty} \alpha_k
u^{*k}(m) 
\right)
e^{\frac{i}{2}<m, \omega>}
e^{i<m, \omega n - g \gamma t>}.
\end{eqnarray}
Thus
\begin{eqnarray}
\label{feq65}
&& F(x_{n-1}(t)-x_n(t)) - F(x_{n}(t)-x_{n+1}(t)) =
\\
\nonumber
&&\quad \quad
\sum_{m \in {\Bbb Z}^{\nu}}
\left(
\sum_{k=1}^{\infty} \alpha_k
u^{*k}(m) 
\right)
\left( - 2i \sin \frac{<m, \omega>}{2} \right)
e^{i<m, \omega n - g \gamma t>}.
\end{eqnarray}
On the other hand, a straight forward calculation shows that
\begin{eqnarray}
\label{feq70}
\ddot{x}_n(t) =
\sum_{m \in {\Bbb Z}^{\nu}}
\tilde{u}(m)
\left( - <m, g>^2 \gamma^2 \right) e^{i<m, \omega n - g \gamma t>}.
\end{eqnarray}
It suffices to verify that the Fourier coefficients in (\ref{feq65})
and (\ref{feq70}) coincide, i.e. to show
\begin{eqnarray}
\label{feq75}
\tilde{u}(m)
\left( - <m, g>^2 \gamma^2 \right) =
\left(
\sum_{k=1}^{\infty} \alpha_k
u^{*k}(m) 
\right)
\left( - 2i \sin \frac{<m, \omega>}{2} \right)
\;\; \mbox{ for all } m \in {\Bbb Z}^{\nu}.
\end{eqnarray}
For $m=0$ both sides of equation (\ref{feq75}) clearly vanish. For
$m \in {\Bbb Z}^{\nu} \setminus \{ 0 \}$ equation (\ref{feq75})
follows from (\ref{feq20}) and (\ref{feq30}). 
\end{proof}

\section{The first induction step}
\label{fis}

The $P$ -- equation is solved by an inductive construction. In this section
we formulate and prove the initial induction statement $({\cal IS})_{j=1}$.

\noindent
{\bf Estimates of lemma \ref{hyp} used in section \ref{fis}:}
\begin{eqnarray}
\label{fis5}
12 D_N \nu 2^{\nu + 1} e^2 \rho_0 < r_{F, b} \\
\label{fis15}
50 \nu \rho_1 N_1 \leq \delta_1\\
\label{fis20}
\frac{s}{(2 N_1 + 4)^{\tau}} \geq 3 \delta_1 \\
\label{fis25}
\delta_1 < 2 \delta_V  
\end{eqnarray}

Induction statement $({\cal IS})_{j=1}$ concerns the following quantities:
$v_1: U(B_{\rho_0} (\lambda^{(0)}), 2 \rho_0) \to X_{2, 1}$, 
${\cal N}^{(1)} \subset {\Bbb R}^{2 \nu}$, $K^{(1)} \subset {\Bbb N}$,
$\lambda_k^{(j)} \in {\Bbb R}^{2 \nu}$ and $I_{C, k}^{(1)} \subset 
(- \pi, \pi)$ for $C \in {\cal C}^{(1)}$, $k \in K^{(1)}$.

\noindent
{\bf Statement $({\cal IS})_{j=1}$:}

\noindent
$\mbox{{\bf (1)}}_{j=1}$
\begin{itemize}
\item[(a)]
$v_1 (\lambda) = 0$ for all $\lambda \in U(B_{\rho_0} (\lambda^{(0)}), 2 
\rho_0)$.
\item[(b)] The map $U(B_{\rho_0} (\lambda^{(0)}), 2  \rho_0) \to X_{2,1}$,
$\lambda \mapsto {\cal P}(v_1(\lambda), \lambda)$ is analytic and
\begin{eqnarray}
\nonumber
\| {\cal P}(v_1(\lambda), \lambda) \|_{2,1} \leq D_W 
(2^{\nu +2} \nu D_N e^2)^2 
|a|^2 \;\; \mbox{ for all } 
\lambda \in U(B_{\rho_0} (\lambda^{(0)}), 2  \rho_0).
\end{eqnarray}
\end{itemize}

\noindent
$\mbox{{\bf (2)}}_{j=1}$
${\cal N}^{(1)} = \bigcup_{k \in K^{(1)}} B_{\rho_1}(\lambda_k^{(1)})$.

\noindent
$\mbox{{\bf (3)}}_{j=1}$
\begin{itemize}
\item[(a)]
For each $C \in {\cal C}^{(1)}$ and $k \in K^{(1)}$ the set $I_{C, k}^{(1)}$
satisfies
\begin{eqnarray}
\nonumber
I_{C, k}^{(1)} \subset
\bigcup_{\mu \in \{0, \pm 2 \omega_i^{(0)} \}}
\{ [\mu] \} + (-\delta_1, \delta_1). 
\end{eqnarray}
Moreover, $I_{C, k}^{(1)}$
is a union of $\# (C \cap {\cal S}) + 1$ 
disjoint open intervals, such that $B(I_{C, k}^{(1)}, 2 \delta_1)$
is again a union of $\# (C \cap {\cal S}) + 1$ disjoint intervals
with $B(I_{C, k}^{(1)}, 2 \delta_1) \subset (- \pi, \pi)$.

\item[(b)]
For $\lambda = (a, \omega) \in U({\cal N}^{(1)}, \rho_1)$ and $m \in 
{\Bbb Z}^{\nu}$ satisfying $0 < |m| \leq 2 N_1$ we have 

\noindent
$\mbox{ dist }( <\omega, m>, 2 \pi {\Bbb Z} ) \geq s d_{\tau, c} 
e^{-\frac{1}{2} |m|^c}$.
\item[(c)]
Let $k \in K^{(1)}$, $\lambda=(a, \omega) \in 
U(B_{\rho_1}(\lambda_k^{(1)}), \rho_1)$, $C \in {\cal C}^{(1)}$,
$m \in {\Bbb Z}^{\nu}$ with $0 < |m| \leq 2 N_1$ and $<m, g> = 0$. 

\noindent
Then $<\omega, m> \in 
U({\Bbb R} \setminus \tilde{I}_{C, k}^{(1)}, \frac{\delta_1}{10})$,
where $\tilde{I}_{C, k}^{(1)} \equiv I_{C, k}^{(1)} + 2 \pi {\Bbb Z}$
(cf. remark \ref{Rova.4}).
 
\item[(d)]
Let $k \in K^{(1)}$, $\lambda=(a, \omega) \in 
U(B_{\rho_1}(\lambda_k^{(1)}), \rho_1)$, $\theta \in {\Bbb C}$,
$C_1, C_2 \in {\cal C}^{(1)}$, $m, n \in {\Bbb Z}^{\nu}$ with
$0 < |n-m| \leq 2 N_1$ and $<m-n, g> = 0$. 

\noindent
Then $\theta + <\omega, m> \in 
U(\tilde{I}_{C_1, k}^{(1)}, \frac{\delta_1}{10})$ implies
$\theta + <\omega, n> \in {\Bbb C} \setminus
U(\tilde{I}_{C_2, k}^{(1)}, \frac{\delta_1}{10})$.
\end{itemize}

\noindent
$\mbox{{\bf (4)}}_{j=1}$
Let $k \in K^{(1)}$, $C \in {\cal C}^{(1)}$, $m \in {\cal S} \cap C$, and
$(\theta, \lambda) \in 
U({\Bbb R} \setminus \tilde{I}_{C, k}^{(1)}, \frac{\delta_1}{10})
\times U(B_{\rho_1}(\lambda_k^{(1)}), \rho_1)$.

\noindent
Then $|V(\theta, \omega)(m)| > \frac{d_V \delta_1}{2}$.

\begin{Lemma}
\label{Lfis}
There exist $v_1$, ${\cal N}^{(1)}$, 
$K^{(1)}$, $\lambda_k^{(1)}$ (for $k \in K^{(1)}$),
$I^{(1)}_{C, k}$ (for $C \in {\cal C}^{(1)}$, $k \in K^{(1)}$) 
such that statement 
$({\cal IS})_{j=1}$ is satisfied.
\end{Lemma}
\begin{proof}
$\mbox{{\bf (1)}}_{j=1}:$
Statement (a) forces the definition 
$v_1(\lambda) := 0$ for all 
$\lambda \in U(B_{\rho_0} (\lambda^{(0)}), 2 \rho_0)$.
To show (b) observe that ${\cal P} (v_1 (\lambda), \lambda) =
P W (\varphi (a))$. 
Definitions (\ref{ova.120}), \ref{DFa1}, and \ref{DFa2} 
together with (\ref{fis5})
imply for $\lambda \in U(B_{\rho_0} (\lambda^{(0)}), 2 \rho_0)$
that
\begin{eqnarray}
\label{fis105}
\nonumber
\| \varphi (a) \|_{2,1} \leq D_N 2 \nu 2^{\nu+1} e^2 |a| 
< 6 D_N \nu 2^{\nu + 1} e^2 \rho_0 <
\frac{r_{F, b}}{2}.
\end{eqnarray}
Lemma \ref{Lenp.1} then implies analyticity of the function
${\cal P}(v_1(\cdot), \cdot)$ on $U(B_{\rho_0} (\lambda^{(0)}), 2 \rho_0)$.
The same lemma \ref{Lenp.1} also yields 
\begin{eqnarray}
\label{fis110}
\nonumber
\| {\cal P} (v_1 (\lambda), \lambda) \|_{2,1} = 
\| P W (\varphi (a)) \|_{2,1} \leq 
\| W (\varphi (a)) \|_{2,1} \leq
D_W (2^{\nu + 2} \nu D_N e^2)^2 |a|^2.
\end{eqnarray}

\noindent
$\mbox{{\bf (2)}}_{j=1}:$ 
We set
\begin{eqnarray}
\label{fis.111}
{\cal N}^{(1)} &:=& B_{\rho_1}(\lambda^{(0)}), \\
\label{fis.112}
K^{(1)} &:=& \{ 1 \}, \\
\label{fis.113}
\lambda_1^{(1)} &:=& \lambda^{(0)}.
\end{eqnarray}
The claim is obvious. 

\noindent
$\mbox{{\bf (3)}}_{j=1}:$
For $C \in {\cal C}^{(1)}$, $k \in K^{(1)}$ we define
\begin{eqnarray}
\label{fis115}
I_{C, k}^{(1)} : = (-\delta_{1}, \delta_1)
\cup
\bigcup_{e_i \in C} 
([-2 \omega_i^{(0)}] - \delta_{1} , [-2 \omega_i^{(0)}] + \delta_1) 
\cup
\bigcup_{-e_i \in C} 
([2 \omega_i^{(0)}] - \delta_{1} , [2 \omega_i^{(0)}] + \delta_1).
\end{eqnarray}
In order to prove statement (a) it suffices to show the following
two conditions.
\begin{eqnarray}
\nonumber
\mbox{ dist}(\mu_1 - \mu_2, 2 \pi {\Bbb Z}) &\geq& 6 \delta_1
\;\; \mbox{ for } \; \mu_1, \mu_2 \in \{ 0, \pm 2 \omega_i^{(0)} \}, \;
\mu_1 \neq \mu_2. \\
\nonumber
\mbox{ dist}(\mu, (2 \pi + 1) {\Bbb Z}) &\geq& 3 \delta_1
\;\; \mbox{ for } \; \mu \in \{ 0, \pm 2 \omega_i^{(0)} \}.
\end{eqnarray}
Both estimates follow from assumption A3 and from $s/(4^{\tau}) \geq 
6 \delta_1$ which in turn is a consequence of (\ref{fis20}).
Statement (b) follows from 
(\ref{fis.111}), 
assumption A3, (\ref{fis15}), (\ref{fis20}) and
(\ref{npc.120}).
Indeed, let $\lambda \in U({\cal N}^{(1)}, \rho_1)$, $0 < |m| \leq 2 N_1$, 
then
\begin{eqnarray}
\label{fis120}
\nonumber
\mbox{ dist }( <\omega, m>, 2 \pi {\Bbb Z} ) &\geq&
\mbox{ dist }( <\omega^{(0)}, m>, 2 \pi {\Bbb Z} ) - \nu (2 \rho_1) (2 N_1) \\
\nonumber
&\geq&
\frac{s}{|m|^{\tau}} - \frac{\delta_1}{10} \\
\nonumber
&\geq& 
\frac{s}{2 |m|^{\tau}} \\
\nonumber
&\geq& s d_{\tau, c} e^{-\frac{1}{2} |m|^c}.
\end{eqnarray}
We prove claim (c) by contradiction. Assume that there exist 
$\lambda \in U({\cal N}^{(1)}, \rho_1)$,
$C \in {\cal C}^{(1)}$ and $m \in {\Bbb Z}^{\nu}$ with 
$0 < |m| \leq 2 N_1$, $<m, g>=0$, such that 
\begin{eqnarray}
\nonumber
<\omega, m> \in
{\Bbb C} \setminus 
U({\Bbb R} \setminus \tilde{I}_{C, 1}^{(1)}, \frac{\delta_1}{10}).
\end{eqnarray}
Again 
(\ref{fis15}) implies $|<\omega, m> - <\omega^{(0)}, m>| < \delta_1/10$ and
therefore 
\begin{eqnarray}
\nonumber
<\omega^{(0)}, m> \in {\Bbb R} \setminus 
({\Bbb R} \setminus \tilde{I}_{C, 1}^{(1)}) = \tilde{I}_{C, 1}^{(1)}.
\end{eqnarray}
From the definition of $\tilde{I}_{C, 1}^{(1)} \equiv
I_{C, 1}^{(1)} + 2 \pi {\Bbb Z}$ and (\ref{fis115})
it follows that there exist $l \in {\Bbb Z}$
and $m_1 \in {\cal S} \cup \{ 0 \}$, such that 
\begin{eqnarray}
\label{fis125}
\nonumber
| <\omega^{(0)}, m> - 2 \pi l - 2 <\omega^{(0)}, m_1> | \; < \; \delta_1
\end{eqnarray}
This implies $\mbox{ dist }(<\omega^{(0)}, m - 2m_1>, 2 \pi {\Bbb Z}) < 
\delta_1$. By assumption A3 and (\ref{fis20}) it then follows that
$m = 2 m_1$. Since $m_1 \in  {\cal S} \cup \{ 0 \}$ the condition 
$0 = <m, g> = <2 m_1, g>$ can only be satisfied, if $m_1 = 0$ and
hence $m=0$. However, this contradicts the assumption $|m| > 0$. 

Claim (d) is also proved by contradiction. Assume there exist
$\lambda \in U({\cal N}^{(1)}, \rho_1)$, $\theta \in {\Bbb C}$,
$C_1, C_2 \in {\cal C}^{(1)}$, $m, n \in {\Bbb Z}^{\nu}$ with
$0 < |m-n| \leq 2 N_1$, $<m-n, g>=0$, and numbers
$\theta_1 \in I_{C_1, 1}^{(1)}$,
$\theta_2 \in I_{C_2, 1}^{(1)}$, 
such that
\begin{eqnarray}
\label{fis130}
\mbox{ dist}(\theta + <\omega, m> - \theta_1, 2 \pi {\Bbb Z})
&<& \frac{\delta_1}{10}, \\
\nonumber
\mbox{ dist}(\theta + <\omega, n> - \theta_2, 2 \pi {\Bbb Z})
&<& \frac{\delta_1}{10}.
\end{eqnarray}
By the definition of the sets $I_{C, k}^{(1)}$ there exist 
$m_1, n_1 \in {\cal S} \cup \{ 0 \}$, such that
\begin{eqnarray}
\nonumber
\mbox{ dist}(\theta_1 - 2<\omega^{(0)}, m_1>, 2 \pi {\Bbb Z}) 
< \delta_ 1 \quad \mbox{ and } \quad 
\mbox{ dist}(\theta_2 - 2<\omega^{(0)}, n_1>, 2 \pi {\Bbb Z}) < \delta_ 1.
\end{eqnarray}
Substituting
$\theta_1$ and $\theta_2$ in (\ref{fis130}) and taking the difference
of the inequalities of (\ref{fis130}) we obtain
\begin{eqnarray}
\label{fis135}
\mbox{ dist}(<\omega, m-n> - 2<\omega^{(0)}, m_1-n_1>, 2 \pi {\Bbb Z}) <
\frac{22}{10} \delta_1.
\end{eqnarray}
Using $|m-n| \leq 2 N_1$ and (\ref{fis15}) we conclude that the distance of
$<\omega^{(0)},m-n-2(m_1-n_1)>$ to some integer multiple of $2 \pi$ is 
less than $3 \delta_1$. Assumption A3 together with 
(\ref{fis20}) imply 
$m-n - 2(m_1 -n_1) = 0$. From $<m-n, g> =0$ we see that $<m_1, g>
=<n_1, g>$. Recall that $m_1, n_1 \in {\cal S} \cup \{ 0 \}$ and therefore
we must have $m_1 = n_1$. Thus $m-n=0$, contradicting the assumption
$|m-n| > 0$.

\noindent
$\mbox{{\bf (4)}}_{j=1}:$
Let $C \in {\cal C}^{(1)}$, $m \in {\cal S} \cap C$ and 
$(\theta, \lambda) \in 
U({\Bbb R} \setminus \tilde{I}_{C, 1}^{(1)}, \frac{\delta_1}{10})
\times U(B_{\rho_1}(\lambda_1^{(1)}), \rho_1)$.
We write $m = \mu e_i$ with $\mu \in \{ -1, 1 \}$ and $i \in \{1, \ldots,
\nu \}$. From definition (\ref{fis115}) it follows that
\begin{eqnarray}
\nonumber
\{ 0, - 2 \mu \omega_i^{(0)} \} + (- \delta_1, \delta_1) + 2 \pi {\Bbb Z}
\subset \tilde{I}^{(1)}_{C, 1}.
\end{eqnarray}
Thus
\begin{eqnarray}
\nonumber
\mbox{ dist }(\theta, 2 \pi {\Bbb Z}) &\geq& 9 \delta_1/10, \\
\nonumber
\mbox{ dist }(\theta, \{- 2 \mu \omega_i^{(0)} \} + 2 \pi {\Bbb Z}) 
&\geq& 9 \delta_1/10.
\end{eqnarray} 
By (\ref{fis15}) we have $|\omega_i - \omega_i^{(0)}| <
\delta_1/3$, implying
\begin{eqnarray}
\label{fis140}
\mbox{ dist }(\theta + \mu \omega_i, \{ \mu \omega_i^{(0)} \} + 2 \pi {\Bbb Z}) 
&>& \frac{1}{2} \delta_1, \\
\nonumber
\mbox{ dist }(\theta + \mu \omega_i, \{ - \mu \omega_i^{(0)} \} 
+ 2 \pi {\Bbb Z}) 
&>& \frac{1}{2} \delta_1.
\end{eqnarray}
Since $\theta + <\omega, m> = \theta + \mu \omega_i$ we conclude from
(\ref{lop.50}), proposition \ref{Plop.1} (a) and (\ref{fis25}) that
\begin{eqnarray}
\label{fis145}
\nonumber
|V(\theta, \omega)(m)| = |V_i (\theta + \mu \omega_i)| > 
\frac{d_V \delta_1}{2}.
\end{eqnarray}
This completes the proof of Lemma \ref{Lfis}.
\end{proof}

\section{The induction statement for $j \geq 2$}
\label{is}

In this section we formulate the induction statements for $j \geq 2$.
They concern the following quantities $v_j: U({\cal N}^{(1)}, \rho_1)
\to X_{1, c}$, ${\cal N}^{(j)} \subset {\cal N}^{(1)}$,
index sets $K^{(j)}$ and $K^{(j-0.5)}$, vectors $\lambda_k^{(j)} \in
{\cal N}^{(1)}$ ($k \in K^{(j)}$), sets $I_{C, k}^{(j)}$ ($k \in K^{(j)}$,
$C \in {\cal C}^{(j)}$) and sets of polynomials 
${\cal POL}^{(j)}_{\tilde{k}} \subset {\cal POL}$ ($\tilde{k} \in
K^{(j-0.5)}$). The proof that such functions exist will stretch
over sections \ref{cv2} -- \ref{comp}.

\noindent
{\bf Statement $({\cal IS})_{j \geq 2}$:}

\noindent
$\mbox{{\bf (1)}}_{j \geq 2}$
\begin{itemize}
\item[(a)]
$v_j: U({\cal N}^{(1)}, \rho_1) \to X_{1, c}$ is a $C^{\infty}$ -- function 
(if $U({\cal N}^{(1)}, \rho_1)$ is considered 
as a subset of ${\Bbb R}^{4 \nu}$).
\item[(b)]
The restriction of the function $v_j$ to the set 
$U({\cal N}^{(j-1)}, \rho_{j-1}/4)$ is analytic.
\item[(c)] The map $\lambda \mapsto {\cal P}(v_j(\lambda), \lambda)$ has an
analytic continuation to $U({\cal N}^{(j-1)}, \rho_{j-1}/4)$ and
\begin{eqnarray}
\nonumber
\| {\cal P}(v_j(\lambda), \lambda) \|_{1/4,c} \leq 
e^{-N_{j-1}^c} |a|^2 \;\;\; \mbox{ for all } 
\lambda \in U({\cal N}^{(j-1)}, \rho_{j-1}/4).
\end{eqnarray}
\item[(d)] For all $\lambda \in U({\cal N}^{(1)}, \rho_1)$:
\begin{eqnarray}
\nonumber
\| v_j(\lambda) - v_{j-1}(\lambda) \|_{1, c} \leq \left\{
\begin{array}{ll}
4 e^6 \nu^2 D_W D_N^2 D_{\tau} N_0^{\tau + 1} |a|^2 &\mbox{ if } j=2, \\
e^{-\frac{1}{8} N_{j-2}^c} |a|^2 &\mbox{ if } j > 2. 
\end{array}
\right.
\end{eqnarray}
\item[(e)] For all $\lambda \in U({\cal N}^{(1)}, \rho_1)$:
\begin{eqnarray}
\nonumber
\|\partial^{\beta}(v_j - v_{j-1})(\lambda) \|_{1, c} \leq \left\{
\begin{array}{ll}
\beta ! 36 e^6 \nu^2 D_W D_N^2 D_{E}(|\beta|_1) N_0^{E_{\rho}(|\beta|_1 - 1)} 
&\mbox{ if } j=2,\\
\beta ! e^{2 \nu} D_{\psi}(|\beta|_1)
\left( \frac{8}{\rho_{j-1}}\right)^{|\beta|_1}
e^{-\frac{1}{8} N_{j-2}^c} &\mbox{ if } j > 2.
\end{array}
\right.
\end{eqnarray}
\item[(f)] For all $\lambda \in U({\cal N}^{(1)}, \rho_1)$
the support of $v_j(\lambda)$ is contained in $B_{N_{j-1}}(0) \setminus
({\cal S} \cup \{ 0 \})$.
\item[(g)]
For all $\lambda \in B({\cal N}^{(1)}, \rho_1)$ and $m \in {\Bbb Z}^{\nu}$
we have $v_j(\lambda)(m) = v_j(\lambda)(-m) \in {\Bbb R}$.
\item[(h)]
Let $\lambda = (a, \omega) \in U({\cal N}^{(1)}, \rho_1)$ with 
$a_i = 0$, for some $1 \leq i \leq \nu$. Then $v_j(\lambda)(m) = 0$, 
if $m_i \neq 0$.
\end{itemize}

\noindent
$\mbox{{\bf (2)}}_{j \geq 2}$
\begin{itemize}
\item[(a)]
${\cal N}^{(j)} = \bigcup_{k \in K^{(j)}} B_{\rho_j}(\lambda_k^{(j)})
\subset {\cal N}^{(j-1)}$.
\item[(b)]
For $\tilde{k} \in K^{(j - 0.5)}$ the following inclusion holds:
\begin{eqnarray}
\nonumber
&&\left( \overline{{\cal N}^{(j-1)}} \setminus 
\overline{{\cal N}^{(j)}} \right)
\cap \overline{B_{\tilde{\rho}_{j-1}} 
\left( \lambda_{\tilde{k}}^{(j-0.5)} \right)} \\
\nonumber
&\subset& 
\bigcup_{\scriptsize
{\scriptstyle
\begin{array}{c}
(p, \vartheta) \in {\cal POL}_{\tilde{k}}^{(j)} 
\\ 2 N_{j-1} < |m| \leq 2 N_j
\end{array}
}
} \!\!\!\!\!\!\!
\vspace{.3cm}
\begin{array}{r}
\{ \lambda \in \overline{B_{\tilde{\rho}_{j-1}} 
( \lambda_{\tilde{k}}^{(j-0.5)})}:
|[<\omega, m>] - \vartheta| \leq \frac{\delta_1}{2} \; \mbox{ and } \\
|p([<\omega, m>] - \vartheta, \lambda)| < (D_K+1) \delta_j
\}.
\end{array}
\end{eqnarray}
\item[(c)] ${\cal POL}_{\tilde{k}}^{(j)} \subset {\cal POL}$ 
and $\# {\cal POL}_{\tilde{k}}^{(j)} \leq 
2 (2 N_{j-1})^{4 \nu} \delta_{j-1}^{-2}$.
\end{itemize}

\noindent
$\mbox{{\bf (3)}}_{j \geq 2}$
Let $C \in {\cal C}^{(j)}$ and $k \in K^{(j)}$. Denote
$C' := C \cap B_{1.5}(0)$, if $j=2$, respectively
$C' := C \cap B_{N_{j-2}}(0)$, if $j>2$ and
$k' := \pi_{j-1}^{(j)} k$. 
\begin{itemize}
\item[(a)]
$I_{C, k}^{(j)} \subset I_{C', k'}^{(j-1)}$ is an open set.
\item[(b)]
For $\lambda = (a, \omega) \in U({\cal N}^{(j)}, \rho_j)$ and $m \in 
{\Bbb Z}^{\nu}$ with $2 N_{j-1} < |m| \leq 2 N_j$:
\begin{eqnarray}
\nonumber
\mbox{ dist }( <\omega, m>, 2 \pi {\Bbb Z} ) \geq s d_{\tau, c} 
e^{-\frac{1}{2} |m|^c}.
\end{eqnarray}
\item[(c)]
For $k \in K^{(j)}$, $\lambda=(a, \omega) \in 
U(B_{\rho_j}(\lambda_k^{(j)}), \rho_j)$, $C \in {\cal C}^{(j)}$,
$m \in {\Bbb Z}^{\nu}$ with $2 N_{j-1} < |m| \leq 2 N_j$ 
and $<m, g> = 0$: 
\begin{eqnarray}
\nonumber
<\omega, m> \in 
U({\Bbb R} \setminus \tilde{I}_{C, k}^{(j)}, \frac{\delta_j}{10}),
\end{eqnarray} 
where $\tilde{I}_{C, k}^{(j)} \equiv I_{C, k}^{(j)}$ (see remark \ref{Rova.4}).
\item[(d)]
For $k \in K^{(j)}$, $\lambda=(a, \omega) \in 
U(B_{\rho_j}(\lambda_k^{(j)}), \rho_j)$, $\theta \in {\Bbb C}$,
$C_1, C_2 \in {\cal C}^{(j)}$, $m, n \in {\Bbb Z}^{\nu}$ with
$2 N_{j-1} < |n-m| \leq 2 N_j$ and $<m-n, g> = 0$: 
\begin{eqnarray}
\nonumber
\theta + <\omega, m> \in 
U(\tilde{I}_{C_1, k}^{(j)}, \frac{\delta_j}{10}) \quad \Longrightarrow
\quad
\theta + <\omega, n> \in {\Bbb C} \setminus
U(\tilde{I}_{C_2, k}^{(j)}, \frac{\delta_j}{10}).
\end{eqnarray}
\end{itemize}

\noindent
$\mbox{{\bf (4)}}_{j \geq 2}$
Let $k \in K^{(j)}$, $C \in {\cal C}^{(j)}$ 
and define the corresponding $C'$, $k'$ as in $\mbox{{\bf (3)}}_{j}$.
There exists a map $(\theta, \lambda) \mapsto G_C^{(j)}(\theta, \lambda)$
which is analytic on  
$U(\tilde{I}_{C', k'}^{(j-1)} \setminus \tilde{I}_{C, k}^{(j)}, \delta_j/10)
\times U(B_{\rho_j}(\lambda_k^{(j)}),\rho_j)$, satisfying
\begin{eqnarray}
\label{is.10}
\| G_C^{(j)} (\theta, \lambda) \|_{\sigma_j, c} \leq M_j
\;\;\; \mbox{ for }
(\theta, \lambda) \in 
U(\tilde{I}_{C', k'}^{(j-1)} \setminus \tilde{I}_{C, k}^{(j)}, \delta_j/10)
\times U(B_{\rho_j}(\lambda_k^{(j)}),\rho_j).
\end{eqnarray}
Furthermore, $G_C^{(j)}(\theta, \lambda)$ 
is the inverse of the matrix $T_C^{(j)}(\theta, \lambda)$
(the restriction of the matrix $T^{(j)}(\theta, \lambda)$
(see (\ref{npc.64})) to the set $C$) for all $(\theta, \lambda) \in
U(\tilde{I}_{C', k'}^{(j-1)} \setminus \tilde{I}_{C, k}^{(j)}, \delta_j/10)
\times U(B_{\rho_j}(\lambda_k^{(j)}),\rho_j) \setminus \tilde{Z}_{C}$.
The entries of $G_C^{(j)} (\theta, \lambda)$ are real for real values
of $(\theta, \lambda)$.

\section{The construction of $v_{j+1}$ for $j=1$}
\label{cv2}
In this section we construct $v_2$ from $v_1$ by a modified
Newton scheme. The induction hypothesis $({\cal IS})(3c)_{j=1}$ and
$({\cal IS})(4)_{j=1}$ will imply that the diagonal parts of the
linearized operators dominate the corresponding off-diagonal parts.
The linearized operators can be inverted by a Neumann series. 
Observe, that the dependence of the estimates $({\cal IS})(1e)_{j=2}$
on the order of the derivative $|\beta|_1$ is better than one could derive
from the analyticity of $v_2$ ( see $({\cal IS})(1b)_{j=2}$) and from statement
$({\cal IS})(1d)_{j=2}$ via Cauchy's formula. The improvement we have gained
for estimates $({\cal IS})(1e)_{j=2}$ will be crucial for obtaining 
lower bounds on the measure of the set of non-resonant parameters in section
\ref{est}.

\noindent
{\bf Estimates of lemma \ref{hyp} used in section \ref{cv2}:}

\begin{eqnarray}
\label{cv2.5}
2^{2 \nu + 4} e^6 \nu^2 D_W D_N^2  \leq e^{\frac{1}{2} N_1^c} \\
\label{cv2.10}
c \leq \frac{A-1}{A} \frac{\log 1.5}{\log 5}\\
\label{cv2.15}
E_{\rho} \geq \tau+2 \\
\label{cv2.20}
6 \nu N_0^{\tau + 1 - E_{\rho}} \leq s \\
\label{cv2.25}
\frac{s}{2 N_0^{\tau}} < \delta_V \\
\label{cv2.30}
\alpha_1 > \frac{d_V s}{2 N_0^{\tau}} \\
\label{cv2.35}
\frac{(\nu+1)^2 \gamma^2}{4 \cosh^2(\nu N_0^{1-E_{\rho}})} - \alpha_1 >
\frac{d_V s}{2 N_0^{\tau}} \\
\label{cv2.40}
12 e^6 \nu^2 D_W D_N^2 D_{\tau} N_0^{\tau + 1} \rho_0 < 1 \\
\label{cv2.45}
(2^{\nu + 2} \nu e^2 D_N+1) 3 \rho_0 < \frac{r_{F, b}}{2} \\
\label{cv2.50}
3 (2^{\nu+2} \nu e^2 D_N + 1) D_N \frac{2 N_0^{\tau - E_{\rho}}}{d_V s} < \frac{1}{2}
\\
\label{cv2.55}
2^{2 \nu + 4}\frac{D_N^2}{d_V s} \leq N_0 \\
\label{cv2.60}
3(2^{\nu + 2} \nu e^2 D_N + 1) \leq N_0 \\
\label{cv2.65}
12 e^3 \nu D_W D_N \leq N_0
\end{eqnarray}

\begin{Lemma}
\label{Lcv2}
Assume that induction statement $({\cal IS})_{j=1}$ holds. Then
there exists a function $v_2: 
U({\cal N}^{(1)}, \rho_1) \to X_{1,c}$ such that induction
statement $({\cal IS}) (\mbox{1})_{j=2}$ is satisfied.
\end{Lemma}

\begin{proof}
By (\ref{npc.400}), (\ref{npc.410}), (\ref{npc.72}) 
there exists a $p_0 \in {\Bbb N}$ such that
$N_1 = N_0 5^{p_0}$. For $1 \leq p \leq p_0 + 1$ define
\begin{eqnarray}
\label{cv2.105}
L_p &:=& N_0 5^{p-1}, \\
\label{cv2.110}
B_p &:=& B_{L_p}(0) \setminus ({\cal S} \cup \{ 0 \}) 
\subset {\Bbb Z}^{\nu}, \\
\label{cv2.115}
\mu_p &:=& 1 + \frac{1}{p},
\\
\label{cv2.120}
r_p &:=& L_p^{-E_{\rho}}.
\end{eqnarray}
For $1 \leq p \leq p_0 + 1$
we construct inductively analytic maps 
$z_p: U(B_{r_p} (\lambda^{(0)}), 2r_p) \to X_{\mu_p, 1}$ satisfying
\begin{itemize}
\item[(i)] For all $\lambda \in U(B_{r_p} (\lambda^{(0)}), 2r_p)$:
\begin{eqnarray}
\nonumber
z_p (\lambda) = 0&& \mbox{ if } p=1,
\\
\nonumber
\|(z_{p} - z_{p-1})(\lambda) \|_{\mu_{p-1}, 1} \leq 
4 e^6 \nu^2 D_W D_N^2 L_p^{\tau + 1} e^{-1.5^{p-1}} |a|^2
&&\mbox{ if }  
2 \leq p \leq p_0+1.
\end{eqnarray}
\item[(ii)] For all $\lambda \in U(B_{r_p} (\lambda^{(0)}), 2r_p)$
the support of $z_p(\lambda)$ is contained in $B_p$.
\item[(iii)]
For all $\lambda \in B(B_{r_p} (\lambda^{(0)}), 2r_p)$ and 
$m \in {\Bbb Z}^{\nu}$
we have $z_p(\lambda)(m) = z_p(\lambda)(-m) \in {\Bbb R}$.
\item[(iv)]
Let $\lambda = (a, \omega) \in U(B_{r_p} (\lambda^{(0)}), 2r_p)$ with 
$a_i = 0$, for some $1 \leq i \leq \nu$. Then $z_p(\lambda)(m) = 0$, 
if $m_i \neq 0$.
\item[(v)] The map $\lambda \mapsto {\cal P}(z_p(\lambda), \lambda)$ has an
analytic continuation to $U(B_{r_p} (\lambda^{(0)}), 2r_p)$ and
\begin{eqnarray}
\nonumber
\| {\cal P}(z_p(\lambda), \lambda) \|_{\mu_p,1} \leq 
2^{2 \nu + 4} e^6 \nu^2 D_W D_N^2 e^{-1.5^p} |a|^2 \;\; \mbox{ for all } 
\lambda \in U(B_{r_p} (\lambda^{(0)}), 2r_p).
\end{eqnarray}
\end{itemize}
Suppose that for $1 \leq p \leq p_0+1$ there exist analytic 
functions $z_p$, satisfying (i)-(v). We verify that 
\begin{eqnarray}
\label{cv2.125}
v_2 := z_{p_0+1}|_{U({\cal N}^{(1)}, \rho_1)}
\end{eqnarray}
satisfies 
$({\cal IS}) (\mbox{1})_{j=2}$.
Observe that (\ref{npc.420}), (\ref{npc.75}) together with the definitions of
$p_0$, $L_p$ and $r_p$ given at the beginning of this proof 
imply $L_{p_0+1} = N_1$ and 
$r_{p_0+1} = \rho_1$. 
Since $\mu_{p_0+1} \geq 1$, $c < 1$ (see (\ref{npc.52}), (\ref{wf.135})) and 
${\cal N}^{(1)} =
B_{\rho_1}(\lambda ^{(0)})$ (see (\ref{fis.111})) 
it is clear that $v_2$ is a well defined map from $U({\cal N}^{(1)}, \rho_1)$
into $X_{1, c}$ 
satisfying properties (a) and (b) of
$({\cal IS}) (\mbox{1})_{j=2}$. Furthermore conditions (f), (g), (h)
follow from properties (ii), (iii), (iv) of
$z_{p_0+1}$. It remains to verify that $v_2$ satisfies
statements (c), (d) and (e) of $({\cal IS}) (\mbox{1})_{j=2}$:

{\em Proof of(c):} Using (\ref{cv2.5}), property (v) of $z_{p_0+1}$ and
(\ref{wf.135}) we only
need to show that 
\begin{eqnarray}
\label{cv2.130}
1.5^{p_0} \geq N_1^c. 
\end{eqnarray}
Recall that $p_0$ was chosen
to satisfy $N_0 5^{p_0} = N_1$. By (\ref{npc.410}) we have 
$N_1 = N_0^A$ and therefore $p_0 \log 5 = (A-1) \log N_0$.
Consequently, (\ref{cv2.130}) can be written in the form 
\begin{eqnarray}
\label{cv2.135}
(A-1) \frac{\log 1.5}{\log 5} \geq Ac.
\end{eqnarray} 
This in turn follows from (\ref{cv2.10}).

{\em Proof of (d):} Since $v_1(\lambda) \equiv 0$ and 
$z_1(\lambda) \equiv 0$
we conclude from (i) and (\ref{cv2.105}) that 
\begin{eqnarray}
\label{cv2.140}
\|(v_2 -v_1)(\lambda)\|_{1, c} \leq \|(z_{p_0+1} - z_1)(\lambda)\|_{1,1} \leq
4 e^6 \nu^2 D_W D_N^2 N_0^{\tau+1} |a|^2 \sum_{p=1}^{p_0}
5^{p(\tau+1)} e^{-1.5^p} 
\end{eqnarray}
for all $\lambda \in U({\cal N}^{(1)}, \rho_1)$.
Claim (d) follows from the definition of the constant $D_{\tau}$ in 
(\ref{npc.110}). 

{\em Proof of (e):} Note that the estimate on $z_p-z_{p-1}$ 
in (i) holds on a complex
$2 r_p$-neighborhood of $B_{r_p}(\lambda^{(0)})$. Furthermore, we
have $r_p \geq \rho_1$ for all $1 \leq p \leq p_0+1$ and therefore
Cauchy's integral formula implies for $\lambda \in U({\cal N}^{(1)}, \rho_1)$,
$1 \leq p \leq p_0$
\begin{eqnarray}
\label{cv2.145}
\nonumber
\| \partial^{\beta} (z_{p+1} -z_{p})(\lambda)\|_{1,1}
&\leq&
\beta ! r_{p+1}^{-|\beta|_1} 4 e^6 \nu^2 D_W D_N^2 L_{p+1}^{\tau+1} 
e^{-1.5^p} (3 r_{p+1})^2 \\
&\leq&
\beta ! 36 e^6 \nu^2 D_W D_N^2 L_{p+1}^{\tau + 1} 
r_{p+1}^{2-|\beta|_1} e^{-1.5^p}.
\end{eqnarray}
From (\ref{cv2.15}) and (\ref{cv2.120}) it follows that
$L_{p+1}^{\tau + 1} < r_{p+1}^{-1}$. A calculation similar to (\ref{cv2.140})
yields for  $\lambda \in U({\cal N}^{(1)}, \rho_1)$
\begin{eqnarray}
\label{cv2.150}
\| \partial^{\beta} (v_2 - v_1)(\lambda) \|_{1,c} \leq
\beta ! 36 e^6 \nu^2 D_W D_N^2 N_0^{E_{\rho}(|\beta|_1 -1)} \sum_{p=1}^{p_0}
5^{p E_{\rho}(|\beta|_1 -1)} e^{-1.5^p}.
\end{eqnarray}
Claim (e) follows from the definition of $D_E(x)$ in 
(\ref{npc.118}).

We have therefore reduced the proof of lemma \ref{Lcv2} to showing that
there exists a finite sequence $(z_p)_{1 \leq p \leq p_0+1}$ of
analytic maps $z_p: U(B_{r_p}(\lambda^{(0)}), 2r_p) \to X_{\mu_p, 1}$
satisfying properties (i)-(v).

\noindent 
{\em Inductive construction of $z_p$:}
 
\noindent
$\underline{p=1:}$ By (i) the choice of $z_1 (\lambda) \equiv 0$ is determined
for all $\lambda \in U(B_{r_1}(\lambda^{(0)}), 2r_1)$. This choice clearly
satisfies also (ii), (iii) and (iv). Since $r_1 = N_0^{-E_{\rho}} = \rho_0$
we have $v_1=z_1$ and property (v) follows from $({\cal IS}) (1)_{j=1}$ (b).

\noindent
$\underline{p \to p+1 \mbox{ for } 1 \leq p \leq p_0:}$
Set
\begin{eqnarray}
\label{cv2.155}
T_p(\lambda) :=
D_v {\cal P} (z_p(\lambda), \lambda)|_{B_{p+1}}.
\end{eqnarray}
We first investigate the invertibility of $T_p(\lambda)$ for 
$\lambda \in U(B_{r_{p+1}}(\lambda^{(0)}), 2r_{p+1})$. To that end
we split 
\begin{eqnarray}
\label{cv2.157}
T_p = D + R,
\end{eqnarray}
where $D(m, n) = V(\omega)(m) \delta_{m, n}$ and 
$R = DW (\varphi(a) + z_p(\lambda))|_{B_{p+1}}$. 
We will show in steps 1 and 2 below that the diagonal matrix $D$ dominates 
$R$ and that we can find the inverse of
$T_p$ by a Neumann series (step 3). In step 4 we define $z_{p+1}$ and
show that $z_{p+1}$ is analytic and satisfies properties (i)-(iv). 
In a final step we verify
condition (v).

\noindent
{\em Step 1: Estimates on $D$.} We prove that for all
$\lambda \in U(B_{r_{p+1}}(\lambda^{(0)}), 2r_{p+1})$ and 
$m \in B_{p+1}$ 
\begin{eqnarray}
\label{cv2.160}
|V(\omega)(m)| > \frac{d_V s}{2 L_{p+1}^{\tau}}.
\end{eqnarray}

{\em Case 1: $1 \leq |<m,g>| \leq \nu$.}

\noindent
Set $l := |<m,g>|$, then $V(\omega)(m) = V_l(< \omega, m>)$.
Since $|\omega - \omega^{(0)}| < 3 r_{p+1}$, it follows from 
assumption A3, $m \pm e_l \neq 0$ (since $m \notin {\cal S}$), and
(\ref{cv2.20}) that
\begin{eqnarray}
\nonumber
\label{cv2.165}
\mbox{ dist }(<\omega, m>, \{ \pm \omega_l^{(0)} \} + 2 \pi {\Bbb Z}) &>&
\mbox{ dist }(<\omega^{(0)}, m \mp e_l>, 2 \pi {\Bbb Z}) - 3 \nu r_{p+1}
L_{p+1} 
\\
\nonumber
&\geq& \frac{s}{L_{p+1}^{\tau}} - 3 \nu r_{p+1} L_{p+1} \geq 
\frac{s}{2 L_{p+1}^{\tau}}.
\end{eqnarray}
Estimate (\ref{cv2.160}) follows from proposition \ref{Plop.1} (a)
together with
(\ref{cv2.25}).

{\em Case 2: $<m, g>=0$.}

\noindent
By definition $V(\omega) = \alpha_1$ and (\ref{cv2.160}) follows
from (\ref{cv2.30}).

{\em Case 3:$|<m,g>| \geq \nu + 1$.}

\noindent
Observe that $|\mbox{ Im}(<\omega, m>)| \leq 2 \nu L_{p+1}^{1 - E_{\rho}}
\leq 2 \nu N_0^{1 - E_{\rho}}$, since $E_{\rho} > 1$ (cf. (\ref{cv2.15})).
Estimate 
(\ref{cv2.160}) then follows from (\ref{cv2.35}).

\noindent
{\em Step 2: Estimates on $R$.}
Using statement (i) of the induction hypothesis, the definition of $D_{\tau}$ 
((\ref{npc.110}), see also (\ref{cv2.140}))
and $\mu_p \leq 2$ we obtain
with (\ref{cv2.40}) and (\ref{cv2.45})
\begin{eqnarray}
\label{cv2.170}
\| \varphi (a) + z_p (\lambda) \|_{\mu_p, 1}
\leq
D_N 2^{\nu + 2} \nu e^2 |a| + 
4 e^6 \nu^2 D_W D_N^2 D_{\tau} N_0^{\tau + 1} |a|^2
< \frac{r_{F, b}}{2}.
\end{eqnarray}
Thus lemma \ref{Lenp.1} can be applied 
and again by (\ref{cv2.40}) we conclude for
$\lambda \in U(B_{r_{p+1}}(\lambda^{(0)}), 2 r_{p+1})$
\begin{eqnarray}
\label{cv2.175}
\| R \|_{\mu_p, 1} \leq 3 (2^{\nu + 2} \nu e^2 D_N + 1) r_{p+1}.
\end{eqnarray}

\noindent
{\em Step 3: Construction of $T_p^{-1}$.}
It follows from (\ref{cv2.160})  that $D^{-1}$ exists and
\begin{eqnarray}
\label{cv2.180}
\| D^{-1} \|_{\mu_p, 1} \leq D_N \frac{2 L_{p+1}^{\tau}}{d_V s}.
\end{eqnarray}
Using (\ref{cv2.175}) and (\ref{cv2.50}) we see that $\| D^{-1} R \|_{\mu_p, 1} < 1/2$
and we can invert $I + D^{-1} R$ by a Neumann series with bound
$\|(I + D^{-1} R)^{-1}\|_{\mu_p, 1} 
\leq  \|I\|_{\mu_p, 1}+ 1 \leq 2 D_N$.
Therefore
$G_p := (I + D^{-1} R)^{-1} D^{-1}$ exists and 
by (\ref{cv2.55})
\begin{eqnarray}
\label{cv2.185}
\| G_p (\lambda) \|_{\mu_p, 1} \leq 4 D_N^2 \frac{L_{p+1}^{\tau}}{d_V s}
\leq 2^{-(2 \nu + 2)} L_{p+1}^{\tau +1}
\end{eqnarray}
for all $\lambda \in U(B_{r_{p+1}}(\lambda^{(0)}), 2r_{p+1})$. 
Note that $G_p$ depends analytically on $\lambda$, since $D^{-1}$ and 
$R$ are analytic in $U(B_{r_{p+1}}(\lambda^{(0)}), 2r_{p+1})$.
Moreover $G_p(\lambda) = T_p^{-1}(\lambda)$ for all $\lambda \in
U(B_{r_{p+1}}(\lambda^{(0)}), 2r_{p+1}) \setminus Z_{B_{p+1}}$ 
(see (\ref{ova.600})).

\noindent
{\em Step 4: Definition of $z_{p+1}$ and proof of properties (i)-(iv).}
We set for $\lambda \in U(B_{r_{p+1}}(\lambda^{(0)}), 2r_{p+1})$
\begin{eqnarray}
\label{cv2.190}
z_{p+1}(\lambda) := z_p(\lambda) - G_p(\lambda) \left[
{\cal P}(z_p(\lambda), \lambda)|_{B_{p+1}} \right].
\end{eqnarray}
Clearly, $z_{p+1}$ is an analytic function on its domain of definition.
Property (i) is an immediate consequence of
induction hypothesis (v), (\ref{cv2.185}) and proposition \ref{Pfa1.2}. 
Property (ii) is obvious,
since $z_{p+1}(m)$ vanishes outside $B_{p+1}$ by construction.
For $\lambda \in B(B_{r_{p+1}}(\lambda^{(0)}), 2 r_{p+1})$ 
we have by induction 
hypothesis that $z_p(\lambda)(m) = z_p(\lambda)(-m) \in {\Bbb R}$, implying
${\cal P}(z_p(\lambda), \lambda) (m) =
{\cal P}(z_p(\lambda), \lambda)(-m) \in {\Bbb R}$ by the definition of 
${\cal P}$ (see (\ref{ova.160})) and the analyticity of 
${\cal P}(z_p(\cdot), \cdot)$. It then follows from proposition 
\ref{Plop.2}(a) and $B_{p+1} = - B_{p+1}$ that $z_{p+1}$ also satisfies
property (iii). To prove (iv) assume that $a_i=0$ for some $1 \leq i \leq 
\nu$. By induction hypothesis $z_p(\lambda)(m) = 0$, if $m_i \neq 0$.
One readily verifies from the definition that 
${\cal P}(z_p(\lambda), \lambda) (m) = 0$ for $m_i \neq 0$. Furthermore,
proposition \ref{Plop.2} (b) shows that the matrix $T_p$ is of block form
$T_p(m, n) = 0$ for $m_i - n_i \neq 0$. Of course, the block form is preserved
under inversion of the matrix. This implies $z_{p+1}(\lambda)(m) = 0$ for
$m_i \neq 0$.

\noindent
{\em Step 5: Proof of property (v).} To prove analyticity of the function
${\cal P}(z_p(\cdot), \cdot)$ across $Z_{B_{p+1}}$ we need to convince
ourselves of the analyticity of $D (z_{p+1} - z_p)$ at the singularities
of the diagonal matrix $D$. From (\ref{cv2.190})
it follows for $\lambda \in 
U(B_{r_{p+1}}(\lambda^{(0)}), 2r_{p+1}) \setminus Z_{B_{p+1}}$ that
\begin{eqnarray}
\label{cv2.197}
D (z_{p+1} - z_p)(\lambda) = 
- {\cal P}(z_p(\lambda), \lambda)|_{B_{p+1}} + R(\lambda) G_p(\lambda)
{\cal P}(z_p(\lambda), \lambda)|_{B_{p+1}},
\end{eqnarray}
which clearly has an analytic continuation across $Z_{B_{p+1}}$.
To obtain an estimate on ${\cal P}(z_{p+1}(\lambda), \lambda)$
we expand
\begin{eqnarray}
\nonumber
{\cal P}(z_{p+1}(\lambda), \lambda) &=&
\left[
{\cal P}(z_{p}(\lambda), \lambda) + T_p (\lambda) (z_{p+1}(\lambda) -
z_p(\lambda))
\right] \\
\nonumber
&&+
\left[
(D_v{\cal P}(z_{p}(\lambda), \lambda) - T_p (\lambda)) (z_{p+1}(\lambda) -
z_p(\lambda))
\right] \\
\nonumber
&&+
\left[
\int_0^1 (1-t) D^2_{vv} {\cal P}
(z_p + t (z_{p+1}-z_p), \lambda)[z_{p+1}-z_p, z_{p+1}-z_p] dt
\right] \\
\nonumber
&=& I + II + III.
\end{eqnarray}
By continuity it suffices to show that the following three
estimates hold for all \\ $\lambda \in 
U(B_{r_{p+1}}(\lambda^{(0)}), 2r_{p+1}) \setminus Z_{B_{p+1}}$.
\begin{eqnarray}
\label{cv2.186}
\| I \|_{\mu_{p+1}, 1} &\leq&
\frac{1}{3} 2^{2\nu +4} e^6 \nu^2 D_W D_N^2 e^{-1.5^{p+1}} |a|^2, \\
\label{cv2.187}
\| II \|_{\mu_{p+1}, 1} &\leq&
\frac{1}{3} 2^{2\nu +4} e^6 \nu^2 D_W D_N^2 e^{-1.5^{p+1}} |a|^2, \\
\label{cv2.188}
\| III \|_{\mu_{p+1}, 1} &\leq&
\frac{1}{3} 2^{2\nu +4} e^6 \nu^2 D_W D_N^2 e^{-1.5^{p+1}} |a|^2.
\end{eqnarray}

From (\ref{cv2.190}), (\ref{wf.130})
it follows that
\begin{eqnarray}
\nonumber
\| I \|_{\mu_{p+1}, 1} &=& 
\| {\cal P}(z_{p}(\lambda), \lambda) -
{\cal P}(z_{p}(\lambda), \lambda)|_{B_{p+1}} \|_{\mu_{p+1}, 1} \\
\nonumber
&\leq&
\| {\cal P}(z_{p}(\lambda), \lambda) \|_{\mu_p, 1} 
e^{-(\mu_p - \mu_{p+1})L_{p+1}}.
\end{eqnarray}
Using the induction hypothesis and (\ref{cv2.115}) it suffices to show that
\begin{eqnarray}
\label{cv2.195}
3 e^{-1.5^p} e^{-\frac{1}{p(p+1)}L_{p+1}} \leq e^{-1.5^{p+1}},
\end{eqnarray}
or equivalently
\begin{eqnarray}
\label{cv2.200}
a_p := \log 3 + \frac{1}{2} 1.5^p \leq \frac{1}{p(p+1)} N_0 5^p =: b_p.
\end{eqnarray}
Inequality (\ref{cv2.200}) follows inductively: 
one verifies easily that $0 < a_1 \leq b_1$ (observe that $N_0 \geq 5$ by
(\ref{npc.400})) and that
\begin{eqnarray}
\nonumber
\frac{a_{p+1}}{a_p} \leq 1.5 < \frac{5}{3} \leq  \frac{b_{p+1}}{b_p}
\;\;\; \mbox{ for } p \in {\Bbb N}.
\end{eqnarray}
Thus (\ref{cv2.186}) is proved.

It follows from the definition of $T_p(\lambda)$ in (\ref{cv2.155})
and from supp$((z_{p+1}-z_p) (\lambda)) \subset B_{p+1}$ that
\begin{eqnarray}
\label{cv2.205}
(II)(m) &=& 0 \;\; \mbox{ for } |m| < L_{p+1}, \\
\label{cv2.210}
(II)(m) &=& \sum_{n \in B_{p+1}} R(m, n) (z_{p+1}-z_p)(n)
\;\; \mbox{ for } |m| \geq L_{p+1}.
\end{eqnarray}
Hence we obtain
\begin{eqnarray}
\label{cv2.215}
\| II \|_{\mu_{p+1}, 1} \leq \| R \|_{\mu_{p}, 1}
\| z_{p+1} - z_p \|_{\mu_p, 1} e^{-(\mu_p - \mu_{p+1})L_{p+1}}.
\end{eqnarray}
Using (\ref{cv2.175}), the already proven property (i), 
(\ref{cv2.15}), and (\ref{cv2.60}) we observe that 
(\ref{cv2.187}) follows from (\ref{cv2.195}).

Using (\ref{cv2.170}) and the estimates leading to (\ref{cv2.170})
we see that both 
\begin{eqnarray}
\nonumber
\| \varphi (a) + z_p (\lambda) \|_{\mu_p, 1} &\leq& \frac{r_{F, b}}{2}
\;\;\; \mbox{ and } \\
\nonumber
\| \varphi (a) + z_{p+1} (\lambda) \|_{\mu_p, 1} &\leq& \frac{r_{F, b}}{2}.
\end{eqnarray}
Therefore
$\| \varphi(a) + z_p(\lambda) + t (z_{p+1}-z_p)(\lambda) \|_{\mu_p, 1}
\leq r_{F,b}/2$ for  all  $\lambda \in 
U(B_{r_{p+1}}(\lambda^{(0)}), 2r_{p+1})$ and $0 \leq t \leq 1$.
We can apply lemma \ref{Lenp.1} and obtain
\begin{eqnarray}
\label{cv2.220}
\| III \|_{\mu_p, 1} \leq D_W \|z_{p+1} - z_p\|_{\mu_p, 1}^2.
\end{eqnarray}
Estimate (\ref{cv2.188}) follows from property (i), and from 
(\ref{cv2.15}), (\ref{cv2.65}).
\end{proof}

\section{The construction of $v_{j+1}$ for $j \geq 2$}
\label{cvj}
In the case $j \geq 2$ we also construct 
$v_{j+1}$ from $v_j$ by a modified Newton scheme. In contrast to the
previous section, the inverse of the linearized operators cannot
be obtained by a Neumann series. However, induction hypothesis
$({\cal IS})(3c)_{1 \leq l \leq j}$ and
$({\cal IS})(4)_{1 \leq l \leq j}$ imply that we can find sufficiently
many local inverse matrices such that the coupling lemma \ref{Lcl.1}
can be applied. We derive estimates on the inverse of the linearized
operators which are sufficient for the 
Newton scheme to work. 

\noindent
{\bf Estimates of lemma \ref{hyp} used in section \ref{cvj}:}

\begin{eqnarray}
\label{cvj.5}
\frac{1}{8} N_1^c (5^c - 1) \geq \log 2 \\
\label{cvj.10}
4 \rho_j \leq \rho_{j-1} \\
\label{cvj.15}
2 \alpha_1 > d_V \delta_1 \\
\label{cvj.20}
\frac{(\nu+1)^2 \gamma^2}{4 \cosh^2 (\nu N_0^{1-E_{\rho}}/2)} - \alpha_1 > 
\frac{d_V \delta_1}{2} \\
\label{cvj.25}
4 D_W e^{-\frac{1}{8} N_{l-1}^c} M_l \leq 1 \;\; \mbox{ for }
2 \leq l \leq j \\
\label{cvj.27}
16 e \nu D_N \rho_1 \leq r_{F, b}\\
\label{cvj.30}
8 e^5 \nu D_W D_N D_{\tau} N_0^{\tau+1} \rho_1 \leq 1 \\
\label{cvj.35}
8 \rho_1 \leq e \nu D_N  \\
\label{cvj.40}
64 e \nu D_W D_N^2 \rho_1 \leq d_V \delta_1 \\
\label{cvj.45}
4 D_N M_l e^{-\frac{1}{4l(l+1)}N_{l-1}^c} \leq 1 \;\; \mbox{ for }
2 \leq l \leq j\\
\label{cvj.50}
M_2  \geq \frac{2}{d_V \delta_1} \\
\label{cvj.55}
4 D_N^2 \leq N_1 \\
\label{cvj.60}
N^{A E_M+1} \leq e^{\frac{1}{16} N^c} \;\; \mbox{ for all }
N \geq N_1 \\
\label{cvj.65}
5^c \leq \frac{19}{18} \\
\label{cvj.70}
3 (4 e \nu D_N)^2 D_W \leq e^{\frac{3}{4}N_1^c} \\
\label{cvj.75}
8 e \nu D_W D_N \rho_1 \leq 1 \\
\label{cvj.80}
3 \leq e^{\frac{1}{8} N_1^c}
\end{eqnarray}

\begin{Lemma}
\label{Lcvj}
Let $j \geq 2$ and suppose that 
$({\cal IS})_{l}$ is satisfied for $1 \leq l \leq j$. 
Then there exists a function $v_{j+1}: 
U({\cal N}^{(1)}, \rho_1) \to X_{1,c}$ satisfying induction statement
$({\cal IS})(1)_{j+1}$.
\end{Lemma}

\begin{proof}
The proof is similar to the proof of lemma \ref{Lcv2}. Again we will
construct $v_{j+1}$ from $v_j$ by a modified Newton scheme. 
The main difference,
however, is that the existence of (and estimates on) the inverse 
of the linearized operator does not follow from the dominance of the diagonal
part. Instead we use the induction hypothesis and the coupling lemma 
\ref{Lcl.1}.

By definitions (\ref{npc.400}), (\ref{npc.410}) and (\ref{npc.72})
there exists a $p_0 \in {\Bbb N}$
such that $N_j = N_{j-1} 5^{p_0}$. For $1 \leq p \leq p_0+1$ define
\begin{eqnarray}
\label{cvj.100}
L_p &:=& N_{j-1} 5^{p-1},\\
\label{cvj.105}
B_p &:=& B_{L_p}(0) \setminus ({\cal S} \cup \{ 0 \}) \subset {\Bbb Z}^{\nu}.
\end{eqnarray}
We construct inductively analytic maps
$z_p: U({\cal N}^{(j)}, \rho_j) \to X_{1, c}$, $1 \leq p \leq p_0+1$, such that
the following holds:
\begin{itemize}
\item[(i)] For all $\lambda \in U({\cal N}^{(j)}, \rho_j)$:
\begin{eqnarray}
\nonumber
z_p (\lambda) = v_{j}(\lambda)&& \mbox{ for } p=1,
\\
\nonumber
\|(z_{p} - z_{p-1})(\lambda) \|_{1, c} \leq 
e^{-\frac{1}{8} L_p^c} |a|^2
&&\mbox{ for }  
2 \leq p \leq p_0+1.
\end{eqnarray}
\item[(ii)] For all $\lambda \in U({\cal N}^{(j)}, \rho_j)$
the support of $z_p(\lambda)$ is contained in $B_p$.
\item[(iii)]
For all $\lambda \in B({\cal N}^{(j)}, \rho_j)$ and 
$m \in {\Bbb Z}^{\nu}$
we have $z_p(\lambda)(m) = z_p(\lambda)(-m) \in {\Bbb R}$.
\item[(iv)]
Let $\lambda = (a, \omega) \in U({\cal N}^{(j)}, \rho_j)$ with 
$a_i = 0$, for some $1 \leq i \leq \nu$. Then $z_p(\lambda)(m) = 0$, 
if $m_i \neq 0$.
\item[(v)]
The map $\lambda \mapsto {\cal P}(z_p(\lambda), \lambda)$ 
has an analytic continuation to
$U({\cal N}^{(j)}, \rho_j)$ and
\begin{eqnarray}
\nonumber
\| {\cal P}(z_p(\lambda), \lambda) \|_{1/4, c} \leq 
e^{-L_p^c} |a|^2 \;\;\; \mbox{ for all } 
\lambda \in U({\cal N}^{(j)}, \rho_j).
\end{eqnarray}
\end{itemize}

We will first complete the proof of Lemma \ref{Lcvj}, assuming that such 
functions $z_p$ exist. Recall the definition of the
function $\psi$ in section \ref{npc} G). For $\epsilon > 0$ set 
\begin{eqnarray}
\label{cvj.115}
\psi_{\epsilon}(x) &:=& \epsilon^{-4 \nu} \psi \left(\frac{x}{\epsilon}\right)
\;\; \mbox{ for } x \in {\Bbb R}^{4 \nu}.
\end{eqnarray}
Furthermore we denote
\begin{eqnarray}
\label{cvj.120}
X^{(j)} &:=& U \left( {\cal N}^{(j)}, \frac{\rho_j}{2} \right) \subset
{\Bbb C}^{2 \nu} \equiv {\Bbb R}^{4 \nu},
\\
\label{cvj.125}
\Psi^{(j)}(x) &:=& \int_{X^{(j)}} \psi_{\rho_j/8}(x-y)dy.
\end{eqnarray}
It is easy to see that the cut-off function
$\Psi^{(j)} \in C^{\infty}({\Bbb R}^{4 \nu})$ 
has the following properties:
\begin{eqnarray}
\label{cvj.130}
\Psi^{(j)}(x) &=& 1 \;\; \mbox{ for } 
x \in U({\cal N}^{(j)}, \rho_j/4), \\
\label{cvj.132}
\mbox{ supp}(\Psi^{(j)}) &\subset&  U({\cal N}^{(j)}, 3 \rho_j/4),
\\
\label{cvj.135}
0 \leq \Psi^{(j)}(x) &\leq& 1 \;\; \mbox{ for } x \in {\Bbb R}^{4 \nu},
\\
\label{cvj.138}
\left|
\partial^{\beta} \Psi^{(j)} (x) \right| &\leq&
\left(
\frac{8}{\rho_j} \right)^{|\beta|_1} D_{\psi}(|\beta|_1) \;\; \mbox{ for }
x \in {\Bbb R}^{4 \nu}.
\end{eqnarray}
We define
\begin{eqnarray}
\label{cvj.140}
v_{j+1}(\lambda) :=
\left\{
\begin{array}{ll}
v_j(\lambda) + (z_{p_0+1}(\lambda) - v_j(\lambda)) \Psi^{(j)}(\lambda) &
\mbox{ for } \lambda \in U({\cal N}^{(j)}, \rho_j), \\
v_j(\lambda)&
\mbox{ for } \lambda \in U({\cal N}^{(1)}, \rho_1) \setminus
U({\cal N}^{(j)}, \rho_j)
\end{array}
\right.
\end{eqnarray}
From (\ref{cvj.132}) we conclude that $v_{j+1}$ is a $C^{\infty}$ function.
Furthermore, by (\ref{cvj.130}), $v_{j+1}$ coincides with $z_{p_0+1}$
on the set $U({\cal N}^{(j)}, \rho_j/4)$ and is therefore analytic on
this set. Thus we have already verified properties (1a) and (1b) of the
induction statement $({\cal IS})_{j+1}$. By definition (see (\ref{cvj.100}))
we have $L_{p_0+1} = N_j$. 
Statement (c) of $({\cal IS})_{j+1}$ is an immediate consequence of 
(\ref{cvj.130}) and property (v) of $z_{p_0+1}$. Note that for 
$\lambda \in U({\cal N}^{(j)}, \rho_j)$
\begin{eqnarray}
\label{cvj.145}
v_{j+1}(\lambda) - v_j(\lambda) &=& \Psi^{(j)}(\lambda) y(\lambda) \quad 
\mbox{ with }\\
\label{cvj.150}
y(\lambda) := \sum_{p=1}^{p_0} z_{p+1}(\lambda) - z_{p}(\lambda).
\end{eqnarray}
The estimate (i) on $z_{p}-z_{p-1}$ together with 
$e^{-\frac{1}{8} L_{p+1}^c} \leq \frac{1}{2} e^{-\frac{1}{8} L_{p}^c}$ for
$1 \leq p \leq p_0$
(a consequence of (\ref{cvj.5})) yields
\begin{eqnarray}
\nonumber
\|z_{p+1}(\lambda) - z_p(\lambda)\|_{1, c} \leq 
2^{-p} e^{-\frac{1}{8} N_{j-1}^c} |a|^2
\end{eqnarray}
and therefore
\begin{eqnarray}
\label{cvj.155}
\| y(\lambda) \|_{1,c} \leq e^{-\frac{1}{8} N_{j-1}^c} |a|^2 \;\;
\mbox{ for } \lambda \in U({\cal N}^{(j)}, \rho_j).
\end{eqnarray}
Since $|\Psi^{(j)}(\lambda)| \leq 1$ (see (\ref{cvj.135})) we 
have proved statement
(d) of $({\cal IS})_{j+1}$. Furthermore, by the analyticity of $y$ in
$U({\cal N}^{(j)}, \rho_j)$ the Cauchy integral formula implies
\begin{eqnarray}
\label{cvj.160}
\| \partial^{\beta} y (\lambda) \|_{1,c} \leq \beta !
\left( \frac{4}{\rho_j} \right)^{|\beta|_1} e^{-\frac{1}{8} N_{j-1}^c} \;\;
\mbox{ for } \lambda \in U({\cal N}^{(j)}, 3 \rho_j/4).
\end{eqnarray}
Applying the Leibniz rule to the right hand side of (\ref{cvj.145})
statement (1e) of $({\cal IS})_{j+1}$ now follows from
(\ref{cvj.132}), (\ref{cvj.138}), (\ref{cvj.145}) and (\ref{cvj.160}).
Statements (1f), (1g) and (1h) of $({\cal IS})_{j+1}$ follow by
(\ref{cvj.132}) and (\ref{cvj.140}) from the corresponding
properties of $v_j$ (see (1f), (1g), (1h) of $({\cal IS})_{j}$)
and of $z_{p_0+1}$ (see (ii), (iii) and (iv) and recall that
$L_{p_0+1} = N_{j}$ by the definition of $p_0$.)

Thus the proof of Lemma \ref{Lcvj} is reduced to constructing
analytic maps $z_p: U({\cal N}^{(j)}, \rho_j) \to X_{1, c}$,
satisfying properties (i)-(v).

\noindent
{\em Inductive construction of $z_p$:}

\underline{$p=1$}: By (i) the choice of $z_1(\lambda) = v_j(\lambda)$,
$\lambda \in U({\cal N}^{(j)}, \rho_j)$, 
is determined. Since $U({\cal N}^{(j)}, \rho_j) \subset
U({\cal N}^{(j-1)}, \rho_{j-1}/4)$ (see (\ref{cvj.10}),
$({\cal IS})(2a)_j$), it follows from
$({\cal IS})_{j}(1b)$ that $z_1$ is analytic. Recall 
furthermore that $L_1$ as defined
in (\ref{cvj.100}) equals $N_{j-1}$ and properties (ii), (iii),
(iv) and (v) are satisfied for $z_1$ by statements
(1f), (1g), (1h) and (1c) of $({\cal IS})_{j}$.

\underline{$p \to p+1$, $1 \leq p \leq p_0$:}
Define
\begin{eqnarray}
\label{cvj.165}
T_{p}(\lambda) := D_v {\cal P}(z_p(\lambda), \lambda)|_{B_{p+1}} \;\;
\mbox{ for } \lambda \in U({\cal N}^{(j)}, \rho_j) \setminus
Z_{B_{p+1}}.
\end{eqnarray}
In order to see that (\ref{cvj.165}) is well defined, we first show that
\begin{eqnarray}
\label{cvj.167}
\| \varphi (a) + z_p (\lambda) \|_{1, c} \leq 
4 e \nu D_N |a| \leq \frac{r_{F, b}}{2} \;\; \mbox{ for }
\lambda \in U({\cal N}^{(j)}, \rho_j).
\end{eqnarray}
Indeed, for $\lambda \in U({\cal N}^{(j)}, \rho_j)$,
\begin{eqnarray}
\nonumber
\| \varphi (a) + z_p (\lambda) \|_{1, c}
&\leq&
\| \varphi (a) \|_{1, c} +
\| v_2 (\lambda) - v_1 (\lambda)\|_{1, c} \\
\nonumber
&+& 
\| v_j (\lambda) - v_2 (\lambda)\|_{1, c}
+
\| z_p (\lambda) - v_j (\lambda)\|_{1, c} \\
\nonumber
&=&
I + II + III + IV,
\end{eqnarray}
where
\begin{eqnarray}
\nonumber
I &\leq&
2 \nu |a| D_N e, \\
\nonumber
II &\leq&
4 e^6 \nu^2 D_W D_N^2 D_{\tau} N_0^{\tau + 1} |a|^2 \leq  e \nu D_N |a|
\;\; (\mbox{ see } ({\cal IS})(1d)_2, (\ref{cvj.30})),
\\
\nonumber
III &\leq&
\sum_{l=3}^j \|v_l - v_{l-1}\|_{1, c} \leq
\sum_{l=3}^j e^{-\frac{1}{8} N_{l-2}^c} |a|^2
\leq 4 e^{-\frac{1}{8} N_{1}^c} \rho_1 |a| \\
\nonumber 
&\leq& \frac{1}{2} e \nu D_N |a| \;\;\; (\mbox{ see } 
({\cal IS})(1d)_l, (\ref{cvj.5}),
(\ref{cvj.35})),\\
IV &\leq&
\sum_{l=2}^p \|z_l - z_{l-1}\|_{1, c} \leq
\sum_{l=2}^p e^{-\frac{1}{8} L_{l}^c} |a|^2
\leq 2 e^{-\frac{1}{8} N_{j-1}^c} \rho_1 |a| \\
\nonumber 
&\leq& \frac{1}{4} e \nu D_N |a| \;\;\; (\mbox{ see } (i), (\ref{cvj.5}),
(\ref{cvj.35})).
\end{eqnarray}
Using in addition (\ref{cvj.27}) we have proved 
(\ref{cvj.167}).

{\em Claim1:}

{\em
\noindent
There exists an analytic map $\lambda \mapsto G_p(\lambda)$ on 
$U({\cal N}^{(j)}, \rho_j)$ satisfying
\begin{eqnarray}
\label{cvj.169}
\| G_p(\lambda) \|_{\sigma_{j+1}, c} \leq L_p^{A E_M + 1} \;\;\;
\mbox{ for } \lambda \in U({\cal N}^{(j)}, \rho_j).
\end{eqnarray}
Moreover, for $\lambda \in U({\cal N}^{(j)}, \rho_j) \setminus Z_{B_{p+1}}$
the matrix $G_p(\lambda)$ is the inverse of $T_p(\lambda)$.
}

{\em Proof of claim 1:}
Let $\lambda \in U({\cal N}^{(j)}, \rho_j)$. Then there
exists a (not necessarily unique) $k_j \in K^{(j)}$ such that
$\lambda \in U(B_{\rho_j}(\lambda_{k_j}^{(j)}), \rho_j)$. Set
$k_l := \pi_l^{(j)} k_j$ for $1 \leq l \leq j-1$. The
definition of the projection $\pi_l^{(j)}$ (section \ref{npc} J) 
together with (\ref{cvj.10}) imply
\begin{eqnarray}
\label{cvj.170}
\lambda \in U \left(B_{\rho_l}(\lambda_{k_l}^{(l)}), \frac{\rho_l}{4} 
\right) \;\;
\mbox{ for } 1 \leq l \leq j-1.
\end{eqnarray}
We set
\begin{eqnarray}
\label{cvj.175}
E_0 := \{ m \in B_{p+1}: |<m, g>| \in \{1, 2, \ldots, \nu \} \}.
\end{eqnarray}
For each $m \in E_0$ we denote by $m' \in
{\Bbb Z}^{\nu}$ the (uniquely defined)
lattice point satisfying $<m', g> = 0$ and $m - m' \in {\cal S}$.
Furthermore, we define for $n \in {\Bbb Z}^{\nu}$
\begin{eqnarray}
\label{cvj.180}
C_l(n) &:=& 
\left\{
\begin{array}{ll}
\left(
B_{p+1} - \{ n \}
\right)
\cap B_{1.5} (0) & \mbox{ if } l = 1,  \\
\left(
B_{p+1} - \{ n \}
\right)
\cap B_{N_{l-1}} (0) & \mbox{ if } 2 \leq l \leq j. 
\end{array}
\right.
\end{eqnarray}
Observe that for $m \in E_0$ with $|m'| > 2 N_{l-1}$, $l \geq 2$, the sets
$C_s(m')$ are elements of ${\cal C}^{(s)}$ for all $1 \leq s \leq l$.
Indeed, $m - m' \in {\cal S} \cap C_s(m')$ by definition. Moreover, it is also
obvious from (\ref{cvj.180}) that $C_s(m') \subset B_{N_{s-1}}(0)$ for
$2 \leq s \leq l$, respectively $C_1(m') \subset B_{1.5}(0)$. To see that
$C_s(m')$ is an interval in the case $2 \leq s \leq l$ we employ 
$|m| \geq 2 N_{l-1}$. Hence the following is a proper definition. 
\begin{eqnarray}
\nonumber
\Gamma_l &:=&
\left\{ m \in E_0: |m'| > 2 N_{l-1} \mbox{ and } 
<\omega, m'> \in 
U
\left(
\tilde{I}_{C_{l-1}(m'), k_{l-1}}^{(l-1)} \setminus
\tilde{I}_{C_{l}(m'), k_{l}}^{(l)}, \frac{\delta_l}{10}
\right)
\right\} \\
\label{cvj.185}
&& \quad \quad \quad \mbox{ for } 2 \leq l \leq j, \\
\label{cvj.190}
\Gamma_1 &:=&
B_{p+1} \setminus
\left(
\bigcup_{l=2}^{j} \Gamma_l
\right).
\end{eqnarray}
The definition of the sets $\Gamma_l$ is motivated
by the following properties:

\begin{itemize}
\item[($\alpha$)]
For all $m \in \Gamma_1$: 
$|V(\omega)(m)| > \frac{d_V \delta_1}{2}$.
\item[($\beta$)]
$\Gamma_l \subset \{ m \in {\Bbb Z}^{\nu}: |m| \geq 2 N_{l-1} \}$
for $2 \leq l \leq j$.
\item[($\gamma$)]
For $2 \leq l \leq j$ and $m \in \Gamma_l$ we have 
\begin{eqnarray}
\label{cvj.195}
\left\| \left( (T_p)(\lambda)|_{C_l(m') + \{m'\}} \right)^{-1} 
\right\|_{\sigma_l, c} \leq 2 D_N M_l
\;\;\; \mbox{ (cf. remark \ref{Rcvj.1} below)}.
\end{eqnarray}
\end{itemize}

{\em
\begin{Remark}
\label{Rcvj.1}
In the proof of claim 1 we have chosen an arbitrary but fixed $\lambda
\in U({\cal N}^{(j)}, \rho_j)$. Since the matrix 
$T_p(\lambda)|_{C_l(m') + \{m'\}}$ does not
have finite entries for $\lambda \in Z_{C_l(m') + \{m'\}}$ we need to
interpret statement ($\gamma$) in this case.
We understand by claim ($\gamma$) the following. There exists an open
neighborhood $U_{\lambda}$ of $\lambda$ and an analytic map
$\mu \mapsto G(\mu)$ defined on $U_{\lambda}$, such that $G(\mu)$ is
the inverse of $(T_p)(\mu)|_{C_l(m') + \{m'\}}$ for $\mu \in 
U_{\lambda} \setminus Z_{C_l(m') + \{m'\}}$. Moreover, $G(\lambda)$
satisfies estimate (\ref{cvj.195}).
\end{Remark}
}

{\em Proof of ($\alpha$):}
Let $m \in \Gamma_1$.

{\em Case 1: $m \in B_{p+1} \setminus E_0$}. Suppose first that $<m, g> = 0$.
Then  $V(\omega)(m) = \alpha_1$ and the claim follows from (\ref{cvj.15}).
Therefore we may assume that $|<m, g>| \geq \nu + 1$. Since
the imaginary part of $<\omega, m>$ is bounded by $\nu \rho_j N_j \leq 
\nu N_0^{1-E_{\rho}}$, the claim follows from (\ref{cvj.20}).

{\em Case 2: $m \in E_0$}. Recall that definition (\ref{cvj.105})
implies ${\cal S} \cap B_{p+1} = \emptyset$. Furthermore, $m'$ has
been chosen such that $m-m' \in {\cal S}$ and hence we have
\begin{eqnarray}
\label{cvj.200}
0 < |m'| \leq N_j.
\end{eqnarray}
We can define $l \in \{1, ..., j\}$ to be the minimal integer such that
$|m'| \leq 2 N_l$ is satisfied. Since
$k_l \in K^{(l)}$, $\lambda \in U(B_{\rho_l}(\lambda_{k_l}^{(l)}), \rho_l)$,
$C_l (m') \in {\cal C}^{(l)}$, 
$<m', g> = 0$ and $0 < |m'| \leq 2 N_1$ for
$l=1$, respectively, $2 N_{l-1} < |m'| \leq 2 N_l$ for $2 \leq l \leq j$,
the inductive statement
$({\cal IS}) (3c)_{l}$ implies that
\begin{eqnarray}
\label{cvj.205}
<\omega , m'> \in U
\left( {\Bbb R} \setminus \tilde{I}_{C_l(m'), k_l}^{(l)}, 
\frac{\delta_l}{10}
\right).
\end{eqnarray}
Observe that the definition of $C_s(m')$ and $k_s$ together with the
induction statement $({\cal IS})(3a)_{s}$ yield 
\begin{eqnarray}
\label{cvj.210}
\tilde{I}_{C_s(m'), k_s}^{(s)} \subset 
\tilde{I}_{C_{s-1}(m'), k_{s-1}}^{(s-1)} \;\; \mbox{ for } 2 \leq s \leq j.
\end{eqnarray}
Thus 
\begin{eqnarray}
\label{cvj.215}
{\Bbb R} \setminus \tilde{I}_{C_l(m'), k_l}^{(l)} =
{\Bbb R} \setminus \tilde{I}_{C_1(m'), k_1}^{(1)} \cup
\bigcup_{s=2}^l \left(
\tilde{I}_{C_{s-1}(m'), k_{s-1}}^{(s-1)} \setminus
\tilde{I}_{C_s(m'), k_s}^{(s)}
\right)
\end{eqnarray}
The definition of $\Gamma_1$ (see (\ref{cvj.190})) 
implies $m \notin \Gamma_s$
for $2 \leq s \leq j$. It follows from (\ref{cvj.205}), (\ref{cvj.215})
and (\ref{cvj.185}) that
\begin{eqnarray}
\label{cvj.220}
<\omega , m'> \in U
\left( {\Bbb R} \setminus \tilde{I}_{C_1(m'), k_1}^{(1)}, 
\frac{\delta_l}{10}
\right).
\end{eqnarray}
Clearly, we have $\delta_l \leq \delta_1$. Induction
statement $({\cal IS}) (4)_1$ ($k \equiv k_1$, $C \equiv C_1(m')$)
implies
\begin{eqnarray}
\label{cvj.225}
|V(<\omega, m'>, \omega)(m-m')| > \frac{d_V \delta_1}{2}.
\end{eqnarray}
Claim ($\alpha$) now follows from 
$V(<\omega, m'>, \omega)(m-m') = V(0, \omega)(m) = V(\omega)(m)$ which is 
a consequence of $<m', g> = 0$ (cf. (\ref{smr2.70}), (\ref{ova.520})).

{\em Proof of ($\beta$):} Obvious from (\ref{cvj.185}) and $|m-m'| = 1$.

{\em Proof of ($\gamma$):} 
Let $m \in \Gamma_l$ for $2 \leq l \leq j$. Recall that $C_l(m') \in
{\cal C}^{(l)}$. Furthermore
\begin{eqnarray}
\label{cvj.230}
(<\omega, m'>, \lambda) \in
U
\left(
\tilde{I}_{C_{l-1}(m'), k_{l-1}}^{(l-1)} \setminus
\tilde{I}_{C_{l}(m'), k_{l}}^{(l)}, \frac{\delta_l}{10}
\right)
\times
U
\left( 
B_{\rho_l}(\lambda_{k_l}^{(l)}), \rho_l
\right).
\end{eqnarray}
Induction statement $({\cal IS}) (4)_l$ provides the existence
of an inverse matrix $G_{C_l(m')}^{(l)}(<\omega, m'>, \lambda)$  which
clearly has an analytic continuation to some neighborhood $U_{\lambda}$ of
$\lambda$.
Using proposition \ref{Plop.3} and $<m', g>=0$ we have found an inverse
matrix of
$T_{C_l(m')+\{m'\}}^{(l)}(0, \lambda)$ (see remark \ref{Rcvj.1}) with
\begin{eqnarray}
\label{cvj.235}
\left\|
\left(
T_{C_l(m')+\{m'\}}^{(l)}
\right)^{-1}(0, \lambda)
\right\|_{\sigma_l, c} \leq M_l.
\end{eqnarray}

According to proposition \ref{Pcr2.1} 
claim ($\gamma$) follows, if we can show that
\begin{eqnarray}
\label{cvj.240}
\left\|
\left(T_p(\lambda) - T^{(l)}(0, \lambda)
\right)|_{C_l(m')+\{m'\}}
\right\|_{\sigma_l, c} \leq \frac{1}{2M_l}.
\end{eqnarray}
Estimate (\ref{cvj.240}), however, is a consequence of
(\ref{cvj.165}), lemma \ref{Lenp.1}, induction 
hypothesis (i), induction statements $({\cal IS}) (1d)_s$ for
$l < s \leq j$, (\ref{cvj.5})
and (\ref{cvj.25}):
\begin{eqnarray}
\nonumber
&&
\left\|
\left(T_p(\lambda) - T^{(l)}(0, \lambda)
\right)|_{C_l(m')+\{m'\}}
\right\|_{\sigma_l, c} \\
\nonumber
&\leq&\| DW (\varphi(a) + z_p(\lambda)) -
DW (\varphi(a) + v_l(\lambda)) \|_{1, c} \\
\nonumber
&\leq& D_W \| z_p(\lambda) - v_l(\lambda) \|_{1, c}
\\
\nonumber
&\leq& D_W
\left(
\sum_{s=l+1}^{j} \|v_{s}(\lambda) - v_{s-1}(\lambda) \|_{1, c}
+
\sum_{q=2}^{p} \|z_q(\lambda) - z_{q-1}(\lambda)\|_{1, c}
\right) \\
\nonumber
&\leq&
D_W
\left(
\sum_{s \geq l+1} e^{-\frac{1}{8} N_{s-2}^c} +
\sum_{q \geq 2} e^{-\frac{1}{8} L_q^c}
\right) \\
\nonumber
& \leq &
2 D_W e^{-\frac{1}{8} N_{l-1}^c} \leq \frac{1}{2 M_l}.
\end{eqnarray}
This completes the proof of statements ($\alpha$) -- ($\gamma$) 
and we return to the proof of 
claim 1.
To this end we apply the coupling lemma \ref{Lcl.1}. 
In the notation of lemma \ref{Lcl.1}
we set
\begin{eqnarray}
\nonumber
\Lambda &\equiv& B_{p+1}, \\
\nonumber
T_p(\lambda) &\equiv& T = D + R, \mbox{ where }\\
\label{cvj.245}
D(m, n) &=& V(\omega)(m) \delta_{m, n}, \\
\nonumber
R(m, n) &=& DW (\varphi(a) + z_p(\lambda))(m, n),\\
\nonumber
\sigma &\equiv& 1,\\
\nonumber
\tilde{\sigma} &\equiv& \sigma_{j+1}.
\end{eqnarray}
Let us first assume that $\lambda \in U({\cal N}^{(j)}, \rho_j) \setminus
Z_{B_{p+1}}$.
Observe from (\ref{cvj.167}) and lemma \ref{Lenp.1} that
\begin{eqnarray}
\label{cvj.250}
\| R \|_{1, c} \leq 8 e \nu D_W D_N \rho_1.
\end{eqnarray}

Next we have to define the
quantities $l_n$, $\mu_n$, $C_n$ and $U(n)$ for each $n \in B_{p+1}$.
By definition (\ref{cvj.190}) we have $B_{p+1} = \bigcup_{s=1}^j \Gamma_s$. 
Therefore we can pick for every $n \in B_{p+1}$ a (not necessarily
unique) integer $s_n \in \{1, \ldots, j\}$ such that $n \in
\Gamma_{s_n}$. In the case $s_n = 1$ we set
\begin{eqnarray}
\label{cvj.247}
\begin{array}{llll}
l_n:=1,&\mu_n:=0,&C_n:=\frac{4 D_N}{d_V \delta_1},&U(n):=\{n\}.
\end{array}
\end{eqnarray}
Recall that for $n \in \Gamma_l$, $l \geq 2$, we have $n \in E_0$
and therefore there exists an unique $n' \in {\Bbb Z}^{\nu}$ satisfying
$<n', g> = 0$ and $n-n' \in {\cal S}$. In the case $s_n \geq 2$
we define
\begin{eqnarray}
\label{cvj.248}
\begin{array}{llll}
l_n:=N_{s_n-1}-1,&\mu_n:=\sigma_{s_n}-\sigma_{j+1},&C_n:=2 D_N M_{s_n},&
U(n):=C_{s_n}(n')+\{n'\}.
\end{array}
\end{eqnarray}

We now verify the hypothesis (\ref{cl.20}) -- (\ref{cl.40}) of 
lemma \ref{Lcl.1}.

{\em Case $s_n=1$:} The set $U(n)$ contains the single point $n$ and
\begin{eqnarray}
\nonumber
| D(n, n) + R(n, n)| \geq 
\frac{d_V \delta_1}{2} - \| R \|_{1,c} \geq 
\frac{d_V \delta_1}{4}
\end{eqnarray}
by ($\alpha$) and (\ref{cvj.40}). Hence (\ref{cl.20}) 
is satisfied.
Condition (\ref{cl.30}) is clear and condition (\ref{cl.40}) 
follows again from (\ref{cvj.40}) and (\ref{cvj.250}).

{\em Case $2 \leq s_n \leq j$:}
Condition (\ref{cl.20}) is satisfied by statement ($\gamma$)
and (\ref{cl.30}) follows from the definition of $C_{s_n}(n')$
(cf. (\ref{cvj.180})) and $|n-n'|=1$. Finally, condition (\ref{cl.40})
is satisfied by (\ref{npc.56}), (\ref{cvj.45}) and (\ref{cvj.75}).

We can apply lemma \ref{Lcl.1} and obtain that the inverse matrix
$T_p(\lambda)^{-1}$ exists satisfying
a matrix $G_{B_{p+1}}$ satisfying
\begin{eqnarray}
\nonumber
\| T_p(\lambda)^{-1} \|_{\sigma_{j+1}, c} \leq  2 D_N
\max \left( \frac{4 D_N}{d_V \delta_1}, 2 D_N M_j \right).
\end{eqnarray}
Keeping (\ref{cvj.50}), (\ref{cvj.55}), (\ref{npc.410}) and 
(\ref{npc.440}) in mind 
we obtain 
\begin{eqnarray}
\label{cvj.255}
\| T_p(\lambda)^{-1} \|_{\sigma_{j+1}, c} \leq
4 D_N^2 M_j \leq N_1 N_{j-1}^{A E_M} \leq L_p^{A E_M + 1}.
\end{eqnarray}
Recall that we have assumed for the construction of $T_p(\lambda)^{-1}$
that $\lambda \in U({\cal N}^{(j)}, \rho_j) \setminus Z_{B_{p+1}}$.
For $\lambda \in U({\cal N}^{(j)}, \rho_j) \cap Z_{B_{p+1}}$ not all
entries of $T_p(\lambda)$ are finite and therefore we have to explain what
we we understand by $T_p(\lambda)^{-1}$ (see e.g. remark \ref{Rcvj.1}).
Observe from the proof of lemma \ref{Lcl.1} that the construction of 
$T_{\Lambda}^{-1}$ only uses $T_{U(n)}^{-1}$ and $R$. The diagonal 
matrix $D$ which contains the singularities does not appear. Note 
furthermore that $T_{U(n)}^{-1}$ is defined (in the sense of remark
\ref{Rcvj.1}) also for $\lambda \in 
U({\cal N}^{(j)}, \rho_j) \cap Z_{B_{p+1}}$ ( see ($\gamma$) and 
remark \ref{Rcvj.1} for $2 \leq s_n \leq j$; in the case $s_n = 1$
simply use $1/\infty \equiv 0$ to define $T_{\{n\}}^{-1}$ in the case
$D(n, n) = \infty$). Thus we may extend the definition of 
$T_p(\lambda)^{-1}$ to the set
$\lambda \in U({\cal N}^{(j)}, \rho_j) \cap Z_{B_{p+1}}$ via 
lemma \ref{Lcl.1}. Estimate (\ref{cvj.255}) holds for this extension, too.
We denote the extension of $T_p(\lambda)^{-1}$ to all of 
$U({\cal N}^{(j)}, \rho_j)$ by $G_p(\lambda)$.

In conclusion, we have constructed $G_p(\lambda)$ for every $\lambda \in 
U({\cal N}^{(j)}, \rho_j)$ such that (\ref{cvj.169}) is satisfied
and $G_p (\lambda) = T_p^{-1}(\lambda)$ for all $\lambda \in 
U({\cal N}^{(j)}, \rho_j) \setminus Z_{B_{p+1}}$. In order to prove claim 1
we still need to show analyticity of $G_p(\cdot)$. 
Fix again $\lambda \in 
U({\cal N}^{(j)}, \rho_j)$. It is clear from the above construction 
that $G_p(\lambda)$ has an analytic continuation $A(\mu)$ to an open 
neighborhood
$U_{\lambda}$ of $\lambda$ (cf. remark \ref{Rcvj.1}), satisfying
$A (\mu) = T_p^{-1}(\mu)$ for all $\mu \in 
U_{\lambda} \setminus Z_{B_{p+1}}$.
Nevertheless, by construction it is not a-priori
clear that $A(\mu) = G_p(\mu)$ for $\mu \in U_{\lambda} \setminus
\{ \lambda \}$ since e.g. the choice of the local neighborhoods
$U(n)$ (see (\ref{cvj.247}), (\ref{cvj.248})) might differ for different
values of $\lambda$. However, we observe that $A(\mu) = G_p(\mu)$ for
all $\mu \in U_{\lambda} \setminus  Z_{B_{p+1}}$ since they are 
both inverse to $T_p(\mu)$. Suppose now that $\lambda'
\in U_{\lambda} \cap Z_{B_{p+1}}$. By the same
reasoning as above there exists an open neighborhood
$U_{\lambda'}$ of $\lambda'$ and an analytic function $B$ on $U_{\lambda'}$
such that $B(\lambda') = G_p(\lambda')$ and 
$B(\mu)  = T_p^{-1}(\mu)$ for all $\mu \in 
U_{\lambda'} \setminus Z_{B_{p+1}}$. Thus $A(\mu) = B(\mu)$ for all
$\mu \in U_{\lambda} \cap U_{\lambda'} \setminus Z_{B_{p+1}}$ and 
since $Z_{B_{p+1}}$ is a set of measure 0 the functions $A$ and $B$
coincide on  $U_{\lambda} \cap U_{\lambda'}$ by analyticity.
This implies in particular that $A(\lambda') = B(\lambda') = 
G_p(\lambda')$. We have shown that $A(\mu) = G_p(\mu)$
for all $\mu \in U_{\lambda}$ which in turn implies the analyticity
of the function $G_p$ at the point $\lambda$. The proof of claim 1 is
complete.

Next we define $z_{p+1}$. 
Set 
\begin{eqnarray}
\label{cvj.260}
\tilde{B}_{p+1} := B_{p+1} \cap B_{3 L_p}(0),
\end{eqnarray}
and denote by $\chi_{\tilde{B}_{p+1}}$ the corresponding characteristic
function on ${\Bbb Z}^{\nu}$.
For 
$\lambda \in U({\cal N}^{(j)}, \rho_{j})$ we define
\begin{eqnarray}
\label{cvj.265}
z_{p+1}(\lambda) := z_p(\lambda) - G_p(\lambda)
\left[
\chi_{\tilde{B}_{p+1}} {\cal P}(z_p(\lambda), \lambda)
\right].
\end{eqnarray}
From claim 1 and induction hypothesis (v) it follows that 
$z_{p+1}$ is an
analytic map from $U({\cal N}^{(j)}, \rho_{j})$ into $X_{1, c}$.
Moreover, property (v) of the induction hypothesis, $\sigma_{j+1} \geq 
1/4$ (see (\ref{npc.56})),
and (\ref{cvj.255}) yield
\begin{eqnarray}
\label{cvj.270}
\| z_{p+1}(\lambda) - z_p(\lambda) \|_{1/4, c} \leq 
L_p^{A E_M+1} e^{-L_p^c} |a|^2.
\end{eqnarray}
Using in addition hypothesis (\ref{cvj.60}), (\ref{cvj.65}) 
and supp$(z_{p+1}) \subset B_{p+1}$ we obtain
\begin{eqnarray}
\label{cvj.275}
\| z_{p+1}(\lambda) - z_p(\lambda) \|_{1, c} \leq L_p^{A E_M + 1} e^{-L_p^c}
e^{\frac{3}{4} L_{p+1}^c} |a|^2 \leq
e^{-L_{p+1}^c(\frac{15}{16} 5^{-c} - \frac{3}{4})} |a|^2 \leq
e^{-\frac{1}{8} L_{p+1}^c} |a|^2,
\end{eqnarray}
proving property (i). Property (ii) is clear by construction.
The proof of properties (iii) and (iv) is similar to the corresponding proof 
in section \ref{cv2} and will not be repeated here. 
It remains
to show (v). Again, analyticity of the function 
${\cal P}(z_{p+1}(\cdot), \cdot)$ can be dealt with in exactly the same way
as in section \ref{cv2}, where one uses that $D G_p = I - R G_p$ is analytic.
We expand
\begin{eqnarray}
\nonumber
{\cal P}(z_{p+1}(\lambda), \lambda) &=&
\left[
{\cal P}(z_p(\lambda), \lambda) - 
\chi_{\tilde{B}_{p+1}} {\cal P}(z_p(\lambda), \lambda) 
\right] \\
&+&
\left[
\chi_{\tilde{B}_{p+1}} 
{\cal P}(z_{p}(\lambda), \lambda) + T_p (\lambda) (z_{p+1}(\lambda) -
z_p(\lambda))
\right] \\
\nonumber
&+&
\left[
(D_v{\cal P}(z_{p}(\lambda), \lambda) - T_p (\lambda)) (z_{p+1}(\lambda) -
z_p(\lambda))
\right] \\
\nonumber
&+&
\left[
\int_0^1 (1-t) D^2_{vv} {\cal P}
(z_p + t (z_{p+1}-z_p), \lambda)[z_{p+1}-z_p, z_{p+1}-z_p] dt
\right] \\
\nonumber
&=& I + II + III +IV.
\end{eqnarray}
Since $II = 0$ by definition (\ref{cvj.265}) it suffices to show that
for $\lambda \in U({\cal N}^{(j)}, \rho_j) \setminus Z_{B_{p+1}}$ the 
following estimates hold:
\begin{eqnarray}
\label{cvj.280}
\| I \|_{1/4, c} &\leq&
\frac{1}{3} e^{-L_{p+1}^c}|a|^2, \\
\label{cvj.285}
\| III \|_{1/4, c} &\leq&
\frac{1}{3} e^{-L_{p+1}^c}|a|^2, \\
\label{cvj.290}
\| IV \|_{1/4, c} &\leq&
\frac{1}{3} e^{-L_{p+1}^c}|a|^2.
\end{eqnarray}

{\em Estimate on (I):}
Since the support of $z_p$ is contained
in $B_{L_p}(0)$, definition (\ref{cvj.260}) implies
\begin{eqnarray}
\label{cvj.295}
{\cal P}(z_p(\lambda), \lambda) - 
\chi_{\tilde{B}_{p+1}} {\cal P}(z_p(\lambda), \lambda)
= W(\varphi(a) + z_p(\lambda)) \chi_{ \{|m| \geq 3 L_p\} }.
\end{eqnarray}
Furthermore we have  
$\|\varphi(a) + z_p(\lambda)\|_{1, c} \leq r_{F, b}/2$ by
(\ref{cvj.167}).
We can therefore apply lemma \ref{Lenp.1}. Using in addition 
(\ref{cvj.65}) and 
(\ref{cvj.70}) we obtain
\begin{eqnarray}
\nonumber
\| I \|_{1/4, c} \leq D_W (4 e \nu D_N)^2 
e^{-\frac{9}{4} L_p^c} |a|^2 \leq \frac{1}{3}
e^{-L_{p+1}^c} |a|^2,
\end{eqnarray}
proving (\ref{cvj.280}).

{\em Estimate on (III):}
Using (\ref{cvj.265}), (\ref{cvj.245}), (\ref{cvj.250}), (\ref{cvj.75}),
(\ref{cvj.255}), induction hypothesis (v), (\ref{cvj.60}) and (\ref{cvj.65})
it follows that
\begin{eqnarray}
\nonumber
&&\| III \|_{1/4, c}
=
\sum_{|n| \geq L_{p+1}, m_1 \in B_{p+1}, |m_2| < 3 L_p}
|R(n,m_1)| |G_p (m_1, m_2)| 
|{\cal P}(z_p, \lambda)(m_2)| w_{1/4, c}(n) \\
\nonumber
&\leq&
\sum_{\ldots}
|R(n,m_1)| 
|G_p (m_1, m_2)| w_{1/4, c}(n-m_2)
|{\cal P}(z_p, \lambda)(m_2)| w_{1/4, c}(m_2) \\
\nonumber
&\leq&
\sum_{\ldots}
|R(n,m_1)| 
|G_p (m_1, m_2)| w_{1/2, c}(n-m_2) e^{-\frac{1}{4}|n-m_2|^c}
|{\cal P}(z_p, \lambda)(m_2)| w_{1/4, c}(m_2) \\
\nonumber
&\leq&
e^{-\frac{1}{4} (2 L_p)^c}
\sum_{\ldots}
|R(n,m_1)| w_{1/2, c}(n-m_1) |G_p (m_1, m_2)| w_{1/2, c}(m_1-m_2) 
|{\cal P}(z_p, \lambda)(m_2)| w_{1/4, c}(m_2)
\\
\nonumber
&\leq&
e^{-\frac{1}{4} (2 L_p)^c}
\|R\|_{1/2, c} \|G_p \|_{1/2, c} 
\| {\cal P}(z_p, \lambda) \|_{1/4, c} \\
\nonumber
&\leq&
e^{-\frac{1}{4} \left( \frac{2}{5} \right)^c L_{p+1}^c}
L_p^{A E_M + 1} e^{-L_p^c} |a|^2 \\
\nonumber
&\leq&
e^{-L_{p+1}^c 
\left(\frac{15}{16} 5^{-c} + \frac{1}{4} \left( \frac{2}{5} \right)^c \right)}
|a|^2
\leq 
e^{-\frac{9}{8} L_{p+1}^c} |a|^2.
\end{eqnarray}
Estimate (\ref{cvj.285}) now follows from (\ref{cvj.80}).

{\em Estimate on (IV):}
Estimate (\ref{cvj.167}) and its proof yield $\|\varphi+z_p\|_{1, c}$,
$\|\varphi+z_{p+1}\|_{1, c} \leq r_{F, b}/2$. Therefore we can apply
lemma \ref{Lenp.1} and obtain via (\ref{cvj.270}), (\ref{cvj.60})
and (\ref{cvj.65})
\begin{eqnarray}
\nonumber
\| IV \|_{1/4, c} &\leq&
D_W \|z_{p+1} - z_p\|_{1/4, c}^2 
\leq 
D_W \left( L_p^{A E_M + 1} e^{-L_p^c} |a|^2 \right)^2 \leq
D_W e^{-L_{p+1}^c \left( \frac{30}{16} 5^{-c}\right)} |a|^4 \\
\nonumber
&\leq&
D_W e^{-\frac{7}{4} L_{p+1}^c} (2 \rho_1)^2 |a|^2.
\end{eqnarray}
Estimate (\ref{cvj.290}) now follows from (\ref{cvj.70}).

\end{proof}

\section{Construction of the set of polynomials 
${\cal POL}_{\tilde{k}}^{(j+1)}$
for $j \geq 1$}
\label{cpol}

In this section we assume $j \geq 1$, that $({\cal IS})_l$ holds for 
$1 \leq l \leq j$ and that $v_{j+1}$ has been constructed, satisfying
$({\cal IS})(1)_{j+1}$. Our goal is to construct inverse
matrices $G_{C}^{(j+1)}(\theta, \lambda)$ of 
$T_{C}^{(j+1)}(\theta, \lambda)$ for suitable values of $(\theta, \lambda)$,
where $C \in {\cal C}^{(j+1)}$ and  
\begin{eqnarray}
\label{cpol.5}
T^{(j+1)}(\theta, \lambda) =
D(\theta, \omega) + DW (\varphi(a) + v_{j+1}(\lambda))
\equiv
D(\theta, \omega) + R(\lambda).
\end{eqnarray}
Note that $T^{(j+1)}$ as defined by (\ref{cpol.5}) is 
\underline{not} projected
onto some subset of ${\Bbb Z}^{\nu}$. In particular, the $P$-projection
has not been applied.

Unfortunately, we will not be able to ensure the invertibility
of $T_C^{(j+1)}(\theta, \lambda)$ for all relevant parameter values
of $(\theta, \lambda)$. Therefore we modify our goal.
We construct polynomials $p$ in the variable $\theta$ with coefficients
depending on $\lambda$, such that the zeros of the polynomials coincide
with those parameter values of $(\theta, \lambda)$ for which $T_C^{(j+1)}
(\theta, \lambda)$ is not invertible. Moreover, away from the zeros
of the polynomials we can estimate the norm of the inverse matrices
$G_{C}^{(j+1)}(\theta, \lambda)$ in terms of $|p(\theta, \lambda)|^{-1}$.
This statement is made precise in lemma \ref{lawp.1} below 
(see also remark \ref{Rova.2}) for a motivation of the polynomials). 
This section contains the heart of the multi-scale analysis. It follows
the ideas introduced by Bourgain in \cite{B1}, \cite{B4}, although we have
chosen a somewhat different presentation.
Before we formulate lemma \ref{lawp.1}, 
we first state the estimates of lemma \ref{hyp} which are used
in this section.

\noindent
{\bf Estimates of lemma \ref{hyp} used in section \ref{cpol}:}

\begin{eqnarray}
\label{conf.2}
8 e^5 \nu D_W D_N D_{\tau} N_0^{\tau+1} \rho_1 \leq 1 \\
\label{conf.3}
8 \rho_1 \leq e \nu D_N  \\
%\label{conf.5}
%\left(
%\{ \mu \} + (-2 \delta_1, 2 \delta_1)
%\right)_{\mu \in \{0, \pm 2 \omega_i^{(0)} \}}
%\mbox{ is a collection of pairwise disjoint subsets of } (-\pi, \pi) \\
\label{conf.7}
16 e \nu D_N \rho_1 \leq r_{F, b}\\
\label{conf.10}
2 \alpha_1 > d_V \delta_1 \\
\label{conf.15}
\frac{\nu \rho_j N_j}{4} + \frac{\delta_j}{100} \leq 
\nu N_0^{1-E_{\rho}}, \\
\label{conf.20}
\frac{(\nu+1)^2 \gamma^2}{4 \cosh^2 (\nu N_0^{1-E_{\rho}}/2)} - \alpha_1 > 
\frac{d_V \delta_1}{2} \\
\label{conf.25}
25 \nu N_j \rho_j \leq \delta_j \\
\label{conf.30}
2 \delta_1 <  \delta_V \\
\label{conf.37}
4 \rho_1 \leq \delta_1\\
\label{conf.35}
\frac{1}{8} N_1^c (5^c - 1) \geq \log 2 \\
\label{conf.40}
4 D_W e^{-\frac{1}{8} N_{l-1}^c} M_l \leq 1 \;\; \mbox{ for }
2 \leq l \leq j \\
\label{conf.45}
64 e \nu D_W D_N^2 \rho_1 \leq d_V \delta_1  \\
\label{conf.50}
2 D_N M_l e^{-\frac{1}{4l(l+1)}N_{l-1}^c}\leq 1  \;\; \mbox{ for }
2 \leq l \leq j \\
\label{conf.55}
M_1  \geq \frac{2}{d_V \delta_1}\\
\label{conf.60}
8 e \nu D_W D_N \rho_1 \leq 1 \\
\label{esf.5}
2 \delta_1 + 2 \rho_1 \leq \delta_V \\
\label{esf.7}
2(\delta_1 + \rho_1) \leq 1 \\
\label{esf.10}
36 B_0! e^6 \nu^2 D_W D_N^2 D_E(B_0) \leq \frac{1}{2} N_0^{E_{\rho}/2} \\
\label{esf.15}
e^{2 \nu} B_0! D_{\psi}(B_0) \left(\frac{8}{\rho_{l-1}}\right)^{B_0} \leq 
e^{\frac{1}{16} N_{l-2}^c}
\;\; \mbox{ for } l \geq 3 \\
\label{esf.20}
2 D_N e \leq \frac{1}{4} N_0^{E_{\rho}/2} \\
\label{esf.25}
\sum_{l=1}^{\infty} e^{- \frac{1}{16} N_{l}^c} \leq \frac{1}{4} \\
\label{esf.30}
D_W \left(
\sum_{p=1}^{B_0 + 1} \frac{1}{p!} p^{B_0}
\right)
\leq N_0^{E_{\rho}/2}
\\
\label{esf.35}
1 + \alpha_1 \leq N_0^{E_{\rho}} \\
\label{esf.40}
E_{\rho} \geq 3\\
\label{esf.45}
d_V \delta_0 < 2 \alpha_1
\\
\label{esf.50}
4 \nu N_1 \rho_1 < \delta_1
\\
\label{esf.55}
\frac{(\nu+1)^2 \gamma^2}{4 \cosh^2 \delta_1} - \alpha_1 > 
\frac{d_V \delta_0}{2}
\\
\label{esf.60}
4 \delta_1 \leq \frac{s}{(N_1 + 3)^{\tau}}
\\
\label{esf.65}
2^{\tau} \delta_0 < s
\\
\label{esf.70}
\delta_0 < 2 \delta_V 
\\
\label{esf.72}
\tau \leq E_{\delta}
\\
\label{esf.75}
N_0^{E_{\rho}} \leq M_0 \delta_0
\\
\label{esf.77}
10 \rho_l \leq \delta_l \;\; \mbox{ for } 2 \leq l \leq j
\\
\label{esf.80}
B_0!
\sum_{p=2}^{\infty}
D_W \left(\frac{8}{\rho_{p}}\right)^{B_0} e^{-\frac{1}{16} N_{p-1}^c} \leq 1
\\
\label{esf.85}
B_0! M_l e^{-\frac{1}{16} N_{l-1}^c} \leq \frac{1}{2}  \;\; \mbox{ for }
2 \leq l \leq j \\
\label{esf.90}
E_{\delta} + 1 \leq q \\
\label{esf.92}
8 e \nu D_V D_W D_N \rho_1^{1/4} \leq \frac{1}{2} \\
\label{esf.95}
2^{B_0} D_V (8 e \nu D_W D_N)^2 D_{1, 3}
\left(\max_{y \geq 0} (1+y)^{3 B_0 q} e^{-\frac{1}{4} |y-1|^c} \right)
M_0^{2 B_0 + 1} \rho_1 \leq \frac{1}{2}\\
\label{esf.99}
D_{1, 5}^4 \rho_1 \leq 1 \\
\label{esf6.300}
8 e \nu D_W D_N \leq N_0 \\
\label{esf6.310}
N_0 \rho_1 \leq \delta_1 \\
\label{esf6.320}
D_{1, 3} M_0 \max_{y \geq 0} \left[
(1 + y)^q e^{-\frac{3}{4} |y-1|^c} 
\right] \leq \delta_1^{-1} \\
\label{awp.1}
B_0 \geq B_1 + B_2\\
\label{awp.2}
400 \nu \sqrt{\rho_1} < 1\\
\label{awp.4}
2 \nu (B_1 + 1) \leq B_2\\
\label{awp.6}
B_1! 2000 \nu^2 D_1 2^{B_1} M_0^{2(B_0 + \nu)} \tilde{\rho}_j <
\left( \frac{\delta_j}{1200 \nu} \right)^{2 \nu} \\
%\label{awp.5}
%320 e \nu (12 \nu)^{2 \nu} 
%B_1! 2^B_1 D_1 M_0^{2 B_0 + 2 \nu} \tilde{\rho}_j \leq
%\delta_j^{2 \nu} \\
\label{awp.8}
(12 \nu)^{2 \nu (B_1 + 1)} \leq B_2!\\
\label{awp.10}
16 \tilde{\rho}_j \leq \rho_j \\
\label{awp.12}
2^{2 \nu} \delta_1 + 2^{4 \nu + 1} \rho_1 \leq \frac{1}{16 \nu}\\
\label{awp.14}
B_1 \geq 8 \nu \\
\label{awp.16}
\delta_1 \leq \frac{1}{4 \nu} \\
\label{awp.18}
(2 \nu)! D_N e^2 2 \left(
\frac{2(D_V+1)}{d_V} \right)^{2 \nu} \delta_1 \leq 1 \\
\label{dpo.5}
2000 (2 \nu + 1) \delta_1 \leq 1 \\
\label{dpoh.1}
\rho_j / \tilde{\rho}_j \in {\Bbb N} \\
\label{dpoh.2}
\tilde{\rho}_j / \rho_{j+1} \in {\Bbb N}
\end{eqnarray}

We now state the main result of this section

\begin{Lemma}
\label{lawp.1}
Let $j \geq 1$. Assume that the 
induction statements $({\cal IS})_l$, $1 \leq l \leq j$ and 
$({\cal IS})(1)_{j+1}$ are satisfied. 
Let $C \in {\cal C}^{(j+1)}$ and $\tilde{k} \in K^{(j+0.5)}$. 
Set $k' := \pi^{(j+0.5)}_j \tilde{k}$ and denote
\begin{eqnarray}
\label{cpol.10}
C' := 
\left\{
\begin{array}{ll}
C \cap B_{N_{j-1}}(0) & \mbox{ for } j \geq 2,\\
C \cap B_{1.5}(0) & \mbox{ for } j = 1.
\end{array}
\right.
\end{eqnarray}
Then $C' \in {\cal C}^{(j)}$. For every $\vartheta \in I_{C', k'}^{(j)}$ at least
one of the two statements $A$ or $B$ is true.

$\underline{\mbox{ Statement } A:}$
Denote by $J$ the connected component of $U(I_{C', k'}^{(j)}, \delta_j/100)$
which contains $\vartheta$.
There exists an analytic map $(\theta, \lambda) \mapsto
G_{C}^{(j+1)} (\theta, \lambda)$ defined on 
$(\theta, \lambda) \in  J \times 
U(\lambda_{\tilde{k}}^{(j+ 0.5)}, 2 \tilde{\rho}_j)$ 
such that
\begin{eqnarray}
\label{cpol.20}
\left\|
G_{C}^{(j+1)} (\theta, \lambda)
\right\|_{\sigma_{j+1}, c} \leq 4 D_N^2 M_j \;\;\;
\mbox{ for } 
(\theta, \lambda) \in J \times 
U(\lambda_{\tilde{k}}^{(j+ 0.5)}, 2 \tilde{\rho}_j).
\end{eqnarray}
Moreover, for 
$(\theta, \lambda) \in J \times 
U(\lambda_{\tilde{k}}^{(j+ 0.5)}, 2 \tilde{\rho}_j) 
\setminus \tilde{Z}_{C}$ the
matrix $G_{C}^{(j+1)} (\theta, \lambda)$ is the inverse of
$T^{(j+1)}_{C}(\theta, \lambda)$.

$\underline{\mbox{ Statement } B:}$
There exists a polynomial
$p_{C, \tilde{k}, \vartheta}$ 
with the following properties. There exist an integer
$1 \leq d \leq 2 \nu$ and analytic functions $b_0, \ldots, b_{d-1}$, 
such that for all $(\theta, \lambda) \in
U(\vartheta, \delta_j/400) \times U(\lambda_{\tilde{k}}^{(j + 0.5)},
2 \tilde{\rho}_j)$
the following holds.
\begin{itemize}
\item[(a)]
$p_{C, \tilde{k}, \vartheta}(\theta, \lambda) = 
\theta^d + \sum_{i<d} b_i(\lambda) \theta^i$.
\item[(b)]
$\left\|
G_{C}^{(j+1)} (\theta, \lambda)
\right\|_{\sigma_{j+1}, c} \leq 4 D_N^2 M_j
\left(
1 + 2 \frac{4 D_N^2 M_j}
{|p_{C, \tilde{k}, \vartheta}(\theta - \vartheta, \lambda)| \delta_1}
\right)$ if $p_{C, \tilde{k}, \vartheta}(\theta - \vartheta, \lambda) \neq 0$.
\item[(c)]
$|b_i(\lambda)| \leq \frac{1}{2d}$.
\item[(d)]
$|\partial^{\beta} b_i(\lambda)| \leq D_2 \left( M_0^{2 B_0 + 2 \nu} 
\right)^{2 |\beta|_1 -1}$ for $1 \leq |\beta|_1 \leq B_1$.
\item[(e)] $b_i(\lambda) \in {\Bbb R}$ for $\lambda \in 
B(\lambda_{\tilde{k}}^{(j + 0.5)}, 2 \tilde{\rho}_j)$ and $0 \leq i \leq
d-1$.
\item[(f)]
$\left|
\partial_{\theta} p_{C, \tilde{k}, \vartheta} (\theta, \lambda)
\right| \leq 1$ for $|\theta| \leq \delta_1/2$.
\end{itemize}
\end{Lemma}

The proof of lemma \ref{lawp.1} stretches over the next three subsections.
In subsection \ref{conf} we construct a function $f_{C, k'}$ defined for
$(\theta, \lambda) \in U(\tilde{I}^{(j)}_{C', k'}, \delta_j/100 ) \times 
U(B_{\rho_j}(\lambda_{k'}^{(j)}), \rho_j/4)$, such that the zeros of $f$ 
coincide with those values of $(\theta, \lambda)$ where 
$T_C^{(j+1)} (\theta, \lambda)$ is not invertible. In principle, one could 
define $f_{C, k'} = \det (T_C^{(j+1)})$. However, we will need estimates
on $f_{C, k'}$ and its derivatives which we cannot prove for the determinant
of a matrix of such large dimensions. Therefore we define $f_{C, k'}$ to
be the determinant of some reduced matrix $b$ which is at most of size
$2 \nu \times 2 \nu$. In subsection \ref{esf} we prove certain 
estimates for the function $f_{C, k'}$ which are sufficient 
to prove lemma \ref{lawp.1} by applying 
the Weierstrass preparation theorem (see subsection \ref{awp}).
In the final subsection \ref{dpo} we define sets of polynomials which
will be important for the completion of the induction step in section 
\ref{comp}. 

\subsection{Construction of the function $f_{C, k'}$}
\label{conf}

Throughout section \ref{conf} we suppose that $j \geq 1$ and that 
induction statements $({\cal IS})_l$, $1 \leq l \leq j$,
$({\cal IS})(1)_{j+1}$ are satisfied. We fix $k' \in K^{(j)}$,
$C \in {\cal C}^{(j+1)}$ and let $C'$ be defined as in 
(\ref{cpol.10}).

\subsubsection{Definition of the singular set $S$}
\label{conf1}
It follows from induction statement ${\cal IS}(3a)_{2 \leq l \leq j}$ that
$\tilde{I}^{(j)}_{C', k'} \subset \tilde{I}^{(1)}_{\tilde{C}, 1}$
for some $\tilde{C} \in {\cal C}^{(1)}$.
Induction statement ${\cal IS}(3a)_{j=1}$ then implies that
for all $\theta \in U(\tilde{I}^{(j)}_{C', k'}, \delta_j/100 )$
precisely one of the following $2 \nu + 1$ cases occurs.
\begin{eqnarray}
\nonumber
\mbox{ Case $I$: }&& 
\mbox{ dist}(\theta, 2 \pi {\Bbb Z} ) < \frac{101}{100} \delta_1 , \\
\label{conf1.5}
\mbox{ Case $II_i$: }&& 
\mbox{ dist}(\theta, \{2 \omega_i^{(0)} \} + 2 \pi {\Bbb Z} ) < 
\frac{101}{100} \delta_1, \;\;  1 \leq i \leq \nu, \\
\nonumber   
\mbox{ Case $III_i$: }&& 
\mbox{ dist}(\theta, \{ - 2 \omega_i^{(0)} \} + 2 \pi {\Bbb Z} ) < 
\frac{101}{100} \delta_1, \;\;  1 \leq i \leq \nu.
\end{eqnarray}
We now define the singular set $S$, depending on the above cases:
\begin{eqnarray}
\nonumber
\mbox{ In Case $I$ we set }&& S := C \cap {\cal S}, \\
\label{conf1.10}
\mbox{ In Case $II_i$ we set }&& S := C \cap \{ - e_i \}, \\
\nonumber
\mbox{ In Case $III_i$ we set }&& S := C \cap \{e_i\}.
\end{eqnarray}

\begin{Remark}
\label{Rconf1.1}
The definition of the set $S$ depends on $\theta$ which is omitted in our
notation. As the definition of $S$ is made in such a way that the diagonal
entries of the linearized operator $T^{(j+1)}(\theta, \lambda)(m,m)$
have small modulus for $m \in S$. More precisely, we can show
\begin{eqnarray}
\label{conf1.12}
|V(\theta, \omega)(m)| \leq 2 D_V \delta_1 \;\;\; \mbox{ for } m \in S.
\end{eqnarray}
Consider first
case $I$. Let $m = \mu e_i \in S$, $\mu \in \{ \pm 1 \}$. By (\ref{lop.50}),
(\ref{ova.520})
$V(\theta, \omega)(m) = V_i(\theta + \mu \omega_i)$. Since 
$|\omega_i - \omega_i^{(0)}| < 2 \rho_1$ we learn from (\ref{conf.37})
that
\begin{eqnarray}
\nonumber
\mbox{ dist}(\theta + \mu \omega_i, \{ \mu \omega_i^{(0)} \} + 2 \pi {\Bbb Z})
< 2 \delta_1.
\end{eqnarray}
Then (\ref{conf.30}) and proposition \ref{Plop.1} imply (\ref{conf1.12}).
Assume now that we are in case $II_i$ and $S \neq \emptyset$. Then
$V(\theta, \omega)(-e_i) = V_i(\theta - \omega_i)$. As above we conclude
\begin{eqnarray}
\nonumber
\mbox{ dist}(\theta - \omega_i, \{ \omega_i^{(0)} \} + 2 \pi {\Bbb Z})
< 2 \delta_1,
\end{eqnarray}
leading again to the estimate (\ref{conf1.12}). Case $III_i$ is similar
to case $II_i$.
\end{Remark}

We set
\begin{eqnarray}
\label{conf1.15}
\Gamma := C \setminus S.
\end{eqnarray}

\begin{Remark}
\label{Rconf1.2}
In order to see that the matrix $T^{(j+1)}(\theta, \lambda)$ is well defined
up to possible singularities on the diagonal we use 
\begin{eqnarray}
\label{cpol.7}
\| \varphi(a) + v_{j+1}(\lambda) \|_{1, c} 
\leq 8 e \nu D_N \rho_1 \leq 
\frac{r_{F, b}}{2} \;\;\; \mbox{ for } \lambda \in U({\cal N}^{(1)}, \rho_1).
\end{eqnarray} 
The proof of estimate (\ref{cpol.7}) is similar to the proof of 
(\ref{cvj.167}) and relies on the induction statement ${\cal IS}(1)_l$ for
$1 \leq l \leq j+1$ and on (\ref{conf.35}), (\ref{conf.2}), (\ref{conf.3}),
(\ref{conf.7}).   
\end{Remark}

\subsubsection{Decomposition of the nonsingular set $\Gamma$}
\label{conf2}
In addition to the assumptions and notation stated in the beginning of
section \ref{conf} we will assume throughout the present subsection
\ref{conf2} that 
$(\theta, \lambda) \in 
U(\tilde{I}^{(j)}_{C', k'}, \delta_j/100) \times
U(B_{\rho_j}(\lambda^{(j)}_{k'}),\rho_j/4)$ is arbitrary but fixed.
For $1 \leq l \leq j$, $n \in {\Bbb Z}^{\nu}$, we define
\begin{eqnarray}
\label{conf2.5}
k_l &:=& \pi^{(j)}_l k', \\
\label{conf2.10}
C_l (n) &:=& 
\left\{
\begin{array}{ll}
(C - \{ n \}) \cap B_{N_{l-1}}(0) & \mbox{ for } 2 \leq l \leq j,\\
(C - \{ n \}) \cap B_{1.5}(0) & \mbox{ for } l = 1.
\end{array}
\right.
\end{eqnarray}
Set
\begin{eqnarray}
\label{conf2.15}
E_0 := \{ m \in \Gamma: |<m, g>| \in \{1, 2, \ldots, \nu \} \}.
\end{eqnarray}
For each $m \in E_0$ we denote by $m' \in
{\Bbb Z}^{\nu}$ the (uniquely defined)
lattice point satisfying $<m', g> = 0$ and $m - m' \in {\cal S}$.
Observe that for every $m \in E_0$, $1 \leq l \leq j$  we have 
$C_l (m') \in {\cal C}^{(l)}$, since $C_l(m')$ is an interval 
(being an intersection of intervals) and $m-m' \in C_l(m')$ guarantees
that $C_l(m') \cap {\cal S}$ is not the empty set.
Next we define subsets of $\Gamma$:
\begin{eqnarray}
\nonumber
\Gamma_l &:=&
\left\{ m \in E_0: |m'| > 2 N_{l-1} \mbox{ and } 
\theta + <\omega, m'> \in 
U
\left(
\tilde{I}_{C_{l-1}(m'), k_{l-1}}^{(l-1)} \setminus
\tilde{I}_{C_{l}(m'), k_{l}}^{(l)}, \frac{\delta_l}{20}
\right)
\right\} \\
\label{conf2.20}
&& \quad \quad \quad \mbox{ for } 2 \leq l \leq j, \\
\label{conf2.25}
\Gamma_1 &:=&
\Gamma \setminus
\left(
\bigcup_{l=2}^{j} \Gamma_l
\right).
\end{eqnarray}
The remainder of this subsection is devoted to proving the following
properties for the sets $\Gamma_l$.
\begin{itemize}
\item[($\alpha$)]
For all $m \in \Gamma_1$: 
$|V(\theta, \omega)(m)| > \frac{d_V \delta_1}{2}$.
\item[($\beta$)]
$\Gamma_l \subset \{ m \in {\Bbb Z}^{\nu}: |m| \geq 2 N_{l-1} \}$
for $2 \leq l \leq j$.
\item[($\gamma$)]
For $2 \leq l \leq j$ and $m \in \Gamma_l$ we have 
\begin{eqnarray}
\label{conf2.30}
\left\| \left( T^{(j+1)}|_{C_l(m') + \{m'\}} \right)^{-1} (\theta, \lambda)
\right\|_{\sigma_l, c} \leq 2 D_N M_l.
\end{eqnarray}
\end{itemize}
For an interpretation of statement ($\gamma$) in the case 
$(\theta, \lambda) \in \tilde{Z}_{C_l(m') + \{m'\}}$ we refer the 
reader to remark \ref{Rcvj.1}.

{\em Proof of ($\alpha$):}
Let $m \in \Gamma_1$.

{\em Case 1: $m \in \Gamma_1 \setminus E_0$}. Suppose first that $<m, g> = 0$.
Then  $V(\theta, \omega)(m) = \alpha_1$ and the claim follows from 
(\ref{conf.10}).
Therefore we may assume that $|<m, g>| \geq \nu + 1$. Since
the imaginary part of $\theta + <\omega, m>$ is bounded by 
$\frac{\nu \rho_j N_j}{4} + \frac{\delta_j}{100} \leq 
\nu N_0^{1-E_{\rho}}$ (see (\ref{conf.15})) the claim follows from (\ref{conf.20}).

{\em Case 2: $m \in E_0 \cap \Gamma_1$ and $m' \neq 0$}.

Recall that the definition of ${\cal C}^{(j+1)}$ (see (\ref{npc.460}))
implies $\Gamma \subset C \subset B_{N_j}(0)$. 
Furthermore, we have assumed that
$m' \neq 0$. Since $|m - m'| = 1$ we obtain
\begin{eqnarray}
\label{conf2.35}
0 < |m'| \leq N_j.
\end{eqnarray}
Thus we can define $l \in \{1, ..., j\}$ to be the minimal integer such that
$|m'| \leq 2 N_l$ is satisfied. 
Recall that 
\begin{eqnarray}
\label{conf2.40}
\theta = \theta + <\omega, 0> \in
U(\tilde{I}^{(j)}_{C_j(0), k'}, \delta_j/100 ) \subset
U(\tilde{I}^{(l)}_{C_l(0), k_l}, \delta_l/10 ).
\end{eqnarray}
Applying induction statement $({\cal IS}) (3d)_l$ it follows that
\begin{eqnarray}
\label{conf2.45}
\theta + <\omega, m'> \in 
{\Bbb C} \setminus
U(\tilde{I}^{(l)}_{C_l(m'), k_l}, \delta_l/10 ) .
\end{eqnarray}
Note that (\ref{conf.25}) implies
\begin{eqnarray}
\label{conf2.50}
| \mbox{ Im}(\theta + <\omega, m'>)| < \frac{\delta_j}{100} +
\frac{\nu N_j \rho_j}{4} \leq \frac{\delta_j}{50}.
\end{eqnarray}
Let $\theta' \in {\Bbb R}$ satisfy $|\theta + <\omega, m'> - \theta'|
\leq \delta_j/50$. In view of (\ref{conf2.45}) it is clear that
$\theta' \in {\Bbb R} \setminus \tilde{I}^{(l)}_{C_l(m'), k_l}$.
Using 
\begin{eqnarray}
\label{conf2.55}
\tilde{I}_{C_s(m'), k_s}^{(s)} \subset 
\tilde{I}_{C_{s-1}(m'), k_{s-1}}^{(s-1)} \;\; \mbox{ for } 2 \leq s \leq j.
\end{eqnarray}
we obtain that
\begin{eqnarray}
\label{conf2.60}
{\Bbb R} \setminus \tilde{I}_{C_l(m'), k_l}^{(l)} =
{\Bbb R} \setminus \tilde{I}_{C_1(m'), k_1}^{(1)} \cup
\bigcup_{s=2}^l \left(
\tilde{I}_{C_{s-1}(m'), k_{s-1}}^{(s-1)} \setminus
\tilde{I}_{C_s(m'), k_s}^{(s)}
\right)
\end{eqnarray}
Since the definition of $\Gamma_1$ (\ref{conf2.25}) 
implies that $m \notin \Gamma_s$
for $2 \leq s \leq j$, it follows from (\ref{conf2.20}) that
\begin{eqnarray}
\label{conf2.65}
\theta' \in  {\Bbb R} \setminus \tilde{I}_{C_1(m'), k_1}^{(1)} 
\end{eqnarray}
and consequently
\begin{eqnarray}
\label{conf2.70}
\theta + <\omega, m'> 
\in U
\left( {\Bbb R} \setminus \tilde{I}_{C_1(m'), k_1}^{(1)}, 
\frac{\delta_j}{50}
\right).
\end{eqnarray}
Since $\delta_j \leq \delta_1$ induction
statement $({\cal IS}) (4)_1$ ($k \equiv k_1$, $C \equiv C_1(m')$)
implies
\begin{eqnarray}
\label{conf2.75}
|V(\theta + <\omega, m'>, \omega)(m-m')| > \frac{d_V \delta_1}{2}.
\end{eqnarray}
Since $<m', g> = 0$ we have $d_V \delta_1 / 2 <
|V(\theta + <\omega, m'>, \omega)(m-m')| = |V(\theta, \omega)(m)|$.

{\em Case 3: $m \in E_0 \cap \Gamma_1$ and $m' = 0$}.

Since $m -m' \in {\cal S}$ and $m \in C$ we have $m \in {\cal S} \cap C$.
Write $m = \mu e_i$, where $\mu \in \{ \pm 1 \}$ and $i \in
\{ 1, \ldots, \nu \}$. We will show by contradiction that
\begin{eqnarray}
\label{conf2.80}
|V(\theta, \omega)(m)| = 
|V_i(\mu \omega_i + \theta)| > \frac{d_V \delta_1}{2}.
\end{eqnarray}
Assume $|V_i(\mu \omega_i + \theta)| \leq d_V \delta_1/2$. Proposition
\ref{Plop.1} (a)
together with (\ref{conf.30})  
imply that either
(\ref{conf2.85}) or (\ref{conf2.90}) below must hold.
\begin{eqnarray}
\label{conf2.85}
\mbox{ dist}(\mu \omega_i + \theta, \{ \mu \omega_i^{(0)} \} + 2 \pi {\Bbb Z})
\leq \frac{\delta_1}{2}.
\end{eqnarray}
\begin{eqnarray}
\label{conf2.90}
\mbox{ dist}(\mu \omega_i + \theta, \{- \mu \omega_i^{(0)} \} + 2 \pi {\Bbb Z})
\leq \frac{\delta_1}{2}.
\end{eqnarray}
Using furthermore that
$|\omega_i - \omega_i^{(0)}| < 2 \rho_1 \leq \delta_1/2$ (see (\ref{conf.37}))
we learn that either 
(\ref{conf2.95}) or (\ref{conf2.100}) below is true.
\begin{eqnarray}
\label{conf2.95}
\mbox{ dist}(\theta, 2 \pi {\Bbb Z})
\leq \delta_1.
\end{eqnarray}
\begin{eqnarray}
\label{conf2.100}
\mbox{ dist}(\theta, \{- 2 \mu \omega_i^{(0)} \} + 2 \pi {\Bbb Z})
\leq \delta_1.
\end{eqnarray}
Suppose (\ref{conf2.95}) holds, then we are in Case $I$ (cf. (\ref{conf1.5}))
and by (\ref{conf1.10}) we have ${\cal S} \cap C = S$. Since 
$m \in  {\cal S} \cap C$ (see beginning of proof of case 3)
this contradicts $m \in \Gamma_1 \subset \Gamma \equiv
C \setminus S$. 
On the other hand, assuming the validity of (\ref{conf2.100}), we have
Case $II_i$ if $\mu = -1$ and Case $III_i$ if $\mu = 1$. In either of these
cases we learn from the definition (\ref{conf1.10}) that $m \in S$,
contradicting the assumption $m \in \Gamma_1$. Hence we have arrived
at a contradiction and (\ref{conf2.80}) holds, completing the proof of 
$(\alpha)$.

{\em Proof of ($\beta$):} Obvious from (\ref{conf2.20}) and $|m-m'| = 1$.

{\em Proof of ($\gamma$):} 
Let $m \in \Gamma_l$ for $2 \leq l \leq j$. Recall that $C_l(m') \in
{\cal C}^{(l)}$. Furthermore
\begin{eqnarray}
\label{conf2.105}
(\theta + <\omega, m'>, \lambda) \in
U
\left(
\tilde{I}_{C_{l-1}(m'), k_{l-1}}^{(l-1)} \setminus
\tilde{I}_{C_{l}(m'), k_{l}}^{(l)}, \frac{\delta_l}{10}
\right)
\times
U
\left( 
B_{\rho_l}(\lambda_{k_l}^{(l)}), \rho_l
\right).
\end{eqnarray}
The induction statement $({\cal IS}) (4)_l$ then implies the 
existence of $G^{(l)}_{C_l(m')}(\theta + <\omega, m'>, \lambda)$
which is the 
inverse of the matrix $T_{C_l(m')}^{(l)}(\theta + <\omega, m'>, \lambda)$
if $(\theta + <\omega, m'>, \lambda) \notin \tilde{Z}_{C_l(m')}$.
Using proposition \ref{Plop.3} and $<m', g>=0$  we obtain
\begin{eqnarray}
\label{conf2.110}
\left\|
\left(
G_{C_l(m')+\{m'\}}^{(l)}
\right)(\theta, \lambda)
\right\|_{\sigma_l, c} \leq M_l.
\end{eqnarray}

In view of proposition \ref{Pcr2.1} claim ($\gamma$) follows, 
if we can show that
\begin{eqnarray}
\label{conf2.115}
\left\|
\left(T^{(j+1)}(\theta, \lambda) - T^{(l)}(\theta, \lambda)
\right)|_{C_l(m')+\{m'\}}
\right\|_{\sigma_l, c} \leq \frac{1}{2M_l}.
\end{eqnarray}
Estimate (\ref{conf2.115}) is a consequence of (\ref{cpol.7}), 
lemma \ref{Lenp.1}, 
induction statements $({\cal IS}) (1d)_s$ for
$l < s \leq j+1$, (\ref{conf.35})
and (\ref{conf.40}):
\begin{eqnarray}
\nonumber&&
\left\|
\left(T^{(j+1)}(\theta, \lambda) - T^{(l)}(\theta, \lambda)
\right)|_{C_l(m')+\{m'\}}
\right\|_{\sigma_l, c} \\
\nonumber
&\leq& \| DW (\varphi(a) + v_{j+1}(\lambda)) -
DW (\varphi(a) + v_l(\lambda)) \|_{1, c} \\
\nonumber 
&\leq& D_W \| v_{j+1}(\lambda) - v_l(\lambda) \|_{1, c}
\\
\nonumber
&\leq& D_W
\sum_{s=l+1}^{j+1} \|v_{s}(\lambda) - v_{s-1}(\lambda) \|_{1, c}
\\
\nonumber
&\leq&
D_W
\sum_{s \geq l+1} e^{-\frac{1}{8} N_{s-2}^c} \\
\nonumber
& \leq &
2 D_W e^{-\frac{1}{8} N_{l-1}^c} \leq \frac{1}{2 M_l}.
\end{eqnarray}

\subsubsection{Construction of $G_{\Gamma}^{(j+1)}(\theta, \lambda)$}
\label{conf3}

\begin{Proposition}
\label{Pconf3.1}
There exists an analytic map $(\theta, \lambda) \mapsto
G_{\Gamma}^{(j+1)} (\theta, \lambda)$ which is defined for \\ 
$(\theta, \lambda) \in U(\tilde{I}^{(j)}_{C', k'}, \delta_j/100 ) \times 
U(B_{\rho_j}(\lambda_{k'}^{(j)}), \rho_j/4)$ such that
\begin{eqnarray}
\label{conf3.80}
\left\|
G_{\Gamma}^{(j+1)} (\theta, \lambda)
\right\|_{\sigma_{j+1}, c} \leq 4 D_N^2 M_j \;\;\;
\mbox{ for } 
(\theta, \lambda) \in U(\tilde{I}^{(j)}_{C', k'}, \delta_j/100 ) \times 
U(B_{\rho_j}(\lambda_{k'}^{(j)}), \rho_j/4).
\end{eqnarray}
Moreover, for 
$(\theta, \lambda) \in U(\tilde{I}^{(j)}_{C', k'}, \delta_j/100 ) \times 
U(B_{\rho_j}(\lambda_{k'}^{(j)}), \rho_j/4) \setminus \tilde{Z}_{\Gamma}$ the
matrix $G_{\Gamma}^{(j+1)} (\theta, \lambda)$ is the inverse of
$T^{(j+1)}_{\Gamma}(\theta, \lambda)$.
\end{Proposition}

\begin{Remark}
\label{Rconf3.1}
Recall from remark \ref{Rconf1.1} that the set $S$ and hence $\Gamma$
depend on the variable $\theta$, casting some doubt on the analyticity
of $G_{\Gamma}^{(j+1)}$.
However, it follows from the definition of the different cases 
(\ref{conf1.5})  and from induction statements $({\cal IS})(3 a)_{1 \leq
l \leq j}$ that $\Gamma$ is constant on connected components of
$U(\tilde{I}^{(j)}_{C', k'}, \delta_j/100 )$ resolving the issue.
\end{Remark}

\begin{proof}
In order to prove proposition \ref{Pconf3.1}
we apply the coupling lemma \ref{Lcl.1}. 
Assume first that $(\theta, \lambda) \notin \tilde{Z}_{\Gamma}$.
In the notation of lemma \ref{Lcl.1} 
\begin{eqnarray}
\nonumber
\Lambda &\equiv& \Gamma, \\
\nonumber
T&\equiv&T^{(j+1)}(\theta, \lambda) = D + R, \mbox{ where by
(\ref{cpol.5}) }\\
\label{conf3.100}
D(m, n) &=& V(\theta, \omega)(m) \delta_{m, n}, \\
\nonumber
R(m, n) &=& DW (\varphi(a) + v_{j+1}(\lambda))(m, n),\\
\nonumber
\sigma &\equiv& 1,\\
\nonumber
\tilde{\sigma} &\equiv& \sigma_{j+1}.
\end{eqnarray}
Next we define the
quantities $l_n$, $\mu_n$, $C_n$ and $U_n$ for each $n \in \Gamma$.
By definition (\ref{conf2.25}) we have $\Gamma = \bigcup_{s=1}^j \Gamma_s$. 
Therefore we can pick for every $n \in \Gamma$ a (not necessarily
unique) integer $s_n \in \{1, \ldots, j\}$ such that $n \in
\Gamma_{s_n}$. In the case $s_n = 1$ we set
\begin{eqnarray}
\label{conf3.105}
\begin{array}{llll}
l_n:=1,&\mu_n:=0,&C_n:=\frac{4 D_N}{d_V \delta_1},&U(n):=\{n\}.
\end{array}
\end{eqnarray}
If $2 \leq s_n \leq j$ it follows from definition (\ref{conf2.20}) that
$n \in E_0$
and therefore there exists an unique $n' \in {\Bbb Z}^{\nu}$ satisfying
$<n', g> = 0$ and $n-n' \in {\cal S}$. We set
\begin{eqnarray}
\label{conf3.110}
\begin{array}{llll}
l_n:=N_{s_n-1}-1,&\mu_n:=\sigma_{s_n}-\sigma_{j+1},&C_n:=2 D_N M_{s_n},&
U(n):=C_{s_n}(n')+\{n'\}.
\end{array}
\end{eqnarray}
From (\ref{cpol.7}) and lemma \ref{Lenp.1} we obtain
\begin{eqnarray}
\label{conf3.120}
\| R \|_{1, c} \leq 8 e \nu D_W D_N \rho_1.
\end{eqnarray}

Next we verify conditions (\ref{cl.20}) -- (\ref{cl.40}) of 
lemma \ref{Lcl.1}.

{\em Case $s_n=1$:} The set $U(n)$ contains the single point $n$ and
\begin{eqnarray}
\nonumber
| D(n, n) + R(n, n)| \geq 
\frac{d_V \delta_1}{2} - \| R \|_{1,c} \geq 
\frac{d_V \delta_1}{4}
\end{eqnarray}
by statement ($\alpha$) in section \ref{conf2}, 
(\ref{conf.45}) and (\ref{conf3.120}). Hence (\ref{cl.20}) is satisfied.
Condition (\ref{cl.30}) is obvious and condition (\ref{cl.40}) 
follows again from (\ref{conf.45}) and (\ref{conf3.120}).

{\em Case $2 \leq s_n \leq j$:}
Condition (\ref{cl.20}) is satisfied by statement ($\gamma$)
in section \ref{conf2} and (\ref{cl.30}) 
follows from the definition of $C_{s_n}(n')$
(cf. (\ref{conf2.10})) and $|n-n'|=1$. 
Note in particular that by definition (\ref{conf2.20}) we have
$|n'| > 2 N_{s_n - 1}$ which guarantees that $C_{s_n}(n') + \{ n' \}$
is not only a subset of $C$ but also a subset of $\Gamma$.
Finally, condition (\ref{cl.40})
is satisfied by (\ref{conf3.120}), (\ref{conf.60}) and (\ref{conf.50}) 
(see also (\ref{npc.56})).
Hence lemma \ref{Lcl.1} can be applied and we obtain
\begin{eqnarray}
\nonumber
\left\|
\left(T_{\Gamma}^{(j+1)}\right)^{-1} (\theta, \lambda)
\right\|_{\sigma_{j+1}, c} \leq 2 D_N 
\max \left( \frac{4 D_N}{d_V \delta_1}, 2 D_N M_j \right)
\end{eqnarray}
for 
$(\theta, \lambda) \in U(\tilde{I}^{(j)}_{C', k'}, \delta_j/100 ) \times 
U(B_{\rho_j}(\lambda_{k'}^{(j)}), \rho_j/4) \setminus \tilde{Z}_{\Gamma}$.
Keeping (\ref{conf.55}) in mind we obtain (\ref{conf3.80}). 
It remains to extend the inverse of $T_{\Gamma}^{(j+1)}$ across the set 
$\tilde{Z}_{\Gamma}$ to a function $G_{\Gamma}^{(j+1)}$ and to prove the
analyticity of $G_{\Gamma}^{(j+1)}$.
The argument, however, 
is very similar to
the one given in the proof of lemma \ref{Lcvj} (see below 
(\ref{cvj.255})) and we do not repeat it here.
\end{proof}

\subsubsection{Definition of $f_{C, k'}$}
\label{conf4}
Consider the following block decomposition of the matrix 
$T_C^{(j+1)}$,
\begin{eqnarray}
\label{conf4.5}
T_C^{(j+1)} (\theta, \lambda) =
\left(
\begin{array}{cc}
T_{\Gamma} & P_1 \\ P_2 & T_S
\end{array}
\right).
\end{eqnarray}
In the case $S \neq \emptyset$ we set
\begin{eqnarray}
\label{conf4.10}
b := T_S - P_2 G_{\Gamma} P_1.
\end{eqnarray}
\begin{Remark}
\label{Rconf4.1}
Remark \ref{Rconf1.1} shows in particular
that $T_S^{(j+1)}(\theta, \lambda)$ has no singularities on the diagonal.
Since $R$ and $G_{\Gamma}^{(j+1)}$ are analytic functions (see proposition
\ref{Pconf3.1}) it follows that (\ref{conf4.10}) defines
an analytic function $b$ on
$U(\tilde{I}^{(j)}_{C', k'}, \delta_j/100) \times
U(B_{\rho_j}(\lambda^{(j)}_{k'}), \rho_j/4)$.
\end{Remark}
We define
\begin{eqnarray}
\label{conf4.15}  
f_{C, k'}(\theta, \lambda) :=
\left\{
\begin{array}{ll}
\det b(\theta, \lambda)&\mbox{ if } S \neq \emptyset, \\
1 &\mbox{ if } S = \emptyset,
\end{array}
\right.
\end{eqnarray}
for $(\theta, \lambda) \in 
U(\tilde{I}^{(j)}_{C', k'}, \delta_j/100) \times
U(B_{\rho_j}(\lambda^{(j)}_{k'}), \rho_j/4)$.
Note that $T_C$ is invertible if $b$ is invertible. Indeed, in case
$b^{-1}$ exists, set
\begin{eqnarray}
\label{conf4.20}
G_C^{(j+1)} :=
\left(
\begin{array}{cc}
G_{\Gamma} + G_{\Gamma} P_1 b^{-1} P_2 G_{\Gamma} & 
- G_{\Gamma} P_1 b^{-1} \\ - b^{-1} P_2 G_{\Gamma} & b^{-1}
\end{array}
\right).
\end{eqnarray}
Then $G_C^{(j+1)}$ is an analytic function on 
$U(\tilde{I}^{(j)}_{C', k'}, \delta_j/100) \times
U(B_{\rho_j}(\lambda^{(j)}_{k'}), \rho_j/4)$ and
\begin{eqnarray}
\nonumber
G_C^{(j+1)} (\theta, \lambda) = \left( T_C^{(j+1)}
\right)^{-1} (\theta, \lambda) \; \mbox{ for } \; (\theta, \lambda) \notin
\tilde{Z}_C.
\end{eqnarray}
In particular, $T_C(\theta, \lambda)$ is invertible
except at the zeros of $f_{C, k'}$.

\subsection{Estimates on $f$}
\label{esf}

Before we state the result of this section we define a real number 
$\theta_{\iota} = \theta_{\iota}(\theta)$ 
depending on the cases defined in (\ref{conf1.5}).
\begin{eqnarray}
\label{esf.2}
\theta_{\iota} :=
\left\{
\begin{array}{ll}
0 & \mbox{ in case } \;I,\\
\left[ 2 \omega_i^{(0)} \right]& \mbox{ in case } \; II_i,\\
\left[- 2 \omega_i^{(0)} \right]& \mbox{ in case } \; III_i.
\end{array}
\right.
\end{eqnarray}

\begin{Lemma}
\label{Lesf.1}
Let $j \geq 1$. Suppose that the 
induction statements $({\cal IS})_l$, $1 \leq l \leq j$ and 
$({\cal IS})(1)_{j+1}$ are satisfied.
For $\theta \in
U(\tilde{I}^{(j)}_{C', k'}, \delta_j/100)$ we set 
\begin{eqnarray}
\label{esf.110}
d \equiv d(\theta) := \# S(\theta).
\end{eqnarray}  
Then we can express 
\begin{eqnarray}
\label{esf.200}
f_{C, k'}(\theta, \lambda) = \tilde{f}(\theta, \lambda)
\cdot \left(
q(\theta, \lambda) + r(\theta, \lambda) \right),
\end{eqnarray}
for $(\theta, \lambda) \in 
U(\tilde{I}^{(j)}_{C', k'}, \delta_j/100) \times
U(B_{\rho_j}(\lambda^{(j)}_{k'}), \rho_j/8)$,
where
\begin{eqnarray}
\label{esf.205}
|\tilde{f}(\theta, \lambda)| &\geq& d_V^{2 \nu} \\
\label{esf.210}
q(\theta, \lambda) &=& ([\theta] - \theta_{\iota})^d + 
\sum_{l < d} \tilde{a}_l (\lambda) ([\theta] - \theta_{\iota})^l,
\\
\label{esf.213}
\left|
\tilde{a}_l (\lambda) \right| &\leq& 2^{2 \nu + 1} \rho_1
\;\;\; \mbox{ for } 0 \leq l < d, \\
\label{esf.215}
\left|
\partial^{\beta} \tilde{a}_l (\lambda) \right| &\leq& 3^{3 \nu} 
\;\;\; \mbox{ for } 0 \leq l < d, \; |\beta|_1 \geq 1,\\
\label{esf.220}
\left|
\partial_{\theta}^{\beta} r(\theta, \lambda)
\right| &\leq& \sqrt{\rho_1},
\;\;\; \mbox{ for } \; 0 \leq \beta \leq B_0,\\
\label{esf.225}
\left|
\partial^{\beta} r(\theta, \lambda)
\right| &\leq&
D_1 M_0^{2 |\beta|_1 + 2 \nu}
\;\;\; \mbox{ for } \; 1 \leq \beta \leq B_0.
\end{eqnarray}
The functions $\tilde{f}$, $r$ and $a_l$ are analytic functions in
their respective domain of definition. 
Furthermore, $\tilde{a}_l(\lambda)$, $r(\theta, \lambda)$ are real if all the
components of the vector 
$(\theta, \lambda)$ are real.
\end{Lemma}

\begin{Remark}
\label{Resf.1}
Recall from remark \ref{Rconf3.1} that $d(\theta)$ is constant on connected 
components of $U(\tilde{I}^{(j)}_{C', k'}, \delta_j/100)$. In the case
that $d(\theta) = 0$, the statement of lemma \ref{Lesf.1} is trivially 
satisfied with $\tilde{f} \equiv 1$, $q \equiv 1$ and $r \equiv 0$.
Therefore we will assume throughout the remainder of section \ref{esf}
that $S \neq \emptyset$.
\end{Remark}

\subsubsection{Definition of $\tilde{f}$, $q$ and $r$}
\label{esf1}

First we define the decomposition 
\begin{eqnarray}
\label{esf1.5}
D_S = \tilde{D} \Lambda.
\end{eqnarray}
{\em Case $I$:}
Let $\mu \in \{ \pm 1\}$, $i \in \{1, \ldots, \nu \}$ and $\mu e_i \in S$.
By proposition \ref{Plop.1} (b) and (\ref{esf.5}) 
\begin{eqnarray}
\label{esf1.10}
D(\theta, \omega) (\mu e_i) =
V_i (\theta + \mu \omega_i) =
\tilde{V}_{i, \mu} ([\theta] + \mu \omega_i) 
([\theta] + \mu (\omega_i - \omega_i^{(0)})).
\end{eqnarray}
On $S$ we define the following diagonal matrices 
\begin{eqnarray}
\label{esf1.15}
\tilde{D}(\theta, \lambda) := \mbox{ diag}
\left( \tilde{V}_{i, \mu} ([\theta] + \mu \omega_i) \right), 
\\
\label{esf1.20}
\Lambda(\theta, \lambda) := \mbox{ diag}
\left( [\theta] + \mu (\omega_i - \omega_i^{(0)}) \right),
\end{eqnarray}
which clearly satisfy (\ref{esf1.5}).

{\em Case $II_i$:} Recall that $S = \{ - e_i \}$ since we assumed
$S \neq \emptyset$. It follows again from proposition \ref{Plop.1} (b), 
(\ref{esf.5}) that (\ref{esf1.5}) holds with 
\begin{eqnarray}
\label{esf1.25}
\tilde{D}(\theta, \lambda)(-e_i, -e_i) := 
\tilde{V}_{i, +} ([\theta - \omega_i] ), \\
\label{esf1.30}
\Lambda(\theta, \lambda)(-e_i, -e_i) :=
([\theta] - \theta_{\iota})  - (\omega_i - \omega_i^{(0)}).
\end{eqnarray}
 
{\em Case $III_i$:} Now $S = \{ e_i \}$. From
proposition \ref{Plop.1} (b) and 
(\ref{esf.5}) we again obtain (\ref{esf1.5}) with 
\begin{eqnarray}
\label{esf1.35}
\tilde{D}(\theta, \lambda)(e_i, e_i) := 
\tilde{V}_{i, -} ([\theta + \omega_i] ), \\
\label{esf1.40}
\Lambda(\theta, \lambda)(e_i, e_i) :=
([\theta] - \theta_{\iota}) + (\omega_i - \omega_i^{(0)}).
\end{eqnarray}

We express the matrix $b$ defined in (\ref{conf4.10}) in the following way.
\begin{eqnarray}
\nonumber
b &=& D_S + R_S - P_2 G_{\Gamma} P_1 =
\tilde{D}
\left(\Lambda + \tilde{D}^{-1}(R_S - P_2 G_{\Gamma} P_1) \right) \\
\label{esf1.45}
&=&
\tilde{D}
(\Lambda + \tilde{R} ), \;\; \mbox{ with } \\
\label{esf1.50}
\tilde{R} &:=& \tilde{D}^{-1}(R_S - P_2 G_{\Gamma} P_1).
\end{eqnarray}
It is clear that the following definitions yield
relation (\ref{esf.200}).
\begin{eqnarray}
\label{esf1.55}
\tilde{f} &:=& \det \tilde{D},
\\
\label{esf1.60}
q &:=& \det \Lambda,
\\ 
\label{esf1.65}
r &:=& \det (\Lambda + \tilde{R} ) - \det \Lambda.
\end{eqnarray}
The following estimates are consequences of the above definitions, 
(\ref{esf.7}) and proposition \ref{Plop.1} (b). For $m \in S$
\begin{eqnarray}
\label{esf1.70}
\left| \partial^{\beta} \Lambda (\theta, \lambda) (m, m) \right|
&\leq& 1 \;\;\; \mbox{ for } |\beta|_1 \geq 0, \\
\label{esf1.75}
\left| \partial^{\beta} \tilde{D}^{-1} (\theta, \lambda) (m, m) \right|
&\leq& D_V \;\;\; \mbox{ for } 0 \leq |\beta|_1 \leq B_0.
\end{eqnarray}

\subsubsection{Proof of properties (\ref{esf.205}) --
(\ref{esf.215}) for $\tilde{f}$ and $q$}
\label{esf2}

Estimate (\ref{esf.205}) follows from (\ref{lop.120}), 
(\ref{esf1.55}), the definition
of $\tilde{D}$ and the fact that $1 \leq d \leq 2\nu$. 
It is clear from the definition that the function $q$
can be written in form (\ref{esf.210}). Moreover, it is 
elementary to verify estimates (\ref{esf.213}), (\ref{esf.215}) for the
coefficients $\tilde{a}_l(\lambda)$ from $|\omega - \omega^{(0)}| <
2 \rho_1 \leq 1$ (see (\ref{esf.7})) and using the Leibniz rule.
The analyticity of $\tilde{f}$, $q$ and $r$ as well as the realness
of $\tilde{a}_l$ and $r$ follow immediately from $({\cal IS})(1b)_{j+1}$
and $({\cal IS})(1g)_{j+1}$. To complete the proof of lemma \ref{Lesf.1} we 
still need to verify (\ref{esf.220})
and (\ref{esf.225}) which will be done in section \ref{esf8}.

\subsubsection{Estimates on $R$}
\label{esf3}
\begin{Proposition}
\label{Pesf3.1}
Let $j \geq 1$. Suppose that the 
induction statements $({\cal IS})_l$, $1 \leq l \leq j$ and 
$({\cal IS})(1)_{j+1}$ are satisfied.
For $1 \leq |\beta|_1 \leq B_0$ and $\lambda \in
U(B_{\rho_j}(\lambda_{k'}^{(j)}), \rho_j/4)$ we obtain
\begin{eqnarray}
\label{esf3.2}
\|R(\lambda)\|_{1, c} &\leq& 8 e \nu D_W D_N \rho_1 \leq 
\frac{d_V \delta_1}{4} \leq 1, \\
\label{esf3.1}
\| \partial^{\beta} R (\lambda) \|_{1,c} &\leq& N_0^{E_{\rho} |\beta|_1}.
\end{eqnarray}
\end{Proposition}

\begin{proof}
Estimate (\ref{esf3.2}) follows from lemma \ref{Lenp.1}, (\ref{cpol.7}), and
(\ref{conf.45}). Assume
now $1 \leq |\beta|_1 \leq B_0$.
Recall that the matrix $R(\lambda)$ is given by 
\begin{eqnarray}
\label{esf3.5}
R(\lambda)(m, n) = \sum_{p \geq 2} p \alpha_p 
\left(
\varphi(a) + v_{j+1}(\lambda)
\right)^{*(p-1)}(m-n)
\end{eqnarray}
(see (\ref{cpol.5}), (\ref{enp.30})).
Using corollary \ref{Ccr1.1}
we obtain the following formula for the derivatives of $R$:
\begin{eqnarray}
\nonumber
(\partial^{\beta} R) y
&=& \sum_{p=1}^{|\beta|_1}
\!\!\!\!\!\!
\sum_{\scriptsize
{\scriptstyle
\begin{array}{c}
\beta_1 + \ldots + \beta_{p} = \beta \\ \beta_i \neq 0 \mbox{ for } 
1 \leq i \leq p
\end{array}
}
} \!\!\!
\frac{1}{p!} \frac{\beta!}{\beta_1! \cdots \beta_p!}
D^{p+1}W(\varphi+v_{j+1})
[\partial^{\beta_1}(\varphi+v_{j+1}), \ldots,
\partial^{\beta_{p}}(\varphi+v_{j+1}), y].\\
\label{esf3.10}
\end{eqnarray}
It follows from the definition of $\varphi$ (\ref{ova.120}),
from the induction statements $({\cal IS})(1e)_{2 \leq l \leq j+1}$ and from 
(\ref{esf.10}) -- (\ref{esf.25}) that for $1 \leq |\beta|_1 \leq B_0$,
$\lambda \in
U(B_{\rho_j}(\lambda_{k'}^{(j)}), \rho_j/4)$
\begin{eqnarray}
\nonumber
\| \partial^{\beta}(\varphi(a) + v_{j+1}(\lambda)) \|_{1, c}
&\leq&
\|\partial^{\beta} \varphi(a) \|_{1, c} +
\sum_{l=2}^{j+1}
\| \partial^{\beta}(v_l - v_{l-1})(\lambda) \|_{1, c}
\\ 
\nonumber
&\leq&
2 D_N e + \frac{1}{2} N_{0}^{E_{\rho}(|\beta|_1 - 0.5)} +
\sum_{l=3}^{j+1} e^{- \frac{1}{16} N_{l-2}^c} \\
\label{esf3.15}
&\leq&
N_{0}^{E_{\rho}(|\beta|_1 - 0.5)}.
\end{eqnarray} 
Estimate (\ref{esf3.1}) is then a consequence of 
(\ref{esf.10}), lemma \ref{Lenp.1}, proposition \ref{Pcr.2} 
and (\ref{esf.30}).
\end{proof}

\subsubsection{Estimates on $G_{\Gamma}$}
\label{esf4}

Recall the definitions of $D_{1,3}$ and $q$ in (\ref{npc.170}) and 
(\ref{npc.78}). 

\begin{Lemma}
\label{Lesf6.1}
Let $j \geq 1$. Suppose that the 
induction statements $({\cal IS})_l$, $1 \leq l \leq j$ and 
$({\cal IS})(1)_{j+1}$ are satisfied. Recall the 
definition of $G_{\Gamma}^{(j+1)}$ in proposition
\ref{Pconf3.1}
For all $0 \leq |\beta|_1 \leq B_0$, 
$(\theta, \lambda) \in U(\tilde{I}^{(j)}_{C', k'}, \delta_j/100 ) \times 
U(B_{\rho_j}(\lambda_{k'}^{(j)}), \rho_j/8)$ and $m$, $n \in \Gamma$
we have
\begin{eqnarray}
\label{esf6.50}
\left|
\partial^{\beta} G_{\Gamma}^{(j+1)}(\theta, \lambda) (m, n)
\right|
w_{1/4, c} (m-n)
\leq D_{1,3} M_0^{2 |\beta|_1 + 1} 
(1 +|n|)^{q(2 |\beta|_1 + 1)}.
\end{eqnarray}
\end{Lemma}

Recall from section \ref{conf3} that $G_{\Gamma}$ was constructed via
the coupling lemma \ref{Lcl.1} from the collection of local inverse
matrices $G_{U(n)}^{(j+1)}$, $n \in \Gamma$. Before we prove lemma 
\ref{Lesf6.1} we first provide estimates on the derivatives of the
matrices $G_{U(n)} \equiv G_{U(n)}^{(j+1)}$.

\begin{Proposition}
\label{Pesf4.1}
Suppose that the assumptions of lemma \ref{Lesf6.1} hold. Let
$n \in \Gamma_1$ (cf. (\ref{conf2.25})). Then
\begin{eqnarray}
\label{esf4.5}
\| \partial^{\beta} G_{U(n)} \|_{1, c}
\leq D_{1,1} M_0^{|\beta|_1 +1} (1 + |n|)^{(E_{\delta}+1)(|\beta|_1 +1)}
\end{eqnarray}
for $(\theta, \lambda) \in 
U(B_{\rho_j}(\lambda_{k'}^{(j)}), \rho_j/8)$,
$0 \leq |\beta|_1 \leq B_0$ (cf. (\ref{npc.150})).
\end{Proposition}

\begin{proof} {\em (Proposition \ref{Pesf4.1}).}

\noindent
Recall that for $n \in \Gamma_1$ we have defined $U(n) = \{ n \}$.
The entry of the corresponding matrix $T_{U(n)}$ can be written in
the form
\begin{eqnarray}
\label{esf4.10}
G_{U(n)} (n, n) =
(T_{U(n)})^{-1}(n,n) = 
\frac{f_n}{g_n},
\end{eqnarray}
where
\begin{eqnarray}
\label{esf4.15}
f_n(\theta, \lambda) &:=&
\left\{
\begin{array}{ll}
1&\mbox{ if } <n, g> = 0, \\
4 \sin^2 \frac{<\omega, n> + \theta}{2}
&\mbox{ if } <n, g> \neq 0, 
\end{array}
\right.
\\
\label{esf4.20}
g_n(\theta, \lambda) &:=&
\left\{
\begin{array}{ll}
\alpha_1 + R(\lambda)(n, n)&\mbox{ if } <n, g> = 0, \\
(\alpha_1 + R(\lambda)(n, n)) 4 \sin^2 \frac{<\omega, n> + \theta}{2}
- <n, g>^2 \gamma^2 &\mbox{ if } <n, g> \neq 0. 
\end{array}
\right.
\end{eqnarray}

{\em Estimates on $f_n$:} Using $2 \sin^2 x = 1 - \cos(2x)$ it is 
easy to verify that
\begin{eqnarray}
\label{esf4.25}
\left|
\partial^{\beta} f_n (\theta, \lambda)
\right| \leq
4 (1 +|n|)^{|\beta|_1}
\end{eqnarray}
for all multi-indices $\beta$.

{\em Estimates on derivatives of $g_n$:}
We will show that for $1 \leq |\beta|_1 \leq B_0$ 
\begin{eqnarray}
\label{esf4.30}
\left|
\partial^{\beta} g_n (\theta, \lambda)
\right| \leq
(2 N_0)^{E_{\rho} |\beta|_1} (1 + |n|)^{|\beta|_1}
\end{eqnarray}
Indeed, it follows from (\ref{esf3.2}), (\ref{esf3.1}),
(\ref{esf4.25}), (\ref{esf.35}),
(\ref{esf.40}) and the product rule that
\begin{eqnarray}
\nonumber
&&\left|
\partial^{\beta}
\left(
(\alpha_1 + R(\lambda)(n, n)) 4 \sin^2 \frac{<\omega, n> + \theta}{2}
\right) \right| \\
\nonumber
&\leq&
(\alpha _1 + 1) 4 (1 + |n|)^{|\beta|_1} + 
\sum_{0 < \alpha \leq \beta}
\left( \begin{array}{c} \beta \\ \alpha \end{array} \right)
N_0^{E_{\rho}|\alpha|_1} 4 (1 + |n|)^{|\beta - \alpha|_1} \\
\nonumber
&\leq& 2^{|\beta|_1 + 2} N_0^{E_{\rho}|\beta|_1} (1 + |n|)^{|\beta|_1}
\leq (2 N_0)^{E_{\rho}|\beta|_1}(1 + |n|)^{|\beta|_1}.
\end{eqnarray}
Note that this calculation contains a proof of (\ref{esf4.30}) in both
cases $<n, g> = 0$ and $<n, g> \neq 0$.

{\em Lower bound on $|g_n|$:}
We shall prove that
\begin{eqnarray}
\label{esf4.35}
|g_n(\theta, \lambda)| \geq \frac{\gamma^2}{4 \alpha_1} 
\frac{d_V \delta_0}{4} (1 + |n|)^{-E_{\delta}}
\end{eqnarray}

Suppose first that $<n, g>= 0$. In this case (\ref{esf4.35})
simply follows from (\ref{esf.45}), $\gamma^2 < \alpha_1$ (see
assumption A2) and from 
\begin{eqnarray}
\label{esf4.37}
|R(\lambda)(n, n)| \leq \frac{d_V \delta_1}{4} <  
\frac{d_V \delta_0}{4} \;\; \mbox{ ( see (\ref{esf3.2})) }.
\end{eqnarray}

Assume $<n ,g> \neq 0$. We introduce the following auxiliary notation
\begin{eqnarray}
\nonumber
x &:=& 4 \sin^2 \frac{<\omega, n> + \theta}{2},\\
\nonumber
y &:=& \alpha_1 + R(\lambda)(n, n), \\
\nonumber
\eta &:=& \frac{d_V \delta_0}{4} (1 +|n|)^{-E_{\delta}}
\end{eqnarray}
We will show below that
\begin{eqnarray}
\label{esf4.40}
|y - \frac{<n, g>^2 \gamma^2}{x}| =
|V(\theta, \omega)(n) + R(\lambda)(n, n)|  
\geq \eta.
\end{eqnarray}
Suppose (\ref{esf4.40}) holds. Then
\begin{eqnarray}
\nonumber
\frac{1}{|y|}
\left|
1 + \frac{<n, g>^2 \gamma^2}{yx - <n, g>^2 \gamma^2}
\right| \leq \frac{1}{\eta}.
\end{eqnarray}
Note that (\ref{esf.45}) and (\ref{esf4.37}) imply
$\eta < y < 2 \alpha_1$. Thus
\begin{eqnarray}
\nonumber
\frac{<n, g>^2 \gamma^2}{|yx - <n, g>^2 \gamma^2|} \leq \frac{y}{\eta} + 1, 
< 2 \frac{y}{\eta}
\end{eqnarray}
implying
\begin{eqnarray}
\nonumber
|g_n(\theta, \lambda)| = |yx - <n, g>^2 \gamma^2| >
\gamma^2 \frac{\eta}{2y}  > \frac{\gamma^2}{4 \alpha_1} \eta.
\end{eqnarray}
Hence the proof of (\ref{esf4.35}) is complete once we 
have verified (\ref{esf4.40}).

{\em Proof of (\ref{esf4.40}):}
Recall from statement ($\alpha$) in section \ref{conf2} that 
\begin{eqnarray}
\label{esf4.43}
|V(\theta, \omega)(n)| > \frac{d_V \delta_1}{2}
\end{eqnarray}
for all $n \in \Gamma_1$. Therefore (\ref{esf4.37}) implies
\begin{eqnarray}
\nonumber
|V(\theta, \omega)(n) + R(\lambda)(n, n)| \geq \frac{d_V \delta_1}{4}
\end{eqnarray}
Since $(1 + |n|)^{-E_{\delta}} 
\leq \delta_1$ for all $|n| \geq N_1$ (see (\ref{npc.430})), we may assume
that $|n| < N_1$ (in addition to the already present assumptions that
$n \in \Gamma_1$ and $<n, g> \neq 0$).
We distinguish the cases $|<n, g>| > \nu$ and $0 < |<n, g>| \leq \nu$.
In the first case we observe that the imaginary part of $<\omega, n>
+ \theta$ is bounded by 
$ \nu 2 \rho_1 N_1 + \delta_1/100 \leq \delta_1$ (use (\ref{esf.50})) 
and thus
(\ref{esf4.40}) follows from (\ref{esf.55}) and (\ref{esf4.37}).
Consider now $|<n, g>=i|$ for $i \in \{ 1, \ldots, \nu \}$. In this case
\begin{eqnarray}
\label{esf4.45}
V(\theta, \omega)(n) = V_i(\theta+<\omega, n>).
\end{eqnarray}
 
Observe that (\ref{esf.50}) together with $|n| \leq N_1$ implies
\begin{eqnarray}
\label{esf4.50}
|<\omega, n> - <\omega^{(0)}, n>| \leq \frac{\delta_1}{2}.
\end{eqnarray}
Next we prove
\begin{eqnarray}
\label{esf4.55}
|V(\theta, \omega)(n)| > \frac{d_V \delta_0}{2}(1 + |n|)^{-\tau}.
\end{eqnarray}
We need to distinguish the $2 \nu + 1$ cases $I$, $II_p$ and $III_p$ (see 
(\ref{conf1.5})).

{\em Case $I$:}
Since $n \in \Gamma_1 \subset \Gamma = C \setminus S$ we have 
$n \pm e_i \neq 0$.
Therefore (\ref{esf4.50}), (\ref{conf1.5}), assumption A3, (\ref{esf.60})
and (\ref{esf.65})
yield
\begin{eqnarray}
\nonumber
\mbox{ dist}\left(
\theta + <\omega, n>, \{ \pm \omega_i^{(0)} \} + 2 \pi {\Bbb Z}
\right)
&\geq&
\frac{s}{|n \mp e_i|^{\tau}} - \frac{1}{2} \delta_1 -
\frac{101}{100} \delta_1 \\
\nonumber
&\geq& 
\frac{s}{2 (1+|n|)^{\tau}} > \frac{\delta_0}{2 (1+|n|)^{\tau}}
\end{eqnarray}
Estimate (\ref{esf4.55}) then follows from 
(\ref{esf4.45}), proposition \ref{Plop.1} and (\ref{esf.70}).

{\em Case $II_p$:}
We first establish that $n + 2 e_p \mp e_i \neq 0$. In fact, assume
that $n + 2 e_p \mp e_i = 0$. Since $|<n, g>|=i$, and $p \in \{1, \ldots,
\nu \}$ this implies that $p = i$ and $n = -e_i$. This, however, contradicts
the fact that $n \notin S$ (since $n \in \Gamma_1$). 
Therefore (\ref{esf4.50}), (\ref{conf1.5}), assumption A3, (\ref{esf.60}),
$|n|+3 \leq 2(1 +|n|)$ and (\ref{esf.65}) imply
\begin{eqnarray}
\nonumber
&&\mbox{ dist}\left(
\theta + <\omega, n>, \{ \pm \omega_i^{(0)} \} + 2 \pi {\Bbb Z}
\right) \\
\nonumber
&\geq&
\mbox{ dist}\left(
\theta + <\omega^{(0)}, n \mp e_i + 2 e_p>, 
\{ 2 \omega_p^{(0)} \} + 2 \pi {\Bbb Z}
\right) - \frac{1}{2} \delta_1
\\
\nonumber
&\geq&
\frac{s}{|n \mp e_i +2 e_p|^{\tau}} - \frac{1}{2} \delta_1 -
\frac{101}{100} \delta_1 \\
\nonumber
&\geq& 
\frac{s}{2 (3+|n|)^{\tau}} > \frac{\delta_0}{2 (1+|n|)^{\tau}}.
\end{eqnarray}
Estimate (\ref{esf4.55}) then follows from 
(\ref{esf4.45}), proposition \ref{Plop.1} and (\ref{esf.70}).

{\em Case $III_p$:}
We first establish that $n - 2 e_p \mp e_i \neq 0$. Assume
that $n - 2 e_p \mp e_i = 0$. Since $|<n, g>|=i$, and $p \in \{1, \ldots,
\nu \}$ this implies that $p = i$ and $n = e_i$. This, however, 
contradicts $n \notin S$. 
Therefore (\ref{esf4.50}), (\ref{conf1.5}), assumption A3, (\ref{esf.60}),
$|n|+3 \leq 2(1 +|n|)$ and (\ref{esf.65}) imply
\begin{eqnarray}
\nonumber
&&\mbox{ dist}\left(
\theta + <\omega, n>, \{ \pm \omega_i^{(0)} \} + 2 \pi {\Bbb Z}
\right) \\
\nonumber
&\geq&
\mbox{ dist}\left(
\theta + <\omega^{(0)}, n \mp e_i - 2 e_p>, 
\{ - 2 \omega_p^{(0)} \} + 2 \pi {\Bbb Z}
\right) - \frac{1}{2} \delta_1
\\
\nonumber
&\geq&
\frac{s}{|n \mp e_i - 2 e_p|^{\tau}} - \frac{1}{2} \delta_1 -
\frac{101}{100} \delta_1 \\
\nonumber
&\geq& 
\frac{s}{2 (3+|n|)^{\tau}} > \frac{\delta_0}{2 (1+|n|)^{\tau}}.
\end{eqnarray}
Estimate (\ref{esf4.55}) again follows from 
(\ref{esf4.45}), proposition \ref{Plop.1} and (\ref{esf.70}).

Thus we have established (\ref{esf4.55}) in all cases. Together with
the estimate (\ref{esf4.43}) we obtain
\begin{eqnarray}
\nonumber
|V(\theta, \omega)(n)| > \max
\left(\frac{d_V \delta_1}{2},
\frac{d_V \delta_0}{2}(1 + |n|)^{-\tau}
\right).
\end{eqnarray}
On the other hand (\ref{esf4.37}) implies
\begin{eqnarray}
\nonumber
|R(\lambda)(n, n)| \leq \frac{1}{2}
\max
\left(\frac{d_V \delta_1}{2},
\frac{d_V \delta_0}{2}(1 + |n|)^{-\tau}
\right).
\end{eqnarray}
The last two inequalities together with 
(\ref{esf.72}) prove (\ref{esf4.40}), 
also in the case $|n| < N_1$. The proof of (\ref{esf4.40}) is complete.
We have so far established estimates (\ref{esf4.25}) -- (\ref{esf4.35})
on the functions $f_n$ and $g_n$. We use them to complete the proof of
proposition \ref{Pesf4.1}.

Applying corollary \ref{Ccr2.1} (b) for the derivatives of quotients yields
\begin{eqnarray}
\nonumber
\left|
\partial^{\beta} \left( \frac{f_n}{g_n} \right) \right|
\leq
\sum_{p=1}^{|\beta|_1 + 1}
\!\!\!\!\!\!
\sum_{\scriptsize
{\scriptstyle
\begin{array}{c}
\beta_1 + \ldots + \beta_p = \beta \\
\beta_r \neq 0 \mbox{ for } 2 \leq r \leq p
\end{array}
}
}
\frac{\beta!}{\beta_1 ! \ldots \beta_p!} |\partial^{\beta_1} f_n|
|g_n|^{-p} |\partial^{\beta_2} g_n| \ldots 
|\partial^{\beta_{p}} g_n| \; 
\end{eqnarray}
for $|\beta|_1 \geq 0$. Using proposition \ref{Pcr.2}
and estimates (\ref{esf4.25}), (\ref{esf4.30}) and
(\ref{esf4.35}) we obtain with (\ref{npc.150})
\begin{eqnarray}
\nonumber
\left|
\partial^{\beta} \left( \frac{f_n}{g_n} \right) \right|
&\leq&
\sum_{p=1}^{|\beta|_1 + 1} p^{|\beta|_1}
\left(\frac{16 \alpha_1}{d_V \gamma^2} \frac{(1+|n|)^{E_{\delta}}}{\delta_0}
\right)^{|\beta|_1 + 1}
(2 N_0)^{E_{\rho} |\beta|_1} 4 (1 + |n|)^{|\beta|_1} \\
\nonumber
&\leq&
D_N^{-1} D_{1,1} \left( \frac{N_0^{E_{\rho}}}{\delta_0} 
(1+|n|)^{E_{\delta}+1} \right)^
{|\beta|_1 + 1}
\end{eqnarray}
for $0 \leq |\beta|_1 \leq B_0$.
Proposition \ref{Pesf4.1} now follows from (\ref{esf.75}).

\end{proof}

Next we investigate $G^{(j+1)}_{U(n)}$ for $n \in \Gamma_l$, $2 \leq l
\leq j$. Recall that (\ref{conf2.20}), (\ref{conf2.15}) imply that
there exists a 
(unique) $n' \in {\Bbb Z}^{\nu}$ with $n-n' \in S$ and $<n', g> = 0$.
%We define 
%\begin{eqnarray}
%\label{esf5.40}
%q:=1+A(E_M + E_{\rho})
%\end{eqnarray}
\begin{Proposition}
\label{Pesf5.1}
Suppose that the assumptions of lemma \ref{Lesf6.1} hold. Let
$n \in \Gamma_l$, $2 \leq l \leq j$. Then
\begin{eqnarray}
\label{esf5.50}
\left\|
\partial^{\beta} G^{(j+1)}_{U(n)}(\theta, \lambda)
\right\|_{\sigma_l, c}
\leq D_{1,2} M_0^{|\beta|_1 + 1} (1+|n|)^{q (|\beta|_1 + 1)}
\end{eqnarray}
for 
$0 \leq |\beta|_1 \leq B_0$ and $(\theta, \lambda) \in
U(B_{\rho_j}(\lambda_{k'}^{(j)}), \rho_j/8)$.
\end{Proposition}

\begin{proof} {\em (Proposition \ref{Pesf5.1})}.
By continuity it suffices to prove estimate (\ref{esf5.50}) for  
$(\theta, \lambda) \in$\\
$U(B_{\rho_j}(\lambda_{k'}^{(j)}), \rho_j/8) \setminus \tilde{Z}_{U(n)}$.
Let $n \in \Gamma_l$, $2 \leq l \leq j$. Recall from (\ref{conf3.110}) that
$U(n) = C_l(n') + \{ n' \}$.
Applying Cauchy's integral formula to induction statement
${\cal (IS)}(4)_l$ we obtain
\begin{eqnarray}
\label{esf5.100}
\left\|
\partial^{\beta}
G_{C_l(n')}^{(l)} (\theta, \lambda) \right\|_{\sigma_l, c}
\leq \beta! M_l 
\left(
\frac{20}{\delta_l}
\right)^{\beta_{\theta}}
\left(
\frac{2}{\rho_l}
\right)^{|\beta_{\lambda}|_1}
\end{eqnarray}
for all $(\theta, \lambda) \in
U(\tilde{I}^{(l-1)}_{C_{l-1}(n'), k_{l-1}} \setminus 
\tilde{I}^{(l)}_{C_{l}(n'), k_{l}}, \delta_l/20) \times
U(B_{\rho_l}(\lambda_{k_l}^{(l)}), \rho_l/2)$ and
$|\beta|_1 \geq 0$.
Observe that $<n', g> = 0$. Proposition \ref{Plop.3} implies
\begin{eqnarray}
\label{esf5.105}
G_{C_l(n')+ \{n'\}}^{(l)} (\theta, \lambda) \cong
G_{C_l(n')}^{(l)} (\theta + <\omega, n'>, \lambda).
\end{eqnarray}
By the definition
of $\Gamma_l$ (see (\ref{conf2.20})) we have
$\theta + <\omega, n'> \in 
U(\tilde{I}^{(l-1)}_{C_{l-1}(n'), k_{l-1}} \setminus 
\tilde{I}^{(l)}_{C_{l}(n'), k_{l}}, \delta_l/20)$. Furthermore 
$\lambda \in U(B_{\rho_j}(\lambda_{k'}^{(j)}), \rho_j/8) \subset
U(B_{\rho_l}(\lambda_{k_l}^{(l)}), \rho_l/2)$. Applying the chain rule
we obtain from (\ref{esf5.100}) and (\ref{esf.77}) 
that for $|\beta|_1 \geq 0$
\begin{eqnarray}
\nonumber
\left\|
\partial^{\beta}
G_{C_l(n') + \{ n' \}}^{(l)} (\theta, \lambda) \right\|_{\sigma_l, c}
&\leq& |\beta|_1! M_l 
\left(
\frac{20}{\delta_l}
\right)^{\beta_{\theta}}
\left( |n'|\frac{20}{\delta_l} +
\frac{2}{\rho_l}
\right)^{|\beta_{\lambda}|_1}
\\
\label{esf5.110}
&\leq& |\beta|_1! M_l
\left[
(1 + |n|) \frac{4}{\rho_l}
\right]^{|\beta|_1}
.
\end{eqnarray}
Finally, in view of proposition \ref{Pcr2.1} 
we need to estimate the derivatives of
$X := T^{(j+1)} - T^{(l)}$. To this end set for $l \leq p \leq j$
\begin{eqnarray}
\label{esf5.115}
X_p &:=&  T^{(p+1)} - T^{(p)} = DW(\varphi + v_{p+1}) -  
DW(\varphi + v_{p})\\
\nonumber
&=&
\left(
\int_0^1 D^2W(\varphi +v_p + t(v_{p+1}-v_p)) dt
\right)
(v_{p+1}-v_p).
\end{eqnarray}
Note that $X_p$ does not depend on $\theta$.
It follows from induction statement $({\cal IS})(1b), (1d)$, $p \geq 2$ 
and lemma \ref{Lenp.1}
that $X_p$ is analytic on $U({\cal N}^{(p)}, \rho_p/4)$ and
\begin{eqnarray}
\nonumber
\| X_p (\lambda) \|_{1, c} \leq D_W e^{-\frac{1}{8} N_{p-1}^c}.
\end{eqnarray}
Applying Cauchy's integral formula and (\ref{esf.80})
we conclude for $\lambda \in 
U(B_{\rho_j}(\lambda_{k'}^{(j)}), \rho_j/8) \subset$ \\ 
$U(B_{\rho_p}(\lambda_{k_p}^{(p)}), \rho_p/8)$, $l \leq p \leq j$, and 
$0 \leq |\beta|_1 \leq B_0$ that
\begin{eqnarray}
\nonumber
\| \partial^{\beta} X (\lambda) \|_{1, c}
\leq  \sum_{p=l}^{j} \| \partial^{\beta} X_p (\lambda) \|_{1, c}
\leq e^{-\frac{1}{16} N_{l-1}^c}
\end{eqnarray}
Estimate (\ref{esf.85}) allows the application of proposition \ref{Pcr2.1},
yielding
\begin{eqnarray}
\label{esf5.120}
\left\|
\partial^{\beta} G^{(j+1)}_{C_l(n') + \{ n' \}}(\theta, \lambda)
\right\|_{\sigma_l, c} \leq
2 B_0! M_l  
\left[
(1 + |n|) \frac{8}{\rho_l}
\right]^{|\beta|_1}
\left(
\sum_{p=1}^{|\beta|_1 + 1} p^{|\beta|_1} D_N^p
\right)
\end{eqnarray}
for $0 \leq |\beta|_1 \leq B_0$.
Finally, statement ($\beta$) in section \ref{conf2} implies
$|n| \geq 2 N_{l-1}$. Thus by (\ref{npc.410}), (\ref{npc.420}),
(\ref{npc.440})
\begin{eqnarray}
\nonumber
\frac{1}{\rho_l} &=& N_l^{E_{\rho}} \leq (1+|n|)^{AE_{\rho}} \\
\nonumber
M_l &=& N_l^{E_M} \leq (1+|n|)^{AE_M}.
\end{eqnarray}
The definition of $q$ (\ref{npc.78}), proposition \ref{Pesf5.1} now follows
from (\ref{esf5.120}) and from the choice of the constant $D_{1,2}$ in
(\ref{npc.160}).
\end{proof}

Before proving lemma \ref{Lesf6.1} we need one more proposition.
\begin{Proposition}
\label{Pesf6.1}
Let $w$, $y> 0$, $t \geq 1$. Then 
\begin{eqnarray}
\nonumber
(t + y)^w e^{-\frac{1}{4} y^c} \leq
t^w \left(
\max_{y \geq 0} (1 + y)^w e^{-\frac{1}{4} y^c} \right).
\end{eqnarray}
\end{Proposition}
The proof of this proposition is trivial and is based on the 
simple observation that $(t+y)^w \leq t^w (1+y)^w$ for $t \geq 1$,
$y$, $w > 0$.
 
\begin{proof} ({\em Lemma \ref{Lesf6.1}}).

We have used the coupling lemma \ref{Lcl.1} to construct $G_{\Gamma}$ from
the collection of local inverse matrices $G^{(j+1)}_{U(n)}$, $n \in \Gamma$.
By (\ref{cl.60}), (\ref{cl.90}), (\ref{cl.100}) 
\begin{eqnarray}
\label{esf6.95}
G_{\Gamma}^{(j+1)} = \tilde{G} (1 + Y)^{-1} 
\end{eqnarray}
where
for $m$, $n \in \Gamma$
\begin{eqnarray}
\nonumber
\tilde{G}(m, n) =
\left\{
\begin{array}{ll}
G_{U(n)}^{(j+1)}(m, n)&\mbox{ if } m \in U(n), \\
0&\mbox{ else },
\end{array}
\right. 
\\
\nonumber
Y(m, n) =
\left\{
\begin{array}{ll}
\sum_{p \in U(n)} R(m, p) \tilde{G}(p, n) &\mbox{ if } 
m \in \Gamma \setminus U(n), \\
0&\mbox{ else }.
\end{array}
\right. 
\end{eqnarray}
It follows from propositions \ref{Pesf4.1}, \ref{Pesf5.1} 
and (\ref{esf.90}) that 
for $D := \max (D_{1,1}, D_{1,2})$ and $n \in \Gamma$
\begin{eqnarray}
\label{esf6.100}
\sum_{m \in \Gamma}
\left|
\partial^{\beta} \tilde{G} (m, n)
\right|
w_{1/2, c} (m-n)
\leq D M_0^{|\beta|_1+1} (1 + |n|)^{q(|\beta|_1+1)}
\end{eqnarray}
for $0 \leq |\beta|_1 \leq B_0$.
Next we estimate the derivatives of the matrix $Y$ up to order 
$0 \leq |\beta|_1 \leq B_0$. To that end recall 
$\| \partial^{\beta} R \|_{1, c} \leq N_0^{E_{\rho} |\beta|_1} \leq 
M_0^{|\beta|_1}$ from proposition \ref{Pesf3.1} and (\ref{esf.75}).
Using in addition (\ref{esf6.100}) we obtain for $n \in \Gamma$, 
$m \in \Gamma \setminus U(n)$
\begin{eqnarray}
\nonumber
&&\left|
\partial^{\beta} Y (m, n)
\right|
w_{1/2, c} (m - n) \\
\nonumber
&\leq&
\sum_{\alpha \leq \beta} 
\sum_{p \in U(n)}
\left( \begin{array}{c} \beta \\ \alpha \end{array} \right)
\left|
\partial^{\beta - \alpha} R (m, p)
\right|
w_{1/2, c} (m - p)
\left|
\partial^{\alpha} \tilde{G} (p, n)
\right|
w_{1/2, c} (p - n)
\\
\nonumber
&\leq&
\sum_{\alpha \leq \beta} 
\sum_{p \in U(n)}
\left( \begin{array}{c} \beta \\ \alpha \end{array} \right)
M_0^{|\beta - \alpha|_1} 
\left|
\partial^{\alpha} \tilde{G} (p, n)
\right|
w_{1/2, c} (p - n)
\\
\nonumber
&\leq&
\sum_{\alpha \leq \beta} 
\left( \begin{array}{c} \beta \\ \alpha \end{array} \right)
D M_0^{|\alpha|_1+1} (1 + |n|)^{q(|\alpha|_1+1)}
M_0^{|\beta - \alpha|_1}
\\
\label{esf6.120}
&\leq&
D (2 M_0)^{|\beta|_1+1} (1 + |n|)^{q(|\beta|_1+1)}.
\end{eqnarray}

Differentiating (\ref{esf6.95}) (see corollary \ref{Ccr2.1} (b)) 
we obtain
\begin{eqnarray}
\label{esf6.135}
\partial^{\beta} G_{\Gamma}^{(j+1)} =
\sum_{p=1}^{|\beta|_1 + 1} (-1)^{p-1}
\!\!\!\!\!\!\!\!\!\!
\sum_{\scriptsize
{\scriptstyle
\begin{array}{c}
\beta_1 + \ldots + \beta_p = \beta \\
\beta_r \neq 0 \mbox{ for } 2 \leq r \leq p
\end{array}
}
}
\!\!\!
\frac{\beta!}{\beta_1 ! \ldots \beta_p!} (\partial^{\beta_1} \tilde{G})
(1+Y)^{-1} (\partial^{\beta_2} Y) \ldots 
(1+Y)^{-1}.
\end{eqnarray}
Furthermore, (\ref{esf6.100}), (\ref{esf6.120}), 
$\| (1 +Y)^{-1} \|_{1/2, c} \leq 2 D_N$ (cf. (\ref{cl.70}))
and lemma \ref{Lwf.1} imply for
$1 \leq p \leq |\beta|_1+1$, $\beta_1 + \ldots + \beta_p = \beta$,
$\beta_r \neq 0$ for  $2 \leq r \leq p$
\begin{eqnarray}
\nonumber
&&\left|
\left[ (\partial^{\beta_1} \tilde{G})
(1+Y)^{-1} (\partial^{\beta_2} Y) \ldots (\partial^{\beta_{p}} Y)
(1+Y)^{-1} 
\right] (m, n)
\right|
w_{1/4, c} (m - n)
\\
\nonumber
&\leq&
w_{1/4, c} (m - n) \!\!\!\!
\sum_{y_1, \ldots, y_{2p-1} \in \Gamma}
(2 D_N)^p D^p (2 M_0)^{|\beta|_1 + p}
\frac{(1 + |y_1|)^{q(|\beta_1|_1+1)} \ldots 
(1 + |y_{2p-1}|)^{q(|\beta_p|_1+1)}}
{
w_{1/2, c}(m - y_1) \ldots w_{1/2, c}(y_{2p-1} - n)
}
\\
\nonumber
&\leq&
(2 D_N D)^p (2 M_0)^{|\beta|_1 + p}
\sum_{y_1, \ldots, y_{2p-1} \in \Gamma}
\left(\frac{w_{1/4, c} (m - n)}
{
w_{1/4, c}(m - y_1) \ldots w_{1/4, c}(y_{2p-1} - n)
}
\right.
\\
\nonumber
&&
\left. \quad \quad \quad
\times \frac{(1 + |y_1|)^{q(|\beta_1|_1+1)} \ldots 
(1 + |y_{2p-1}|)^{q(|\beta_p|_1+1)}}
{e^{\frac{1}{4}(|m - y_1|^c + \ldots + |y_{2p-1} - n|^c)}}
\right)
\\
\nonumber
&\leq&
(2 D_N D)^{|\beta|_1 + 1} (2 M_0)^{2 |\beta|_1 + 1}
\left(
\max_{y_1, \ldots, y_{2p-1} \in \Gamma}
\frac{(1 + |y_1|)^{q(|\beta_1|_1+1)} \ldots 
(1 + |y_{2p-1}|)^{q(|\beta_p|_1+1)}}
{e^{\frac{1}{4}(|m - y_1|^c + \ldots + |y_{2p-1} - n|^c)}}
\right) \\
\label{esf6.140}
\end{eqnarray}
Denote by $MAX$ the maximum which appears in the previous line
(\ref{esf6.140}). We claim
\begin{eqnarray}
\label{esf6.150}
MAX \leq (1 + |n|)^{q(2 |\beta|_1 + 1)} \left(
\max_{y \geq 0} (1+y)^{q(2 |\beta|_1 + 1)} e^{-\frac{1}{4} y^c} \right).
\end{eqnarray}
Indeed, since $\Gamma$ is finite, the maximum is achieved, say at
lattice points
$y_1^*, \ldots, y_{2p-1}^* \in \Gamma$. Set 
$\Delta := \max \{ |y_1^*|, \ldots, |y_{2p-1}^*| \}$. 
In case $\Delta \leq |n|$ estimate (\ref{esf6.150}) is obvious.
If $\Delta > |n|$ we see that  $MAX$ is bounded by 
\begin{eqnarray}
\nonumber
(1 + \Delta)^{q(2 |\beta|_1 + 1)}
e^{-\frac{1}{4}(\Delta - |n|)^c}.
\end{eqnarray}
Using proposition \ref{Pesf6.1} with $w \equiv q(2 |\beta|_1 + 1)$,
$y \equiv \Delta - |n|$ and $t \equiv 1 + |n|$ we conclude
that (\ref{esf6.150}) is satisfied.
Using (\ref{esf6.135}), (\ref{esf6.140}), (\ref{esf6.150}),
proposition \ref{Pcr.2} and
(\ref{npc.170}) we have finally proved lemma \ref{Lesf6.1}.
\end{proof}

\begin{Corollary}
\label{Cesf6.1}
Let $j \geq 1$. Suppose that the 
induction statements $({\cal IS})_l$, $1 \leq l \leq j$ and 
$({\cal IS})(1)_{j+1}$ are satisfied. Then
the matrix valued function $b$ defined in (\ref{conf4.10}) is 
analytic on 
$U(\tilde{I}^{(j)}_{C', k'}, \delta_j/100) \times
U(B_{\rho_j}(\lambda_{k'}^{(j)}), \rho_j/8)$ and the entries of $b$ are
uniformly bounded by
\begin{eqnarray}
\label{esf6.400}
|b(\theta, \lambda)(m, n)| \leq 2(D_V + 1)\delta_1
\;\;\; \mbox{ for } m, n \in S.
\end{eqnarray} 
\end{Corollary}

\begin{proof}
The analyticity of $b$ was already stated in remark \ref{Rconf4.1}.
By definition 
\begin{eqnarray}
\nonumber
b(m, n) &=&
V(\theta, \omega)(m) \delta_{m, n} +
R(m, n) -
\sum_{y_1, y_2 \in \Gamma} R(m, y_1) G_{\Gamma}^{(j+1)}(y_1, y_2) R(y_2, n)\\
\label{esf6.450}
&=& I + II + III.
\end{eqnarray}
It follows from remark
\ref{Rconf1.1} that $| I | \leq 2 D_V \delta_1$. Furthermore estimate
(\ref{esf3.2}) together with conditions (\ref{esf6.300}) and 
(\ref{esf6.310}) imply $|II| \leq \delta_1$. Finally, by lemma \ref{Lesf6.1},
proposition \ref{Pesf3.1}, lemma \ref{Lwf.1}
and (\ref{esf6.300}) -- (\ref{esf6.320})
\begin{eqnarray}
\nonumber
|III| &\leq&
\sum_{y_1, y_2 \in \Gamma} |R(m, y_1) 
G_{\Gamma}^{(j+1)}(y_1, y_2) R(y_2, n)| \\
&\leq&
\delta_1^2 
\sum_{y_1, y_2 \in \Gamma} D_{1, 3} M_0 
\frac{(1 + |y_2|)^q e^{- \frac{3}{4} |y_2 - n|^c}}
{w_{\frac{1}{4}, c}(m - y_1) w_{\frac{1}{4}, c}(y_1 - y_2) 
w_{\frac{1}{4}, c}(y_2 - n)} \\
&\leq&
\delta_1^2 D_{1, 3} M_0 \max_{y_2 \in \Gamma} \left[
(1 + |y_2|)^q e^{-\frac{3}{4} |y_2 - n|^c} 
\right] \leq \delta_1.
\end{eqnarray}
\end{proof}

\subsubsection{Estimates on $\tilde{R}$}
\label{esf7}
Recall the definition of $\tilde{R}$ in (\ref{esf1.50}).
\begin{Proposition}
\label{Pesf7.1}
Let $j \geq 1$. Suppose that the 
induction statements $({\cal IS})_l$, $1 \leq l \leq j$ and 
$({\cal IS})(1)_{j+1}$ are satisfied. Then
for all $0 \leq |\beta|_1 \leq B_0$, 
$(\theta, \lambda) \in U(\tilde{I}^{(j)}_{C', k'}, \delta_j/100 ) \times 
U(B_{\rho_j}(\lambda_{k'}^{(j)}), \rho_j/8)$ and $n$, $m \in S$
\begin{eqnarray}
\label{esf7.50}
\left|
\partial^{\beta}
\tilde{R} (n, m)
\right|
&\leq&
D_{1,4}
M_0^{2 |\beta|_1 + 1}, \\
\label{esf7.55}
\left|
\partial_{\theta}^{\beta}
\tilde{R} (n, m)
\right|
&\leq&
\rho_1^{3/4}.
\end{eqnarray}
\end{Proposition}

\begin{proof}
The product rule implies
\begin{eqnarray}
\label{esf7.100}
\partial^{\beta} \tilde{R} (m, n)
&=&
\sum_{\alpha \leq \beta}
\left( \begin{array}{c} \beta \\ \alpha \end{array} \right)
(\partial^{\alpha} \tilde{D}^{-1})(m, m) (\partial^{\beta - \alpha} R)(m, n) \\
\nonumber
&-& \!\!\!\!\!\!\!\!\!\!\!\!\!\!\!\!
\sum_{\scriptsize{\scriptstyle
\begin{array}{c}
\alpha_1 + \ldots + \alpha_4 = \beta \\
y_1, y_2 \in \Gamma
\end{array}
}
}
\!\!\!\!\!\!\frac{\beta!}{\alpha_1! \ldots \alpha_4!}
(\partial^{\alpha_1} \tilde{D}^{-1})(m, m)
(\partial^{\alpha_2} R)(m, y_1)
(\partial^{\alpha_3} G_{\Gamma}^{(j+1)}) (y_1, y_2) 
(\partial^{\alpha_4} R ) (y_2, n) \\
\nonumber
&=& S_1 - S_2
\end{eqnarray}
It follows from (\ref{esf1.75}), proposition \ref{Pesf3.1} and (\ref{esf.75})
that
\begin{eqnarray}
\label{esf7.110}
|S_1| \leq 2 ^{|\beta|_1} D_V N_0^{E_{\rho} |\beta|_1} \leq 
2 ^{|\beta|_1} D_V M_0^{|\beta|_1}.
\end{eqnarray}
Using in addition (\ref{esf6.50}), lemma \ref{Lwf.1} and proposition 
\ref{Pcr.2}
\begin{eqnarray}
\nonumber
|S_2| &\leq& 
D_V  \!\!\!\!
\sum_{\alpha_1 + \ldots + \alpha_4 = \beta}
\frac{\beta!}{\alpha_1! \ldots \alpha_4!}
M_0^{|\alpha_2|_1 + |\alpha_4|_1}
\sum_{y_1, y_2 \in \Gamma}
D_{1, 3} M_0^{2 |\alpha_3|_1 + 1} 
\frac{(1 + |y_2|)^{q(2 |\alpha_3|_1 + 1)} e^{-\frac{1}{4} |y_2 - n|^c}}
{w_{1/4, c}(m-y_1) \ldots w_{1/4,c}(y_2-n)} \\
\label{esf7.120}
&\leq&
D_V D_{1, 3} 4^{|\beta|_1} M_0^{2 |\beta|_1 + 1} 
\left(
\max_{y \geq 0}
(1 + y)^{q(2 |\beta|_1 + 1)} e^{-\frac{1}{4} |y-1|^c}
\right).
\end{eqnarray}
Claim (\ref{esf7.50}) follows from 
(\ref{esf7.100}), (\ref{esf7.110}), (\ref{esf7.120}) and (\ref{npc.180}).

The proof of claim (\ref{esf7.55}) is similar but easier since the
matrix $R$ does not depend on $\theta$. Indeed,
\begin{eqnarray}
\label{esf7.130}
\partial_{\theta}^{\beta} \tilde{R} (m, n)
&=&
(\partial_{\theta}^{\beta} \tilde{D}^{-1})(m, m) R(m, n) \\
\nonumber
&-&
\sum_{\scriptsize{\scriptstyle
\begin{array}{c}
\alpha \leq \beta \\
y_1, y_2 \in \Gamma
\end{array}
}
}
\left( \begin{array}{c} \beta \\ \alpha \end{array} \right) 
(\partial_{\theta}^{\alpha} \tilde{D}^{-1})(m, m)
R(m, y_1)
(\partial_{\theta}^{\beta - \alpha} G_{\Gamma}^{(j+1)}) (y_1, y_2) 
R (y_2, n) \\
\nonumber
&=& S_1 - S_2
\end{eqnarray}
Using similar estimates as above it follows that
\begin{eqnarray}
\label{esf7.140}
|S_1| 
\leq D_V 8 e \nu D_W D_N \rho_1
\end{eqnarray}
and 
\begin{eqnarray}
\label{esf7.150}
|S_2| &\leq& 
2^{|\beta|_1} D_V 
(8 e \nu D_W D_N \rho_1)^2
D_{1, 3} M_0^{2 |\beta|_1 + 1} 
\left(
\max_{y \geq 0}
(1 + y)^{q(2 |\beta|_1 + 1)} e^{-\frac{1}{4} |y-1|^c}
\right).
\end{eqnarray}
Claim (\ref{esf7.55}) follows from 
(\ref{esf7.130}), (\ref{esf7.140}), (\ref{esf7.150}), and 
(\ref{esf.92}), (\ref{esf.95}).

\end{proof}

\subsubsection{Completion of the proof of lemma \ref{Lesf.1}}
\label{esf8}

Recall that we have already proved lemma \ref{Lesf.1} in subsections \ref{esf1}
and \ref{esf2} except for 
estimates (\ref{esf.220})
and (\ref{esf.225}) 
We will complete the proof 
by showing
\begin{eqnarray}
\label{esf8.55}
\left|
\partial_{\theta}^{\beta}
r(\theta, \lambda)
\right|
&\leq&
D_{1, 5} \rho_1^{3/4}, \\
\label{esf8.50}
\left|
\partial^{\beta}
r(\theta, \lambda)
\right|
&\leq&
D_{1}
M_0^{2 |\beta|_1 + 2 \nu}
\end{eqnarray}
for all
$0 \leq |\beta|_1 \leq B_0$ and 
$(\theta, \lambda) \in U(\tilde{I}^{(j)}_{C', k'}, \delta_j/100 ) \times 
U(B_{\rho_j}(\lambda_{k'}^{(j)}), \rho_j/8)$
and by using (\ref{esf.99}). 
The function $r$ was defined in (\ref{esf1.65}).
The proof of (\ref{esf8.55}) and (\ref{esf8.50})
is a straight forward calculation. Indeed, using the representation of the
determinant as a sum of products of the entries 
and the fact that $\Lambda$ is a 
diagonal matrix, we can write
\begin{eqnarray}
\label{esf8.100}
r = \sum_{\sigma \in \Sigma} \prod_{i=1}^{d} a_{i, \sigma},
\end{eqnarray}
where $d$ denotes the cardinality of the singular set $S$ ($d \leq 2 \nu$),
$\Sigma$ contains less than $2^d (d!)$ elements, $a_{i, \sigma}$ denotes
(up to a sign)
an entry of either $\Lambda$ or $\tilde{R}$ and for each $\sigma \in 
\Sigma$ at least one of the terms $a_{i, \sigma}$, $1 \leq i \leq d$
is an entry of $\tilde{R}$. Representation (\ref{esf8.100}), 
together with (\ref{esf7.50}), (\ref{esf7.55}), (\ref{esf1.70}), 
(\ref{npc.200}), (\ref{npc.190}) yield estimates
(\ref{esf8.50}) and (\ref{esf8.55}).

\subsection{Proof of lemma \ref{lawp.1} -- 
Application of the Weierstrass Preparation Theorem}
\label{awp}

\begin{proof}
Let $\vartheta \in I_{C', k'}^{(j)}$. Define $S$ as in (\ref{conf1.10}),
where $\theta$ is replaced by $\vartheta$ in order to determine the
case in (\ref{conf1.5}). Suppose first that $S = \emptyset$. Then
statement A of lemma \ref{lawp.1} holds by proposition \ref{Pconf3.1}
(cf. remark \ref{Rconf3.1}) and (\ref{awp.10}). For the remainder of the proof we assume that
$S \neq \emptyset$. Define $q$ and $r$ through lemma \ref{Lesf.1}, by
writing
\begin{eqnarray}
\label{awp.30}
f_{C, k'}(\theta, \lambda) = \tilde{f}(\theta, \lambda)
(q(\theta, \lambda) + r(\theta, \lambda)),
\end{eqnarray}
such that (\ref{esf.205}) -- (\ref{esf.225}) are satisfied.  
We apply the Weierstrass preparation theorem in the form of lemma
\ref{Lwpt.1} to the function
\begin{eqnarray}
\label{awp.20}
f(z, \lambda) := q(z+ \vartheta, \lambda) + r(z + \vartheta, \lambda).
\end{eqnarray}
In the notation of lemma \ref{Lwpt.1} we set 
$d \equiv \#S$, $\delta \equiv \delta_j/100$, $\rho \equiv 2 \tilde{\rho}_j$,
$\epsilon \equiv \sqrt{\rho_1}$, $B_1$, $B_2$ as in (\ref{npc.90}), 
(\ref{npc.100}),
$C^* \equiv D_1 M_0^{2(B_0 + \nu)}$, $C \equiv 30 C^* B_1! 2^{B_1}$,
$\lambda_0 \equiv \lambda_{\tilde{k}}^{(j + 0.5)}$.

Assumptions (\ref{wpt.5}) -- (\ref{wpt.50}) of lemma \ref{Lwpt.1}
are satisfied by (\ref{awp.2}) -- (\ref{awp.8}).
The analyticity of $f$ as defined in (\ref{awp.20}) 
follows from lemma \ref{Lesf.1} since
$\vartheta \in I_{C', k'}^{(j)}$, $\lambda_{\tilde{k}}^{(j + 0.5)} \in
B_{\rho_j}(\lambda_{k'}^{(j)})$ and $2 \tilde{\rho}_j \leq \rho_j/4$ by
(\ref{awp.10}). Furthermore $f$ is of the form (\ref{wpt.70}), where
(\ref{wpt.100}) follows from (\ref{esf.220}), (\ref{esf.225}) and (\ref{awp.1}).
To verify (\ref{wpt.80}) and (\ref{wpt.90}) observe that 
\begin{eqnarray}
\nonumber
q(z + \vartheta, \lambda) = 
(z + (\vartheta - \theta_{\iota}))^d + \sum_{0 \leq l < d}
\tilde{a}_l(\lambda) (z + (\vartheta - \theta_{\iota}))^l 
= z^d  + \sum_{0 \leq l < d} a_l (\lambda) z^l,
\end{eqnarray}
with
\begin{eqnarray}
\nonumber
a_l(\lambda) =
\left( \begin{array}{c} d \\ l \end{array} \right) 
(\vartheta - \theta_{\iota})^{d-l} +
\sum_{l \leq r < d}
\tilde{a}_r(\lambda)
\left( \begin{array}{c} r \\ l \end{array} \right) 
(\vartheta - \theta_{\iota})^{r-l} \;\;\; \mbox{ for } 0 \leq l < d.
\end{eqnarray}
Since $\vartheta \in I_{C', k'}^{(j)}$ the definition of $\theta_{\iota}$
implies $|\vartheta - \theta_{\iota}| < \delta_1 \leq \frac{1}{2}$ (see
also (\ref{awp.16})). Using (\ref{esf.213}), (\ref{esf.215}), (\ref{awp.12})
and (\ref{awp.14}) we obtain for $0 \leq l < d$
\begin{eqnarray}
\nonumber
|a_l(\lambda)| \leq 2^{2 \nu} \delta_1 + 2^{4 \nu + 1} \rho_1 \leq 
\frac{1}{8d}.
\end{eqnarray}
and
\begin{eqnarray}
\nonumber
\left|
\partial^{\beta} a_l (\lambda)
\right| \leq 3^{3 \nu} 2^{2 \nu} < 2^{8 \nu} < 2^{B_1} < C.
\end{eqnarray}
The hypothesis of lemma \ref{Lwpt.1} are therefore satisfied
and we obtain functions $Q$ and $b_i$ satisfying (\ref{wpt.120}) --
(\ref{wpt.160}). Set
\begin{eqnarray}
\nonumber
p_{C, \tilde{k}, \vartheta} (\theta, \lambda) :=
\theta^d + \sum_{i=0}^{d-1} b_i(\lambda) \theta^i.
\end{eqnarray}
Statement (a) is then trivially satisfied.
Statements (c) -- (e) follow from (\ref{wpt.140}) -- (\ref{wpt.160}) and
(\ref{npc.210}). Claim (f) is certainly true in the case $d=1$
and follows for $d \geq 2$ from (\ref{wpt.140}) and 
(\ref{awp.16}) via
\begin{eqnarray}
\nonumber
\left|
\partial_{\theta} p_{C, \tilde{k}, \vartheta} (\theta, \lambda)
\right| \leq d \sum_{k=1}^{d-1} \left( \frac{\delta_1}{2} \right)^k +
\frac{1}{2d} \leq d \delta_1 + \frac{1}{2d} \leq 1.
\end{eqnarray} 
It remains to prove (b). 
It follows from (\ref{awp.30}), (\ref{esf.205}), (\ref{awp.20}), 
(\ref{wpt.120}) and (\ref{wpt.130}) that
\begin{eqnarray}
\label{awp.40}
|f_{C, k'}(\theta, \lambda)| \geq \frac{1}{2}
d_V^{2 \nu}
|p_{C, \tilde{k}, \vartheta}(\theta - \vartheta, \lambda)|
\end{eqnarray}
for all $(\theta, \lambda) \in 
U(\vartheta, \delta_j/400) \times
U(\lambda_{\tilde{k}}^{(j+0.5)}, 2 \tilde{\rho}_j)$. Using 
the cofactor matrix 
to represent the inverse of $b(\theta, \lambda)$
we obtain from corollary 
\ref{Cesf6.1} 
\begin{eqnarray}
\label{awp.50}
\| b^{-1}(\theta, \lambda)\|_{1, c} \leq d D_N e^2
\frac{(d-1)! [2(D_V+1) \delta_1]^{d-1}}{|f_{C, k'}(\theta, \lambda)|}
\end{eqnarray}
for all $(\theta, \lambda) \in 
U(\vartheta, \delta_j/400) \times
U(\lambda_{\tilde{k}}^{(j+0.5)}, 2 \tilde{\rho}_j)$ with 
$f_{C, k'}(\theta, \lambda) \neq 0$. Using hypothesis (\ref{awp.18})
one arrives at
\begin{eqnarray}
\nonumber
\| b^{-1}(\theta, \lambda)\|_{1, c} \leq 
\frac{1}{|p_{C, \tilde{k}, \vartheta}(\theta - \vartheta, \lambda)| \delta_1}
\end{eqnarray}
for all $(\theta, \lambda) \in 
U(\vartheta, \delta_j/400) \times
U(\lambda_{\tilde{k}}^{(j+0.5)}, 2 \tilde{\rho}_j)$ with 
$p_{C, \tilde{k}, \vartheta}(\theta - \vartheta, \lambda) \neq 0$.
Claim (b) then follows from (\ref{conf4.20}), proposition \ref{Pconf3.1}
and (\ref{esf3.2}).
\end{proof}

\subsection{Definition of ${\cal POL}_{\tilde{k}}^{(j+1)}$}
\label{dpo}
In this section we again assume that $j \geq 1$ 
and that 
induction statements $({\cal IS})_l$, $1 \leq l \leq j$ ,
$({\cal IS})(1)_{j+1}$ are satisfied.
Let $C' \in {\cal C}^{(j)}$ and $k' \in K^{(j)}$.
Using (\ref{dpo.5}) and 
\begin{eqnarray}
\nonumber
I^{(j)}_{C', k'} \subset 
(-\delta_1, \delta_1) \cup \bigcup_{i=1}^{\nu} 
([2 \omega_i^{(0)}] - \delta_1, [2 \omega_i^{(0)}] + \delta_1) 
\cup \bigcup_{i=1}^{\nu} 
([- 2 \omega_i^{(0)}] - \delta_1, [- 2 \omega_i^{(0)}] + \delta_1)
\end{eqnarray}
(cf. (\ref{fis115})) it is clear that we can find
a set $\Theta_{C', k'} \subset I^{(j)}_{C', k'}$ satisfying
\begin{eqnarray}
\label{dpo.50}
\Theta_{C', k'} \mbox{ contains at most }  \delta_j^{-1} 
\mbox{ elements }.&& \\
\label{dpo.55}
I^{(j)}_{C', k'} \subset \bigcup_{\vartheta \in \Theta_{C', k'}}
B_{\delta_j/800}(\vartheta).&&
\end{eqnarray}

We are ready to define the set of polynomials ${\cal POL}^{(j+1)}_{\tilde{k}}$,
where we use the notation $p \ominus q$ introduced in section \ref{res}, 
definition \ref{Dres.1}.
\begin{Definition}
\label{Ddpo.1}
Let $\tilde{k} \in K^{(j+0.5)}$. Set
$k' := \pi^{(j+0.5)}_j \tilde{k}$. For $C \in {\cal C}^{(j+1)}$ 
we define $C'$ as in (\ref{cpol.10}).
\begin{eqnarray}
\label{dpo.100}
{\cal POL}_{\tilde{k}}^{(j+1)}
:=
{\cal PI}_{\tilde{k}}^{(j+1)} \cup
{\cal PII}_{\tilde{k}}^{(j+1)} \cup
{\cal PIII}_{\tilde{k}}^{(j+1)},
\end{eqnarray}
where
\begin{eqnarray}
\label{dpo.105}
{\cal PI}_{\tilde{k}}^{(j+1)} &:=&
\{ (q, 0) \} \mbox{ with } q(\theta, \lambda) := \theta, \\
\label{dpo.110}
{\cal PII}_{\tilde{k}}^{(j+1)} &:=& 
\{ (p_{C, \tilde{k}, \vartheta}, \vartheta): C \in {\cal C}^{(j+1)},
\vartheta \in \Theta_{C', k'} \},
\\
\nonumber
{\cal PIII}_{\tilde{k}}^{(j+1)} &:=&
\{ (p_{C_1, \tilde{k}, \vartheta_1} \ominus p_{C_2, \tilde{k}, \vartheta_2}, 
[\vartheta_1 - \vartheta_2]): C_1, C_2 \in {\cal C}^{(j+1)},
\vartheta_1 \in \Theta_{C_1', k'}, \vartheta_2 \in \Theta_{C_2', k'} \}. \\
\label{dpo.115} 
\end{eqnarray}
\end{Definition}
The following proposition is a consequence of
lemma \ref{lawp.1}, lemma \ref{Lres.1}, (\ref{npc.490}), (\ref{npc.480}),
proposition \ref{Pcr.2}, (\ref{npc.460}), and (\ref{dpo.50}).

\begin{Proposition}
\label{Pdpo.1}
Let $j \geq 1$, $\tilde{k} \in K^{(j+0.5)}$ 
and assume that induction statements $({\cal IS})_{1 \leq l \leq j}$,
$({\cal IS})(1)_{j+1}$ hold. Then
${\cal POL}_{\tilde{k}}^{(j+1)} \subset {\cal POL}$ and
${\cal POL}_{\tilde{k}}^{(j+1)}$ contains at most $1 + (2 N_j)^{2 \nu} \delta_j^{-1} +
(2 N_j)^{4 \nu} \delta_j^{-2}$ elements. 
\end{Proposition}

\section{Completion of the induction step $j \to j+1$}
\label{comp}

In this section we will use the polynomials constructed in section
\ref{cpol} to define sets $K^{(j+1)}$, ${\cal N}^{(j+1)}$ and 
$I_{C, k}^{(j+1)}$ such that induction statements 
$({\cal IS})_{j+1}(2)-(4)$ are satisfied. We first state the estimates of 
lemma \ref{hyp} used in this section.

\noindent
{\bf Estimates of lemma \ref{hyp} used in section \ref{comp}:}

\begin{eqnarray}
\label{comp.30}
D_K \geq 2 \\
\label{comp.35}
40 \nu \rho_{j+1} N_{j+1} < \delta_{j+1} \\
\label{comp.40}
4 \nu D_P (4 \nu^2 N_{j+1} + 6 M_0^{4(B_0 + \nu)} ) \rho_{j+1} 
\leq \delta_{j+1} \\
\label{comp.45}
5 \delta_1 < \frac{s}{8^{\tau}}\\
\label{comp.55}
\delta_{j+1} \geq s d_{\tau, c} e^{-\frac{1}{2} N_j^c}\\
\label{comp.60}
4 D_N^2 M_j (1 + 16 D_N^2 M_j \delta_{j+1}^{-1} \delta_1^{-1}) \leq M_{j+1} \\
\label{comp.70}
80 \delta_{j+1} \leq \delta_j
\end{eqnarray}

\begin{Lemma}
\label{Lcomp.1}
Let $j \geq 1$. Assume that induction statements 
$({\cal IS})_l$, $1 \leq l \leq j$ ,
$({\cal IS})(1)_{j+1}$ are satisfied.
Then there exist sets $K^{(j+0.5)}$, $K^{(j+1)}$, ${\cal N}^{(j+1)}$ and 
$I_{C, k}^{(j+1)}$ (for $C \in {\cal C}^{(j+1)}$, $k \in K^{(j+1)}$),
such that induction statements 
$({\cal IS})_{j+1}(2)-(4)$ are satisfied.
\end{Lemma}

\begin{proof}
We define $K^{(j+0.5)}$ to be the set of indices which is generated
by the cube decomposition of the set ${\cal N}^{(j)}$ into sub-cubes
of radius $\tilde{\rho}_j$ (see section \ref{npc} J, (\ref{dpoh.1}) ).
We decompose this set ${\cal N}^{(j)}$ further into sub-cubes
of radius $\rho_{j+1}$ (see (\ref{dpoh.2})), generating a set of
indices $\hat{K}^{(j+1)}$ and a corresponding set of cube midpoints
$\{ \lambda_k^{(j+1)}: k \in \hat{K}^{(j+1)} \}$.
We obtain
\begin{eqnarray}
\nonumber
\overline{{\cal N}^{(j)}} &=& \bigcup_{k \in \hat{K}^{(j+1)}}
\overline{B_{\rho_{j+1}}(\lambda_k^{(j+1)})},\\
\nonumber
B_{\rho_{j+1}}(\lambda_{k_1}^{(j+1)}) \cap
B_{\rho_{j+1}}(\lambda_{k_2}^{(j+1)}) &=& \emptyset \;\;
\mbox{ for } k_1, k_2 \in \hat{K}^{(j+1)}, \; k_1 \neq k_2. 
\end{eqnarray} 
For $k \in \hat{K}^{(j+1)}$ we denote by 
$\tilde{k} := \pi^{(j+1)}_{j+0.5} k \in K^{(j+0.5)}$. Recall the definition
of the constant $D_K$ in (\ref{npc.115}) and set
\begin{eqnarray}
\label{comp.100}
K^{(j+1)} &:=& \hat{K}^{(j+1)} \setminus
\bigcup_{
2 N_j < |m| \leq 2 N_{j+1}
}  K_{m}^{(j+1)} \;\; \mbox{ where } \\
\nonumber
K_m^{(j+1)} &:=&
\{ k \in \hat{K}^{(j+1)}: \exists (p, \vartheta) \in
{\cal POL}_{\tilde{k}}^{(j+1)} \mbox{ with }
|p([<\omega_k^{(j+1)}, m>] - \vartheta, \lambda_k^{(j+1)} )| 
< D_K \delta_ {j+1} \\
\nonumber
&& \quad \quad \quad \quad \quad
\mbox{ and } | [<\omega_k^{(j+1)}, m>] - \vartheta | < \frac{\delta_1}{5}
\}
.
\end{eqnarray}
Furthermore, we define
\begin{eqnarray}
\label{comp.110}
{\cal N}^{(j+1)} := \bigcup_{k \in K^{(j+1)}} B_{\rho_{j+1}}
(\lambda_k^{(j+1)}).
\end{eqnarray} 

For $C \in {\cal C}^{(j+1)}$, $k \in K^{(j+1)}$ we denote $C'$ as
defined in (\ref{cpol.10}), $\tilde{k}:=
\pi^{(j+1)}_{j+0.5} k$ and $k' := \pi^{(j+1)}_{j}k$. Set
\begin{eqnarray}
\nonumber
I_{C, k}^{(j+1)} :=
\{
\theta \in I^{(j)}_{C', k'} &:&
\exists \vartheta \in \Theta_{C', k'} \mbox{ such that }
p_{C, \tilde{k}, \vartheta} \mbox{ exists }, |\theta - \vartheta| < 
\frac{\delta_1}{20}, \\
\label{comp.120}
&& \quad \quad \mbox{ and } 
|p_{C, \tilde{k}, \vartheta}(\theta- \vartheta, \lambda_k^{(j+1)})| < \delta_{j+1}
\}.
\end{eqnarray}
Recall that $p_{C, \tilde{k}, \vartheta}$ exists if and only if the 
set $S$ is not empty. The set $S$ was defined in (\ref{conf1.10}). In order
to determine the case we replace $\theta$ in (\ref{conf1.5}) by $\vartheta$. 

We now show that the inductive statements
$({\cal IS})_{j+1}(2)-(4)$ are satisfied. 

{\em $(2a)_{j+1}$:}
The statement is obvious from the definition (\ref{comp.110}) of
${\cal N}^{(j+1)}$.

{\em $(2b)_{j+1}$:}
Fix $\tilde{k} \in K^{(j+0.5)}$. Using definitions (\ref{comp.110}), 
(\ref{comp.100}) we see
\begin{eqnarray}
\nonumber
\left( \overline{{\cal N}^{(j)}} \setminus 
\overline{{\cal N}^{(j+1)}} \right)
\cap 
\overline{B_{\tilde{\rho}_j} \left(\lambda_{\tilde{k}}^{(j+0.5)}\right)}
\subset
\bigcup_{
2 N_j < |m| \leq 2 N_{j+1}
} 
\bigcup_{
\scriptsize
\begin{array}{cc}
k \in K_m^{(j+1)} \\
\pi^{(j+1)}_{j+0.5} k = \tilde{k}
\end{array}
}
\overline{B_{\rho_{j+1}}(\lambda_k^{(j+1)})}.
\end{eqnarray}
This means that for $\lambda \in
\left( 
\overline{{\cal N}^{(j)}} \setminus 
\overline{{\cal N}^{(j+1)}} \right)
\cap 
\overline{B_{\tilde{\rho}_j} \left(\lambda_{\tilde{k}}^{(j+0.5)}\right)}$.
there exist $m \in {\Bbb Z}^{\nu}$ with $2 N_j < |m| \leq 2 N_{j+1}$,
$k \in \hat{K}^{(j+1)}$ with $\pi^{(j+1)}_{j+0.5} k = \tilde{k}$
and $(p, \vartheta) \in {\cal POL}_{\tilde{k}}^{(j+1)}$ such that
\begin{eqnarray}
\label{comp.300}
|\lambda - \lambda_k^{(j+1)}| &\leq& \rho_{j+1}, \\
\label{comp.310}
|p([<\omega_k^{(j+1)}, m>]-\vartheta, \lambda_k^{(j+1)})| 
&<& D_K \delta_{j+1},\\
\label{comp.320}
|[<\omega_k^{(j+1)}, m>]-\vartheta| &<& \frac{\delta_1}{5}.
\end{eqnarray}
Recall from
proposition \ref{Pdpo.1} that $p \in {\cal POL}$ and therefore
the first order derivatives of $p$ are bounded by
\begin{eqnarray}
\label{comp.324}
|\partial_{\theta} p| &\leq& 8 \nu^2 D_P, \\
\label{comp.326}
|\partial_{\lambda_i} p| &\leq& 2 D_p M_0^{4(B_0 + \nu)}
\end{eqnarray}
on the set $\{|\theta| \leq 1/2\} \times 
B_{\tilde{\rho}_j} \left(\lambda_{\tilde{k}}^{(j+0.5)}\right)$.
Using in addition (\ref{comp.40}), induction statement 
$({\cal IS})(2b)_{j+1}$
follows, as soon as we have established the estimates
\begin{eqnarray}
\label{comp.330}
\left|
[<\omega, m>] - [<\omega_k^{(j+1)}, m>]
\right| \leq 2 \nu N_{j+1} \rho_{j+1} \leq \frac{\delta_1}{5}.
\end{eqnarray}
The latter inequality in (\ref{comp.330}) 
is a consequence of condition (\ref{comp.35}).
The first inequality of (\ref{comp.330}) seems to be a direct consequence
of (\ref{comp.300}). However, since the function $[ \cdot ]$ is 
discontinuous at points with real parts in the set 
$(2 {\Bbb Z} +1) \pi$ we need to 
establish that
\begin{eqnarray}
\label{comp.340}
\mbox{ dist}(\vartheta, (2 {\Bbb Z} + 1)\pi) \geq \frac{\delta_1}{2}.
\end{eqnarray}
Observe that the definitions of ${\cal POL}_{\tilde{k}}^{(j+1)}$
and $\Theta_{C', k'}^{(j+1)}$ in section \ref{dpo} imply that
\begin{eqnarray}
\nonumber
\mbox{ dist}(\vartheta, 
\{ 0 \} \cup \{ \pm 2 \omega_j^{(0)} : 1 \leq j \leq \nu \}
\cup \{ \pm 2 \omega_i^{(0)} \pm 2 \omega_j^{(0)}: 1 \leq i, j \leq \nu \}
+ 2 \pi {\Bbb Z}) 
\leq 2 \delta_1.
\end{eqnarray}
Furthermore, assumption A3 together with (\ref{comp.45}) imply
that
\begin{eqnarray}
\nonumber
\mbox{ dist}(<\omega^{(0)}, n>, 2 \pi {\Bbb Z}) > 5 \delta_1
\;\;\; \mbox{ for all } \; 0 < |n| \leq 8,
\end{eqnarray}
yielding (\ref{comp.340}) by contradiction.

{\em $(2c)_{j+1}$:}
The statement follows from proposition \ref{Pdpo.1}.

{\em $(3a)_{j+1}$:}
The statement is an immediate consequence of definition (\ref{comp.120}).

{\em $(3b)_{j+1}$:}
Let $\lambda \in U({\cal N}^{(j+1)}, \rho_{j+1})$ and
$m \in {\Bbb Z}^{\nu}$ with $2 N_j < |m| \leq 2 N_{j+1}$. 
Then there exists a
$k \in K^{(j+1)}$ such that $\lambda \in U_{2 \rho_{j+1}}
(\lambda_k^{(j+1)})$. Set $\tilde{k} := \pi^{(j+1)}_{j+0.5} k$.
Since  ${\cal PI}_{\tilde{k}}^{(j+1)}
\subset {\cal POL}_{\tilde{k}}^{(j+1)}$ (see (\ref{dpo.100})) the definition
of $K^{(j+1)}$ in (\ref{comp.100}) implies
\begin{eqnarray}
\nonumber
|[<\omega_k^{(j+1)}, m>]| \geq D_K \delta_{j+1}.
\end{eqnarray}
Using in addition that $|\omega - \omega_k^{(j+1)}| < 2 \rho_{j+1}$ 
together with (\ref{comp.35}) 
we obtain
\begin{eqnarray}
\nonumber
\mbox{ dist}(<\omega, m>, 2 \pi {\Bbb Z}) =
|[<\omega, m>]| \geq (D_K - 1) \delta_{j+1}.
\end{eqnarray}
Statement $(3b)_{j+1}$ now follows from (\ref{comp.30}) and (\ref{comp.55}).

{\em $(3c)_{j+1}$:}
Let $k \in K^{(j+1)}$, $\lambda \in 
U(B_{\rho_{j+1}}(\lambda_k^{(j+1)}), \rho_{j+1})$,
$C \in {\cal C}^{(j+1)}$ and $m \in {\Bbb Z}^{\nu}$
with
$2 N_j < |m| \leq 2 N_{j+1}$.
From (\ref{comp.30}), (\ref{comp.100}), (\ref{comp.120}) and
${\cal PII}_{\tilde{k}}^{(j+1)} \subset
{\cal POL}_{\tilde{k}}^{(j+1)}$ (see (\ref{dpo.100})) it follows 
by contradiction that
\begin{eqnarray}
\label{comp.130}
[<\omega_k^{(j+1)}, m>] \in (-\pi, \pi] \setminus I_{C, k}^{(j+1)}.
\end{eqnarray}
This implies 
\begin{eqnarray}
\label{comp.135}
<\omega_k^{(j+1)}, m> \in {\Bbb R} \setminus \tilde{I}_{C, k}^{(j+1)}.
\end{eqnarray}
The claim now follows from (\ref{comp.35}).

{\em $(3d)_{j+1}$:}
Proof by contradiction:
Let $k \in K^{(j+1)}$, $\lambda \in 
U(B_{\rho_{j+1}}(\lambda_k^{(j+1)}), \rho_{j+1})$,
$\theta \in {\Bbb C}$, $C_1$, $C_2 \in {\cal C}^{(j+1)}$, $m$, 
$n \in {\Bbb Z}^{\nu}$ with $2 N_j < |n-m| \leq 2 N_{j+1}$ and assume that
\begin{eqnarray}
\nonumber
\theta + <\omega, m> &\in& 
U \left( \tilde{I}_{C_1, k}^{(j+1)}, \frac{\delta_{j+1}}{10} \right)
\;\; \mbox{ and } \\
\nonumber
\theta + <\omega, n> &\in& 
U \left( \tilde{I}_{C_2, k}^{(j+1)}, \frac{\delta_{j+1}}{10} \right).
\end{eqnarray}
Since $B(I_{C, 1}^{(1)}, \delta_1) \subset (-\pi, \pi)$ for any 
$C \in {\cal C}^{(1)}$ 
(see $({\cal IS}) (3a)_{j=1}$) we conclude that
\begin{eqnarray}
\nonumber
[\theta + <\omega, m>] &\in& 
U \left( I_{C_1, k}^{(j+1)}, \frac{\delta_{j+1}}{10} \right)
\;\; \mbox{ and } \\
\nonumber
[\theta + <\omega, n>] &\in& 
U \left( I_{C_2, k}^{(j+1)}, \frac{\delta_{j+1}}{10} \right).
\end{eqnarray}
Set $\tilde{k} := \pi^{(j+1)}_{j+0.5} k$, $k' := \pi^{(j+1)}_{j} k$
and $C_1'$, $C_2'$ according to (\ref{cpol.10}).

Using statement (f) in lemma \ref{lawp.1} and
definition (\ref{comp.120}) there exist $\vartheta_1
\in \Theta_{C_1', k'}$, $\vartheta_2 \in \Theta_{C_2', k'}$ such that
\begin{eqnarray}
\label{comp.410}
\left|
p_{C_1, \tilde{k}, \vartheta_1}
\left(
[\theta + <\omega, m>] - \vartheta_1, \lambda_k^{(j+1)}
\right)
\right| < \frac{11}{10} \delta_{j+1}, \\
\label{comp.420}
\left|
p_{C_2, \tilde{k}, \vartheta_2}
\left(
[\theta + <\omega, n>] - \vartheta_2, \lambda_k^{(j+1)}
\right)
\right| < \frac{11}{10} \delta_{j+1}
\end{eqnarray}
and $|[\theta + <\omega, m>] - \vartheta_1| < \delta_1/14$, 
$|[\theta + <\omega, n>] - \vartheta_2| < \delta_1/14$; here we have 
used that (\ref{comp.70}) yields $\delta_1/20 + \delta_{j+1}/10 \leq
\delta_1/14$.
Introduce the auxiliary variables
\begin{eqnarray}
\nonumber
x := [\theta + <\omega, m>] - \vartheta_1\;, \;\;\;
y := [\theta + <\omega, n>] - \vartheta_2.
\end{eqnarray}
We calculate
\begin{eqnarray}
\label{comp.400}
x - y = [x-y] = [<\omega, m-n> - (\vartheta_1 - \vartheta_2)]
= [<\omega, m-n>] - [\vartheta_1 - \vartheta_2],
\end{eqnarray}
where the last equality is justified by $|x-y| < \delta_1/7$ and by
dist$([\vartheta_1 - \vartheta_2], (2 {\Bbb Z} + 1) \pi) \geq \delta_1/2$
(see (\ref{comp.340}) above with 
$\vartheta = [\vartheta_1 - \vartheta_2]$). Define 
\begin{eqnarray}
\nonumber
(q, \vartheta) :=
(p_{C_1, \tilde{k}, \vartheta_1} \ominus  p_{C_2, \tilde{k}, \vartheta_2},
[\vartheta_1 - \vartheta_2]).
\end{eqnarray}
Since ${\cal PIII}_{\tilde{k}}^{(j+1)} \subset 
{\cal POL}_{\tilde{k}}^{(j+1)}$ we observe that
$(q, \vartheta) \in {\cal POL}_{\tilde{k}}^{(j+1)}$.
Using the notation of lemma \ref{Lres.1} one obtains
\begin{eqnarray}
\nonumber
q(x-y) = R_1(x, y) p_{C_1, \tilde{k}, \vartheta_1}(x, \lambda_k^{(j+1)}) +
R_2(x, y) p_{C_2, \tilde{k}, \vartheta_2}(y, \lambda_k^{(j+1)}).
\end{eqnarray}
From (\ref{comp.410}), (\ref{comp.420}), (\ref{comp.400}),
lemma \ref{Lres.1}, lemma \ref{lawp.1} B (c) and the definition of $D_K$  (\ref{npc.115}) it follows that
\begin{eqnarray}
\label{comp.436}
| q([<\omega, m-n>] - \vartheta, \lambda_k^{(j+1)} ) |
\leq (D_K - \nu^2 D_P) \delta_{j+1}.
\end{eqnarray}
Repeating the derivation of (\ref{comp.330}) in the proof of statement 
$(2b)_{j+1}$ we obtain
\begin{eqnarray}
\label{comp.438}
|[<\omega, m-n>] - [ <\omega^{(j+1)}_k, m-n>] | \leq 4 \nu N_{j+1} \rho_{j+1}
\leq \frac{\delta_{j+1}}{10} \quad \mbox{(cf. (\ref{comp.35}))}.
\end{eqnarray}
Since $|\partial_{\theta} q| \leq 8 \nu^2 D_P$ (cf. (\ref{comp.324}))
estimates (\ref{comp.436}) and (\ref{comp.438}) yield
\begin{eqnarray}
\label{comp.440}
| q([<\omega^{(j+1)}_k, m-n>] - \vartheta, \lambda_k^{(j+1)} ) |
< D_K \delta_{j+1}.
\end{eqnarray}
In addition, (\ref{comp.400}), (\ref{comp.35}), and (\ref{comp.70}) imply
\begin{eqnarray}
\nonumber
|[<\omega^{(j+1)}_k, m-n>] - \vartheta| &\leq&
|x - y|  + |[<\omega, m-n>] - [ <\omega^{(j+1)}_k, m-n>]| \\
\nonumber
&\leq&
\frac{\delta_1}{7} + \frac{\delta_{j+1}}{10} < \frac{\delta_1}{5}.
\end{eqnarray}
By definition (\ref{comp.100}) we conclude $k \in 
K^{(j+1)}_{m-n}$, contradicting the assumption that $k \in K^{(j+1)}$.

{\em $(4)_{j+1}$:}
Let $k \in K^{(j+1)}$, $C \in {\cal C}^{(j+1)}$. Set 
$k': = \pi_j^{(j+1)} k$, $\tilde{k} := \pi^{(j+1)}_{j+0.5} k$ 
and define $C'$ as
in (\ref{cpol.10}).
Let $(\theta, \lambda) \in 
U(\tilde{I}^{(j)}_{C', k'} \setminus \tilde{I}^{(j+1)}_{C, k},
\delta_{j+1}/10) \times 
U(B_{\rho_{j+1}}(\lambda_k^{(j+1)}), \rho_{j+1})$.
Recall from the construction of the inverse matrix $G_C^{(j+1)}
(\theta, \lambda)$ in section \ref{cpol} that we need to distinguish
the cases whether the singular set $S$ defined in (\ref{conf1.10})
is empty or not. If $S = \emptyset$ then statement A of lemma
\ref{lawp.1} holds and estimate (\ref{is.10}) follows from (\ref{cpol.20}),
(\ref{comp.60}), and (\ref{comp.70}). Assume now that $S \neq \emptyset$.
There exists 
$\tilde{\theta} \in 
I^{(j)}_{C', k'} \setminus I^{(j+1)}_{C, k}$ satisfying 
\begin{eqnarray}
\label{comp.170}
|[\theta] - \tilde{\theta}| < \delta_{j+1}/10 \leq \delta_j/800 
\end{eqnarray}
(see (\ref{comp.70})) and a 
$\vartheta \in \Theta_{C', k'}$ such that
$|\tilde{\theta} - \vartheta| < \delta_j/800$ (see (\ref{dpo.55})). 
Hence
\begin{eqnarray}
\label{comp.200}
|[\theta] - \vartheta| < \frac{\delta_j}{400}.
\end{eqnarray}
Note that $(p_{C, \tilde{k}, \vartheta}, \vartheta) \in
{\cal PII}_{\tilde{k}}^{(j+1)} \subset
{\cal POL}_{\tilde{k}}^{(j+1)}$ and since 
$\tilde{\theta} \notin I_{C, k}^{(j+1)}$ definition
(\ref{comp.120}) implies
\begin{eqnarray}
\nonumber
|p_{C, \tilde{k}, \vartheta}(\tilde{\theta} - \vartheta, \lambda_k^{(j+1)})|
\geq \delta_{j+1}.
\end{eqnarray}
Using statement B (f) of lemma \ref{lawp.1} and (\ref{comp.170}) we conclude
\begin{eqnarray}
\nonumber
|p_{C, \tilde{k}, \vartheta}([\theta] - \vartheta, \lambda_k^{(j+1)})| >
\frac{9}{10} \delta_{j+1}.
\end{eqnarray}
The estimate on the first derivatives of $p_{C, \tilde{k}, \vartheta}$ 
with respect
to $\lambda$ (see (\ref{comp.326})), $|\lambda - \lambda_k^{(j+1)}|
< 2 \rho_{j+1}$ and (\ref{comp.40}) imply
\begin{eqnarray}
\label{comp.220}
|p_{C, \tilde{k}, \vartheta}([\theta] - \vartheta, \lambda)| >
\frac{1}{2} \delta_{j+1}.
\end{eqnarray}
Statement B (b) of lemma 
\ref{lawp.1} together with (\ref{comp.60})
yield the estimate (\ref{is.10}).
It follows from the constructions in section \ref{cpol} that $G_C^{(j+1)}$
is analytic and represents the inverse of $T_C^{(j+1)}$ except on the
set $\tilde{Z}_C$. Realness of the entries of $G_C^{(j+1)}$ follow from
realness of the matrix $T_C^{(j+1)}$, which in turn is a consequence
of $({\cal IS})(1g)_{j+1}$ and proposition \ref{Plop.2}. 
\end{proof}

\section{Solution of the ${\cal P}$ -- equation}
\label{sp}

In this section we use the induction statements proved in sections \ref{fis}
-- \ref{comp} to produce a solution of the $P$ - equation and to show
a few properties of this solution.

\noindent
{\bf Estimate of lemma \ref{hyp} used in section \ref{sp}:}

\begin{eqnarray}
\label{sp.50}
e^{2 \nu -6} D_{\psi}(B_1 +1)
\left(
\sum_{j=1}^{\infty} 
\left( \frac{8}{\rho_{j+1}} \right)^{B_1 + 1} e^{-\frac{1}{8} N_j^c}
\right) \leq 1.
\end{eqnarray}

Define 
\begin{eqnarray}
\label{sp.10}
{\cal N}^{\infty} &:=& \bigcap_{j=1}^{\infty} 
\overline{{\cal N}^{(j)}}, \\
\label{sp.17}
v(\lambda) &:=& \lim_{j \to \infty} v_j(\lambda) \;\;\; \mbox{ for }
\lambda \in {\cal N}^{(1)} \;\;(= B_{\rho_1}(\lambda^{(0)}), 
\mbox{ see (\ref{fis.111})) }. 
\end{eqnarray}
The existence of the limit $v$ and the proof of
theorem \ref{TP} below is a consequence of the induction statements 
$({\cal IS})(1)_{j \geq 1}$, $({\cal IS})(3b)_{j \geq 1}$
and of the definition of $D_3$ in (\ref{npc.220})

\begin{Theorem}
\label{TP}
The function
\begin{eqnarray}
\nonumber
v : U(B_{\rho_1}(\lambda^{(0)}), \rho_1) \to X_{1, c},
\end{eqnarray}
defined in (\ref{sp.17}) exists and satisfies
\begin{itemize}
\item[(a)]
$v$ is a $C^{\infty}$ -- function.
\item[(b)]
${\cal P}(v(\lambda), \lambda) = 0$ for all $\lambda \in {\cal N}^{\infty}$.
\item[(c)]
$\|v(\lambda)\|_{1, c} \leq D_3 N_0^{\tau+1} |a|^2$ for all $\lambda \in 
U(B_{\rho_1}(\lambda^{(0)}), \rho_1)$.
\item[(d)] 
$\|\partial^{\beta}(\varphi(a) + v(\lambda))\|_{1, c} \leq D_3 N_0^{E_{\rho}(|\beta|_1 -1)}$ 
for all $\lambda \in 
U(B_{\rho_1}(\lambda^{(0)}), \rho_1)$ and $1 \leq |\beta|_1 \leq B_1 + 1$.
\item[(e)]
For all $\lambda \in U(B_{\rho_1}(\lambda^{(0)}), \rho_1)$ the support
of $v(\lambda)$ is contained in ${\Bbb Z}^{\nu} \setminus
({\cal S} \cup \{ 0 \})$.
\item[(f)]
For all $\lambda \in B(B_{\rho_1}(\lambda^{(0)}), \rho_1)$ and
$m \in {\Bbb Z}^{\nu}$ we have
$v(\lambda)(m) = v(\lambda)(-m) \in {\Bbb R}$.
\item[(g)]
Let $\lambda =(a, \omega) \in U(B_{\rho_1}(\lambda^{(0)}), \rho_1)$
with $a_i = 0$ for some $1 \leq i \leq \nu$. Then 
$v(\lambda)(m) = 0$ if $m_i \neq 0$.
\item[(h)]
Let $\lambda = (a, \omega) 
\in {\cal N}^{\infty}$ and $m \in {\Bbb Z}^{\nu}
\setminus \{ 0 \}$. Then
dist$(< \omega, m>, 2 \pi {\Bbb Z}) \geq s d_{\tau, c} e^{-\frac{1}{2} |m|^c}$.
\end{itemize}
\end{Theorem}

\newpage
\noindent
\begin{center}
{\huge \bf Chapter III} \vspace{1.2cm}\\
{\huge \bf The Bifurcation Equation} \vspace{2cm}
\end{center}

In this chapter we solve the bifurcation equation (section \ref{solq}), 
provide lower bounds for
the measure of the set of non-resonant parameters (section \ref{est}), 
and prove our main result, theorem \ref{Tsmr3.1}, in section \ref{pmt}.

\section{Solution of the ${\cal Q}$-equation}
\label{solq}
In this section we will first show by a standard implicit function theorem
that the bifurcation equation ($Q$ -- equation) can be solved for $\omega$.
Observe that the bifurcation is degenerate. However, we will use
property (g) of the solution of the $P$ -- equation (see theorem \ref{TP})
to factor out the degeneracy. In the remaining part of the section we
derive estimates on the derivatives of the implicitly defined function $\omega$.

\noindent
{\bf Estimates of lemma \ref{hyp} used in section \ref{solq}:}

\begin{eqnarray}
\label{est3.10}
4 \nu^2 D_{4, 1} N_0^{2 E_{\rho}} \rho_1  \leq d_V \\
\label{est3.15}
\rho_1 < \delta_V \\
\label{est3.20}
2 \nu D_{4, 1} N_0^{E_{\rho}} \tilde{\rho}_1 \leq \rho_1 d_V
\end{eqnarray}

\begin{Theorem}
\label{Tsolq.1}
Assume that $v$ is defined
by (\ref{sp.17}). Then there exists an unique  
$C^{\infty}$ -- function 
\begin{eqnarray}
\label{solq.5}
\omega: {\Bbb R}^{\nu} \supset B_{\tilde{\rho}_1}(0) \to 
B_{\rho_1}(\omega^{(0)}) \subset  {\Bbb R}^{\nu},
\end{eqnarray}
such that for $\lambda(a) := (a , \omega(a))$ we have
\begin{eqnarray}
\label{solq.10}
{\cal Q}(v(\lambda(a)), \lambda(a)) = 0 \;\; 
\mbox{ for all } a \in B_{\tilde{\rho}_1}(0).
\end{eqnarray}
Furthermore, the derivatives of $\omega$ up to order $B_1$ can be estimated
by
\begin{eqnarray}
\label{solq.15}
\left| \partial^{\alpha} \omega (a) \right| \leq
\left( D_4 N_0^{B_1 E_{\rho}} \right)^{2 |\alpha|_1 - 1} \;\;\;
\mbox{ for } 1 \leq |\alpha|_1 \leq B_1, \; |a| < \tilde{\rho}_1.
\end{eqnarray} 
\end{Theorem}

The proof of theorem \ref{Tsolq.1} proceeds in several steps. 
In section \ref{solq1} we define a map $\tilde{h}: 
B_{\rho_1}(\lambda^{(0)}) \to {\Bbb R}^{\nu}$ such that
$\tilde{h}(\lambda) = 0$ implies 
${\cal Q}(v(\lambda), \lambda) = 0$.
In section \ref{solq2} we show $\tilde{h}(\lambda^{(0)}) = 0$
and the invertibility of the derivative $D_{\omega} \tilde{h}
(\lambda^{(0)})$. Consequently we can apply the implicit function 
theorem and obtain a function $\omega$ defined on some neighborhood
$B_{\tilde{\rho}}(0)$
of the origin in ${\Bbb R}^{\nu}$ satisfying $\tilde{h}(a, \omega(a))=0$
for $a \in B_{\tilde{\rho}}(0)$. We proceed in section \ref{est3} to
show various estimates for the functions $\tilde{h}$, $\omega$ and their 
derivatives which imply $\tilde{\rho} \geq \tilde{\rho}_1$ and 
the estimate (\ref{solq.15}). This proves
theorem \ref{Tsolq.1}. Moreover, the estimates of proposition
\ref{Pest3.1} will also be used in the subsequent section \ref{est}
to obtain lower bounds on the measure of the set 
$\{ a: (a, \omega(a)) \in {\cal N}^{\infty} \}$. 

\subsection{Symmetry reductions of the ${\cal Q}$-equation}
\label{solq1}

We define 
\begin{eqnarray} 
\label{solq1.10}
h&:& B_{\rho_1}(\lambda^{(0)}) \ni (a, \omega) \mapsto
h(a, \omega) \in {\Bbb R}^{\nu}, \\ 
\nonumber
&&h_j(a, \omega) := {\cal Q}(v(a, \omega), (a, \omega)) (e_j) \;\;
\mbox{ for } 1 \leq j \leq \nu.
\end{eqnarray}
Indeed, we know from theorem \ref{TP} (a), (f) that all $h_j$ are real-valued 
$C^{\infty}$-functions.
Furthermore, 
\begin{eqnarray}
\nonumber
{\cal Q}(v(a, \omega), (a, \omega)) (e_j) =
{\cal Q}(v(a, \omega), (a, \omega)) (-e_j) \;\; \mbox{ for all }
(a, \omega) \in B_{\rho_1}(\lambda^{(0)}), 1 \leq j \leq \nu.
\end{eqnarray}
Therefore it suffices to determine the zeros of $h$ in order to
solve the $Q$-equation.

Statement (g) of theorem \ref{TP} implies
\begin{eqnarray}
\label{solq1.30}
h_j (a, \omega) = 0 \;\; \mbox{ if } a_j = 0.
\end{eqnarray}
Set $\hat{a}_j := a - a_j e_j$ the vector which equals $a$ except at the
$j$-th component, which is set to be zero. Then (\ref{solq1.30}) and 
the fundamental theorem of
calculus yield
\begin{eqnarray}
\label{solq1.40}
h_j (a, \omega) = a_j \left( \int_0^1
\frac{\partial h_j}{\partial a_j} (\hat{a}_j + t a_j e_j, \omega) dt
\right)
\end{eqnarray}
Therefore we can define a $C^{\infty}$-function
\begin{eqnarray}
\label{solq1.45}
\tilde{h}&:& B_{\rho_1}(\lambda^{(0)}) \ni (a, \omega) \mapsto
\tilde{h}(a, \omega) \in {\Bbb R}^{\nu}, \\ 
\nonumber
&&\tilde{h}_j(a, \omega) := \int_0^1
\frac{\partial h_j}{\partial a_j} (\hat{a}_j + t a_j e_j, \omega) dt \;\;
\mbox{ for } 1 \leq j \leq \nu,
\end{eqnarray}
satisfying
\begin{eqnarray}
\label{solq1.50}
h_j (a, \omega) = a_j \tilde{h}_j (a, \omega).
\end{eqnarray}
Note that the zero-set of $\tilde{h}$ is contained in the zero-set of
$h$.
  
\subsection{Towards a solution of the reduced equations}
\label{solq2}

In this subsection we will show that 
\begin{eqnarray}
\label{solq2.5}
&\tilde{h}(\lambda^{(0)}) = 0,& \\
\label{solq2.10}
&(D_{\omega} \tilde{h})(\lambda^{(0)}) \;\; \mbox{ is invertible }.&
\end{eqnarray}
It follows then from the implicit function theorem that there
exists a $\tilde{\rho} > 0$ and a function $\omega$ 
such that $\tilde{h}(a, \omega(a)) = 0$ for all
$a \in B_{\tilde{\rho}}(0)$. However, since we need to show that
$\tilde{\rho}$ can be chosen to equal $\tilde{\rho}_1$ we will defer
the application of the implicit function theorem to section \ref{est3}
(see proposition \ref{Pest3.1} below) where we provide the necessary
estimates on the function $\tilde{h}$ and its derivatives to guarantee
that $\omega$ can be defined on all of $B_{\tilde{\rho}_1}(0)$. 

We introduce the following notation
\begin{eqnarray}
\label{solq2.20}
u(a, \omega) &:=& \varphi(a) + v(a, \omega) \;\; \mbox{ for }
(a, \omega) \in B_{\rho_1}(\lambda^{(0)}) , \\
\label{solq2.25}
l(a, \omega) &:=& (l_j(a ,\omega))_{1 \leq j \leq \nu} \;\; 
\mbox{ where }
l_j (a, \omega) := a_j
\left( \alpha_1 - \frac{(j \gamma)^2}{4 \sin^2 (\frac{\omega_j}{2})} 
\right), \\
\label{solq2.30}
n(a, \omega) &:=& (n_j(a ,\omega))_{1 \leq j \leq \nu} \;\; 
\mbox{ where }
n_j (a, \omega) := W(u(a, \omega))(e_j).
\end{eqnarray}
Recall from (\ref{ova.150}) and (\ref{solq1.10}) that 
\begin{eqnarray}
\label{solq2.35}
h = l + n.
\end{eqnarray}
It follows from statement (g) of theorem \ref{TP} 
that $v(a, \omega) = 0$ for $a = 0$.
We conclude that
\begin{eqnarray}
\label{solq2.40}
\partial_{\omega}^{\beta} u(0, \omega) = 0 \;\; \mbox{ for all } 
\omega \in B_{\rho_1}(\omega^{(0)}), |\beta|_1 \geq 0.
\end{eqnarray} 
Since $W$ only contains term of order $\geq 2$ it is clear 
(see lemma \ref{Lenp.1}) that
\begin{eqnarray}
\label{solq2.45}
DW(0) = 0.
\end{eqnarray}
The following statements are easy to verify.
\begin{eqnarray}
\nonumber
\tilde{h}_j (\lambda^{(0)}) &=&  
\frac{\partial h_j}{\partial a_j} (\lambda^{(0)}), \\
\nonumber
\frac{\partial \tilde{h}_j}{\partial \omega_k} (\lambda^{(0)}) &=&
\frac{\partial^2 h_j}{\partial a_j \partial \omega_k} (\lambda^{(0)}), \\ 
\nonumber
\frac{\partial l_j}{\partial a_j} (\lambda^{(0)}) &=& 
V_j (\omega_j^{(0)}) = 0, \\
\nonumber
\frac{\partial^2 l_j}{\partial a_j \partial \omega_k} (\lambda^{(0)})
&=&  \delta_{j, k} V_j'(\omega_j^{(0)}) =
\delta_{j, k} \frac{(j \gamma)^2 \cos (\frac{\omega_j^{(0)}}{2})}
{4 \sin^3(\frac{\omega_j^{(0)}}{2})} \neq 0, \\
\nonumber
\frac{\partial n_j}{\partial a_j} (\lambda^{(0)}) &=&
DW(0)[\partial_{a_j} u](e_j) = 0 \;\; \mbox{ (see (\ref{solq2.45})) }, \\
\nonumber
\frac{\partial^2 n_j}{\partial a_j \partial \omega_k} (\lambda^{(0)}) 
&=&
DW(0)[\partial_{a_j} \partial_{\omega_k} u](e_j) + 
D^2 W(0)[\partial_{a_j} u, \partial_{\omega_k} u](e_j) = 0 \;\;
\mbox{ by (\ref{solq2.45}), (\ref{solq2.40}) }.
\end{eqnarray}
Claims (\ref{solq2.5}) and (\ref{solq2.10}) 
follow readily from these observations (see 
(\ref{solq2.35})).
For later convenience we define
\begin{eqnarray}
\label{solq2.105}
\Lambda &:=& \mbox{ diag}(\Lambda_j),  \quad \mbox{ where }\\
\label{solq2.100}
\Lambda_j &=& V_j'(\omega_j^{(0)}) = \tilde{V}_{j, +} (\omega_j^{(0)}) =
\frac{(j \gamma)^2 \cos (\frac{\omega_j^{(0)}}{2})} 
{4 \sin^3(\frac{\omega_j^{(0)}}{2})}
\geq d_V  \;\;
\mbox{ (cf. (\ref{smr2.60}), proposition \ref{Plop.1}). }
\end{eqnarray}
Then 
\begin{eqnarray}
\label{solq2.110}
D_{\omega}\tilde{h}(\lambda^{(0)}) = \Lambda. 
\end{eqnarray}

\subsection{Estimates on the functions $\tilde{h}$ and $\omega$}
\label{est3}

\begin{Proposition}
\label{Pest3.1}
\begin{itemize}
\item[(a)] For all $\lambda \in B_{\rho_1}(\lambda^{(0)})$
\begin{eqnarray}
\label{est3.100}
\left|
\partial^{\beta} \tilde{h} (\lambda) 
\right| &\leq& 
D_{4, 1} N_0^{E_{\rho}|\beta|_1} \;\;\; \mbox{ for all }
0 \leq |\beta|_1 \leq B_1, \\
\label{est3.145}
\left\|
\left( D_{\omega} \tilde{h} \right)^{-1} (\lambda) 
\right\| 
&\leq&  
\frac{2}{ d_V}.
\end{eqnarray}
\item[(b)]
There exists an unique $C^{\infty}$ -- function $\omega: 
B_{\tilde{\rho}_1}(0) \to B_{\rho_1}(\omega^{(0)})$ such that
$\tilde{h}(a, \omega(a)) = 0$ for all $a \in B_{\tilde{\rho}_1}(0)$.
\item[(c)]
\begin{eqnarray}
\label{est3.50}
\left| \partial^{\alpha} \omega (a) \right| \leq
\left( D_4 N_0^{B_1 E_{\rho}} \right)^{2 |\alpha|_1 - 1} \;\;\;
\mbox{ for } 1 \leq |\alpha|_1 \leq B_1, \; |a| < \tilde{\rho}_1.
\end{eqnarray} 
\end{itemize}
\end{Proposition}

\begin{proof}

\noindent
{\em (a)}
It follows from (\ref{solq2.25}), proposition \ref{Plop.1} (a), and
(\ref{est3.15}) that 
\begin{eqnarray}
\label{est3.110}
\left|
\partial^{\beta} l (\lambda) 
\right| \leq 
D_{V} \;\;\; \mbox{ for all }
\lambda \in B_{\rho_1}(\lambda^{(0)}), \; 0 \leq |\beta|_1 \leq B_1 + 1.
\end{eqnarray}
In order to estimate the derivatives on the function $n$ (see (\ref{solq2.30}))
we use proposition \ref{Pcr.1} and obtain
for $|\beta|_1 \geq 1$
\begin{eqnarray}
\label{est3.120}
\partial^{\beta}_{\lambda} W(u(\lambda)) =
\sum_{p=1}^{|\beta|_1} \frac{1}{p!}
\sum_{\scriptsize
{\scriptstyle
\begin{array}{c}
\beta_1 + \ldots + \beta_p = \beta \\
\beta_i \neq 0 \mbox{ for } 1 \leq i \leq p
\end{array}
}
}
\frac{\beta !}{\beta_1 ! \ldots \beta_p !}
D^pW(u)[\partial^{\beta_1}u, \ldots, \partial^{\beta_p} u].
\end{eqnarray}
Statement (d) of theorem \ref{TP} together with 
(\ref{npc.230}) yield
\begin{eqnarray}
\label{est3.140}
\left|
\partial^{\beta} n (\lambda) 
\right| \leq 
D_{4, 2} N_0^{E_{\rho}(|\beta|_1 - 1)}
\;\;\; \mbox{ for all }
\lambda \in B_{\rho_1}(\lambda^{(0)}), \; 1 \leq |\beta|_1 \leq B_1 + 1.
\end{eqnarray}
Since $h = l + n$ (see (\ref{solq2.35})) 
claim (\ref{est3.100}) follows from
(\ref{solq1.45}), (\ref{est3.110}), (\ref{est3.140}) and (\ref{npc.250}).
In order to prove (\ref{est3.145}) we first observe from 
(\ref{solq2.100}) and (\ref{solq2.110})
that
\begin{eqnarray}
\label{est3.180}
\left\|
\left( D_{\omega} \tilde{h} \right)^{-1} (\lambda^{(0)})
\right\| \leq \frac{1}{d_V},
\end{eqnarray}
where $\| \cdot \|$ denotes the operator norm for linear maps
$({\Bbb R}^{\nu}, | \cdot |) \to ({\Bbb R}^{\nu}, | \cdot |)$. 
Using (\ref{est3.10}) and (\ref{est3.100}) we obtain
\begin{eqnarray}
\label{est3.190}
\left\|
\left( D_{\omega} \tilde{h} \right) (\lambda) - 
\left( D_{\omega} \tilde{h} \right) (\lambda^{(0)})
\right\| \leq  \frac{d_V}{2} \;\;\; \mbox{ for all } 
\lambda \in B_{\rho_1}(\lambda^{(0)}).
\end{eqnarray}
This proves (\ref{est3.145}).

\noindent
{\em (b)}
It suffices to show that for every $a \in B_{\tilde{\rho}_1}(0)$, the map
\begin{eqnarray}
\label{est3.200}
T_{a}: B_{\rho_1}(\omega^{(0)}) \to {\Bbb R}^{\nu}; \;\;
T_a(\omega) = \omega - (D_{\omega} \tilde{h})^{-1}(\lambda^{(0)}) 
\tilde{h}(a, \omega),
\end{eqnarray}
defines a contraction on $B_{\rho_1}(\omega^{(0)})$. Indeed,
by the Banach fixed point theorem this would imply that for every
$a \in B_{\tilde{\rho}_1}(0)$ there exists an unique 
$\omega(a) \in B_{\rho_1}(\omega^{(0)})$ satisfying 
$\tilde{h}(a, \omega(a)) = 0$. The smoothness of the such defined function
$\omega$ would then follow from a standard implicit function theorem (cf.
(\ref{est3.145})). 

We will now show that $T_a$ is a contraction on
$B_{\rho_1}(\omega^{(0)})$ by verifying the estimates
(\ref{est3.210}), (\ref{est3.220}) below for all
$a \in B_{\tilde{\rho}_1}(0)$.
\begin{eqnarray}
\label{est3.210}
|T_a(\omega_2) - T_a(\omega_1)| &\leq& \frac{1}{2} |\omega_2 - \omega_1|, \\
\label{est3.220}
|T_a(\omega^{(0)}) - \omega^{(0)}| &\leq& \frac{\rho_1}{2}.
\end{eqnarray}
In order to verify (\ref{est3.210}), observe that
\begin{eqnarray}
\nonumber
T_a(\omega_2) - T_a(\omega_1) =
\Lambda^{-1}
\left( 
\int_0^1 D_{\omega} \tilde{h}(\lambda^{(0)}) - 
D_{\omega} \tilde{h}(a, \omega_1 + t(\omega_2 - \omega_1)) dt
\right)
(\omega_2 - \omega_1),
\end{eqnarray}
and use (\ref{est3.180}), (\ref{est3.190}). Estimate (\ref{est3.220})
follows from (\ref{est3.180}), (\ref{solq2.5}), (\ref{est3.100}) and 
(\ref{est3.20}).

\noindent
{\em (c)}
We will prove (\ref{est3.50}) by induction with respect to 
$k \equiv |\alpha|_1$.

$\underline{k=1:}$ Implicit differentiation of $\tilde{h}(a, \omega(a)) = 0$
with respect to $a_j$ gives
\begin{eqnarray}
\nonumber
\partial_{a_j} \omega = 
- (D_{\omega} \tilde{h})^{-1} \partial_{a_j} \tilde{h}.
\end{eqnarray} 
The claim follows from (\ref{est3.145}), 
(\ref{est3.100}) and $D_4 \geq \frac{2}{d_V} D_{4, 1}$ (see (\ref{npc.260})).

$\underline{k-1 \to k:}$
Let $2 \leq k \leq B_1$ and let $\alpha$ be a multi-index with 
$|\alpha|_1 = k$.
Using proposition \ref{Pcr4.1} we obtain 
\begin{eqnarray}
\label{est3.150}
0 = \frac{\partial^{\alpha} \tilde{h}}{\partial a^{\alpha}} (a, \omega(a))
= \left( D_{\omega} \tilde{h} \right) \partial^{\alpha} \omega +
\sum_{s \in S_0(\alpha)}
\left( \partial^{\gamma^{(s)}} \tilde{h}  \right)
\prod_{l=1}^{l^{(s)}} \partial^{\alpha_l^{(s)}} \omega_{i_l^{(s)}}, 
\end{eqnarray}
where $\# S_0(\alpha) \leq (\nu + k)^k$, and for all $s \in S_0(\alpha)$
we have $|\gamma^{(s)}|_1 \leq k$,
$1 \leq |\alpha_l^{(s)}|_1 \leq k-1$,
$\sum_{1 \leq l \leq l^{(s)}} |\alpha_l^{(s)}|_1 \leq k$.
Define for $p \in {\Bbb N}$
\begin{eqnarray}
\nonumber
C(p) := \left( D_4 N_0^{B_1 E_{\rho}} \right)^{2p-1}.
\end{eqnarray}
Using (\ref{est3.150}), (\ref{est3.145}), (\ref{est3.100}), 
the induction hypothesis, and (\ref{npc.260})
we obtain
\begin{eqnarray}
\nonumber
|\partial^{\alpha} \omega|
&\leq&
\frac{2}{d_V} (\nu + k)^k D_{4, 1} N_0^{E_{\rho} k}
\left( \max_{s \in S_0(\alpha)} \prod_{l=1}^{l^{(s)}}
C(|\alpha_l^{(s)}|_1) \right) \\
\nonumber
&\leq&
D_4 N_0^{B_1 E_{\rho}}
\left( \max_{1 \leq p_i \leq k-1, \; p_1 + \ldots +p_l = k} 
C(p_1) \cdot \ldots C(p_l) \right) \\
\nonumber
&\leq&
D_4 N_0^{B_1 E_{\rho}} \left( D_4 N_0^{B_1 E_{\rho}} \right)^{2k-2}
\\
\nonumber
&\leq& C(k).
\end{eqnarray}
We have completed the proof of proposition \ref{Pest3.1} and have thus
also shown theorem \ref{Tsolq.1}.
\end{proof}

\section{Estimates on the measure of the set of resonant parameters} 
\label{est}

In this section we use the non-degeneracy assumption A4 to show that
the function $\omega$ introduced in the previous section is not
degenerate to second order. We use this fact to derive lower bounds on the
measure of the set
$\{ a \in B_{\tilde{\rho}_1}(0) : (a, \omega(a)) \in {\cal N}^{\infty} \}$
(see lemma \ref{Lest4.1}).

\noindent
{\bf Estimates of lemma \ref{hyp} used in section \ref{est}:}

\begin{eqnarray}
\label{est1.200}
2 (|\alpha_2| + 1) \frac{80 N_1^{\tau}}{d_V^2 s} \leq 
e^{\frac{15}{16} N_1^c} \\
\label{est1.210}
128 e^{2 \nu} D_{\psi}(2) N^{2 A E_{\rho}} \leq \frac{1}{80(|\alpha_2| + 1)d_V}
e^{\frac{1}{16} N^c} \;\; \mbox{ for all } N \geq N_1 \\
\label{est1.220}
\frac{1}{16} N_1^c (N_1^{(A-1)c} - 1) \geq \log 2
\\
\label{est2.10}
4 \sqrt{\nu} e^{-\frac{1}{16} N_1^c} \leq d_{min} \\
\label{est2.20}
4 \nu \sqrt{\nu}
\left(D_4 N_0^{B_1 E_{\rho}}\right)^5 \tilde{\rho}_1 \leq d_{min}
\\
\label{est3.470}
\nu \left( D_4 N_0^{B_1 E_{\rho}} \right)^3 \sqrt{\rho_1} \leq 1 \\
\label{est3.475}
B_1 \geq 8 \nu^2 \\
\label{est3.480}
\left(8 \nu^{2.5}\right)^{4 \nu^2} 
\left( D_4 N_0^{B_1 E_{\rho}} \right)^{12 \nu^2}  
\sqrt{\rho_1} \leq \frac{1}{4} \left(\frac{d_{min}}{4}\right)^{4 \nu^2}\\
\label{est3.490}
D_P (9 \nu^2)^{11 \nu^2} [(5 \nu^2)!] (D_4 \nu)^{16 \nu^2}
N_0^{16 \nu^2 B_1 E_{\rho} + 32 \nu^2 (B_0 + \nu) E_M - A} \leq
\frac{1}{2} 
\left( \frac{d_{min}}{4} \right)^{4 \nu^2}
\\
\label{est4.2}
4 e^{-\frac{1}{16} N_1^c} \leq \tilde{d}_{min}\\
\label{est4.4}
4 \nu \left(D_4 N_{0}^{B_1 E_{\rho}}\right)^5 \tilde{\rho}_1
\leq \tilde{d}_{min}\\
\label{est4.6}
\nu^2 \left(D_4 N_{0}^{B_1 E_{\rho}}\right)^3 \tilde{\rho}_1
\leq 1\\
\label{est4.7}
4 \nu \tilde{\rho}_1 \leq 1\\
\label{est4.8}
(\nu+2)B_1 \sqrt{\frac{8(D_K+1)}{d_{min}}} 5^{\nu} 2^{7 \nu + 1} \leq N_2
\\
\label{est4.9}
\delta_{j+1}^{\frac{1}{B_1}} (\tilde{\rho}_j)^{\nu-1}
N_{j+1}^{\nu + 2} N_j^{4 \nu} \delta_j^{-2} \leq
2^{-j} \tilde{\kappa} \tilde{\rho}_j^{\nu}
\end{eqnarray}

\subsection{Second order Taylor expansion for the function $\omega$.}
\label{est1}

\begin{Proposition}
\label{Pest1.1}
Let $\omega$ be the function as defined in theorem \ref{Tsolq.1}. Then
\begin{eqnarray}
\label{est1.5}
\mbox{ (a) }&& \frac{\partial \omega}{\partial a_j}(0) = 0 \;\;
\mbox{ for } 1 \leq j \leq \nu, \\
\label{est1.10}
\mbox{ (b) }&&  
\left| \frac{\partial^2 \omega}{\partial a_j \partial a_k}(0) -
\delta_{j, k} \Omega^{(j)} \right| \leq e^{-\frac{1}{16} N_1^c}
\;\; \mbox{ for } \;\; 
1 \leq j, k \leq \nu, 
\end{eqnarray}
where 
\begin{eqnarray}
\label{est1.15}
\Omega_l^{(j)} =
\left\{
\begin{array}{ll}
6 \Lambda_l^{-1} \left(
2 \alpha_2^2 \frac{1}{V(\omega^{(0)})(2 e_l)} - 3 \alpha_3
\right) & \mbox{ if } l=j \\
4 \Lambda_l^{-1} \left(
2 \alpha_2^2 \left(
\frac{1}{V(\omega^{(0)})(e_l+e_j)} 
+ \frac{1}{V(\omega^{(0)})(e_l-e_j)}
\right)
- 3 \alpha_3
\right) & \mbox{ if } l \neq j
\end{array}
\right.
\end{eqnarray}
\end{Proposition}

\begin{proof}
{\em (a):}
Statement (c) of theorem \ref{TP} implies
\begin{eqnarray}
\label{est1.50}
\partial_{a_j} v(0, \omega)=0 \;\; \mbox{ for all } \omega \in 
B_{\rho_1}(\omega^{(0)}).
\end{eqnarray}
Denote for $1 \leq j \leq \nu$ 
\begin{eqnarray}
\label{est1.60}
\tilde{e}_j : = \delta_{e_j} + \delta_{-e_j} \in X_{1, c}.
\end{eqnarray}
The definition of $\varphi(a)$ (see (\ref{ova.120})) and $u(a, \omega)$ 
(see (\ref{solq2.20})) implies
\begin{eqnarray}
\label{est1.70}
\partial_{a_j} u(0, \omega) = \tilde{e}_j \;\; \mbox{ for all } \omega \in 
B_{\rho_1}(\omega^{(0)}).
\end{eqnarray}
Using (\ref{solq1.45}) and (\ref{solq2.25}) -- (\ref{solq2.45}) we obtain
\begin{eqnarray}
\nonumber
\frac{\partial \tilde{h}_k}{\partial a_j}(0, \omega^{(0)}) &=&
\frac{\partial^2 h_k}{\partial a_k \partial a_j}(0, \omega^{(0)}) \\
\nonumber
&=&
\frac{\partial^2 n_k}{\partial a_k \partial a_j}(0, \omega^{(0)}) \\
\nonumber
&=& (D^2 W)(0)[\partial_{a_k} u, \partial_{a_j} u] (e_k) \\
\label{est1.80}
&=&
2 \alpha_2 (\tilde{e}_k * \tilde{e}_j)(e_k) = 0.
\end{eqnarray}
Differentiating $\tilde{h}(a, \omega(a))=0$ with respect to $a_j$,
$1 \leq j \leq \nu$, and using (\ref{solq2.110}), (\ref{est1.80}) yields
\begin{eqnarray}
\label{est1.90}
\frac{\partial \omega}{\partial a_j}(0) = - \Lambda^{-1}
\frac{\partial \tilde{h}}{\partial a_j}(0, \omega^{(0)}) = 0.
\end{eqnarray}

{\em (b):} Differentiating $\tilde{h}(a, \omega(a))=0$ with respect to $a_j$
and $a_k$,
$1 \leq j, k \leq \nu$ and using the already proven claim (a) 
together with (\ref{solq2.110}) leads to
\begin{eqnarray}
\label{est1.100}
\frac{\partial^2 \omega}{\partial a_k \partial a_j}(0) = - \Lambda^{-1}
\frac{\partial^2 \tilde{h}}{\partial a_k \partial a_j}(0, \omega^{(0)}).
\end{eqnarray}
It follows from (\ref{solq1.45}),  (\ref{solq2.25}) -- (\ref{solq2.45}),
(\ref{solq2.20}) (\ref{ova.120}) and (\ref{est1.70}) that
\begin{eqnarray}
\label{est1.110}
\frac{\partial^2 \tilde{h}_l}{\partial a_k \partial a_j}(0, \omega^{(0)}) &=&
\frac{\partial^3 n_l}{\partial a_l \partial a_k \partial a_j}
(0, \omega^{(0)}) \\
\nonumber
&=& \left(DW(0)[\partial^3_{a_l, a_k, a_j}u] + 
(D^2 W)(0)[\partial^2_{a_l, a_k} u, \partial_{a_j} u] +
(D^2 W)(0)[\partial^2_{a_l, a_j} u, \partial_{a_k} u] \right. \\
\nonumber
&& +
\left.
(D^2 W)(0)[\partial^2_{a_k, a_j} u, \partial_{a_l} u] +
(D^3 W)(0)[\partial_{a_l} u, \partial_{a_k} u, \partial_{a_j} u]
\right) (e_l) \\
\nonumber
&=&
2 \alpha_2
\left(
(\partial^2_{a_l, a_k} v) * \tilde{e}_j +
(\partial^2_{a_l, a_j} v) * \tilde{e}_k +
(\partial^2_{a_k, a_j} v) * \tilde{e}_l
\right)(e_l) 
+ 6 \alpha_3
\left(
\tilde{e}_l * \tilde{e}_k * \tilde{e}_j
\right) (e_l) 
\end{eqnarray}
Next we compute $(\partial^2_{a_k, a_j} v) ( \lambda^{(0)})$.
Induction statement $({\cal IS})(1c)_{j=2}$ and Cauchy's integral formula
give
\begin{eqnarray}
\label{est1.115}
\left\|
\partial^{\beta} {\cal P}(v_2(\cdot), \cdot)(\lambda^{(0)})
\right\|_{1/4, c} \leq 2 e^{- N_1^c} \;\;\; \mbox{ for }
|\beta|_1 = 2.
\end{eqnarray}
Recalling induction statements $({\cal IS})(1d)_{j=2}$ and
$({\cal IS})(1a)_{j=1}$ it is clear that 
$v_2(\lambda^{(0)}) = 0$ and also 
$(\partial_{a_j} v_2)(\lambda^{(0)}) = 0$ for $1 \leq j \leq \nu$.
Hence,
\begin{eqnarray}
\nonumber
\frac{\partial^2 {\cal P}(v_2(\cdot), \cdot)}{\partial a_k \partial a_j}
(\lambda^{(0)}) =
P \left[
D(\omega^{(0)}) 
\frac{\partial^2 v_2}{\partial a_k \partial a_j} (\lambda^{(0)}) +
2 \alpha_2 (\tilde{e}_k * \tilde{e}_j)
\right].
\end{eqnarray}
Recall furthermore from (\ref{cv2.160}) that  
\begin{eqnarray}
\label{est1.117}
| V(\omega^{(0)})(m)| > \frac{d_V s}{2 N_1^{\tau}} \;\; \mbox{ for }
m \in B_{N_1}(0) \setminus ({\cal S} \cup \{ 0 \} ).
\end{eqnarray}
Employing induction statement
$({\cal IS})(1f)_{j=2}$, (\ref{est1.200}), (\ref{est1.115}) and 
(\ref{est1.117}) we arrive at
\begin{eqnarray}
\label{est1.120}
\left\|
\frac{\partial^2 v_2}{\partial a_k \partial a_j} (\lambda^{(0)}) +
2 \alpha_2 D^{-1}(\omega^{(0)}) P (\tilde{e}_k * \tilde{e}_j)
\right\|_{1/4, c} \leq \frac{1}{40(|\alpha_2| + 1)d_V} e^{-\frac{1}{16}N_1^c}.
\end{eqnarray}
Induction statements $({\cal IS})(1e)_{j \geq 2}$ together
with (\ref{est1.210}) and (\ref{est1.220}) show that
\begin{eqnarray}
\label{est1.130}
\left\|
\frac{\partial^2 v}{\partial a_k \partial a_j} (\lambda^{(0)}) -
\frac{\partial^2 v_2}{\partial a_k \partial a_j} (\lambda^{(0)})
\right\|_{1/4, c} \leq \frac{1}{40(|\alpha_2| + 1)d_V} e^{-\frac{1}{16}N_1^c}.
\end{eqnarray}
From (\ref{est1.110}), (\ref{est1.120}) and (\ref{est1.130}) we learn
\begin{eqnarray}
\nonumber
|
\partial^2_{a_k, a_j} \tilde{h}_l (\lambda^{(0)})&+&
4 \alpha_2^2 \left([D^{-1}P(\tilde{e}_l * \tilde{e}_k)]*\tilde{e}_j +
[D^{-1}P(\tilde{e}_l * \tilde{e}_j)]*\tilde{e}_k +
[D^{-1}P(\tilde{e}_k * \tilde{e}_j)]*\tilde{e}_l
\right)(e_l)\\
\label{est1.150}
&-& 6 \alpha_3 (
\tilde{e}_l * \tilde{e}_k * \tilde{e}_j
)(e_l)| \leq \frac{1}{d_V} e^{-\frac{1}{16} N_1^c}.
\end{eqnarray}
Statement (b) of proposition \ref{Pest1.1} now follows from 
(\ref{solq2.100}), (\ref{est1.100}) and (\ref{est1.150}) 
by an explicit calculation.
\end{proof}

\subsection{Application of the non-degeneracy condition}
\label{est2}

\begin{Proposition}
\label{Pest2.1}
For every $m \in {\Bbb Z}^{\nu}
\setminus \{ 0 \}$ there exists a $\sigma \in \{ 1, \ldots, \nu \}$ 
such that for all $a_0 \in {\Bbb R}^{\nu}$
\begin{eqnarray}
\label{est2.50}
\left|
\left( \frac{d}{dt} \right)^2 <\omega(a_0 + t e_{\sigma}), m> \right| \geq
\frac{d_{min}}{2} |m|_2 \;\;\;\; \mbox{ for all $t$ with } 
|a_0 + e_{\sigma} t| < \tilde{\rho}_1.
\end{eqnarray}
\end{Proposition}

\begin{proof}
Recall that the $j$-th row of the matrix $\Omega$, defined in (\ref{smr2.90}),
is given by the vector $\Omega^{(j)}$ of proposition \ref{Pest1.1} 
\begin{eqnarray}
\label{est2.100}
\Omega_{j, l} = \Omega_l^{(j)} \quad \mbox{ (cf. (\ref{npc.70}))} 
\end{eqnarray}
Recall further the definition of $d_{min}$ in (\ref{npc.130}),
\begin{eqnarray}
\label{est2.110}
d_{min} = 
\min_{\xi \in {\Bbb R}^{\nu}: |\xi|_2 = 1} |\Omega \xi| > 0.
\end{eqnarray}
We infer from statement (c) of proposition \ref{Pest3.1} and (\ref{est2.20})
that
\begin{eqnarray}
\label{est2.120}
\left|
\frac{\partial^2 \omega}{\partial a_j^2}(a) -
\frac{\partial^2 \omega}{\partial a_j^2}(0) 
\right|_2 \leq \frac{d_{min}}{4} \;\;\; \mbox{ for all }
|a| < \tilde{\rho}_1, \; 1 \leq j \leq \nu.
\end{eqnarray}
Observe further that condition (\ref{est2.10}) together with statement
(b) of proposition \ref{Pest1.1} imply
\begin{eqnarray}
\label{est2.125}
\left|
\frac{\partial^2 \omega}{\partial a_j^2}(0) - \Omega^{(j)} \right|_2
\leq \frac{d_{min}}{4} \;\;\; \mbox{ for all } 1 \leq j \leq \nu.
\end{eqnarray}
Fix $m \in {\Bbb Z}^{\nu} \setminus \{ 0 \}$. By (\ref{est2.110}) there
exists $\sigma \in \{ 1, \ldots, \nu \}$ with 
\begin{eqnarray}
\label{est2.130}
\left| \left( \Omega \frac{m}{|m|_2} \right)_{\sigma} 
\right| \geq d_{min}.
\end{eqnarray}
Then for $|a_0 + t e_{\sigma}| < \tilde{\rho}_1$:
\begin{eqnarray}
\nonumber
\left|
\left( \frac{d}{dt} \right)^2 <\omega(a_0 + t e_{\sigma}), m> \right|
&=&
\left| <(\partial_{a_{\sigma}}^2 \omega)(a_0 + t e_{\sigma}), \frac{m}{|m|_2}> \right|
|m|_2 \\
\nonumber
&\geq&
\left(
\left| <\Omega^{(\sigma)}, \frac{m}{|m|_2}> \right| -
\frac{d_{min}}{2} \right)
|m|_2 \\
\nonumber
&=&
\left(
\left| \left( \Omega \frac{m}{|m|_2}\right)_{\sigma} \right| -
\frac{d_{min}}{2} \right)
|m|_2 \\
\nonumber
&\geq&
\frac{d_{min}}{2} |m|_2.
\end{eqnarray}
\end{proof}

We will now discuss an important consequence of proposition \ref{Pest2.1}.
To that end, fix 
$(p, \vartheta) \in {\cal POL}$ with
\begin{eqnarray}
\label{est3.400}
p(\theta, \lambda) = \theta^d + \sum_{k<d} b_k(\lambda) \theta^k,
\end{eqnarray}
Let $l \in {\Bbb Z}$, $m \in {\Bbb Z}^{\nu} \setminus \{ 0 \}$. Choose
$\sigma \in \{1, \ldots, \nu \}$ as a function of $m$ as 
described in proposition \ref{Pest2.1}. Furthermore, we set
\begin{eqnarray}
\label{est3.410}
a(t) &:=& a_0 + e_{\sigma} t \;\;\; \mbox{ for $t$ with } |a(t)| < 
\tilde{\rho}_1, \\
\label{est3.420}
\theta(t) &:=& <\omega(a(t)), m> - \vartheta - 2 \pi l, \\
\label{est3.430}
\lambda(t) &:=& (a(t), \omega(a(t))), \\
\label{est3.440}
f(t) &:=& (\theta(t), \lambda(t)), \\
\label{est3.450}
g(t) &:=& p(f(t)).
\end{eqnarray}
\begin{Proposition}
\label{Pest3.2}
For $|m| \geq 2 N_1$ the following estimate (\ref{est3.500})
holds for all $t \in 
{\Bbb R}$ satisfying $|a(t)| < \tilde{\rho}_1$ and $|\theta(t)| < 1$:
\begin{eqnarray}
\label{est3.500}
\left|
\left( \frac{d}{dt} \right)^{2d}g(t)\right| \geq (2d)!
\left( \frac{d_{min}}{8} |m|_2 \right)^d .
\end{eqnarray}
\end{Proposition}
\begin{proof}
According to proposition \ref{Pcr3.1} (see also (\ref{cr3.10}))
\begin{eqnarray}
\label{est3.510}
g^{[2d]} = \sum_{1 \leq |\alpha|_1 \leq 2d}
p^{[\alpha]}
\sum_{\scriptsize
{\scriptstyle
\begin{array}{c}
\delta^{i} \in {\Bbb N}^{\alpha_i} \\
\sum_{i} |\delta^i|_1 = 2d
\end{array}
}
}
\prod_{i=1}^{2 \nu + 1}
\prod_{j=1}^{\alpha_i}
f_i^{[\delta^{i}_j]}
\end{eqnarray}
From (\ref{est3.400}) it is obvious that
\begin{eqnarray}
\label{est3.520}
p^{[\beta]}(\theta, \lambda) =
\left\{
\begin{array}{ll}
0 & \mbox{ if } \beta_1 > d, \\
1 & \mbox{ if } \beta = (d, 0), \\
0 & \mbox{ if } \beta_1 =d \mbox{ and } |\beta|_1 > d. 
\end{array}
\right.
\end{eqnarray}
We use (\ref{est3.520}) to express the sum in (\ref{est3.510}) 
\begin{eqnarray}
\nonumber
g^{[2d]} &=& p^{[(d, 0)]} (\theta^{[2]})^d \\
\nonumber
&&+ p^{[(d, 0)]} 
\sum_{\scriptsize
{\scriptstyle
\begin{array}{c}
k_1 + \ldots + k_d = 2d , k_j \in {\Bbb N}\\
\exists i \in \{1, \ldots, d\}: k_i = 1
\end{array}
}
}
\theta^{[k_1]} \cdots \theta^{[k_d]}
\\
\nonumber
&&+
\sum_{\scriptsize
{\scriptstyle
\begin{array}{c}
1 \leq |\alpha|_1 \leq 2d \\
\alpha_1 < d
\end{array}
}
}
p^{[\alpha]}
\sum_{\scriptsize
{\scriptstyle
\begin{array}{c}
\delta^{i} \in {\Bbb N}^{\alpha_i} \\
\sum_{i} |\delta^i|_1 = 2d
\end{array}
}
}
\prod_{i=1}^{2 \nu + 1}
\prod_{j=1}^{\alpha_i}
f_i^{[\delta^{i}_j]}
\\
\label{est3.530}
&=&I + II + III. 
\end{eqnarray}
It follows from proposition \ref{Pest2.1}, (\ref{est3.410}), 
(\ref{est3.420}), (\ref{est3.520}) and the choice of $\sigma$ that 
\begin{eqnarray}
\label{est3.540}
|I| \geq \left( \frac{d_{min}}{4} |m|_2 \right)^d.
\end{eqnarray}
In order to estimate the second term, we observe from 
propositions \ref{Pest1.1}, \ref{Pest3.1} and from (\ref{est3.470})
that for $|\beta|_1 = 1$ and $|a|< \tilde{\rho}_1$ 
\begin{eqnarray}
\label{est3.550}
|\partial^{\beta} \omega (a)| =
|\partial^{\beta} \omega (a) - \partial^{\beta} \omega (0)|
\leq \nu \left( D_4 N_0^{B_1 E_{\rho}} \right)^3 \rho_1 \leq \sqrt{\rho_1}.
\end{eqnarray}
This implies immediately that
\begin{eqnarray}
\label{est3.560}
|\theta^{[1]}(t)| \leq \sqrt{\nu \rho_1} |m|_2.
\end{eqnarray}
On the other hand, it follows from (\ref{est3.50}), (\ref{est3.475}) 
that for $1 \leq k \leq 2d \leq 8 \nu^2$ (see (\ref{npc.490}))
\begin{eqnarray}
\label{est3.570}
|\theta^{[k]}(t)| \leq  
\sqrt{\nu} \left( D_4 N_0^{B_1 E_{\rho}} \right)^{2k-1} |m|_2.
\end{eqnarray}
Expression II is a sum with less than 
$(2d)^d \leq (8 \nu^2)^{4 \nu^2}$ terms. 
Together with (\ref{est3.480}) we obtain
\begin{eqnarray}
\label{est3.580}
|II| \leq \frac{1}{4}\left( \frac{d_{min}}{4} |m|_2 \right)^d.
\end{eqnarray}

Expression III can be written as a sum of products with at most
$(2 d +1)^{2\nu + 1} (2d)^{2d} \leq 
(9 \nu^2)^{11 \nu^2}$ summands. Moreover, (\ref{npc.490}) implies
\begin{eqnarray}
\nonumber
\left| p^{[\alpha]} \right|
\leq D_P [(d+1)!] M_0^{32 \nu^2 (B_0 + \nu)} \leq
D_P(5 \nu^2)! N_0^{32 \nu^2 (B_0 + \nu) E_M}.
\end{eqnarray}
The assumption that $|m| \geq 2 N_1$ yields
\begin{eqnarray}
\nonumber
\prod_{i=1}^{2 \nu + 1}
\prod_{j=1}^{\alpha_i}
f_i^{[\delta^{i}_j]} &\leq&
(\sqrt{\nu} |m|_2)^{d-1} (D_4 N_0^{B_1 E_{\rho}})^{4d} \\
\nonumber
&\leq& \frac{1}{2 N_1} \nu^{2 \nu^2}
(D_4 N_0^{B_1 E_{\rho}})^{16 \nu^2} |m|_2^d.
\end{eqnarray}
Finally, using (\ref{est3.490}) we obtain
\begin{eqnarray}
\label{est3.680}
|III| \leq \frac{1}{4}\left( \frac{d_{min}}{4} |m|_2 \right)^d.
\end{eqnarray} 
\end{proof}

\subsection{The excision procedure}
\label{est4}

Recall from theorem \ref{TP} that we have solved the ${\cal P}$ --
equation only for those values of the parameter $\lambda$  
which satisfy $\lambda \in {\cal N}^{\infty}$. 
The ${\cal Q}$ -- equation is solved if $\lambda$ is of the 
form $\lambda = (a, \omega(a))$. Combining both observations,
we obtain a solution of our original equation
(\ref{ova.80}) for those values of the parameter $a$ which lie in the set
\begin{eqnarray}
\label{est4.10}
{\cal M}^{\infty} := \{ a \in B_{\tilde{\rho}_1}(0):
(a , \omega(a)) \in {\cal N}^{\infty} \}. 
\end{eqnarray}
The goal of this section 
is to show that the set of bad parameters for which we have no solutions 
is relatively small.
More precisely, we will show that
\begin{Lemma}
\label{Lest4.1} 
$\mbox{ vol}( B_{\tilde{\rho}_1}(0) \setminus {\cal M}^{\infty} )
\leq \tilde{\kappa} \tilde{\rho}_1^{\nu}$.
\end{Lemma}
Recall that $\tilde{\kappa}$ was defined in (\ref{npc.270}).
In order to prove lemma \ref{Lest4.1} we first recall 
the definition of ${\cal N}^{\infty}$ in (\ref{sp.10}) as
the intersection of a sequence of nested sets. Thus
\begin{eqnarray}
\label{est4.20}
{\cal M}^{\infty} &=& \bigcap_{j=1}^{\infty} {\cal M}^{(j)},
\;\;\; \mbox{ where }\\
\label{est4.30}
{\cal M}^{(j)} &:=& \{ a \in B_{\tilde{\rho}_1}(0):
(a , \omega(a)) \in \overline{{\cal N}^{(j)}} \}. 
\end{eqnarray}
${\cal M}^{(j)}$ is again a sequence of nested sets with
${\cal M}^{(1)} = B_{\tilde{\rho}_1}(0)$. Therefore
lemma \ref{Lest4.1} is a consequence of the following
proposition \ref{Pest4.1}.  
\begin{Proposition}
\label{Pest4.1}
Let $j \geq 1$. Then 
$\mbox{ vol}( {\cal M}^{(j)} \setminus {\cal M}^{(j+1)} )
\leq 2^{-j} \tilde{\kappa} \tilde{\rho}_1^{\nu}$.
\end{Proposition}

\begin{proof} ({\em Proposition \ref{Pest4.1}}). 
Fix $j \in {\Bbb N}$.
Since $\tilde{\rho}_1/\tilde{\rho}_j \in {\Bbb N}$ (by definitions
(\ref{npc.72}) -- (\ref{npc.75}), (\ref{npc.400}) --
(\ref{npc.425}) and the choice of $n_0$ in Lemma 
\ref{hyp}) there exists a cube decomposition of the set
$B_{\tilde{\rho}_1}(0)$ into cubes of radius $\tilde{\rho}_j$, i.e.
\begin{eqnarray}
\nonumber
\overline{B_{\tilde{\rho}_1}(0)} &=& 
\bigcup_{y \in Y} 
\overline{B_{\tilde{\rho}_j}(a^{(y)})}, \\
\nonumber
B_{\tilde{\rho}_j}(a^{(y_1)}) \cap B_{\tilde{\rho}_j}(a^{(y_2)}) &=&
\emptyset \;\; \mbox{ for } y_1, y_2 \in Y, \; y_1 \neq y_2.
\end{eqnarray}
To prove proposition \ref{Pest4.1} it suffices to show that for 
each $y \in Y$
\begin{eqnarray}
\label{est4.100}
\mbox{ vol}( ({\cal M}^{(j)} \setminus {\cal M}^{(j+1)}) \cap
B_{\tilde{\rho}_j}(a^{(y)}) )
\leq 2^{-j} \tilde{\kappa} \tilde{\rho}_j^{\nu}.
\end{eqnarray}
Fix $y \in Y$. 
%For $\tilde{k} \in K^{(j+0.5)}$ we write
%$\lambda_{\tilde{k}}^{(j+0.5)} = (a_{\tilde{k}}^{(j+0.5)}, 
%\omega_{\tilde{k}}^{(j+0.5)})$. 
Define
\begin{eqnarray}
\label{est4.110}
K_y := \{ \tilde{k} \in K^{(j+0.5)}: 
(a, \omega(a)) \in B_{\tilde{\rho}_j}(\lambda_{\tilde{k}}^{(j+0.5)})
\; \mbox{ for some } a \in B_{\tilde{\rho}_j}(a^{(y)})
\}.
\end{eqnarray}
In view of induction statement $({\cal IS})(2b)_{j+1}$ the relation
\begin{eqnarray}
\label{est4.120}
a \in ({\cal M}^{(j)} \setminus {\cal M}^{(j+1)}) \cap 
B_{\tilde{\rho}_j}(a^{(y)})
\end{eqnarray}
implies that
there exist $\tilde{k} \in K_y$, $(p, \vartheta) \in 
{\cal POL}^{(j+1)}_{\tilde{k}}$, $m \in {\Bbb Z}^{\nu}$, $l \in {\Bbb Z}$
such that the following holds
\begin{eqnarray}
\label{est4.130}
(a, \omega(a)) &\in& 
\overline{B_{\tilde{\rho}_j}(\lambda_{\tilde{k}}^{(j+0.5)})}, \\
\label{est4.140}
2 N_j < &|m|& \leq  2 N_{j+1}, \\
\label{est4.150}
|<\omega(a), m> - \vartheta - 2 \pi l| &<& 1, \\
\label{est4.160}
|p(<\omega(a), m> - \vartheta - 2 \pi l, (a, \omega(a)) )| &<& (D_K + 1) 
\delta_{j+1}.
\end{eqnarray}
Fix $\tilde{k} \in K_y$, $(p, \vartheta) \in 
{\cal POL}^{(j+1)}_{\tilde{k}}$, $m \in {\Bbb Z}^{\nu}$, $l \in {\Bbb Z}$
and set 
\begin{eqnarray}
\label{est4.170}
{\cal B} := \{a \in B_{\tilde{\rho}_j}(a^{(y)}):
(\ref{est4.130}) - (\ref{est4.160}) \;\mbox{ are satisfied } \}.
\end{eqnarray}
In order to estimate the (Lebesgue-) measure of the set
${\cal B}$ we use Fubini's theorem. Choose $\sigma \in
\{1, \ldots, \nu \}$ as a function of $m$ according to proposition
\ref{Pest2.1}. Denote by $H$ the hyperplane in ${\Bbb R}^{\nu}$
which contains the point $a^{(y)}$ and which is perpendicular to the 
vector $e_{\sigma}$. For $a_0 \in H \cap B_{\tilde{\rho}_j}(a^{(y)})$
we define 
\begin{eqnarray}
\label{est4.180}
{\cal B}_{a_0} := \{ t \in {\Bbb R}:
a_0 + t e_{\sigma} \in {\cal B} \}.
\end{eqnarray}
According to the notation we introduced in (\ref{est3.410}) -- 
(\ref{est3.450}) let
\begin{eqnarray}
\nonumber
a(t) &=& a_0 + t e_{\sigma},\\
\nonumber
\theta(t) &=& < \omega( a(t) ), m > - \vartheta - 
2 \pi l.
\end{eqnarray}
Recall from proposition \ref{Pest2.1} that $\theta''$ does not change 
sign. This implies immediately that all values of $t$ for which
condition (\ref{est4.150}) is satisfied form a set which is the union
of at most two intervals (including the possibility that the set is empty).
Observe that proposition \ref{Pest1.1}(b), proposition \ref{Pest3.1}(c)
together with assumption A4, (\ref{npc.140}), (\ref{est4.2}) and (\ref{est4.4}) 
imply 
\begin{eqnarray}
\nonumber
\left| \frac{d^2}{dt^2} \omega_l(a(t)) \right| \geq \frac{\tilde{d}_{min}}{2}
\end{eqnarray}
for all $1 \leq l \leq \nu$ and $|a(t)| \leq \tilde{\rho}_1$. Hence
the set of those values of $t$ for which condition (\ref{est4.130}) is 
satisfied is the union of at most $\nu + 1$ intervals
(including the possibility of the empty set). 
Consequently, the values of $t$ for which both conditions (\ref{est4.130})
and (\ref{est4.150}) are satisfied form a set which is the union of
at most $\nu + 2$ intervals. In each of these intervals we can apply
the excision lemma \ref{Lel.1} and proposition 
\ref{Pest3.2} to estimate the measure of the set ${\cal B}_{a_0}$. 
Using in addition that the degree $d$ of the polynomial $p$ is bounded
by $4 \nu^2 = B_1/2$ (see (\ref{npc.490}), (\ref{npc.90})) we 
obtain
\begin{eqnarray}
\label{est4.200}
\mbox{ vol}({\cal B}_{a_0}) \leq 
(\nu + 2)  2 B_1 \sqrt{\frac{8 (D_K + 1)}{d_{min} |m|_2}} 
\delta_{j+1}^{\frac{1}{B_1}}.
\end{eqnarray}
To obtain an estimate on ${\cal B}$ we integrate (\ref{est4.200})
with respect to $a_0$ over the set $H \cap B_{\tilde{\rho}_j}(a^{(y)})$,
yielding 
\begin{eqnarray}
\label{est4.210}
\mbox{ vol}({\cal B}) \leq 
(\nu + 2)  2 B_1 \sqrt{\frac{8 (D_K + 1)}{d_{min}}} 
\delta_{j+1}^{\frac{1}{B_1}} (2 \tilde{\rho}_j)^{\nu - 1}.
\end{eqnarray}
Recall that we have fixed 
$\tilde{k} \in K_y$, $(p, \vartheta) \in 
{\cal POL}^{(j+1)}_{\tilde{k}}$, $m \in {\Bbb Z}^{\nu}$, $l \in {\Bbb Z}$
for the definition of the set ${\cal B}$.
In order to complete the proof we need to determine for how many values of
$\tilde{k}$, $(p, \vartheta)$, $m$ and $l$ the corresponding set ${\cal B}$ 
is not empty.
\begin{itemize}
\item[$\bullet$]
$\# m:$ condition (\ref{est4.140}) implies that the number of possible
lattice points is bounded by $(5 N_{j+1})^{\nu}$.
\item[$\bullet$]
$\# (p, \vartheta):$ Recall from induction statement $({\cal IS})(2c)_{j+1}$
that the set ${\cal POL}^{(j+1)}_{\tilde{k}}$ has at most
$2 (2 N_j)^{4 \nu} \delta_j^{-2}$ elements.
\item[$\bullet$]
$\# \tilde{k}:$ Using proposition \ref{Pest1.1}(a), 
proposition \ref{Pest3.1}(c) and (\ref{est4.6}) we learn that
\begin{eqnarray}
\label{est4.250}
\left| \frac{\partial \omega}{\partial a_l}(a) \right| 
\leq \frac{1}{\nu}
\;\;\; \mbox{ for all } 1 \leq l \leq \nu, \; |a| \leq \tilde{\rho}_1.
\end{eqnarray}
This implies that $|\omega(a) - \omega(a^{(y)})| \leq |a - a^{(y)}|$
which in turn shows that the set $K_y$ contains at most $2^{2 \nu}$ elements.
\item[$\bullet$]
$\# l:$ It follows 
from (\ref{est4.250}) that for $a_1$, $a_2 \in B_{\tilde{\rho}_j}(a^{(y)})$
and $|m| \leq 2N_{j+1}$ we have 
\begin{eqnarray}
\nonumber
|<\omega(a_1) , m> - <\omega(a_2), m>| \leq 4 \nu N_{j+1} \tilde{\rho}_j.
\end{eqnarray}
Using (\ref{est4.7}) we conclude that the number of integers
$l$, for which (\ref{est4.150}) is satisfied for some fixed $\vartheta$ and
$|m| \leq 2 N_{j+1}$ is bounded above by
\begin{eqnarray}
\nonumber
\frac{N_{j+1} + 2}{2 \pi} + 1 \leq N_{j+1}.
\end{eqnarray}
\end{itemize}
These estimates, together with (\ref{est4.210}), (\ref{est4.8}) and
(\ref{est4.9}) prove (\ref{est4.100}), completing the proof of proposition
\ref{Pest4.1}.
\end{proof}

\section{Proof of the main theorem \ref{Tsmr3.1}}
\label{pmt}
Our main result essentially follows from theorems \ref{TP}, \ref{Tsolq.1}
and lemma \ref{Lest4.1} via proposition \ref{Pfeq}.
\noindent
{\bf Estimates of lemma \ref{hyp} used in section \ref{pmt}:}
\begin{eqnarray}
\label{pmt.10}
D_3 N_0^{\tau+1} \tilde{\rho}_1^2 + 2 \nu e D_N \tilde{\rho}_1 < r_{F, b}\\
\label{pmt.20}
\nu \left( D_4 N_0^{B_1 E_{\rho}} \right)^3 \tilde{\rho}_1 < 1\\
\label{pmt.30}
\sin \frac{\omega^{(0)}_1 - \nu \tilde{\rho}_1}{2} > \tilde{\rho}_1
\end{eqnarray}

\begin{Notation}
\label{Npmt.1}
For every $z \in {\Bbb C}^{\nu}$ we denote its polar coordinates by
$(a, \zeta) \in ({\Bbb R}_{+} \cup \{ 0 \})^{\nu} 
\times [0, 2 \pi)^{\nu}$ which satisfy
\begin{eqnarray}
\nonumber
(z_1, \ldots, z_{\nu}) = (a_1 e^{i \zeta_1}, \ldots, a_{\nu} 
e^{i \zeta_{\nu}}).
\end{eqnarray}
\end{Notation}

Set
\begin{eqnarray}
\label{pmt.50}
\rho &:=& \tilde{\rho}_1, \\
\label{pmt.60}
{\cal Z} &:=& \{ z \in {\Bbb C}^{\nu}: 
a \in {\cal M}^{\infty} \}, \\
\label{pmt.70}
\omega && \mbox{ as defined in theorem \ref{Tsolq.1} }, \\
\label{pmt.80}
u(z)(m) &:=& ( \; \varphi(a)(m) + 
v(a, \omega(a))(m) \; ) e^{i <\zeta, m>},\\
\label{pmt.90}
\chi(\xi, z) &:=& \sum_{m \in {\Bbb Z}^{\nu} \setminus \{ 0 \} }
\frac{u(z)(m)}{- 2i \sin \frac{<\omega(a), m>}{2}} e^{i <m, \xi>}.
\end{eqnarray}
Note that the ambiguity in the definition of the polar coordinates
for $z=0$ plays no role in the definition (\ref{pmt.80}) since
$\varphi(0) + v(\lambda^{(0)}) = 0$.
We show now that the functions $x_n$ as defined through 
(\ref{smr3.10}) define a solution of system (\ref{smr.10}).
In order to prove this we fix $z \in {\cal Z}$ and
use proposition \ref{Pfeq} with $s_0 :=
s d_{\tau, c}$, $\omega := \omega(a)$ and $u := u(z)$.
We need to verify conditions (\ref{feq10}) -- (\ref{feq20}).
\begin{itemize}
\item[$\bullet$]
(\ref{feq10}): follows from theorem \ref{TP}(c), (\ref{ova.120}), 
and (\ref{pmt.10}).
\item[$\bullet$]
(\ref{feq15}): follows from theorem \ref{TP}(e) and (\ref{ova.120}).
\item[$\bullet$]
(\ref{feq17}): Using theorem \ref{TP}(f) and (\ref{ova.120})
we see that 
\begin{eqnarray}
\nonumber
u(-m) = u(z)(-m) = u(a)(-m)e^{i <\zeta, - m>} = 
\overline{u(a)(m)e^{i <\zeta, m>}} = \overline{u(z)(m)} = \overline{u(m)}.
\end{eqnarray}
\item[$\bullet$]
(\ref{feq25}): Recall that $a \in {\cal M}^{\infty}$ implies 
by definition (\ref{est4.10}) that $(a, \omega(a)) \in {\cal N}^{\infty}$ and 
the claim follows from theorem \ref{TP}(h) and the definition of
$s_0$ given above.
\item[$\bullet$]
(\ref{feq20}): Using (\ref{feq25}) condition (\ref{feq20}) is 
equivalent to satisfying $D(\omega)u(m) + W(u)(m) = 0$ for all
$m \in {\Bbb Z}^{\nu} \setminus \{ 0 \}$. Since $(a, \omega(a))
\in {\cal N}^{\infty}$ we recall from theorem \ref{TP}(b) and
theorem \ref{Tsolq.1} that
\begin{eqnarray}
\nonumber
D(\omega)u(a)(m) + W(u(a))(m) = 0 \;\;\; \mbox{ for all }
\; m \in {\Bbb Z}^{\nu} \setminus \{ 0 \}.
\end{eqnarray}
Observe that $u(z)(m) = u(a)(m)e^{i<\zeta, m>}$. It follows from
the definition of $W$ and property (\ref{fa1.57}) of convolutions that
$W(u(z))(m) = W(u(a))(m) e^{i<\zeta, m>}$. This proves (\ref{feq20}).
\end{itemize}
We have seen that for every $z \in {\cal Z}$ the functions $x_n$ 
defined through (\ref{smr3.10}) and (\ref{pmt.90}) solve system
(\ref{smr.10}) and are of type (\ref{smr.20}).
 
It remains to prove assertions (\ref{smr3.20}) -- (\ref{smr3.50}).
Claim (\ref{smr3.20}) follows from the definition of $\tilde{\kappa}$
in (\ref{npc.270}), from lemma \ref{Lest4.1}, and from proposition 
\ref{Pnpc.1}. 

In order to verify 
(\ref{smr3.30}) we first observe that proposition \ref{Pest3.1}(c)
together with proposition \ref{Pest1.1}(a) and (\ref{pmt.20}) imply
that $|\partial_{a_j} \omega(a)| < 1$ for all $|a| < \rho$. Using in addition
(\ref{pmt.30}) we conclude that the function
\begin{eqnarray}
\label{pmt.110}
a_j \mapsto \frac{a_j}{\sin\frac{\omega_j(a)}{2}}
\end{eqnarray}
is strictly increasing for all $1 \leq j \leq \nu$. Suppose now that 
$\chi(\cdot, z) = \chi(\cdot, \tilde{z})$. Evaluating the $e_j$ 
Fourier coefficients we obtain from theorem \ref{TP}(e) and 
(\ref{ova.120})
\begin{eqnarray}
\nonumber
\frac{z_j}{2 i \sin\frac{\omega_j(a)}{2}} =
\frac{\tilde{z}_j}{2 i \sin \frac{\omega_j(\tilde{a})}{2}}
\end{eqnarray}
for all $1 \leq j \leq \nu$.  Comparing moduli and arguments of these
complex numbers and using the strict monotonicity of the functions
defined in (\ref{pmt.110}) we conclude that $a_j = \tilde{a}_j$, 
$1 \leq j \leq \nu$, and therefore $z=\tilde{z}$.

Claim (\ref{smr3.40}) follows immediately from proposition \ref{Pest1.1} (a).
and proposition \ref{Pest3.1} (c)

Finally, claim (\ref{smr3.50}) is a consequence of statements (c) and (h) of
theorem \ref{TP} and from 
\begin{eqnarray}
\nonumber
\frac{z_j}{- 2 i \sin \frac{\omega_j(a)}{2}} - z_j' =
\frac{z_j}{- 2 i \sin \frac{\omega_j(a)}{2}} -
\frac{z_j}{- 2 i \sin \frac{\omega_j^{(0)}}{2}} =
{\cal O}(|\omega(a) - \omega^{(0)}|) = {\cal O}(|z|^2) \;\; 
\mbox{ by } (\ref{smr3.40}).
\end{eqnarray}
The proof of theorem \ref{Tsmr3.1} is complete. 

\newpage
\setcounter{page}{121}
\noindent
\begin{center}
{\huge \bf Appendix} \vspace{2cm}
\end{center}

\section{The assumptions of theorem \ref{Tsmr3.1}}
\label{ass}

In this section we study the assumptions of our main result which were
formulated in section \ref{smr2}. 
The first lemma shows that the entries of the matrix $\Omega$ as 
given in (\ref{smr2.90}) are well defined 
if assumptions A1 and A2 are satisfied. Then we prove that for a
generic set of force laws and sufficiently low frequencies we
may satisfy assumptions A1 -- A4 by varying the averaged lattice
spacing $b$. This result is formulated in theorem \ref{Tass.1}.
A description of the exceptional set of force laws
which is not covered by theorem \ref{Tass.1} is given in lemma \ref{Lass.3}.  

\begin{Lemma}
\label{Lass.1}
Let $F: {\Bbb R} \to {\Bbb R}$, $b \in {\Bbb R}$ and $\gamma \in {\Bbb R}_+$
be given 
such that assumptions A1 and A2 of section \ref{smr2} are satisfied.
Let $\nu$, $\omega^{(0)}$, $\Lambda_j$ and $V(\omega)(m)$ be defined as 
in (\ref{smr2.20}) -- (\ref{smr2.40}), (\ref{smr2.60}) -- (\ref{smr2.70}). 
Then
\begin{eqnarray}
\label{ass.10}
\Lambda_j &\in& {\Bbb R} \setminus \{ 0 \} \;\;\; 
\mbox{ for } 1 \leq j \leq \nu, \\
\label{ass.20}
V(\omega^{(0)})(e_j + e_l)
&\in& {\Bbb R} \setminus \{ 0 \} \;\;\;
\mbox{ for } 1 \leq j, l \leq \nu, \\
\label{ass.30}
V(\omega^{(0)})(e_j - e_l)
&\in& {\Bbb R} \setminus \{ 0 \} \;\;\;
\mbox{ for } 1 \leq j \neq l \leq \nu.
\end{eqnarray}
\end{Lemma}

\begin{proof}
Statement (\ref{ass.10}) is obvious from $\omega^{(0)}_j/2 \in
(0, \pi/2)$ for $1 \leq j \leq \nu$. Furthermore, for $1 \leq j, l \leq \nu$
we observe that $\frac{1}{2}<\omega^{(0)}, e_j + e_l> \in
(0, \pi)$ and $V(\omega^{(0)})(e_j + e_l) \in {\Bbb R}$. To see 
$V(\omega^{(0)})(e_j + e_l) \neq 0$ we use the inequality 
\begin{eqnarray}
\label{ass.40}
\sin (\alpha + \beta) < \sin \alpha + \sin \beta \;\;\;
\mbox{ for } \alpha, \beta \in (0, \pi/2).
\end{eqnarray}
Definition (\ref{smr2.30}) implies
\begin{eqnarray}
\label{ass.45}
j \gamma = 2 \sqrt{F'(-b)} \sin \frac{\omega^{(0)}_j}{2}, \;\;\;
l \gamma = 2 \sqrt{F'(-b)} \sin \frac{\omega^{(0)}_l}{2}.
\end{eqnarray}
Adding the two equations of (\ref{ass.45}) and using (\ref{ass.40})
yields
\begin{eqnarray}
\label{ass.50}
(l + j)\gamma > 
2 \sqrt{F'(-b)} \sin \frac{\omega^{(0)}_j + \omega^{(0)}_l}{2} > 0.
\end{eqnarray}
This implies $V(\omega^{(0)})(e_j + e_l) < 0$ by definition
(\ref{smr2.70}).

In order to prove (\ref{ass.30}), let $1 \leq j \neq l \leq \nu$.
Since $V(\omega)(m) = V(\omega)(-m)$ it suffices to consider the case
$j > l$. Again, $V(\omega^{(0)})(e_j - e_l) \in {\Bbb R}$ follows from
$\frac{1}{2}<\omega^{(0)}, e_j - e_l> \in
(0, \pi)$. Furthermore, (\ref{ass.40}) and (\ref{ass.45}) yield
\begin{eqnarray}
\label{ass.60}
0 < (j-l)\gamma &=& 
2 \sqrt{F'(-b)} \left(\sin \frac{\omega^{(0)}_j}{2} - 
\sin \frac{\omega^{(0)}_l}{2}\right)  
< 2 \sqrt{F'(-b)} \sin \frac{\omega^{(0)}_j - \omega^{(0)}_l}{2}.
\end{eqnarray}
Thus $V(\omega^{(0)})(e_j - e_l) > 0$.
\end{proof}

\begin{Lemma}
\label{Lass.2}
Let $\nu \in {\Bbb N}$, $\nu \geq 2$. Define
\begin{eqnarray}
\nonumber
\eta(x) := (\eta_1, \ldots, \eta_{\nu})(x), \;\;
\eta_j (x) := 2 \arcsin(j x)  
\;\;\; \mbox{ for }
|x| < 1/\nu.
\end{eqnarray}
For every $x_0 \in B_{1/\nu}(0)$ there exist positive constants $r$, 
$\tau$
and a (Lebesgue --) zero set $N \subset {\Bbb R}$ such that for every
$x \in B_{r}(x_0)  \setminus N$ there exists a real number $s > 0$ with
\begin{eqnarray}
\label{ass.70}
\mbox{ dist}(<m, \eta(x)>, 2 \pi {\Bbb Z}) > \frac{s}{|m|^{\tau}}
\;\;\; \mbox{ for all } m \in {\Bbb Z}^{\nu} \setminus \{ 0 \}.
\end{eqnarray}
\end{Lemma}
\begin{proof}
The proof proceeds in four steps:

\noindent
{\bf Step 1:}

Claim 1: 
{\em The functions $\eta_1', \ldots, \eta_{\nu}'$ are linearly
independent.}

\noindent
{\em Proof of claim 1:}
Denote by $\eta^{(k)}_j := (d/dx)^k \eta_j$ the $k$-th derivative of 
$\eta_j$. We define the $\nu \times \nu$ matrix
\begin{eqnarray}
\label{ass.100}
M := \left( \eta^{(2k-1)}_j (0) \right)_{1 \leq j, k \leq \nu}.
\end{eqnarray}
Note that $\eta^{(2k-1)}_j (0) = j^{2k-1} \eta^{(2k-1)}_1 (0)$ and 
therefore
\begin{eqnarray}
\label{ass.110}
\det M &=& 
\left(
\prod_{k=1}^{\nu}
\eta^{(2k-1)}_1 (0)
\right) \det \tilde{M}, \;\;\; \mbox{ where } \\
\nonumber
\tilde{M} := \left( j^{2k-1} \right)_{1 \leq j, k \leq \nu}.
\end{eqnarray}
Expanding
\begin{eqnarray}
\nonumber
\arcsin'(x) = (1 - x^2)^{-\frac{1}{2}} = 
\sum_{j \geq 0} \left( 
\begin{array}{c}
-\frac{1}{2}\\j \end{array} \right)
(-x^2)^j,
\end{eqnarray}
it is clear that $\eta^{(2k-1)}_1 (0) \neq 0$ for all $k \in {\Bbb N}$. 
In order to prove
$\det M \neq 0$ it suffices to verify that $\det \tilde{M} \neq 0$.
Suppose $\det \tilde{M} = 0$. Then the columns of 
$\tilde{M}$ are linearly dependent and consequently
there exists a nonzero vector
$(\alpha_1, \ldots, \alpha_k)$ with $\sum_{k=1}^{\nu} \alpha_k j^{2k -1} = 0$
for all $1 \leq j \leq \nu$. Hence the nonzero polynomial
$p(z) := \sum_{k=1}^{\nu} \alpha_k z^{2k -1}$ has $2\nu$ zeros
$\pm 1$, $\pm 2, \ldots, \pm \nu$ yielding a contradiction.
We have thus established that $\det M \neq 0$. This 
implies that the
functions $\eta_j'$ are linearly independent. Indeed, linear dependence 
of the vectors $\eta_j'$ would translate to linear 
dependence of the columns of $M$.

\noindent
{\bf Step 2:}

Claim 2: 
{\em
There exist positive constants $r$, $\delta$ and $K$,
such that for every $m \in {\Bbb Z}^{\nu} \setminus \{ 0 \}$ there
exists $k_m \in {\Bbb N}$, $k_m \leq K$ satisfying
\begin{eqnarray}
\label{ass.150}
\left|
< \frac{m}{|m|_2}, \eta^{(k_m)}(x) > \right| > k_m ! \delta^{k_m}
\;\; \mbox{ for all } x \in B_r(x_0)
\end{eqnarray}
and $|x_0| + r < 1/\nu$.
}

\noindent
{\em Proof of claim 2:}
Let $\xi \in S^{\nu -1}$, where $S^{\nu - 1}=\{ \xi \in {\Bbb R}^{\nu}:
|\xi|_2 = 1 \}$ denotes the Euclidean unit sphere in ${\Bbb R}^{\nu}$.
The function $f_{\xi}(x) := < \xi, \eta'(x) >$ is analytic and does not vanish
identically on $B_{1/\nu}(0)$ by claim 1. 
Hence there exists a $k_{\xi} \in {\Bbb N}$ such that
$< \xi, \eta^{(k_{\xi})}(x_0) > \neq 0$. By continuity we find positive
constants $r_{\xi}$, $\delta_{\xi}$ and a open neighborhood
$\xi \in V_{\xi} \subset S^{\nu - 1}$ such that
\begin{eqnarray}
\label{ass.160}
\left| < \zeta, \eta^{(k_{\xi})}(x) > \right| > k_{\xi}! \delta_{\xi}^
{k_{\xi}} \;\;\; \mbox{ for all } \zeta \in V_{\xi}, \; x \in 
B_{r_{\xi}}(x_0).
\end{eqnarray}
The sets $(V_{\xi})_{\xi \in S^{\nu - 1}}$ form an open cover of the
unit sphere. Claim 2 now follows from (\ref{ass.160}) and from the
compactness of $S^{\nu - 1}$.

\noindent
{\bf Step 3:}

Claim 3:  
{\em
Let $r$, $K$ be as determined in claim 2. 
For $0< t \leq 1$ we define the sets
\begin{eqnarray}
\label{ass.190}
{\cal B}_{t} :=
\left\{
x \in B_r(x_0): \exists m \in {\Bbb Z}^{\nu} \setminus \{ 0 \}
\; \mbox{ with } \;
\mbox{ dist}(<m, \eta(x)>, 2 \pi {\Bbb Z}) < t |m|_2^{-K(\nu+1)}  
\right\}.
\end{eqnarray}
There exists a constant $C > 0$ such that 
\begin{eqnarray}
\nonumber
\lambda( {\cal B}_{t} ) \leq C t^{1/K} 
\;\;\; \mbox{ for all } \; 0 < t \leq 1,
\end{eqnarray}
where $\lambda$ denotes the Lebesgue measure on ${\Bbb R}$.
}

\noindent
{\em Proof of claim 3:}
Denote by $C_1 := \sup_{x \in B_r(x_0)} |\eta'(x)|_2$ which is finite
since $|x_0| + r < 1/\nu$ (see claim 2).
Let $m \in {\Bbb Z}^{\nu} \setminus \{ 0 \}$ and set 
$g_m(x) := <m, \eta(x)>$ for $x \in B_r(x_0)$. From claim 2 it follows that
\begin{eqnarray}
\label{ass.200}
\left|
g_m^{(k_m)}(x)
\right| > k_m! \delta^{k_m} |m|_2.
\end{eqnarray}
According to lemma \ref{Lel.1} this implies for $l \in {\Bbb Z}$ and
\begin{eqnarray}
\nonumber
{\cal B}_{t, m, l}  :=
\left\{
x \in B_r(x_0): |g_m(x) - 2 \pi l| < t |m|_2^{-K(\nu+1)}  
\right\}
\end{eqnarray}
the following estimate
\begin{eqnarray}
\nonumber
\lambda
\left( {\cal B}_{t, m, l} 
\right) \leq 2 k_m \delta^{-1} |m|_2^{-\frac{1}{k_m}} 
(t |m|_2^{-K(\nu+1)} )^{\frac{1}{k_m}} 
\leq 2 K \frac{t^{\frac{1}{K}}}{\delta |m|_2^{\nu + 1 + \frac{1}{K}}}.
\end{eqnarray}
It follows from the definition of $C_1$ above that the 
range $\{ g_m(x): x\in B_r(x_0) \}$ is an interval of length 
$\leq 2 r C_1 |m|_2$. This implies that the number
of integers $l \in {\Bbb Z}$ for which ${\cal B}_{t, m, l}$
is not the empty set is bounded above by $1 + (2 r C_1 |m|_2 + 2t)/(2\pi)$.
The Lebesgue measure of the set
\begin{eqnarray}
\nonumber
{\cal B}_{t, m}  :=
\left\{
x \in B_r(x_0): \mbox{ dist}(g_m(x), 2 \pi {\Bbb Z}) < t |m|_2^{-K(\nu+1)}  
\right\}.
\end{eqnarray}
is therefore bounded by
\begin{eqnarray}
\label{ass.250}
\lambda
\left( {\cal B}_{t, m} 
\right) \leq (2 r C_1 + 3) |m|_2 \frac{2 K}{\delta} t^{\frac{1}{K}}
|m|_2^{-(\nu + 1 + \frac{1}{K})}.
\end{eqnarray}
Since ${\cal B}_t = \bigcup_{m \neq 0} {\cal B}_{t, m}$ claim 3
is an immediate consequence of (\ref{ass.250}).

\noindent
{\bf Step 4:}
 
\noindent
{\em Proof of the lemma.} Let $r$, $K$ be defined as in claim 2 and set
$\tau := K(\nu + 1)$. Furthermore we define 
$N := \bigcap_{n \geq 1} {\cal B}_{1/n}$ as intersection of countably many
nested sets. It follows from claim 3 that $\lambda (N) = 0$. On the 
other hand for $x \in B_r(x_0) \setminus N$ there exists an integer
$n \in {\Bbb N}$ such that $x \in B_r (x_0) \setminus {\cal B}_{1/n}$.
By definition (\ref{ass.190}) it is clear that (\ref{ass.70}) is satisfied
for $s = 1/(2n)$.
\end{proof}

In order to investigate the non-degeneracy condition, i.e. assumption A4,
we introduce some more notation. For $\nu \in {\Bbb N} \setminus
\{ 1 \}$ let $\Omega_1^{(\nu)}$ denote the $\nu \times \nu$ matrix
\begin{eqnarray}
\label{ass.400}
(\Omega_1^{(\nu)})_{i, j} := \delta_{i, j} + 2.
\end{eqnarray}
It follows from (\ref{smr2.20}) -- (\ref{smr2.40}), (\ref{smr2.70}),
(\ref{smr2.90}), and lemma \ref{Lass.1} 
that there exists a real analytic two parameter family
of symmetric $\nu \times \nu$ matrices $\Omega_2^{(\nu)} (t, \gamma)$
defined on 
\begin{eqnarray}
\label{ass.405}
{\cal D}_2^{(\nu)} :=
\{ (t, \gamma) \in {\Bbb R}^2: t > 0, \;\; \nu \gamma < 2 \sqrt{t}
< (\nu + 1) \gamma \},
\end{eqnarray}
such that $\Omega^{(\nu)}$ (see (\ref{smr2.90})) can be written in the form
\begin{eqnarray}
\label{ass.410}
\Omega^{(\nu)} = 
\left(
F''(-b)^2 \Omega_2^{(\nu)} (F'(-b), \gamma) - F'''(-b) \Omega_1^{(\nu)}
\right) \cdot \mbox{ diag}(\Lambda_l^{-1}).
\end{eqnarray}
Set
\begin{eqnarray}
\label{ass.420}
{\cal D}_3^{(\nu)} := {\Bbb R} \times {\cal D}_2^{(\nu)},
\end{eqnarray}
For $1 \leq j, l \leq \nu$ define real analytic functions 
\begin{eqnarray}
\label{ass.423}
G_{j, l}^{(\nu)} :  {\cal D}_3^{(\nu)} \to {\Bbb R}; \;\;\;
G_{j, l}^{(\nu)} (u, t, \gamma) := \frac{u^2}{2 + \delta_{j, l}} 
\left(
\Omega_2^{(\nu)} (t, \gamma)
\right)_{j, l}.
\end{eqnarray}
Observe that 
\begin{eqnarray}
\label{ass.427}
\Omega_{j, l}^{(\nu)} = 0  \Leftrightarrow
F'''(-b) = G_{j, l}^{(\nu)} (F''(-b), F'(-b), \gamma).
\end{eqnarray}
Next we turn to the condition $\det \Omega^{(\nu)} \neq 0$ in assumption
A4.
 
It is easy to see that 
$\det \Omega_1^{(\nu)} = 2 \nu + 1$ (use for example that
$\Omega_1^{(\nu)}$ is of the form identity + rank 1 matrix). Hence 
$\Omega_1^{(\nu)}$ is invertible. We define
for $(u, t, \gamma) \in {\cal D}_3^{(\nu)}$
\begin{eqnarray}
\label{ass.430}
\Gamma^{(\nu)} (u, t, \gamma) := u^2 
\Omega_2^{(\nu)} (t, \gamma) \left( \Omega_1^{(\nu)} \right)^{-1},
\end{eqnarray}
which defines a real analytic family of symmetric $\nu \times \nu$
matrices. Hence there exist real analytic functions
$\lambda_j^{(\nu)} : {\cal D}_3^{(\nu)} \to {\Bbb R}$, $1 \leq j \leq \nu$,
not necessarily pairwise distinct, which denote the eigenvalues of 
$\Gamma^{(\nu)}$ (see \cite{Kat}). Observe that (\ref{ass.410}) and
(\ref{ass.430})
imply
\begin{eqnarray}
\label{ass.440}
\det \Omega^{(\nu)} = 0
\Leftrightarrow
\det \left( \Gamma^{(\nu)} (F''(-b), F'(-b), \gamma) - F'''(-b) \right) = 0
\end{eqnarray}
and hence
\begin{eqnarray}
\label{ass.450}
\det \Omega^{(\nu)} = 0
\Leftrightarrow
F'''(-b) \in \{\lambda_j^{(\nu)}(F''(-b), F'(-b), \gamma): 1 \leq j \leq \nu \}.
\end{eqnarray}
Since it is our goal to avoid that either an entry of $\Omega^{(\nu)}$ 
vanishes or that
$\det \Omega^{(\nu)}$ vanishes we exclude 
those analytic force laws $F$ which, for some $1 \leq j, l \leq \nu$,
satisfy one of the differential equations
\begin{eqnarray}
\label{ass.477}
F'''(-b) &=& G_{j, l}^{(\nu)}(F''(-b), F'(-b), \gamma), \\
\label{ass.480}
F'''(-b) &=& \lambda_j^{(\nu)}(F''(-b), F'(-b), \gamma),
\end{eqnarray}
in a neighborhood of $-b_0$. Furthermore, for reasons which will become clear
in the proof of theorem \ref{Tass.1} below we also exclude the 
case $F'$ is constant. Note that in this case $F$ is linear and the 
quasi-periodic wave solutions of (\ref{smr.10}) are easily constructed
(see section \ref{smr1}). We define for $b_0 \in {\Bbb R}$, $r_0 > 0$ 
\begin{eqnarray}
\label{ass.470}
{\cal A}_{b_0, r_0} :=
\left\{
F: B_{r_0}(- b_0) \to {\Bbb R} \mbox{ real analytic }: F'(-b_0) > 0 \right\}.
\end{eqnarray}
For $\nu \in {\Bbb N} \setminus \{ 1 \}$ we say that $F$ is 
$\nu$ -- degenerate, if there 
exists $\gamma > 0$ with
\begin{eqnarray}
\nonumber
\nu \gamma < 2 \sqrt{F'(-b_0)} < (\nu + 1) \gamma
\end{eqnarray}
and $F$ satisfies one of the differential
equations (\ref{ass.477}), (\ref{ass.480}) in some neighborhood of $-b_0$.
Set
\begin{eqnarray}
\label{ass.490}
{\cal F}_{b_0, r_0}^{(\nu)} := 
{\cal A}_{b_0, r_0} \setminus
\left(
\{ F \in {\cal A}_{b_0, r_0}: F \mbox{ linear } \} 
\cup \{ F \in {\cal A}_{b_0, r_0}: F \mbox{ is 
$\nu$ -- degenerate} \}
\right)
\end{eqnarray}

\begin{Lemma}
\label{Lass.3}
Let $\nu \in {\Bbb N} \setminus \{ 1 \}$, $b_0 \in {\Bbb R}$, and 
$r_0 > 0$.
\begin{itemize}
\item[(a)]
The set ${\cal F}_{b_0, r_0}^{(\nu)}$ is an open and dense subset of 
${\cal A}_{b_0, r_0}$
with respect to the topology generated by the sup-norm. 
\item[(b)]
The exceptional set
${\cal A}_{b_0, r_0} \setminus {\cal F}_{b_0, r_0}^{(\nu)}$ 
consists of all linear functions
in ${\cal A}_{b_0, r_0}$ together with a finite collection of four-parameter
families of functions, each of which can be parameterized by
$F(-b_0)$, $F'(-b_0)$, $F''(-b_0)$ and $\gamma$.
\end{itemize}
\end{Lemma}
\begin{proof}
The proof of (b) follows from definition (\ref{ass.490}). Claim (a) 
is a consequence of (b).
\end{proof}

We now state the main result of this section. It makes the notion precise
that {\em assumptions A1 -- A4 are generically satisfied}.

\begin{Theorem}
\label{Tass.1}
Let $\nu \in {\Bbb N} \setminus \{ 1 \}$, $b_0 \in {\Bbb R}$,
$r_0 > 0$, $F \in {\cal F}_{b_0, r_0}^{(\nu)}$, and $\gamma > 0$ such
that 
\begin{eqnarray}
\label{ass.300}
\nu \gamma < 2 \sqrt{F'(-b_0)} < (\nu + 1) \gamma.
\end{eqnarray}
Then
there exist $\beta > 0$ and a (Lebesgue --) zero set $M \subset {\Bbb R}$ 
such that assumptions A1 -- A4 of section \ref{smr2}
are satisfied for all $b \in B_{\beta}(b_0) \setminus M$.
The set ${\cal F}_{b_0, r_0}^{(\nu)}$ is an open 
and dense subset of ${\cal A}_{b_0, r_0}$
with respect to the topology generated by the sup-norm.
\end{Theorem}

\begin{proof}
Assumptions A1 and A2 are satisfied for
$F$, $b_0$, $\gamma$ by the hypothesis
of theorem \ref{Tass.1}. Clearly there exists a 
$\beta_1 > 0$ such that 
\begin{eqnarray}
\label{ass.305}
\nu \gamma < 2 \sqrt{F'(-b)} < (\nu + 1) \gamma
\end{eqnarray}
for all $b \in B_{\beta_1}(b_0)$ and therefore 
assumptions A1 and A2 are satisfied for
$F$, $b$, $\gamma$ if $b \in B_{\beta_1}(b_0)$. 

To verify the diophantine condition (assumption A3) set
\begin{eqnarray}
\label{ass.320} 
G(b) := \frac{\gamma}{2 \sqrt{F'(- b)}} \;\;\; \mbox{ for } 
b \in B_{\beta_1}(b_0)
\end{eqnarray}
and define $x_0 := G(b_0) \in B_{1/\nu}(0)$ (by (\ref{ass.300})).
We apply lemma \ref{Lass.2} and obtain positive constants $r$, $\tau$ and
a zero-set $N$ such that for every $x \in B_r(x_0) \setminus N$ there
exists a constant $s > 0$ such that (\ref{ass.70}) holds.
Choose $0 < \beta \leq \beta_1$ such that 
$G(B_{\beta}(b_0)) \subset B_r(x_0)$ and define
\begin{eqnarray}
\label{ass.330}
M_1 := \{ b \in B_{\beta}(b_0) : G(b) \in N \}.
\end{eqnarray}
By (\ref{ass.305}) and the definition of ${\cal F}_{b_0, r_0}^{(\nu)}$ 
(in (\ref{ass.490}))
we observe that $G$ is an analytic, non constant function 
on $B_{\beta}(b_0)$. Thus
$M_1$ is again a set of zero Lebesgue measure.
For all $b \in B_{\beta}(b_0) \setminus M_1$ we have 
$G(b) \in B_r(x_0) \setminus N$ and hence there
exists a constant $s > 0$ such that (\ref{ass.70}) holds.
Assumption A3 is satisfied for all 
$b \in B_{\beta}(b_0) \setminus M_1$.

In order to investigate assumption A4 recall that the choice of $\beta$
implies that (\ref{ass.305}) is true for all $b \in B_{\beta}(b_0)$.
This implies that $(F'(-b), \gamma) \in {\cal D}_2^{(\nu)}$ for all
$b \in B_{\beta}(b_0)$.
Since $F$ is not $\nu$ -- degenerate and $F$,
$\lambda_j^{(\nu)}$, $G_{j, l}^{(\nu)}$
$1 \leq j, l \leq \nu$ are real analytic functions the set
\begin{eqnarray}
\nonumber
M_2 :=
\left\{
b \in B_{\beta}(b_0):
F'''(-b) \right. &\in& 
\{\lambda_j^{(\nu)}(F''(-b), F'(-b), \gamma): 1 \leq j \leq \nu \} 
\\
\nonumber
&\cup& \left.
\{G_{j, l}^{(\nu)}(F''(-b), F'(-b), \gamma): 1 \leq j, l \leq \nu \} \;
\right\}
\end{eqnarray}
consists of at most countably many points. By (\ref{ass.450}) 
and (\ref{ass.427}) it is clear that no entry of $\Omega^{(\nu)}$ vanishes
and that
$\det \Omega^{(\nu)} \neq 0$ for $b \in B_{\beta}(b_0) \setminus M_2$.

Set $M := M_1 \cup M_2$. Assumptions A1 -- A4 hold for all 
$b \in B_{\beta}(b_0) \setminus M$. 
The genericity of the set 
${\cal F}_{b_0}$ has already been proved in lemma \ref{Lass.3}.
\end{proof}

\section{Weight functions}
\label{wf}

In this section we prove proposition \ref{Pwf.1} which implies the useful
property
(\ref{wf.120}) for the families of weight functions introduced in
\ref{DFa1}. The proof of proposition \ref{Pwf.1} 
can be found in \cite{Kri} and we repeat it
here for the sake of completeness. 
%Similar classes of weights were already used by
%DeLatte \cite{DeL}, Craig -- Wayne \cite{CW2},  Bourgain \cite{B4}
%and Kriecherbauer \cite{Kri} in the context of small divisor problems.

Let $\nu \in {\Bbb N}$ and $0 < x \leq 1$. 
Recall that $|\cdot|$ denotes the maximum norm on ${\Bbb Z}^{\nu}$.
We define
\begin{eqnarray}
\label{wf.10}
D_{\nu, x} :=
\left\{
\begin{array}{ll}
\sup_{n \in {\Bbb Z}^{\nu}, \sigma \geq 1/4}
\sum_{m \in {\Bbb Z}^{\nu}} 
\frac{e^{\sigma |n|^x}}{e^{\sigma |m|^x} e^{\sigma |n-m|^x}} &
\mbox{ if } \; 0 < x < 1,\\
\sup_{n \in {\Bbb Z}^{\nu}, \sigma \geq 1/4}
\sum_{m \in {\Bbb Z}^{\nu}} 
\frac{(1 +|n|)^{\nu+1} e^{\sigma |n|}}
{(1 + |m|)^{\nu+1}e^{\sigma |m|}
(1 + |n-m|)^{\nu+1}e^{\sigma |n-m|}} &
\mbox{ if } \; x=1.
\end{array}
\right.
\end{eqnarray}

\begin{Proposition}
\label{Pwf.1}
Let $\nu \in {\Bbb N}$ and $0 < x \leq 1$. Then $1 \leq D_{\nu, x} < \infty$.
\end{Proposition}
\begin{proof}
The claim $D_{\nu, x} \geq 1$ is obviously satisfied.

{\em Case $0 < x < 1$:}

\noindent
For every $n \in {\mathbb Z}^{\nu}$, 
we divide the sum in the definition of $D_{\nu, x}$
into four parts, depending on $n$, and show that each of the partial
sums is bounded independently of $n$. Set
\begin{eqnarray}
\nonumber
A_n &:=& \left\{ m \in {\mathbb Z}^{\nu}: \; |m| < \frac{2}{3} |n| \right\},
\\
\nonumber
B_n &:=&
\left\{ m \in {\mathbb Z}^{\nu} \setminus A_n:  
|m-n| < \frac{2}{3} |n| \right\}, \\
\nonumber
C_n &:=& \left\{ m \in {\mathbb Z}^{\nu} \setminus (A_n \cup B_n):  
|m| < \frac{5}{3} |n| \right\}, \\
\nonumber
D_n &:=&
{\mathbb Z}^{\nu} \setminus (A_n \cup B_n \cup C_n).
\end{eqnarray}
To estimate the sum over $A_n$, observe that
\begin{eqnarray}
\label{wf.20}
\sum_{m \in A_n} e^{-\sigma(|m|^x + |n-m|^x - |n|^x)} \leq
\sum_{k=0}^{\lfloor 2|n|/3 \rfloor} 2 \nu (2 k + 1)^{\nu-1} 
e^{-\sigma(k^x + (|n|-k)^x - |n|^x)}
\end{eqnarray}
and 
\begin{eqnarray}
\label{wf.30}
a^x + b^x - (a+b)^x = x \int_0^a s^{x-1} 
\left(1-\left(\frac{s}{s+b}\right)^{1-x} \right)ds
\geq a^x \left(1-\left(\frac{a}{a+b}\right)^{1-x} \right) 
\nonumber
\end{eqnarray}
for  $a \geq 0$, $b > 0$.
Furthermore, there exists a constant $C_1 (\nu, x)$ such that
\begin{eqnarray}
\label{wf.40}
2 \nu (2 k + 1)^{\nu-1} e^{- \sigma k^x 
\left(1-\left(\frac{2}{3}\right)^{1-x} \right)} \leq
\frac{C_1(\nu, x)}{(k+1)^2} \quad \mbox{ for all } \sigma \geq \frac{1}{4},
\; k \geq 0,
\end{eqnarray}
and we can bound (\ref{wf.20}) by
$2 C_1(\nu, x)$.
By translation, the same estimate holds for the sum over $B_n$.
Clearly, there exists another constant $C_2(\nu, x)$, such
that
\begin{eqnarray}
\label{wf.50}
\sum_{m \in C_n} e^{-\sigma(|m|^x + |n-m|^x - |n|^x)} \leq
(\frac{10}{3} |n| +1)^{\nu} e^{-\frac{\sigma}{3} |n|^x} \leq
C_2(\nu, x) \;\;\; \mbox{ for all } \sigma \geq \frac{1}{4}, \;\;
n \in {\mathbb Z}^{\nu}.
\end{eqnarray}
Finally, for $m \in D_n$, one easily derives 
$|m|^x + |n-m|^x - |n|^x \geq |m|^x \left(1 - (3/5)^x
\right)$, which in turn implies a uniform bound on the sum over $D_n$.

{\em Case $x=1$:} 

\noindent
Let $n \in {\Bbb Z}^{\nu}$. 
We use the same decomposition of ${\Bbb Z}^{\nu}$ into sets
$A_n$, $B_n$, $C_n$ and $D_n$ as above.
%\begin{eqnarray}
%\label{wf.60}
%\left(\frac{1+|n|}{(1+|m|)(1+|n-m|)}\right)^{\nu+1} 
%e^{-\sigma(|m| + |n-m| - |n|)}
%\end{eqnarray} 
In order to estimate the sum over $A_n$ and $B_n$,
we use
\begin{eqnarray}
\label{wf.70}
\sum_{|m| < \frac{2|n|}{3}} \left(\frac{1+|n|}{(1+|m|)(1+|n-m|)}\right)
^{\nu +1} 
< \sum_{m \in {\mathbb Z}^{\nu}}
\left(\frac{3}{1+|m|}\right)^{\nu + 1} < \infty.
\end{eqnarray}
For $C_n$ we use
\begin{eqnarray}
\label{wf.80}
\frac{1+|n|}{(1+|m|)(1+|n-m|)} \leq 1 \;\;\;
\mbox{ for all } m,n \in {\mathbb Z}^{\nu},
\end{eqnarray}
and
\begin{eqnarray}
\label{wf.90}
\sum_{m \in C_n} e^{-\sigma(|m| + |n-m| - |n|)} \leq
(\frac{10}{3} |n| +1)^{\nu} e^{-\frac{\sigma}{3} |n|}, 
\end{eqnarray}
which has an $n$- independent bound.
Finally, to estimate the sum over $D_n$, we use (\ref{wf.80}) and
\begin{eqnarray}
\label{wf.100}
\sum_{m \in D_n} e^{-\sigma(|m| + |n-m| - |n|)} &\leq&
\sum_{m \in {\Bbb Z}^{\nu}} 
e^{-\frac{2 \sigma}{5} |m|}.
\end{eqnarray}
\end{proof}

\section{Properties of the nonlinear part}
\label{enp}

In this section we provide estimates for the nonlinear part of the 
equation. Let $\nu \in {\Bbb N}$ and recall that  $c = 0.01$ (see 
(\ref{npc.52})). For $\sigma \geq 1/4$ 
we consider sequence spaces $X_{\sigma, c}$, $X_{\sigma, 1}$ 
on ${\Bbb Z}^{\nu}$
and corresponding matrix spaces
${\cal L}_{\sigma, c}$, ${\cal L}_{\sigma, 1}$. These spaces were defined in
\ref{DFa2} (see also \ref{DFa1}, (\ref{fa1.20}), (\ref{fa1.60}) and
(\ref{npc.54})).

\begin{Lemma}
\label{Lenp.1}
Suppose that
\begin{eqnarray}
\label{enp.10}
\sum_{k=0}^{\infty} \alpha_k y^k
\end{eqnarray}
defines a real power series with positive radius of convergence
$0 < r \leq \infty$. There exists a positive constant $D_W > 0$ such that
the following holds. 
For $x \in \{1, c\}$ and $\sigma \geq 1/4$ the map defined by
\begin{eqnarray}
\label{enp.20}
W: \{u \in X_{\sigma, x}: \|u\|_{\sigma, x} < r \} &\to& X_{\sigma, x},
\\
\nonumber  
W(u) &:=& \sum_{k=2}^{\infty} \alpha_k u^{*k}
\end{eqnarray}
is analytic and statements (a) -- (c) are true.
\begin{itemize}
\item[(a)]
For $j \in {\Bbb N}$, $u, y_1, \ldots, y_j \in X_{\sigma, x}$ with
$\| u \|_{\sigma, x} < r$,
\begin{eqnarray}
\label{enp.30}
(D^j W)(u)[y_1, \ldots, y_j] =
\sum_{k \geq \max(2,j)} 
\frac{k!}{(k-j)!} \alpha_k u^{*(k-j)}*y_1 * \ldots *y_j.
\end{eqnarray}

\item[(b)]
For $j \in {\Bbb N}$, $2 \leq j \leq 22 \nu^3 + 1$ (cf. 
(\ref{npc.80})), $y_1$, $\ldots$, $y_{j-1}$, 
$u \in X_{\sigma, x}$ with
$\| u \|_{\sigma, x} \leq  \min(1, r)/2$
\begin{eqnarray}
\nonumber
\| W(u) \|_{\sigma, x} &\leq& D_W \|u \|_{\sigma, x}^2 \\
\label{enp.40} 
\| (DW)(u) \|_{\sigma, x} &\leq& D_W \|u \|_{\sigma, x} \\
\nonumber
\| D^j W(u)[y_1, \ldots, y_{j-1}, \cdot] \|_{\sigma, x} &\leq&
D_W \| y_1 \|_{\sigma, x} \cdot \; \ldots \; \cdot \|y_{j-1}\|_{\sigma, x}
%\left\|
%\sum_{k \geq \max(2,j)} 
%\frac{k!}{(k-j)!} \alpha_k u^{*(k-j)}  
%\right\| &\leq& j! \left( 2/\rho \right)^j C \;\;
%\mbox{ for } j \geq 2.
\end{eqnarray}

\item[(c)] 
For $L > 0$, $u \in X_{1, c}$ with $\| u \|_{1, c} \leq 
\min(1, r)/2$ and
supp$(u) \subset B_L(0)$ (i.e. $u(m) = 0$ for $|m| \geq L$)
\begin{eqnarray} 
\label{enp.50}
\left\|
W(u) \chi_{ \{ |m| \geq 3L \} } 
\right\|_{1/4, c}
\leq D_W e^{-\frac{9}{4} L^c} \|u\|^2_{1, c}.
\end{eqnarray}
\end{itemize}
\end{Lemma}

\begin{proof}
The analyticity of $W$ and statement (a) are obvious. In view of
proposition \ref{Pfa1.1} the estimates of (b) follow from standard
arguments for power series. In order to prove (c) we first observe that
for $k \in {\Bbb N}$ and for numbers $L \geq a_1 \geq \ldots \geq a_k \geq 0$
with $a_1 + \ldots + a_k \geq 3L$ the following holds
\begin{eqnarray}
\nonumber
a_1^c + \ldots + a_k^c - 3L^c &=&
c
\left(
\int_{0}^{a_1} t^{c-1}dt + \ldots + \int_{0}^{a_k} t^{c-1}dt - 3
\int_{0}^{L} t^{c-1}dt
\right) \\
\nonumber
&=& c
\left[
\left(
\sum_{j=4}^{k} \int_{0}^{a_j} t^{c-1}dt 
\right)
-
\left(
\int_{a_1}^{L} t^{c-1}dt + \int_{a_2}^{L} t^{c-1}dt + 
\int_{a_3}^{L} t^{c-1}dt
\right)
\right] \\
\label{enp.60}
&=& c [ I - II ] \geq 0.
\end{eqnarray}
Indeed, the last inequality, or equivalently $ I \geq II$, is easy to see:
firstly, the total length of the domain of integration in $I$ equals to
$a_4 + \ldots + a_k \geq 3 L - a_1 -a_2 -a_3$ and is hence bigger or equal
to the total length of the integration domain in $II$. Secondly,
by the ordering of the $a_j$'s and the monotonicity of $t^{c-1}$ it is
clear that the integrands in $I$ are pointwise bigger than 
$a_4^{c-1}$ whereas the integrands in $II$ are pointwise bounded above
by $a_3^{c-1} \leq a_4^{c-1}$. Hence (\ref{enp.60}) is established.
Using in addition lemma \ref{Lwf.1} (cf. (\ref{fa1.10})) we conclude
\begin{eqnarray}
&&\nonumber
\left\|
W(u) \chi_{ \{ |m| \geq 3L \} } 
\right\|_{1/4, c} \\
&\leq&
\sum_{k \geq 2} |\alpha_k|
\sum_{\scriptsize
\begin{array}{c}
s_1 + \ldots + s_k = m \\
|m| \geq 3 L
\end{array}
}
|u(s_1)| \cdot \ldots \cdot |u(s_k)| w_{1/4, c} (m)\\
\nonumber
&\leq&
\sum_{k \geq 2} |\alpha_k|
\sum_{\scriptsize
\begin{array}{c}
s_1 + \ldots + s_k = m \\
|m| \geq 3 L
\end{array}
}
|u(s_1)| w_{1/4, c}(s_1) \cdot \ldots \cdot |u(s_k)| w_{1/4, c} (s_k)\\
\nonumber
&\leq&
\sum_{k \geq 2} |\alpha_k|
\sum_{\scriptsize
{\scriptstyle
\begin{array}{c}
s_1 + \ldots + s_k = m \\
|m| \geq 3 L
\end{array}
}
}
|u(s_1)| w_{1, c}(s_1) \cdot \ldots \cdot |u(s_k)| w_{1, c} (s_k)
e^{-\frac{3}{4} (|s_1|^c + \ldots + |s_k|^c)}\\
\nonumber
&\leq&
e^{-\frac{3}{4} 3 L^c}
\sum_{k \geq 2} |\alpha_k|
\sum_{\scriptsize
{\scriptstyle
\begin{array}{c}
s_1 + \ldots + s_k = m \\
m \in {\Bbb Z}^{\nu}
\end{array}
} 
}
|u(s_1)| w_{1, c}(s_1) \cdot \ldots \cdot |u(s_k)| w_{1, c} (s_k)
\\
\nonumber
&\leq&
e^{-\frac{3}{4} 3 L^c}
\sum_{k \geq 2} |\alpha_k|
\| u \|_{1, c}^k.
\end{eqnarray}
This establishes claim (c).
\end{proof}

\section{Properties of the linearized operators}
\label{lop}

In this section we prove various properties of the matrices
$T^{(j)}(\lambda)$, $T^{(j)}(\theta, \lambda)$ which are used in the  
analysis of chapter II. We provide estimates on the diagonal entries
(proposition \ref{Plop.1}) as well as certain symmetry properties
(propositions \ref{Plop.2} and \ref{Plop.3}).

\begin{Proposition}
\label{Plop.1}
Recall the definitions of $\nu$ (\ref{smr2.20}), $\omega_j^{(0)}$ 
(\ref{smr2.30}), and of the function 
$V_j$, $j \in {\Bbb N}$ (\ref{lop.50}).
There exist constants $0 < d_V, \delta_V \leq 1 \leq D_V$ 
such that the following holds.
\begin{itemize}
\item[(a)]  For all $1 \leq j \leq \nu$, 
$0 < \delta < \delta_V$ and $\vartheta \in {\Bbb C}$
\begin{eqnarray}
\nonumber
\mbox{ dist}(\vartheta, \{\pm w_j^{(0)} \} + 2 \pi {\Bbb Z}) &>& \delta
\; \Rightarrow \; |V_j(\vartheta)| > d_V \delta, \\
\nonumber
\mbox{ dist}(\vartheta, \{\pm w_j^{(0)} \} + 2 \pi {\Bbb Z}) &\leq& \delta
\; \Rightarrow \; |V_j(\vartheta)| \leq D_V \delta, \\ 
\nonumber
|(d/d\vartheta)^k V_j (\vartheta)| &\leq& D_V
\;\;\; \mbox{ for } 0 \leq k \leq 1+ 8 \nu^2, \;
|\vartheta  \mp \omega_j^{(0)}| < \delta_V.
\end{eqnarray}  

\item[(b)] For $1 \leq j \leq \nu$ there exist smooth functions
$\tilde{V}_{j, +}: B_{\delta_V}(\omega_j^{(0)}) \to {\Bbb R}$ and
$\tilde{V}_{j, -}: B_{\delta_V}(-\omega_j^{(0)}) \to {\Bbb R}$, such that
\begin{eqnarray}
\label{lop.100}
V_j(\vartheta) &=& \tilde{V}_{j, +} (\vartheta) (\vartheta - \omega_j^{(0)}) 
\;\;\; \mbox{ for } \vartheta \in B_{\delta_V}(\omega_j^{(0)}), \\
\label{lop.110}
V_j(\vartheta) &=& \tilde{V}_{j, -} (\vartheta) (\vartheta + \omega_j^{(0)}) 
\;\;\; \mbox{ for } \vartheta \in B_{\delta_V}(-\omega_j^{(0)}).
\end{eqnarray}
Furthermore 
\begin{eqnarray}
\label{lop.120}
|\tilde{V}_{j, \pm}(\vartheta) | &\geq& d_V 
\;\;\; \mbox{ for } \vartheta \in B_{\delta_V}(\pm \omega_j^{(0)}), \\
\label{lop.130}
|(d/d\vartheta)^k (1/\tilde{V}_{j, \pm})(\vartheta)| &\leq& D_V
\;\;\; \mbox{ for } 0 \leq k \leq 22 \nu^3, \;
\vartheta \in B_{\delta_V}(\pm \omega_j^{(0)}).
\end{eqnarray}

\end{itemize}
\end{Proposition}

\begin{proof}
Recall that $\omega_j^{(0)} \in (0, \pi)$ by assumption A2 
(cf. (\ref{smr2.30})).
Therefore $V_j'(\pm \omega_j^{(0)}) \neq 0$ and there exist $\delta_V > 0$,
$d_{V, 1} > 0$ such that $\omega_1^{(0)} - \delta_V > 0$ and
\begin{eqnarray}
\label{lop.150}
|V_j (\vartheta)| \geq d_{V, 1} |\vartheta - \omega_j^{(0)}|
\;\;\; \mbox{ for all }\; \vartheta \in U_{\delta_V}(\omega_j^{(0)}), \\
\label{lop.160}
|V_j (\vartheta)| \geq d_{V, 1} |\vartheta + \omega_j^{(0)}|
\;\;\; \mbox{ for all }\; \vartheta \in U_{\delta_V}(-\omega_j^{(0)}).
\end{eqnarray}
Observe further that the zeros of $V_j$ form a set
$\{\pm w_j^{(0)} \} + 2 \pi {\Bbb Z}$. Since $V_j$ is a 
$2 \pi$ -- periodic function and $\lim_{|Im z| \to \infty} V_j(z) = \alpha_1
> 0$ we conclude 
\begin{eqnarray}
\label{lop.170}
m := \inf \{ |V_j(\vartheta)| \; : \; \vartheta \in {\Bbb C},
\mbox{ dist}(\vartheta, \{\pm w_j^{(0)} \} + 2 \pi {\Bbb Z}) \geq \delta_V
\} > 0
\end{eqnarray}
Setting
\begin{eqnarray}
\label{lop.180}
d_V := \frac{1}{2} \min(d_{V, 1}, m/\delta_V, 2)
\end{eqnarray}
we have satisfied the first estimate in (a). 
The second and third estimates in (a) are true if
\begin{eqnarray}
\label{lop.183}
D_V \geq \max_{1 \leq j \leq \nu, 0 \leq k \leq 8 \nu^2 + 1} \;\;
\sup_{\vartheta \in 
U_{\delta_V}(\pm \omega_j^{(0)})} |(d/d\vartheta)^k V_j(\vartheta)|.
\end{eqnarray}
The right hand side of (\ref{lop.183}) is finite since the distance
of $U_{\delta_V}(\pm \omega_j^{(0)})$ to the set $2 \pi {\Bbb Z}$
is positive by the choice of $\delta_V$
($\omega_1^{(0)} - \delta_V > 0)$. 

To prove (b) observe first 
that (\ref{lop.100}) is satisfied
with
\begin{eqnarray}
\label{lop.190}
\tilde{V}_{j, +}(\vartheta) = \int_{0}^{1}
V_j'(\omega_j^{(0)} + s(\vartheta-\omega_j^{(0)})) ds.
\end{eqnarray}
The estimate (\ref{lop.120}) follows from (\ref{lop.100}), 
(\ref{lop.150}) and (\ref{lop.180}). The existence of the constant
$D_V$ follows from the estimate (\ref{lop.120}), since $\delta_V$ is chosen
such that $\omega_1^{(0)} - \delta_V > 0$
(and consequently dist$(B_{\delta_V}(\omega_j^{(0)}), 2 \pi {\Bbb Z}) > 0$). 
The proof of the corresponding
statements for $V_{j, -}$ is similar and will not be repeated here. 

\end{proof}

Recall the definition of the nonlinear map $W$ in section \ref{enp} and
of the functions $V(\omega)$, $V(\theta, \omega)$ in 
(\ref{smr2.70}), (\ref{ova.520}). For 
$u \in X_{\sigma, x}$, $\sigma \geq 1/4$, $x \in \{0.01, 1\}$ and
$\|u\|_{\sigma, x} < r_{F, b}$ (cf. (\ref{npc.60})) we may define
\begin{eqnarray}
\label{lop.2}
T(\omega)(m, n) &=& V(\omega)(m) \delta_{m, n} +
DW(u) (m, n), \\
\label{lop.3}
T(\theta, \omega)(m, n) &=& V(\theta, \omega)(m) \delta_{m, n} +
DW(u) (m, n).
\end{eqnarray}
By lemma \ref{Lenp.1} (a) we can express
\begin{eqnarray}
\label{lop.5}
DW(u)(m, n) = \left( \sum_{k \geq 2} k \alpha_k u^{*(k-1)} \right)(m-n).
\end{eqnarray}

Using (\ref{lop.5}) and the definition 
of the convolution (see (\ref{fa1.40})), 
the following two propositions are easy to verify.

\begin{Proposition}
\label{Plop.2}
Let $T(\omega)$ be defined as in (\ref{lop.2}).
\begin{itemize}
\item[(a)]
Suppose that $\omega \in {\Bbb R}^{\nu}$ and $u$ satisfies
$u(m) = u(-m) \in {\Bbb R}$ for all 
$m \in {\Bbb Z}^{\nu}$. Then
\begin{eqnarray}
\label{lop.200}
T(\omega)(m, n) = T(\omega)(-m, -n) 
\in {\Bbb R} \cup \{ \infty \}
\;\;\; \mbox{ for all } \; m, n \in {\Bbb Z}^{\nu}.
\end{eqnarray} 
\item[(b)]
Let $1 \leq j \leq \nu$ and $\omega \in {\Bbb C}^{\nu}$. Suppose that
$u(m) = 0$ for all $m \in {\Bbb Z}^{\nu}$ with $m_j \neq 0$. Then
\begin{eqnarray}
\label{lop.210}
T(\omega)(m, n) = 0 \;\;\; \mbox{ if } m_j \neq n_j.
\end{eqnarray}
\end{itemize}
\end{Proposition}

\begin{Proposition}
\label{Plop.3}
Let $\theta \in {\Bbb C}$, 
$\omega \in {\Bbb C}^{\nu}$ and let $T(\theta, \omega)$ be defined
as in (\ref{lop.3}). For $k \in {\Bbb Z}^{\nu}$ with
$<k, g> = 0$ the following translation property holds.
\begin{eqnarray}
\label{lop.300}
T(\theta, \omega)(m + k, n + k) = 
T(\theta + <\omega, k>, \omega)(m, n).
\end{eqnarray}
\end{Proposition}

\section{The coupling lemma}
\label{cl}

In this section we state and prove a version of the coupling lemma which
first appeared in \cite{Kri}. The main advantage of this coupling
lemma over other versions is that one may couple local inverse
matrices of all length scales in one step. The proof of the coupling
lemma \ref{Lcl.1} proceeds by constructing a parametrix for the 
inverse matrix. The columns of the parametrix are given
by the corresponding columns of the local inverse matrices.

Recall the definitions of $\nu$ (\ref{smr2.20}), $c = 0.01$ (\ref{npc.52}), of
the weight functions $w_{\sigma, c}$ on ${\Bbb Z}^{\nu}$ 
(definition \ref{DFa1}), and of the corresponding matrix spaces
${\cal L}_{\sigma, c}$ (definition \ref{DFa2}).

\begin{Lemma}
\label{Lcl.1}
Let $\sigma \geq \tilde{\sigma} \geq 1/4$ and let
$\Lambda$ be a finite subset of ${\Bbb Z}^{\nu}$.
Suppose that 
\begin{eqnarray}
\label{cl.10}
T(m, n) = D(m, n) + R(m, n), \;\;\; m, n \in \Lambda,
\end{eqnarray}
where $D$ is a diagonal matrix. Suppose further, that for each
$n \in \Lambda$ there exist numbers $l_n > 0$, 
$0 \leq \mu_n \leq \sigma - \tilde{\sigma}$,
$C_n > 0$ and a set $U(n) \subset \Lambda$, $n \in U(n)$, satisfying
\begin{eqnarray}
\label{cl.20}
&&T_{U(n)}^{-1} \mbox{ exists and } 
\| T_{U(n)}^{-1} \|_{\tilde{\sigma} + \mu_n, c} \leq C_n, \\
\label{cl.30}
&&\mbox{ dist }(n, \Lambda \setminus U(n)) \geq l_n, \\
\label{cl.40}
&&C_n e^{- \mu_n l_n^c} \|R\|_{\sigma, c} \leq \frac{1}{2}.
\end{eqnarray}
Then the matrix $T_{\Lambda}$ is invertible with
\begin{eqnarray}
\label{cl.50}
\| T_{\Lambda}^{-1} \|_{\tilde{\sigma}, c} \leq 
(1 + w_{\tilde{\sigma}, c}(0)) C, \quad \mbox{ where } \;
C := \sup_{n \in \Lambda} C_n.
\end{eqnarray}
\end{Lemma}

\begin{proof}
The proof of lemma \ref{Lcl.1} proceeds in three steps.
First we define matrices $P$ and $K$ satisfying 
\begin{eqnarray}
\label{cl.60}
\left(T_{\Lambda} \right) P = I + K,
\end{eqnarray}
where $I$ denotes the identity matrix restricted to the sub-lattice
$\Lambda$.
In step 2 we show that 
\begin{eqnarray}
\label{cl.70}
\| K \|_{\tilde{\sigma}, c} \leq \frac{1}{2}.
\end{eqnarray}
Since the norm is sub-multiplicative (see proposition \ref{Pfa1.2} and
lemma \ref{Lwf.1})
the inverse of $I+K$ can be expanded in
a Neumann series with  $\| (I+K)^{-1} \|_{\tilde{\sigma}, c} \leq
\| I\|_{\tilde{\sigma}, c} + 1 = 1 + w_{\tilde{\sigma}, c}(0)$.
Hence we have constructed a right inverse of $T_{\Lambda}$ and, as
$\Lambda$ is finite, the inverse matrix
$T_{\Lambda}^{-1}$ exists. Claim (\ref{cl.50}) then follows from
\begin{eqnarray}
\label{cl.80}
\| P \|_{\tilde{\sigma}, c} \leq C,
\end{eqnarray}
which is derived in step 3. 

\noindent
{\bf Step 1:} Proof of (\ref{cl.60}).
Define 
\begin{eqnarray}
\label{cl.90}
P(m, n) &:=&
\left\{
\begin{array}{ll}
T_{U(n)}^{-1}(m,n)&\mbox{for } n \in \Lambda, m \in U(n),\\
0&\mbox{ else },
\end{array}
\right. \\
\label{cl.100}
K(m, n) &:=&
\left\{
\begin{array}{ll}
\sum_{p \in U(n)} R(m, p) P(p,n) &\mbox{for } n \in \Lambda, 
m \in \Lambda \setminus U(n),\\
0&\mbox{ else }.
\end{array}
\right.
\end{eqnarray}
Relation (\ref{cl.60}) then follows from
\begin{eqnarray}
\nonumber 
\sum_{p \in \Lambda} \!\!  \left(T_{\Lambda}\right)(m, p)
P(p, n) =
\left\{
\begin{array}{ll}
\sum_{p \in U(n)}  T(m,p) T_{U(n)}^{-1}(p,n)= \delta_{m, n}
&\mbox{for } n \in \Lambda, \;\; m \in U(n), \\
\sum_{p \in U(n)} R(m, p) P(p, n)  = K(m, n)
&\mbox{for } n \in \Lambda, \;\; m \in \Lambda \setminus U(n). 
\end{array}
\right.
\end{eqnarray}

\noindent
{\bf Step 2:} Proof of (\ref{cl.70}). Fix $n \in \Lambda$. Using
lemma \ref{Lwf.1}, (\ref{cl.30}) we obtain
\begin{eqnarray}
\label{cl.120}
&&\sum_{m \in \Lambda} w_{\tilde{\sigma}, c}(m-n) |K(m, n)| 
\leq
\sum_{m \in \Lambda \setminus U(n)} 
w_{\tilde{\sigma}+\mu_n, c}(m - n) |K(m, n)| e^{-\mu_n l_n^c} \nonumber \\
&\leq&
\sum_{m \in \Lambda \setminus U(n)} \sum_{p \in U(n)} 
w_{\tilde{\sigma}+\mu_n, c}(m-p) |R(m, p)| 
w_{\tilde{\sigma}+\mu_n, c}(p-n) 
|T_{U(n)}^{-1}(p,n)| 
e^{-\mu_n l_n^c} \nonumber \\
\nonumber
&\leq&
\sum_{p \in U(n)} \| R \|_{\tilde{\sigma} + \mu_n, c} \;
w_{\tilde{\sigma}+\mu_n, c}(p-n) \;
|T_{U(n)}^{-1}(p, n)| \;
e^{-\mu_n l_n^c} \\
\nonumber
&\leq& 
\| R \|_{\tilde{\sigma} + \mu_n, c} \;
\| T_{U(n)}^{-1} \|_{\tilde{\sigma} + \mu_n, c} \; e^{-\mu_n l_n^c} \leq \frac{1}{2}.
\end{eqnarray}
 
\noindent
{\bf Step 3:} Proof of (\ref{cl.80}). Fix $n \in \Lambda$. Clearly,
\begin{eqnarray}
\label{cl.130}
\sum_{m \in \Lambda} w_{\tilde{\sigma}, c}(m - n) |P(m, n)| \leq 
\sum_{m \in U(n)} w_{\tilde{\sigma}+\mu_n, c}(m -n) 
|T_{U(n)}^{-1}(m,n)| \leq
\| T_{U(n)}^{-1} \|_{\tilde{\sigma} + \mu_n, c} \leq C_n \leq C. 
\nonumber
\end{eqnarray}

\end{proof}

\section{A version of the Weierstrass Preparation Theorem}
\label{wpt}

In our analysis of chapter II we need a version of 
the Weierstrass preparation theorem which provides estimates
on the derivatives of the coefficients of the resulting polynomials.
Such a version was stated and proved by Bourgain in \cite{B4}.
For the sake of completeness we reproduce Bourgain's proof in this section.

We begin with a simple application of the Banach fixed point theorem.
\begin{Proposition}
\label{Pwpt.1}
Let $(X, |\cdot|)$ be a Banach space, $x_0 \in X$, $\eta > 0$ 
and denote $K_{\eta}(x_0) := \{x \in X: |x-x_0| \leq \eta\}$.
Assume furthermore, that $F:K_{\eta}(x_0) \to X$ is a $C^1$ map with
\begin{eqnarray}
\label{wpt.2}
\|DF(x) - I\| \leq 1/5, \;\;\; \mbox{ for all } x \in K_{\eta}(x_0),
\end{eqnarray} 
where $\|\cdot \|$
denotes the operator norm and $I$ is the identity map.

If $|F(x_0)| \leq 2 \eta/5$ then there exists a unique $y \in K_{\eta}(x_0)$
with $F(y) = 0$. Furthermore, $|y - x_0| \leq \frac{5}{2} |F(x_0)|$.
\end{Proposition}
\begin{proof}
It follows from (\ref{wpt.2}) that $DF(x_0)$ is invertible with
$\| DF(x_0)^{-1} \| \leq 5/4$. We define
\begin{eqnarray}
\nonumber
T : K_{\eta}(x_0) \to X \; ; \;\;\; T(x) = x - DF(x_0)^{-1}F(x)
\end{eqnarray}
Note that the zeros of $F$ are precisely the fixed points of $T$.
For $x_1$, $x_2 \in K_{\eta}(x_0)$ we obtain 
\begin{eqnarray}
\nonumber
T(x_1) - T(x_2) &=& x_1 - x_2 + 
DF(x_0)^{-1} 
\left(
\int_{0}^{1} DF(x_1 + t(x_2 - x_1))dt
\right) (x_2 - x_1) \\
\nonumber
&=&
DF(x_0)^{-1} 
\left(
\int_{0}^{1} [DF(x_1 + t(x_2 - x_1)) - DF(x_0)]dt
\right) (x_2 - x_1).
\end{eqnarray}
From (\ref{wpt.2}) we conclude that $\|DF(x_1 + t(x_2 - x_1)) - DF(x_0)\|
\leq 2/5$ and therefore 
\begin{eqnarray}
\label{wpt.4}
|T(x_1) - T(x_2)| \leq \frac{5}{4} \cdot \frac{2}{5} |x_2 - x_1| = 
\frac{1}{2} |x_1 - x_2|.
\end{eqnarray}
Furthermore, for all $x \in K_{\eta}(x_0)$ we have
\begin{eqnarray}
\nonumber
|T(x) - x_0| \leq |T(x) - T(x_0)| + |T(x_0) - x_0| \leq
\frac{1}{2} |x - x_0| + \frac{5}{4} |F(x_0)| < \frac{\eta}{2} + 
\frac{5}{4} \cdot \frac{2 \eta}{5} = \eta.
\end{eqnarray}
This shows that $T$ is a contraction on $K_{\eta}(x_0)$ and has thus
an unique fixed point $y$. Finally,
\begin{eqnarray}
\nonumber
|y - x_0| = |T(y) - x_0| \leq |T(y) - T(x_0)| + |T(x_0) - x_0| 
\leq \frac{1}{2} |y - x_0| + \frac{5}{4} |F(x_0)|, 
\end{eqnarray}  
from which $|y- x_0| \leq \frac{5}{2} |F(x_0)|$ follows.
\end{proof}

\begin{Lemma}
\label{Lwpt.1}
Let $d$, $\nu$, $B_1$, $B_2 \in {\Bbb N}$, $\lambda_0 \in {\Bbb R}^{2 \nu}$,  
and let $\delta$, $\rho$, $\epsilon$, 
$C$, $C^*$ be positive constants satisfying
\begin{eqnarray}
\label{wpt.5}
\delta \leq 1 \\
\label{wpt.7}
C^{*} \geq 1 \\
\label{wpt.10}
200 d \epsilon < 1 \\
\label{wpt.20}
30 C^{*} (B_1!) 2^{B_1} \leq C \\
\label{wpt.30}
d(B_1 + 1) \leq B_2 \\
\label{wpt.40}
16 \nu C d \rho < \left(\frac{\delta}{6d}\right)^d \\
\label{wpt.50}
(6d)^{d(B_1+1)} \leq B_2! 
\end{eqnarray}
Assume further that 
\begin{eqnarray}
\nonumber
f: \{ z \in {\Bbb C}: |z| , \delta \} \times
\{\lambda \in {\Bbb C}^{2 \nu}: |\lambda - \lambda_0|< \rho \} \to {\Bbb C}
\end{eqnarray}
is an analytic function of the form
\begin{eqnarray}
\label{wpt.70}
f(z, \lambda) = z^d + \sum_{0 \leq j < d} a_j (\lambda) z^j + r(z, \lambda),
\end{eqnarray}
where for $z \in U_{\delta}(0)$, $\lambda \in U_{\rho}(\lambda_0)$,
\begin{eqnarray}
\label{wpt.80}
|a_j(\lambda)| &\leq& \frac{1}{8d} \;\;\;
\mbox{ for all } 0 \leq j < d, \\
\label{wpt.90}
\left|\partial^{\beta} a_j (\lambda)
\right| &\leq& C \;\;\;
\mbox{ for all } 1 \leq |\beta|_1 \leq B_1, \;
0 \leq j < d, \\
\label{wpt.100}
\left|
\partial^{\beta} r (z, \lambda)
\right|
&\leq&
\left\{
\begin{array}{llll}
\epsilon&\mbox{ for }& |\beta_{\lambda}|_1=0,
&0 \leq |\beta_{z}| \leq B_2,
\\
C^{*}&\mbox{ for }& 1 \leq |\beta_{\lambda}|_1 \leq B_1 ,
&0 \leq |\beta_{z}| \leq B_2.
\end{array}
\right.
\end{eqnarray}
Suppose that $a_j$, $r$ are real functions, i.e. 
$a_j(\lambda)$, $r(z, \lambda) \in {\Bbb R}$ for 
$z \in U_{\delta}(0) \cap {\Bbb R}$, 
$\lambda \in U_{\rho}(\lambda_0) \cap {\Bbb R}^{2 \nu}$.

Then there exist functions
\begin{eqnarray}
\nonumber
Q &:& U_{\delta/4}(0) \times U_{\rho}(\lambda_0) \to {\Bbb C}, \\
\nonumber
b_j &:& U_{\rho}(\lambda_0)
\to {\Bbb C}, \;\;
0 \leq j < d,
\end{eqnarray}
such that for all $z \in U_{\delta/4}(0)$, $\lambda \in U_{\rho}(\lambda_0)$ 
and $0 \leq j < d$ the following holds:
\begin{eqnarray}
\label{wpt.120}
f(z, \lambda) &=& [ 1 + Q(z, \lambda)] 
\left(
z^d + \sum_{0 \leq j < d} b_j (\lambda) z^j
\right), \\
\label{wpt.130}
|Q(z, \lambda)| &\leq& \frac{1}{10}, \\
\label{wpt.140}
|b_j(\lambda)| &\leq& \frac{1}{2d}, \\
\label{wpt.150}
\left| \partial^{\beta} b_j (\lambda) \right|
&\leq&
\left( \frac{5}{2} C \right)^{2 |\beta|_1 - 1} (d + B_1)^{|\beta|_1^2}
\;\;\; \mbox{ for } 1 \leq |\beta|_1 \leq B_1, \\
\label{wpt.160}
b_j(\lambda) &\in& {\Bbb R} \;\;\; \mbox{ for }
\lambda \in U_{\rho}(\lambda_0) \cap {\Bbb R}^{2 \nu}.
\end{eqnarray}

\end{Lemma}

\begin{proof}
Without loss of generality we assume for the proof that $\lambda_0 = 0$.
Recall the basic idea to prove the Weierstrass preparation theorem.
Denote for $b \in {\Bbb C}^d$ the polynomial 
\begin{eqnarray}
\label{wpt.200}
p_b(z) = z^d + \sum_{0 \leq j < d} b_j z^j.
\end{eqnarray}
Computing $(p_b(s) -p_b(z))/(s-z)$ one easily derives the following
formula
\begin{eqnarray}
\label{wpt.210}
\frac{1}{s - z} &=&
\frac{1}{p_b(s) (s-z)} p_b(z) + \sum_{k=0}^{d-1} \frac{q_k(b, s)}{p_b(s)} z^k,
\;\;\; \mbox{ where }
\\
\label{wpt.220}
q_k(b, s) &:=& s^{d-1-k} + \sum_{l=k+1}^{d-1} b_l s^{l-1-k}.
\end{eqnarray}
Let $0 < \alpha < \delta$.
By the Cauchy integral formula
\begin{eqnarray}
\label{wpt.225}
r(z, \lambda) =
\frac{1}{2 \pi i} \oint_{|s| = \alpha} \frac{r(s, \lambda)}{s-z}ds
\end{eqnarray}
for $|z| < \alpha$.
Using representation (\ref{wpt.210}) we obtain
\begin{eqnarray}
\label{wpt.230}
f(z, \lambda) &=& 
\left(
1 + \oint_{|s| = \alpha} \frac{r(s, \lambda)}{p_b(s) (s-z)}\frac{ds}{2 \pi i}
\right) p_b (z)
\\
\nonumber
&&+
\sum_{k=0}^{d-1}
\left(a_k(\lambda) - b_k + 
\oint_{|s| = \alpha} \frac{r(s, \lambda) q_k(b, s)}
{p_b(s)}\frac{ds}{2 \pi i}
\right) z^k. 
\end{eqnarray}
To prove the Weierstrass preparation theorem one needs to show
that $b = b(\lambda)$ can be chosen such that the second line 
in (\ref{wpt.230}) vanishes.

For our purposes it is convenient to split the remainder term $r$ into two
parts and apply the above described procedure twice with different choices of
$\alpha$. According to Taylor's formula we decompose
\begin{eqnarray}
\label{wpt.250}
r(z, \lambda) &=& \sum_{j=0}^{B_2-1} c_j(\lambda) z^j
+ \tilde{r}(z, \lambda),
\;\;\; \mbox{ with } \\
\label{wpt.260}
c_j(\lambda) &:=& \frac{1}{j!} \partial_z^{j} r (0, \lambda), \\
\label{wpt.270}
\tilde{r}(z, \lambda) &:=&
\frac{1}{(B_2 - 1)!} \int_0^{z} 
\left( \partial_z^{B_2} r \right) (s, \lambda) (z - s)^{B_2 - 1} ds
\end{eqnarray} 
and denote
\begin{eqnarray}
\label{wpt.280}
f(z, \lambda) = z^d + \sum_{0 \leq j < d} a_j (\lambda) z^j + 
\sum_{j=0}^{B_2-1} c_j(\lambda) z^j
+ \tilde{r}(z, \lambda) \equiv \tilde{f}(z, \lambda) + \tilde{r}(z, \lambda).
\end{eqnarray}
From hypothesis (\ref{wpt.100}) the following estimates are immediate.
\begin{eqnarray}
\label{wpt.290}
|c_j(\lambda)| \leq \frac{1}{j!} \epsilon&&\mbox{ for } \lambda \in
U_{\rho}(0), \\
\label{wpt.292}
|\partial^{\beta} c_j(\lambda) | \leq \frac{1}{j!} C^*&&  
\mbox{ for } \lambda \in U_{\rho}(0), \;\; 1 \leq |\beta|_1 \leq B_1, \\
\label{wpt.294}
|\tilde{r}(z, \lambda)| \leq \frac{1}{B_2!} \epsilon |z|^{B_2} 
&&\mbox{ for } z \in U_{\delta}(0), \;\; \lambda \in U_{\rho}(0), \\
\label{wpt.296}
|\partial^{\beta}_{\lambda} \tilde{r}(z, \lambda)| 
\leq \frac{1}{B_2!} C^* |z|^{B_2} 
&&\mbox{ for } z \in U_{\delta}(0), \;\; \lambda \in U_{\rho}(0),
\;\; 1 \leq |\beta|_1 \leq B_1.
\end{eqnarray}

In a first step we apply the above described procedure to the 
auxiliary function $\tilde{f}$ (see (\ref{wpt.280})).
For $b \in {\Bbb C}^{d}$, $|b| < 1/(2d)$ and 
$s \in {\Bbb C}$, $|s|=1$ we have $|p_b(s)| \geq 1/2$. We can therefore
apply (\ref{wpt.210}) and obtain for
$(b, z, \lambda) \in U_{1/(2d)}(0) \times U_{\delta}(0) \times
U_{\rho}(0) \subset {\Bbb C}^d \times {\Bbb C} \times {\Bbb C}^{2\nu}$
\begin{eqnarray}
\label{wpt.300}
\tilde{f}(z, \lambda) =
(1 + \tilde{Q}(b, z, \lambda)) p_b(z) +
\sum_{k=0}^{d-1}
(a_k(\lambda) - b_k + \tilde{R}_k(b, \lambda) )z^k,
\end{eqnarray}
where
\begin{eqnarray}
\label{wpt.310}
\tilde{Q}(b, z, \lambda) &:=& \sum_{j=0}^{B_2-1} c_j(\lambda)
\oint_{|s|=1} \frac{s^j}{p_b(s)(s-z)} \frac{ds}{2 \pi i}, \\
\label{wpt.320}
\tilde{R}_k(b, \lambda) &:=& \sum_{j=0}^{B_2-1} c_j(\lambda)
\oint_{|s|=1} \frac{q_k(b, s) s^j}{p_b(s)} \frac{ds}{2 \pi i},
\;\; \mbox{ for } 0 \leq k < d.
\end{eqnarray}
We will now determine $b = \tilde{b}(\lambda)$ in such a way that
$a_k(\lambda) - b_k - \tilde{R}(b, \lambda) = 0$ for all $0 \leq k < d$. 
To achieve this we define
\begin{eqnarray}
\nonumber
\tilde{G}&:& U_{1/(2d)}(0) \times U_{\rho}(0) \to {\Bbb C}^d; \;\;\;
\tilde{G} = (
\tilde{G}_0, \ldots \tilde{G}_{d-1}) \;\;\; \mbox{ with }\\
\label{wpt.330}
\tilde{G}_k(b, \lambda) &:=& b_k - a_k(\lambda) - \tilde{R}_k(b, \lambda).
\;\;\; \mbox{ for } 0 \leq k < d.
\end{eqnarray}

\noindent
{\em
{\bf Claim 1:}
For every $\lambda \in U_{\rho}(0)$ there exists an
unique $\tilde{b}(\lambda) \in U_{\frac{3}{8d}}(0)$ such that
$\tilde{G}(\tilde{b}(\lambda), \lambda) = 0$. Moreover, the function
$\lambda \mapsto \tilde{b}(\lambda)$ is analytic and the 
first order derivatives are bounded by
$|\partial^{\beta} \tilde{b} (\lambda)| \leq 2 C$ for all 
$\lambda \in U_{\rho}(0)$, $|\beta|_1 = 1$.
}

{\em Proof of claim 1:}
The definitions of $p_b$ and $q_k$ (see (\ref{wpt.200}), (\ref{wpt.220})) 
yield the following
estimates for $|b| < 1/(2d)$ and $|s|=1$.
\begin{eqnarray}
\label{wpt.350}
\left| \frac{q_k(b, s)s^j}{p_b(s)} \right| &\leq& 4, 
\;\;\; 0 \leq k < d, \;\; j \in {\Bbb N}_0, \\
\label{wpt.360}
\left| \partial_{b_l} \frac{q_k(b, s)s^j}{p_b(s)} \right| &\leq& 10,
\;\;\; 0 \leq k, l < d, \;\; j \in {\Bbb N}_0.
\end{eqnarray}
These estimates together with (\ref{wpt.290}), (\ref{wpt.10}),
(\ref{wpt.90}), (\ref{wpt.292}) imply for $|b| < 1/(2d)$,
$|\lambda| < \rho$, $0 \leq k, l < d$ (denoting $a:=(a_j)_{0 \leq j < d}$)
\begin{eqnarray}
\label{wpt.370}
|\tilde{G}_k(a(\lambda), \lambda)|
&=& |\tilde{R}_k(a(\lambda), \lambda)| 
\leq \sum_{j \geq 0} \frac{\epsilon}{j!} 4
\leq 12 \epsilon \leq \frac{1}{15d}, \\
\label{wpt.380}
\left| \partial_{b_l} \tilde{G}_k (b, \lambda)  - \delta_{k, l} \right|
&\leq& \sum_{j \geq 0} \frac{\epsilon}{j!} 10 \leq 30 \epsilon 
< \frac{1}{6d},  \\
\label{wpt.390}
\left| \partial_{\lambda}^{\beta} \tilde{G}_k (b, \lambda)  \right|
&\leq& C + \sum_{j \geq 0} \frac{C^*}{j!} 4 \leq C + 12 C^{*},
\;\;\mbox{ for } 1 \leq |\beta|_1 \leq B_1. 
\end{eqnarray}
We apply proposition \ref{Pwpt.1} to obtain $\tilde{b}(\lambda)$. 
Choose $(X, |\cdot|)$ to be ${\Bbb C}^d$
together with the maximum norm. We fix  
$\lambda \in B_{\rho}(0)$ and set $x_0 \equiv a(\lambda)$, $\eta \equiv
\frac{1}{4d}$ and
\begin{eqnarray}
\nonumber
F(x) := \tilde{G}(x , \lambda) \;\;\; \mbox{ for} \;\;
x \in K_{1/(4d)}(a(\lambda))
\subset U_{1/(2d)}(0) \;\;\; (cf. (\ref{wpt.80})).
\end{eqnarray}
Note that it follows from (\ref{wpt.380}) that $\|DF(x) - I\|= 
\|D_x \tilde{G}(x, \lambda) - I \| \leq 30 \epsilon d < 
\frac{1}{6}$. Inequality
(\ref{wpt.370}) implies 
that $|F(x_0)| \leq \frac{1}{15d} < \frac{2}{5} \eta$. Hence there
exists an unique $\tilde{b}(\lambda)$ satisfying $F(\tilde{b}(\lambda))=0$.
Furthermore $|\tilde{b}(\lambda) - a(\lambda)| \leq 2.5 \frac{1}{15d} <
\frac{1}{4d}$. Since $|a(\lambda)| \leq 1/(8d)$ by (\ref{wpt.80}) 
we conclude that
$\tilde{b}(\lambda) \in U_{\frac{3}{8d}}(0)$.

The analyticity of $\tilde{b}$ as a function of $\lambda$ follows 
from the analyticity of $\tilde{G}$ by
a standard implicit function theorem. Moreover, 
differentiating $G(\tilde{b}(\lambda), \lambda) = 0$ with respect to 
$\lambda_l$, $ 1 \leq l \leq 2 \nu$, and using
$\| D_b \tilde{G} - I \| \leq 1/5$, (\ref{wpt.390}), and (\ref{wpt.20}) 
we obtain 
\begin{eqnarray}
\nonumber
\left| \partial_{\lambda_l} \tilde{b} \right| \leq
\left| (D_b \tilde{G})^{-1} \partial_{\lambda_l} \tilde{G} \right| \leq 
\frac{5}{4} (C + 12 C^*) < 2 C.
\end{eqnarray}
This completes the proof of claim 1.

We have constructed a representation for $\tilde{f}$ in the desired form
\begin{eqnarray}
\nonumber
\tilde{f}(z, \lambda) = (1 + Q(\tilde{b}(\lambda), z, \lambda))
p_{\tilde{b}(\lambda)}(z).
\end{eqnarray}
Since $f = \tilde{f} + \tilde{r}$ we now apply (\ref{wpt.225}),
(\ref{wpt.210}) with $r$ being replaced by the smaller remainder
term $\tilde{r}$. First we choose the radius for the contour of integration.

\noindent
{\em
{\bf Claim 2:}
There exists $\frac{\delta}{2} < \alpha < \delta$ such that
\begin{eqnarray}
\label{wpt.396}
\left| p_{b}(s) \right| \geq \frac{1}{2}
\left(\frac{\delta}{6d} \right)^{d} \;\;\;
\mbox{ for } |s|= \alpha,\;\; \mbox{ and } \;\;  
b \in {\cal B} :=
U_{\frac{1}{2d} \left(\frac{\delta}{6d} \right)^{d}}(\tilde{b}(0)).
\end{eqnarray}
For all $\lambda \in U_{\rho}(0)$ the following estimate holds
\begin{eqnarray}
\label{wpt.398}
|\tilde{b}(\lambda) - \tilde{b}(0)| < \frac{1}{4d}
\left(\frac{\delta}{6d} \right)^{d}.
\end{eqnarray}
}

{\em Proof of claim 2:}
Estimate (\ref{wpt.398}) follows from claim 1 and (\ref{wpt.40}).
Denote the zeros of $p_{\tilde{b}(0)}$ by $\xi_1, \ldots, \xi_d$. 
Choose $\alpha \in (\delta/2, \delta)$ such that
$|\alpha - |\xi_l|| > \frac{\delta}{6d}$ for $1 \leq l \leq d$.
For $s \in {\Bbb C}$ with $|s| = \alpha$ and $b \in {\cal B}$ 
it follows from (\ref{wpt.5}) that 
\begin{eqnarray}
\nonumber
\left| p_{b}(s) \right| \geq 
\left| p_{\tilde{b}(0)}(s) \right| -  \frac{1}{2}
\left(\frac{\delta}{6d} \right)^{d}\geq
\prod_{l=1}^{d} ||s|-|\xi_l|| - 
\frac{1}{2} \left(\frac{\delta}{6d} \right)^{d}
\geq \frac{1}{2} \left(\frac{\delta}{6d} \right)^{d}.
\end{eqnarray}
This proves claim 2.

Note that $b \in {\cal B}$, (\ref{wpt.5}) and claim 1 imply
\begin{eqnarray}
\nonumber
|b| \leq |\tilde{b}(0)| + \frac{1}{12 d} < \frac{3}{8d} + \frac{1}{12 d} <
\frac{1}{2d}.
\end{eqnarray}
We conclude ${\cal B} \subset U_{1/(2d)}(0)$. 
For $z \in U_{\alpha}(0)$, $\lambda \in U_{\rho}(0)$ and 
$b \in {\cal B}$ we obtain the
following representation of the function $f$
\begin{eqnarray}
\label{wpt.400}
f(z, \lambda) =
(1 + Q(b, z, \lambda)) p_b(z) + 
\sum_{k=0}^{d-1} (a_k(\lambda) - b_k + R_k(b, \lambda) ) z^k,
\end{eqnarray}
where
\begin{eqnarray}
\label{wpt.410}
Q(b, z, \lambda) &:=& \tilde{Q}(b, z, \lambda) +
\oint_{|s|=\alpha} \frac{\tilde{r}(s, \lambda)}{p_b(s)(s-z)} 
\frac{ds}{2 \pi i}, \\
\label{wpt.420}
R_k(b, \lambda) &:=&
\tilde{R}_k(b, \lambda) + 
\oint_{|s|=\alpha} \frac{q_k(b, s) \tilde{r}(s, \lambda)}{p_b(s)} 
\frac{ds}{2 \pi i},
\;\; \mbox{ for } 0 \leq k < d.
\end{eqnarray}
We define
\begin{eqnarray}
\nonumber
G&:& {\cal B} \times U_{\rho}(0) \to {\Bbb C}^d; \;\;\;
G = (
G_0, \ldots, G_{d-1}) \;\;\;\mbox{ with } \\
\label{wpt.430}
G_k(b, \lambda) &:=& b_k - a_k(\lambda) - R_k(b, \lambda)
\;\; \mbox{ for } 0 \leq k < d.
\end{eqnarray}

\noindent
{\em
{\bf Claim 3:}
There exists an analytic function
$b:U_{\rho}(0) \to {\cal B}$ satisfying
$G(b(\lambda), \lambda) = 0$.}

{\em Proof of claim 3:}
Using $\tilde{G}(\tilde{b}(\lambda), \lambda) = 0$ 
we can express $G_k$ by
\begin{eqnarray}
\nonumber
G_k(b, \lambda) = (b_k - \tilde{b}_k (\lambda)) +
[\tilde{R}_k(\tilde{b}(\lambda), \lambda) - 
\tilde{R}_k(b, \lambda) ] +
[\tilde{R}_k(b, \lambda) - R_k(b, \lambda)].
\end{eqnarray}
Fix $\lambda \in U_{\rho}(0)$. We introduce the new variable
\begin{eqnarray}
\label{wpt.450}
x := b - \tilde{b}(\lambda).
\end{eqnarray}
By (\ref{wpt.396}) and (\ref{wpt.398}) it suffices to show that there
exists $x \in \overline{ 
U_{\frac{1}{4d} \left(\frac{\delta}{6d} \right)^{d}}(0)} \subset
{\Bbb C}^d$ solving $F(x) = 0$, where $F = (F_0, \ldots, F_{d-1})$ is 
defined by
\begin{eqnarray}
\label{wpt.460}
F_k(x) := x_k -
[\tilde{R}_k(x+ \tilde{b}(\lambda), \lambda) - 
\tilde{R}_k(\tilde{b}(\lambda), \lambda)]
-
[
R_k(x+ \tilde{b}(\lambda), \lambda) -
\tilde{R}_k(x+ \tilde{b}(\lambda), \lambda)
] 
\end{eqnarray}
for $0 \leq k < d$.
To obtain the existence of the zero of $F$ we again apply proposition
\ref{Pwpt.1}.
In the notation of proposition \ref{Pwpt.1} $(X, |\cdot|)$ is given by ${\Bbb C}^d$ together with
the maximum norm, $x_0 \equiv 0$ and $\eta \equiv
\frac{1}{4d} \left(\frac{\delta}{6d} \right)^{d}$. We need to derive
the estimates on $|F(0)|$ and $\|DF(x) - I\|$. For 
$0 \leq k, l < d$ the above definitions yield
\begin{eqnarray}  
\nonumber
F_k(0) &=& - \oint_{|s| = \alpha}
\frac{\tilde{r}(s, \lambda) q_k(\tilde{b}(\lambda), s)}
{p_{\tilde{b}(\lambda)}(s)}
\frac{ds}{2 \pi i}, \\
\nonumber
\partial_{x_l} F_k (x) - \delta_{l, k} &=&
- \sum_{j=0}^{B_2-1} c_j(\lambda) \oint_{|s|=1}
s^j \partial_{x_l}
\left(
\frac{q_k(x + \tilde{b}(\lambda), s)}
{p_{x + \tilde{b}(\lambda)}(s)}
\right)
\frac{ds}{2 \pi i}\\
\nonumber
&-&
\oint_{|s|=\alpha}
\tilde{r}(s, \lambda) \partial_{x_l}
\left(
\frac{q_k(x + \tilde{b}(\lambda), s)}
{P_{x + \tilde{b}(\lambda)}(s)}
\right)
\frac{ds}{2 \pi i}.
\end{eqnarray}
For $|s| = \alpha$ and 
$|x| \leq \frac{1}{4d} \left(\frac{\delta}{6d} \right)^{d}$ we have
$|x + \tilde{b}(\lambda) - \tilde{b}(0)| 
< \frac{1}{2d} \left(\frac{\delta}{6d} \right)^{d}$, i.e.
$x + \tilde{b}(\lambda) \in {\cal B} \subset U_{1/(2d)}(0)$. Then
claim 2 implies for $0 \leq l < d$ (use also $\alpha < \delta \leq 1$ by 
(\ref{wpt.5}))
\begin{eqnarray}
\nonumber
\left|
\partial_{x_l}
\left(
\frac{q_k(x + \tilde{b}(\lambda), s)}
{P_{x + \tilde{b}(\lambda)}(s)}
\right)
\right| \leq 8 \left(\frac{6d}{\delta}\right)^{2d}  + 
2 \left(\frac{6d}{\delta}\right)^{d} \leq
10 \left(\frac{6d}{\delta}\right)^{2d}.
\end{eqnarray}
Together with (\ref{wpt.294}), (\ref{wpt.296}), (\ref{wpt.5}), 
(\ref{wpt.30}), (\ref{wpt.50}), (\ref{wpt.10}), and (\ref{wpt.360}) 
we conclude
\begin{eqnarray}
\nonumber
|F(0)| &\leq& 4 \left(\frac{6d}{\delta}\right)^{d} \frac{\epsilon}{B_2!}
\delta^{B_2}  = \frac{4 (6d)^d}{B_2!} \epsilon \delta^{B_2 - d} \\
\label{wpt.500}
&\leq& 
4 \epsilon \left(\frac{\delta}{6d}\right)^{d B_1}
\leq \frac{1}{50 d} \left(\frac{\delta}{6d}\right)^{d}
< \frac{2 \eta}{5}, \\
\nonumber
\| DF(x) - I \| &\leq&
d \left(
\sum_{j=0}^{B_2-1} \frac{10 \epsilon}{j!} + 
\frac{\epsilon}{B_2!} \delta^{B_2}   
10 \left(\frac{6d}{\delta}\right)^{2d}
\right)
\\
\label{wpt.510}
&\leq&
\epsilon d \left(30 + 10 \left(\frac{\delta}{6d}\right)^{d (B_1-1)} \right)
\leq 40 \epsilon d \leq
\frac{1}{5}.
\end{eqnarray}
Proposition \ref{Pwpt.1} can hence be applied and we obtain $y \in 
{\Bbb C}^d$, $|y| < \eta = \frac{1}{4d} \left(\frac{\delta}{6d} \right)^{d}$,
solving $F(y) = 0$. Setting $b(\lambda) := y + \tilde{b}(\lambda)$
it is clear from (\ref{wpt.450}), (\ref{wpt.460}) 
that $G(b(\lambda), \lambda) = 0$. So far
we have defined $b(\lambda)$ pointwise 
for each $\lambda \in U_{\rho}(0)$. The
analytic dependence of $b$ on $\lambda$ follows from 
the analyticity of $G$ by a standard implicit
function theorem.
This completes the proof of claim 3.

By a slight abuse of notation we set 
\begin{eqnarray}
\label{wpt.550}
Q(z, \lambda) := Q(b(\lambda), z, \lambda),
\end{eqnarray}
with the function $b$ as defined in claim 3. From 
(\ref{wpt.400}) and claim 3 it is clear that
we have found a representation of $f$ of the desired form 
(\ref{wpt.120}), satisfying
(\ref{wpt.140}). 
Observe that the realness condition on $a_j$ and $r$ imply
that $G(b, \lambda) \in {\Bbb R}^d$ for $b \in {\cal B} \cap
{\Bbb R}^d$, $\lambda \in U_{\rho}(0) \cap {\Bbb R}^{2 \nu}$,
proving (\ref{wpt.160}). 
In the remaining two steps we will verify (\ref{wpt.130}) and 
(\ref{wpt.150}).

{\em
{\bf Claim 4:}
$|Q(z, \lambda)| \leq 1/10$ for all $|z| < \delta/4$ and 
$|\lambda| < \rho$.
}

{\em Proof of claim 4:}
It follows from (\ref{wpt.410}), (\ref{wpt.310}), (\ref{wpt.290}),
(\ref{wpt.294}), $\delta/2 < \alpha < \delta$, (\ref{wpt.396}),
(\ref{wpt.30}), (\ref{wpt.50}), (\ref{wpt.10}) 
that for $|z| < \delta/4$, $|\lambda| < \rho$
\begin{eqnarray}
\nonumber
|Q(z, \lambda)| &\leq& \sum_{j=0}^{B_2 -1} \frac{\epsilon}{j!} 4 +
\frac{\epsilon}{B_2!} \delta^{B_2} 2 \left(\frac{6d}{\delta}\right)^{d}
\frac{4}{\delta} \cdot \alpha\\
\label{wpt.570} 
&\leq& \epsilon 
\left( 
12 + \frac{8 (6d)^d}{B_2!} \delta^{B_2 - d}
\right)\\
\nonumber
&\leq&
\epsilon 
\left( 
12 + 8\left(\frac{\delta}{6d}\right)^{d B_1}
\right)
\leq 20 \epsilon < \frac{1}{10}.
\end{eqnarray}
Claim 4 is verified.

{\em
{\bf Claim 5:} For $0 \leq j < d$, $|\lambda| < \rho$ and 
$1 \leq |\beta|_1 \leq B_1$ the following holds:
\begin{eqnarray}
\nonumber
\left| \partial^{\beta} b_j (\lambda) \right|
\leq
\left( \frac{5}{2} C \right)^{2 |\beta|_1 - 1} (d + B_1)^{|\beta|_1^2}.
\end{eqnarray}
}

{\em Proof of claim 5:}
From the definition of $p_b$ and $q_k$ in (\ref{wpt.200}) and
(\ref{wpt.220}) we conclude that for $|s| \leq 1$
and $|b| \leq 1/(2d)$
\begin{eqnarray}
\nonumber
\left|
\partial^{\beta}_b
\frac{q_k(b, s)}{p_b(s)}
\right|
&\leq&
|q_k(b, s)|
\left|  
\partial^{\beta}_b
\frac{1}{p_b(s)}
\right| +
\sum_{k < l < d, \beta_l > 0} 
\beta_l
\left|
\partial^{\beta - e_l}_b
\frac{1}{p_b(s)}
\right| \\
\nonumber
&\leq&
2 |\beta|_1 !
\frac{1}{|p_b(s)|^{|\beta|_1 + 1}} +
|\beta|_1 ! \frac{1}{|p_b(s)|^{|\beta|_1}} \leq
\frac{|\beta|_1 !}{|p_b(s)|^{|\beta|_1}}
\left(
\frac{2}{|p_b(s)|} + 1
\right).
\end{eqnarray}
Using in addition (\ref{wpt.420}), (\ref{wpt.320}), (\ref{wpt.290}) --
(\ref{wpt.296}),(\ref{wpt.7}), 
(\ref{wpt.396}), (\ref{wpt.30}), (\ref{wpt.50}), and (\ref{wpt.20})
this implies for multi-indices $\beta = (\beta_b, \beta_{\lambda})$ with
$1 \leq |\beta|_1 \leq B_1$, $|\lambda| < \rho$
and $|b| \leq 1/(2d)$ that
\begin{eqnarray}
\nonumber
\left| (\partial^{\beta} R)(b, \lambda)
\right| &\leq&
\sum_{j=0}^{B_2 -1} 
\left| (\partial^{\beta_{\lambda}} c_j)(\lambda) \right|
\oint_{|s|=1}
|\beta_b|_1 ! |p_b(s)|^{-|\beta_b|_1}
\left(
\frac{2}{|p_b(s)|} + 1
\right)
\frac{d|s|}{2 \pi} \\
\nonumber
&+&
\oint_{|s|=\alpha}
\left| (\partial_{\lambda}^{\beta_{\lambda}} \tilde{r})(s, \lambda) \right|
|\beta_b|_1 ! |p_b(s)|^{-|\beta_b|_1}
\left(
\frac{2}{|p_b(s)|} + 1
\right)
\frac{d|s|}{2 \pi} \\ 
\nonumber
&\leq&
3 C^{*} (B_1!) 2^{B_1} 5 + \frac{\delta^{B_2} C^*}{B_2!} (B_1!) 
3 \alpha \left(2  \left(\frac{6d}{\delta}\right)^d \right)^{B_1 + 1} \\
\label{wpt.700}
&\leq&
C^{*} B_1! \left( 15 \cdot 2^{B_1} + 3 \cdot 2^{B_1 + 1} \right) 
\leq
C.
\end{eqnarray}
For $p \in {\Bbb N}$ set 
\begin{eqnarray}
\label{wpt.710}
C(p) := (2.5 C)^{2 p -1} (d +B_1)^{p^2}.
\end{eqnarray}
We conclude the proof by showing that for $1 \leq |\beta|_1 \leq B_1$ we have
\begin{eqnarray}
\label{wpt.720}
|\partial^{\beta}b(\lambda)| \leq C(|\beta|_1) \;\;\;
\mbox{ for } \lambda \in U_{\rho}(0).
\end{eqnarray}
The proof of (\ref{wpt.720}) proceeds by induction on $p \equiv |\beta|_1$.

$\underline{p=1}:$ Differentiating
\begin{eqnarray}
\label{wpt.730}
b(\lambda) - a(\lambda) - R(b(\lambda), \lambda) = 0
\end{eqnarray} 
with respect to $\lambda_l$ ,$1 \leq l \leq 2 \nu$ yields
\begin{eqnarray}
\nonumber
(I - (D_b R)) (\partial_{\lambda_l} b) =
\partial_{\lambda_l} a + \partial_{\lambda_l} R.
\end{eqnarray}
Since $D_b R = I - D_b G = I - DF$ we learn from (\ref{wpt.510})
that
$\| D_b R (b(\lambda), \lambda) \| \leq 1/5$ for all 
$\lambda \in U_{\rho}(0)$. Using in addition (\ref{wpt.700}) and
(\ref{wpt.90}) we conclude
\begin{eqnarray}
\label{wpt.750}
|\partial_{\lambda_l} b| \leq \frac{5}{4} (C + C) = 2.5 C \leq C(1).
\end{eqnarray}  

$\underline{p-1 \to p:}$ Let $\beta$ be a multi-index of order $p$ with 
$2 \leq p \leq B_1$. Applying $\partial^{\beta}$ to 
(\ref{wpt.730}) we obtain with proposition \ref{Pcr4.1} ($d_1 = 2 \nu$,
$d_2 = d_3 = d$)
\begin{eqnarray}
\label{wpt.760}
(I - (D_b R)) (\partial^{\beta} b) = 
\partial^{\beta} a + 
\sum_{s \in S_0(\beta)} 
\left( \partial^{\gamma^{(s)}} R \right)
\prod_{l=1}^{l^{(s)}}
\partial^{\beta_l^{(s)}} b_{i_l^{(s)}},
\end{eqnarray}
where $\# S_0(\beta) \leq (d + B_1)^p$, 
and for all $s \in S_0(\beta)$ we have
$1 \leq |\gamma^{(s)}|_1 \leq p$, $1 \leq |\beta_l^{(s)}|_1 \leq p-1$ and
$\sum_{1 \leq l \leq l^{(s)}} |\beta_l^{(s)}|_1 \leq p$.
Using again that $\| D_b R \| \leq 1/5$ conclude with (\ref{wpt.90}),
(\ref{wpt.700})
\begin{eqnarray}
\label{wpt.800}
|\partial^{\beta} b| \leq \frac{5}{4}
\left[
C + (d + B_1)^p C \cdot
\left(
\max_{s \in S_0(\beta)}
\prod_{l=1}^{l^{(s)}}
C(|\beta_l^{(s)}|_1) 
\right)
\right].
\end{eqnarray}
To estimate the maximum appearing in (\ref{wpt.800}) observe that
for $p_1 + \ldots + p_k = p$ with $1 \leq p_l \leq p-1$ for $1 \leq l \leq k$
we have
\begin{eqnarray} 
\nonumber
C(p_1) \cdot \ldots \cdot C(p_k) \leq (2.5 C)^{2p -2} (d+B_1)^{p^2 -2(p-1)}.
\end{eqnarray}
Since $C \geq 1$ (see (\ref{wpt.7}), (\ref{wpt.20})) we find 
\begin{eqnarray}
\nonumber
|\partial^{\beta} b| &\leq& \frac{5}{4} C
\left(
1 + (d + B_1)^p (2.5C)^{2p -2} (d+B_1)^{p^2 -2(p-1)}
\right) \\
\label{wpt.820}
&\leq&  \frac{5}{2} C 
(d + B_1)^p (2.5 C)^{2p -2} (d+B_1)^{p^2 -2(p-1)} \leq C(p).
\end{eqnarray}
The induction is complete. Claim 5 and hence lemma \ref{Lwpt.1} 
are proven.

\end{proof}

\section{An excision lemma}
\label{el}
Our estimates on the measure of the set of resonant parameters
use the following elementary lemma.
 
\begin{Lemma}
\label{Lel.1}
Let $k \in {\Bbb N}$, let $J \subset {\Bbb R}$ be an 
interval and assume that
$g: J \to {\Bbb R}$ is $k$-times continuously differentiable.
Suppose further that there exists a constant $\delta > 0$ with
\begin{eqnarray}
\label{el.10}
\left| g^{(k)} (x) \right| \geq k! \delta^k \;\;\;
\mbox{ for all } x \in J, 
\end{eqnarray}
where $g^{(k)}$ denotes the $k$-th derivative of $g$.
For $t > 0$ define $I_t := \{ x \in J: |g(x)| < t \}$. 

Then the set $I_t$ is a union of at most $k$ intervals of total length
$|I_t| \leq \frac{2k}{\delta} t^{\frac{1}{k}}$.
\end{Lemma}

\begin{proof}
We prove the lemma inductively by showing that for each $0 \leq l \leq k$
the following holds:
\begin{itemize}
\item[(a)]
The function $g^{(k-l)}$ has at most $l$ zeros.
\item[(b)]
The set $I^{(l)} :=
\{ x \in J: |g^{(k-l)}(x)| < \frac{k!}{l!} \delta^{k-l} t^{l/k} \}$ 
is a union of at most $l$ intervals of total length 
$\leq \frac{2l}{\delta} t^{\frac{1}{k}}$.
\end{itemize}

$\underline{l = 0:}$ Statements (a) and (b) follow from the hypothesis of the
lemma.

$\underline{l-1 \to l \mbox{ for } 1 \leq l \leq k:}$ 
Without loss of generality we 
may assume that $J \neq \emptyset$.
Denote $x_0 := \inf J$ and $x_{l} := \sup J$.
By induction hypothesis (a)
we can find points $x_1, \ldots, x_{l-1} \in J$ such that the corresponding 
intervals $J_s := (x_{s-1}, x_{s})$, $1 \leq s \leq l$ contain no zero
of $g^{(k-l+1)}$. Hence $g^{(k-l)}$ is a strictly monotone function in 
each of the intervals $J \cap \overline{J_s}$, $1 \leq s \leq l$.
This proves (a). Furthermore, for each $1 \leq s \leq l$, the set
\begin{eqnarray}
\nonumber
I_s^{(l)} := 
\{ x \in J \cap \overline{J_s}: 
|g^{(k-l)}(x)| < \frac{k!}{l!} \delta^{k-l} t^{l/k} \}
\end{eqnarray}
is a (possibly empty) interval.
Hence $I^{(l)}=\bigcup_{s=1}^l I_s^{(l)}$ 
is a union of at most $l$ intervals.
In order to estimate the length of $I_s^{(l)}$ we use
the lower bound on $|g^{(k-l+1)}|$ outside $I^{(l-1)}$. We obtain
\begin{eqnarray}
\nonumber
\left| I_s^{(l)} \right| &\leq&
\frac{2 \frac{k!}{l!} \delta^{k-l} t^{l/k} }
{\frac{k!}{(l-1)!} \delta^{k-l+1} t^{(l-1)/k}}
+
\left|
J_s \cap I^{(l-1)}
\right| \\
\label{el.50}
&\leq&
\frac{2 t^{1/k}}{l \delta} + 
\left|
J_s \cap I^{(l-1)}
\right|
\end{eqnarray} 
Using (\ref{el.50}) and the induction hypothesis the total length of $I^{(l)}$
can be bounded by
\begin{eqnarray}
\nonumber
\left| I^{(l)} \right| \leq 
\frac{2 t^{1/k}}{\delta} + 
\left|
I^{(l-1)}
\right| \leq \frac{2l}{\delta} t^{\frac{1}{k}}.
\end{eqnarray}
\end{proof}

\section{The construction of $p \ominus q$}
\label{res}

In this section we present a resultant type construction. More precisely, for
given polynomials
\begin{eqnarray}
\label{res.10}
p(z) &=& z^{d_1} + \sum_{j < d_1} a_j z^j \in {\Bbb C}[z], \\ 
\label{res.20}
q(z) &=& z^{d_2} + \sum_{j < d_2} b_j z^j \in {\Bbb C}[z],
\end{eqnarray}
of degree $d_1$, $d_2 \geq 1$ we will construct a polynomial $r$ of 
degree $d_1 d_2$ with
\begin{eqnarray}
\label{res.30}
r(x-y) = R_1(x, y) p(x) + R_2(x, y) q(y),
\end{eqnarray}
where $R_1$, $R_2 \in {\Bbb C}[x,y]$ are polynomials in two variables.
Furthermore we will obtain some information on how the coefficients
of $r$, $R_1$ and $R_2$ depend on the coefficients of $p$ and $q$.
We start with a few elementary definitions and observations.

Given $p$ and $q$ as in (\ref{res.10}), (\ref{res.20}). Denote 
\begin{eqnarray}
\label{res.40}
K := \{ \pm a_j: 0 \leq j < d_1 \} \cup \{ \pm b_j: 0 \leq j < d_2 \}.
\end{eqnarray}
For $k$, $l \in {\Bbb N}_0$ we define the following finite subsets of 
${\Bbb C}$
\begin{eqnarray}
\label{res.50}
{\cal T}(k, l) :=
\left\{
f \in {\Bbb C}:
f = \sum_{j=1}^{l'} \prod_{i=1}^k x_{i, j} ; \; \mbox { for } 
x_{i, j} \in \{0, \pm 1\} \cup K \mbox{ and } 0 \leq l' \leq l
\right\},
\end{eqnarray}
where we use the standard definition
that an empty sum has value $0$ and an empty product has value $1$.
The following two properties are immediate from definition (\ref{res.50}).
For $k$, $l$, $m$, $n \in {\Bbb N}_0$
\begin{eqnarray}
\label{res.60}
{\cal T}(k, l) &\subset& {\cal T}(m, n) \;\;\; \mbox{ if }
k \leq m \; \mbox{ and } l \leq n, \\
\label{res.70}
{\cal T}(k, l) \cdot {\cal T}(m, n) &\subset& {\cal T}(k+m, ln),
\end{eqnarray}
where 
$
{\cal T}(k, l) \cdot {\cal T}(m, n) = \{ f \cdot g:
f \in {\cal T}(k, l), g \in {\cal T}(m, n) \}
$.

\begin{Proposition}
\label{Pres.1}
Let $p$, $q \in {\Bbb C}[z]$ be given as in (\ref{res.10}), (\ref{res.20})
with degrees $d_1$, $d_2 \geq 1$ and let $s \in {\Bbb N}_0$.
Then there exist unique polynomials $g_s$, $h_s$, $v_s$, $w_s$ of degrees
$\deg(g_s) \leq s$, $\deg(h_s) \leq s$, $\deg(v_s) < d_1$, $\deg(w_s) < d_2$ and
\begin{eqnarray}
\label{res.80}
z^s&=&
g_s(z) p(z) + v_s(z),
\\
z^s&=&
h_s(z) q(z) + w_s(z).
\label{res.90}
\end{eqnarray}
The coefficients of $g_s$, $h_s$, $v_s$, $w_s$ all lie in the set 
${\cal T}(s, 2^s)$.
\end{Proposition}

\begin{proof}
We only need to proof the result for the polynomial $p$. The unique
existence of $g_s$ and $v_s$ satisfying (\ref{res.80}) and 
$\deg(g_s) \leq s$, $\deg(v_s) < d_1$ is obvious.
Denote
\begin{eqnarray}
\nonumber
g_s(z) &=& \sum_{k=0}^s g_k^{(s)} z^k, \\
\nonumber
v_s(z) &=& \sum_{k=0}^{d_1 -1} v_k^{(s)} z^k.
\end{eqnarray}
We show that $g_k^{(s)}$, $v_k^{(s)} \in {\cal T}(s, 2^s)$ by an induction
in $s$.

$\underline{s=0:}$ Obviously, $g_0^{(0)} = 0$, $v_0^{(0)} = 1$ and
$v_k^{(0)} = 0$ for $1 \leq k < d_1$. By (\ref{res.50}) 
all coefficients lie in ${\cal T}(0, 1)$.

$\underline{s \to s+1:}$
The induction hypothesis yields
\begin{eqnarray}
\nonumber
z^{s+1} =
\left(
\sum_{k=0}^s g_k^{(s)} z^{k+1}
\right) p(z) + 
\sum_{k=0}^{d_1 -2} v_k^{(s)} z^{k+1}
+ 
v_{d_1 -1}^{(s)} z^{d_1}.  
\end{eqnarray}
Applying (\ref{res.10}) to the last term in the sum we obtain
\begin{eqnarray}
\nonumber
g_{k}^{(s+1)} &=&
\left\{
\begin{array}{ll}
v_{d_1-1}^{(s)}&\mbox{ if } k=0,\\
g_{k-1}^{(s)} & \mbox{ if } 1 \leq k \leq s+1,
\end{array}
\right.
\\ 
\nonumber
v_{k}^{(s+1)} &=&
\left\{
\begin{array}{ll}
- a_0 v_{d_1-1}^{(s)}&\mbox{ if } k=0,\\
v_{k-1}^{(s)} - a_k v_{d_1 - 1}^{(s)} & \mbox{ if } 1 \leq k < d_1.
\end{array}
\right.
\end{eqnarray}
The claim $g_k^{(s+1)}$, $v_k^{(s+1)} \in {\cal T}(s+1, 2^{s+1})$ now follows
from the induction hypothesis, (\ref{res.60}) and from 
definition (\ref{res.50}).
\end{proof}

\begin{Lemma}
\label{Lres.1}
Let $p$, $q \in {\Bbb C}[z]$ be given as in (\ref{res.10}), (\ref{res.20})
with degrees $d_1$, $d_2 \geq 1$.
Then there exist polynomials $r \in {\Bbb C}[z]$, $R_1$, $R_2 \in
{\Bbb C}[x, y]$ satisfying
\begin{eqnarray}
\label{res.95}
&&
r(x-y) = R_1(x, y) p(x) + R_2(x, y) q(y),
\\
\label{res.100}
&&
\deg(r) = d_1 d_2, \;\;\; \mbox{ and the leading coefficient of 
$r$ equals $1$ },
\\
\label{res.105}
&&\deg(R_1) < d_1 d_2, \deg(R_2) < d_1 d_2,
\\
\label{res.110}
&&
\mbox{ the coefficients of $r$ lie in } {\cal T}\left(d_1 + d_2, 
2^{d_1+d_2}[(d_1 d_2)!] \right),
\\
\label{res.120}
&&
\mbox{ the coefficients of $R_1$, $R_2$ lie in } 
{\cal T}(\tilde{d}, 4^{\tilde{d}} \tilde{d}!), \; \mbox{ where }
\tilde{d} := (d_1 + 1)(d_2 + 1).
\end{eqnarray}
\end{Lemma}

\begin{proof}
Set $V := {\Bbb C}[x, y] /(p(x) {\Bbb C}[x, y] + q(y) {\Bbb C}[x, y])$.
The quotient $V$ is a vector space of (complex) dimension $d_1 d_2$ with
basis $e_{i, j} := x^i y^j$, $0 \leq i < d_1$, $0 \leq j < d_2$. 
We define the linear map
\begin{eqnarray}
\label{res.130}
A: V \to V, \;\; (Ag)(x, y) = (x - y) g(x, y).
\end{eqnarray}
We set $r$ to be the characteristic polynomial of A
\begin{eqnarray}
\label{res.140}
r(z) := \det(z - A).
\end{eqnarray}
Then $r$ is a polynomial of degree $d_1 d_2$ with leading coefficient $1$
and the theorem of 
Cayley-Hamilton implies that $r(A) e_{0, 0} = 0$ in $V$. This means
\begin{eqnarray}
\label{res.150}
r(x-y) \in p(x) {\Bbb C}[x, y] + q(y) {\Bbb C}[x, y].
\end{eqnarray}
Hence there exist $R_1$, $R_2 \in {\Bbb C}[x, y]$ satisfying (\ref{res.95})
and (\ref{res.100}). Observe that $R_1$ and $R_2$ are not
uniquely defined. However, in order to show (\ref{res.105}) and 
(\ref{res.120}) we will now make a 
special choice for $R_1$ and $R_2$ using the 
notation introduced in proposition
\ref{Pres.1}. Write
\begin{eqnarray}
\label{res.160}
r(z) = \sum_{s=0}^{d_1 d_2} c_s z^s = z^{d_1 d_2} + 
\sum_{s < d_1 d_2} c_s z^s
\end{eqnarray}
and note that
\begin{eqnarray}
\nonumber
(x-y)^s &=& \sum_{l=0}^{s}
\left( \begin{array}{c} s \\ l \end{array} \right)
(-1)^{s-l} x^l y^{s-l} \\
\nonumber
&=&
\left(
\sum_{l=0}^{s}
\left( \begin{array}{c} s \\ l \end{array} \right)
(-1)^{s-l} g_l(x) y^{s-l}
\right) p(x) \\
\nonumber
&& +
\left(
\sum_{l=0}^{s}
\left( \begin{array}{c} s \\ l \end{array} \right)
(-1)^{s-l} v_l(x) h_{s-l}(y)
\right) q(y) \\
\nonumber
&&+
\sum_{l=0}^{s}
\left( \begin{array}{c} s \\ l \end{array} \right)
(-1)^{s-l} v_l(x) w_{s-l}(y).
\end{eqnarray}
It follows from $\deg(v_l) < d_1$, $\deg(w_{s-l}) < d_2$ and 
from (\ref{res.150}) that
\begin{eqnarray}
\nonumber
\sum_{s=0}^{d_1 d_2} c_s
\sum_{l=0}^{s}
\left( \begin{array}{c} s \\ l \end{array} \right)
(-1)^{s-l} v_l(x) w_{s-l}(y) = 0.
\end{eqnarray}
Thus we obtain a representation of the form (\ref{res.95}), satisfying
(\ref{res.100}) and (\ref{res.105}) by setting
\begin{eqnarray}
\label{res.170}
R_1(x, y) &=&
\sum_{s=0}^{d_1 d_2} 
\sum_{l=0}^{s} c_s 
\left( \begin{array}{c} s \\ l \end{array} \right)
(-1)^{s-l} g_l(x) y^{s-l}, \\
\label{res.180}
R_2(x, y) &=&
\sum_{s=0}^{d_1 d_2} 
\sum_{l=0}^{s} c_s 
\left( \begin{array}{c} s \\ l \end{array} \right)
(-1)^{s-l} v_l(x) h_{s-l}(y).
\end{eqnarray}
In order to prove (\ref{res.110}) we determine the entries of the matrix $z-A$
when expressed with respect to the basis $(e_{i, j})$. We compute
$(z - A) e_{i, j} =$
\begin{eqnarray}
\nonumber
ze_{i, j} - e_{i+1, j} + e_{i, j+1}&&\mbox{ if } i < d_1 -1,
j < d_2 -1, \\
\nonumber
(z+ a_{d_1-1}) e_{i, j} + \sum_{l=0}^{d_1 - 2} a_l e_{l, j} + e_{i, j+1}
&&\mbox{ if } i = d_1 -1, j < d_2 -1, \\
\nonumber
(z - b_{d_2-1}) e_{i, j} - e_{i+1, j} - \sum_{l=0}^{d_2 - 2} b_l e_{i, l}
&&\mbox{ if } i < d_1 -1, j = d_2 -1, \\
\nonumber
(z+ a_{d_1-1}- b_{d_2-1}) e_{i, j} + \sum_{l=0}^{d_1 - 2} a_l e_{l, j}
- \sum_{l=0}^{d_2 - 2} b_l e_{i, l} 
&&\mbox{ if } i = d_1 -1, j = d_2 -1.
\end{eqnarray}
Observe that all entries of $z - A$ are contained in the set
\begin{eqnarray}
\nonumber
\{0, 1, -1 \} \cup K \cup 
\{ z, z+ a_{d_1-1}, z - b_{d_2-1}, z+ a_{d_1-1}- b_{d_2-1} \}.
\end{eqnarray}
Expanding the determinant of $z - A$ in a sum of products and multiplying
out the diagonal entries $z+ a_{d_1-1}$, $z - b_{d_2-1}$ and 
$z+ a_{d_1-1}- b_{d_2-1}$ one sees easily that the coefficients $c_s$,
$0 \leq s < d_1 d_2$ of
the characteristic polynomial $r$ satisfy
\begin{eqnarray}
\label{res.190}
c_s \in {\cal T}(d_1 + d_2 - 1, 3 \cdot 2^{d_1 + d_2 - 2} [(d_1 d_2)!] ).
\end{eqnarray}
This proves (\ref{res.110}) (see (\ref{res.60})).
Finally, using (\ref{res.170}), (\ref{res.180}), (\ref{res.190}),
proposition \ref{Pres.1}, (\ref{res.70}), and (\ref{res.60}) 
we conclude that the 
coefficients of $R_1$, $R_2$ lie in the set
\begin{eqnarray}
\nonumber
{\cal T}
\left( d_1 d_2 + d_1 + d_2 - 1,
(d_1 d_2 + 1) 2^{d_1 d_2} \cdot 3 \cdot 2^{d_1 + d_2 -2} [(d_1 d_2)!] 2^{d_1 d_2}  
\right),
\end{eqnarray}
proving (\ref{res.120}).
\end{proof}

\begin{Definition}
\label{Dres.1}
Given $p$, $q \in {\Bbb C}[z]$ as in (\ref{res.10}), (\ref{res.20}) with
$d_1$, $d_2 \geq 1$. Then we denote by $p \ominus q$ the polynomial $r$
as defined through (\ref{res.130}), (\ref{res.140}) in the proof of lemma
\ref{Lres.1}.
\end{Definition}

\begin{Remark}
\label{Rres.1}
The following example shows that for given polynomials $p$ and $q$ (as in
(\ref{res.10}), (\ref{res.20})) the polynomial $r$ is not necessarily uniquely
defined by (\ref{res.95}) and (\ref{res.100}).
In fact, set $p(z) = q(z) = z^2$ and define $A$ as in (\ref{res.130}). One
checks easily that $A^3 = 0$ and hence for $r(z) = z^4 + c z^3$ a
representation of
the form (\ref{res.95}) can be found for all $c \in {\Bbb C}$.
\end{Remark}

\section{Remarks on the higher order chain rule}
\label{cr}

In this section we collect those specialized versions of the higher order
chain rule which are used in the paper. We formulate them for
functions which are defined on open subsets of ${\Bbb R}^d$ and which
satisfy some finite regularity assumptions. Of course, the formulae
also hold in the case of analytic functions defined on open
subsets of ${\Bbb C}^d$.

\begin{Proposition}
\label{Pcr.1}
Let $U \subset {\Bbb R}^d$ be an open set and let $X$ and $Y$ be Banach
spaces. Assume that $V \subset X$ is an open set and that $f: U \to
V$ and $g:V \to Y$ are $C^s$ -- functions with $s \in {\Bbb N}$.
Then $h := g \circ f : U \to Y$ is also $s$-times continuously
differentiable. For multi-indices $\beta$ with 
$1 \leq |\beta|_1 \leq s$ and $\lambda \in U$, 
\begin{eqnarray}
\label{cr.10}
\partial^{\beta} h(\lambda) =
\sum_{k=1}^{|\beta|_1} \frac{1}{k!}
\sum_{\scriptsize
{\scriptstyle
\begin{array}{c}
\alpha_1 + \ldots + \alpha_{k} = \beta \\
\alpha_i \neq 0 \mbox{ for } 1 \leq i \leq k
\end{array}
}
}
\frac{\beta!}{\alpha_1! \ldots \alpha_k!}
(D^{k}g)(f(\lambda))
[\partial^{\alpha_1}f(\lambda), \ldots, \partial^{\alpha_{k}}f(\lambda)].
\end{eqnarray}
\end{Proposition}
\begin{proof}
We prove (\ref{cr.10}) by calculating the corresponding term in the
Taylor-expansion of $h$. Observe first that for $\lambda$, $\mu \in U$,
\begin{eqnarray}
\label{cr.20}
f(\mu) - f(\lambda) = \sum_{1 \leq |\alpha| \leq s} \frac{1}{\alpha!}
(\partial^{\alpha}f)(\lambda) (\mu - \lambda)^{\alpha} + 
o(|\mu - \lambda|^s).
\end{eqnarray}
Furthermore,
\begin{eqnarray}
\label{cr.30}
g(v) - g(x) = \sum_{k=1}^{s} \frac{1}{k!}
(D^k g)(x)[v-x, \ldots, v-x] + o(|v-x|_{X}^s).
\end{eqnarray}
From (\ref{cr.20}) and (\ref{cr.30}) we obtain
\begin{eqnarray}
\nonumber
&&h(\mu) - h(\lambda)
\\
\nonumber
&=& \sum_{k=1}^{s} \frac{1}{k!}
(D^k g)(f(\lambda))
\left[
\sum_{1 \leq |\alpha| \leq s} \frac{1}{\alpha!}
(\partial^{\alpha}f)(\lambda) (\mu - \lambda)^{\alpha}, 
\ldots, 
\sum_{1 \leq |\alpha| \leq s} \frac{1}{\alpha!}
(\partial^{\alpha}f)(\lambda) (\mu - \lambda)^{\alpha}
\right] 
\\
\nonumber
&&+ o(|\mu - \lambda|^s)
\\
\nonumber
&=&
\sum_{1 \leq |\beta|_1 \leq s} c_{\beta}(\lambda) (\mu - \lambda)^{\beta}
+ o(|\mu - \lambda|^s), \;\;\;\; \mbox{ where} \\
\nonumber
c_{\beta}(\lambda) &:=& 
\sum_{k=1}^{|\beta|_1} \frac{1}{k!}
\sum_{\scriptsize
{\scriptstyle
\begin{array}{c}
\alpha_1 + \ldots + \alpha_{k} = \beta \\
\alpha_i \neq 0 \mbox{ for } 1 \leq i \leq k
\end{array}
}
}
\frac{1}{\alpha_1! \ldots \alpha_k!}
(D^{k}g)(f(\lambda))
[\partial^{\alpha_1}f(\lambda), \ldots, \partial^{\alpha_{k}}f(\lambda)].
\end{eqnarray}
This proves (\ref{cr.10}).
\end{proof}

To estimate the number of terms in (\ref{cr.10}) the following identity
is useful.
\begin{Proposition}
\label{Pcr.2}
Let $k$, $r \in {\Bbb N}$ and let $\beta \in {\Bbb N}_0^r$ be a multi-index.
Then
\begin{eqnarray}
\label{cr.100}
\sum_{
\alpha_1 + \ldots + \alpha_{k} = \beta
}
\frac{\beta!}{\alpha_1! \ldots \alpha_k!} = k^{|\beta|_1}.
\end{eqnarray}
\end{Proposition}

\begin{proof}
Applying the formula for multinomials we obtain for $k$, $s \in {\Bbb N}$
\begin{eqnarray}
\label{cr.150}
k^s = (1 + \ldots + 1)^s = \sum_{a_1 + \ldots + a_k = s}
\frac{s!}{a_1! \cdot \ldots \cdot a_k!}.
\end{eqnarray}
Furthermore, $k^{|\beta|_1} = k^{\beta_1} \cdot \ldots \cdot
k^{\beta_r}$. Applying (\ref{cr.150}) to each of the the factors
$k^{\beta_i}$, $1 \leq i \leq r$, yields (\ref{cr.100}).
\end{proof}

\subsection{Higher derivatives of $DW$}
\label{cr1}

We use the notation of section \ref{enp}.
The following is an immediate consequence of proposition \ref{Pcr.1}. 
\begin{Corollary}
\label{Ccr1.1}
Let $W$ be defined as in section \ref{enp} with respect to some real
power series $\sum \alpha_k y^k$ with positive radius of convergence
$0 < r \leq \infty$. Assume further that 
$u: U \to X$ is a $C^{\infty}$ function, where $U$ is some open 
subset of ${\Bbb R}^{d}$, $(X, | \cdot |)$ is a Banach space,
and $| u (\lambda) | < r$
for all $\lambda \in U$. Set $R(\lambda) := DW(u(\lambda))$, then
for $v \in X$
\begin{eqnarray}
\label{cr1.10}
(\partial^{\beta}_{\lambda} R)v =
\sum_{k=1}^{|\beta|_1} 
\sum_{\scriptsize
{\scriptstyle
\begin{array}{c}
\beta_1 + \ldots + \beta_{k} = \beta \\
\beta_i \neq 0 \; \mbox{ for } \; 1 \leq i \leq k
\end{array}
}
} 
\frac{1}{k!} 
\frac{\beta!}{\beta_1! \ldots \beta_k!}
(D^{k+1} W)(u)
[\partial^{\beta_1}u, \ldots, \partial^{\beta_{k}}u, v].
\end{eqnarray}
\end{Corollary}

\subsection{Derivatives of inverse matrices}
\label{cr2}

\begin{Corollary}
\label{Ccr2.1}
Let ${\cal B}$ be a Banach algebra with unity.  
Let $U \subset {\Bbb R}^d$ be an open set and assume that $G$, $H : U \to 
{\cal B}$ are $C^{s}$ -- maps, 
such that $H^{-1}(\lambda)$ exists in ${\cal B}$ for
all $\lambda \in U$. Then for $1 \leq |\beta|_1 \leq s$,
\begin{itemize}
\item[(a)]
\begin{eqnarray}
\nonumber
\partial^{\beta}(H^{-1}) =
\sum_{k=1}^{|\beta|_1} (-1)^{k}
\!\!\!\!\!\!\!\!\!\!\!\!\!\!
\sum_{\scriptsize
{\scriptstyle
\begin{array}{c}
\alpha_1 + \ldots + \alpha_{k} = \beta \\
\alpha_i \neq 0 \mbox{ for } 1 \leq i \leq k
\end{array}
}
}\!\!\!\!\!\!\!
\frac{\beta!}{\alpha_1! \ldots \alpha_k!}
H^{-1} (\partial^{\alpha_1}H) H^{-1} \ldots H^{-1} 
(\partial^{\alpha_k} H) H^{-1}
\end{eqnarray}

\item[(b)]
\begin{eqnarray}
\nonumber
\partial^{\beta}(G H^{-1}) =
\sum_{k=1}^{|\beta|_1+1} (-1)^{k-1}
\!\!\!\!\!\!\!\!\!\!\!\!\!\!\!
\sum_{\scriptsize
{\scriptstyle
\begin{array}{c}
\alpha_1 + \ldots + \alpha_{k} = \beta \\
\alpha_i \neq 0 \mbox{ for } 2 \leq i \leq k
\end{array}
}
}\!\!\!\!\!\!\!
\frac{\beta!}{\alpha_1! \ldots \alpha_k!}
(\partial^{\alpha_1} G)
H^{-1} (\partial^{\alpha_2}H) H^{-1} \ldots H^{-1} 
(\partial^{\alpha_k} H) H^{-1}
\end{eqnarray}
\end{itemize}
\end{Corollary}
\begin{proof}
To prove (a) it is convenient to repeat the proof of Proposition \ref{Pcr.1}
where one replaces (\ref{cr.30}) by the Neumann series
\begin{eqnarray}
\nonumber
V^{-1} - X^{-1} = -X^{-1}(V-X)X^{-1} + X^{-1}(V-X)X^{-1}(V-X)X^{-1} \mp
\ldots .
\end{eqnarray}
Claim (b) follows from (a) via the Leibniz rule.
\end{proof}

\begin{Proposition}
\label{Pcr2.1}
Let $({\cal B}, \| \cdot \|)$ 
be a Banach algebra with unity $e$ and let $U \subset
{\Bbb R}^d$ be an open set. Assume that
$A$, $X: U \to {\cal B}$ are $C^s$ -- maps, $s \in {\Bbb N}_0$, 
and $A^{-1}(\lambda)$ exists
for all $\lambda \in U$. Suppose further that there exist constants
$\epsilon > 0$, $M \geq 1$, $C \geq 1$, satisfying 
\begin{eqnarray}
\label{cr2.100}
\| \partial^{\beta} A^{-1}(\lambda) \| &\leq& M C^{|\beta|_1}
\;\;\; \mbox{ for } \lambda \in U, \; 0 \leq |\beta|_1 \leq s,
\\
\label{cr2.110}
\| \partial^{\beta} X(\lambda) \| &\leq& \epsilon C^{|\beta|_1}
\;\;\; \mbox{ for } \lambda \in U, \; 0 \leq |\beta|_1 \leq s,
\\
\label{cr2.120}
\epsilon M &\leq& \frac{1}{2}.
\end{eqnarray}
Then $(A+X)(\lambda)$ is invertible in ${\cal B}$ for all $\lambda \in U$ and
\begin{eqnarray}
\label{cr2.130}
\| \partial^{\beta} (A+X)^{-1}(\lambda) \| \leq 
2 \left(
\sum_{k=1}^{|\beta|_1+1} k^{|\beta|_1} \| e \|^k
\right)
M (2C)^{|\beta|_1}
\;\;\; \mbox{ for } \lambda \in U, \; 0 \leq |\beta|_1 \leq s.
\end{eqnarray}
\end{Proposition}

\begin{proof}
From condition (\ref{cr2.120}) it follows that
$\|X A^{-1}\| \leq 1/2$. Therefore we can invert $S := e + X A^{-1}$ by
a Neumann series. Keeping in mind that $\|e\| \geq 1$ (by sub-multiplicativity
of the norm) we obtain  that
\begin{eqnarray}
\label{cr2.140}
\|S^{-1}\| \leq \|e\| + 1 \leq 2 \|e\|.
\end{eqnarray}
Furthermore, we can express
\begin{eqnarray}
\label{cr2.150}
(A + X)^{-1} = A^{-1} S^{-1}. 
\end{eqnarray}
We have thus proved (\ref{cr2.130}) for 
$\beta = 0$. In order to estimate derivatives of $(A+X)^{-1}$
we use the 
Leibniz rule, (\ref{cr2.100}) -- (\ref{cr2.120}), 
proposition \ref{Pcr.2} and find for $1 \leq |\beta|_1 \leq s$,
\begin{eqnarray}
\nonumber
\| \partial^{\beta} S \| &\leq& 
\sum_{\alpha \leq \beta} 
\frac{\beta!}{\alpha! (\beta- \alpha)!} \| \partial^{\alpha} X \|
\| \partial^{\beta - \alpha} A^{-1} \| \\
\nonumber
&\leq&
\sum_{\alpha \leq \beta} 
\frac{\beta!}{\alpha! (\beta- \alpha)!} \epsilon M C^{|\beta|_1} \\
\label{cr2.200} 
&\leq&
\frac{1}{2} (2 C)^{|\beta|_1}.
\end{eqnarray}
We apply (\ref{cr2.150}), Corollary \ref{Ccr2.1},
(\ref{cr2.100}), (\ref{cr2.140}), (\ref{cr2.200}) and proposition
\ref{Pcr.2} to obtain
for $1 \leq |\beta|_1 \leq s$ the following estimates.
\begin{eqnarray}
\nonumber
\| \partial^{\beta} (A+X)^{-1} \| &=& 
\| \partial^{\beta} (A^{-1} S^{-1}) \| \\
\nonumber
&\leq&
\sum_{k=1}^{|\beta|_1+1} 
\!\!\!\!\!\!\!\!\!
\sum_{\scriptsize
{\scriptstyle
\begin{array}{c}
\alpha_1 + \ldots + \alpha_{k} = \beta \\
\alpha_i \neq 0 \mbox{ for } 2 \leq i \leq k
\end{array}
}
}
\frac{\beta!}{\alpha_1! \ldots \alpha_k!}
M C^{|\alpha_1|_1} (2 \|e\|)^k \left( \frac{1}{2} \right)^{k-1} 
(2 C)^{|\alpha_2|_1 + \ldots + |\alpha_k|_1} \\
\nonumber
&\leq&
\sum_{k=1}^{|\beta|_1+1} 
\!\!\!\!\!\!\!
\sum_{\scriptsize
{\scriptstyle
\begin{array}{c}
\alpha_1 + \ldots + \alpha_{k} = \beta \\
\alpha_i \neq 0 \mbox{ for } 2 \leq i \leq k
\end{array}
}
}
\frac{\beta!}{\alpha_1! \ldots \alpha_k!}
M (2C)^{|\beta|_1} 2 \|e\|^k \\
\nonumber
&\leq&
2 \sum_{k=1}^{|\beta|_1+1} k^{|\beta|_1} \|e\|^k M (2C)^{|\beta|_1}.
\end{eqnarray}
This completes the proof of Proposition \ref{Pcr2.1}.
\end{proof}

\subsection{The chain rule in a special case}
\label{cr3}
The goal of this section is to make the multi-linear form $D^k g$
which appears in formula (\ref{cr.10}) more explicit in a special case.
Before we state our proposition we introduce the following notation.
For $f: U \to X$  
smooth ($U \subset {\Bbb R}^d$ open,
$X$ a Banach space), and a multi-index $\beta \in {\Bbb N}_0^d$ we denote
\begin{eqnarray}
\label{cr3.10} 
f^{[\beta]} := \frac{1}{\beta!} \partial^{\beta} f.
\end{eqnarray} 
\begin{Proposition}
\label{Pcr3.1}
Let $U \subset {\Bbb R}$ and $V \subset {\Bbb R}^d$ be open sets
and let $f: U \to V$, $p: V \to {\Bbb R}$ be $C^s$ -- functions.
Then the composition $g = p \circ f$ is also a $C^s$ -- function
with
\begin{eqnarray}
\label{cr3.20}
g^{[s]} = \sum_{1 \leq |\alpha|_1 \leq s}
(p^{[\alpha]} \circ f)
\sum_{\scriptsize
{\scriptstyle
\begin{array}{c}
\delta^{i} \in {\Bbb N}^{\alpha_i} \mbox{ for } 1 \leq i \leq d \\
\sum_{i} |\delta^i|_1 = s
\end{array}
}
}
\prod_{i=1}^{d}
\prod_{j=1}^{\alpha_i}
f_i^{[\delta^{i}_j]}.
\end{eqnarray}
\end{Proposition}

\begin{Remark}
\label{Rcr3.1}
{\em
In formula (\ref{cr3.20}) we need to clarify the notation
$\delta^{i} \in {\Bbb N}^{\alpha_i}$, if $\alpha_i = 0$.
In this case we understand $|\delta^i|_1 = 0$ and use the standard
convention for the empty product
\begin{eqnarray}
\nonumber
\prod_{j=1}^{\alpha_i}
f_i^{[\delta^{i}_j]} = 1.
\end{eqnarray}
}
\end{Remark}

\begin{proof}
One could derive proposition \ref{Pcr3.1} from proposition
\ref{Pcr.1}. However, it is more convenient to give a direct proof via
Taylor expansion. In fact, proposition \ref{Pcr3.1} is a consequence of
the expansions (\ref{cr3.50}) and (\ref{cr3.60}) below. For
$\tau$, $t \in U$ 
\begin{eqnarray}
\nonumber
p(f(\tau)) - p(f(t)) &=&
\sum_{1 \leq |\alpha|_1 \leq s}
p^{[\alpha]}(f(t)) (f(\tau) - f(t))^{\alpha} + 
o(|f(\tau) - f(t)|^s) \\
\label{cr3.50}
&=& \sum_{1 \leq |\alpha|_1 \leq s}
p^{[\alpha]}(f(t))
\prod_{i=1}^d (f_i(\tau) - f_i(t))^{\alpha_i} + o(|\tau - t|^s),\\
\label{cr3.60}
f_i(\tau) - f_i(t) &=&
\sum_{p=1}^s f_i^{[p]}(\tau - t)^p + o(|\tau - t|^s).
\end{eqnarray}
\end{proof}

\subsection{Differentiating implicitly defined functions}
\label{cr4}

\begin{Proposition}
\label{Pcr4.1}
Let $k$, $d_1$, $d_2$, $d_3 \in {\Bbb N}$ and 
let $U \subset {\Bbb R}^{d_1}$, $V \subset {\Bbb R}^{d_2}$ be open sets.
Suppose further that
$f: U \to V$, $ x \mapsto f(x)$ and $g:V \times U \to {\Bbb R}^{d_3}$,
$(y, x) \mapsto g(y, x)$
are $C^{k}$ -- functions. Then $G : U \to {\Bbb R}^{d_3}$, 
$G(x) := g(f(x), x)$ is again a $C^k$ -- function and  the derivatives
$\partial^{\beta} G$, $1 \leq |\beta|_1 \leq k$, 
can be written in the following form. 
\begin{eqnarray}
\label{cr4.10}
\partial^{\beta} G =
\sum_{i=1}^{d_2} \left( \partial_{y_i} g \right) \partial^{\beta} f_i
+ 
\sum_{s \in S_{0}(\beta)} 
\left( \partial^{\gamma^{(s)}} g \right)
\cdot \prod_{l=1}^{l^{(s)}} 
\partial^{\alpha_l^{(s)}} f_{i_l^{(s)}},  
\end{eqnarray}
where
\begin{eqnarray}
\label{cr4.20}
S_0(\beta) \; \mbox{ is a set of cardinality } \; \# S_0(\beta) \leq 
(d_2 + |\beta|_1)^{|\beta|_1} - d_2, \\
\label{cr4.30}
1 \leq |\gamma^{(s)}|_1 \leq |\beta|_1 
\;\;\; \mbox{ for } \; s \in S_0(\beta),
\\
\label{cr4.40}
0 \leq l^{(s)} \leq |\beta|_1 \;\;\; \mbox{ for } \; s \in S_0(\beta),
\\
\label{cr4.50}
1 \leq |\alpha_l^{(s)}|_1 \leq |\beta|_1 - 1
\;\;\; \mbox{ for } \; s \in S_0(\beta) \; \mbox{ and } \; 
1 \leq l \leq l^{(s)},
\\
\label{cr4.60}
\sum_{l=1}^{l^{(s)}} |\alpha_l^{(s)}|_1 \leq |\beta|_1 
\;\;\; \mbox{ for } \; s \in S_0(\beta).
\end{eqnarray}
\end{Proposition}

\begin{proof}
We prove proposition \ref{Pcr4.1} by induction on $|\beta|_1$. The case
$|\beta|_1 = 1$ is straight forward to verify. Suppose that (\ref{cr4.10})
holds for some $1 \leq |\beta|_1 \leq k-1$ and denote $\partial^{\beta} G = A(\beta)
+ B(\beta)$,
where $A(\beta)$ and $B(\beta)$ denote the sums in 
(\ref{cr4.10}). Let $1 \leq m \leq d_1$.
Then $\partial^{\beta + e_m} G = \partial^{e_m} A(\beta) 
+ \partial^{e_m} B(\beta)$ with
\begin{eqnarray}
\nonumber
\partial^{e_m} A(\beta) &=&
\left(
\sum_{i=1}^{d_2} \left( \partial_{y_i} g \right) \partial^{\beta +e_m} f_i
\right) +
\sum_{i=1}^{d_2}
\left(
\sum_{j=1}^{d_2} \left( \partial_{y_j} \partial_{y_i} g \right) 
\partial^{e_m} f_j + \left( \partial_{x_m} \partial_{y_i} g \right)
\right)
\partial^{\beta} f_i \\
\nonumber
&=& I + II.
\end{eqnarray}
Note that $I = A(\beta + e_m)$, whereas $II$ contributes $d_2 (d_2 + 1)$ terms
to $B(\beta+e_m)$ which satisfy conditions (\ref{cr4.30}) -- (\ref{cr4.60})
with respect to $\beta + e_m$ (note that $|\beta + e_m|_1 \geq 2$). 
It is not difficult to see that
$\partial^{e_m} B(\beta)$ is a sum of at most 
$\# S_0 (\beta) \cdot (d_2 + 1 + |\beta|_1)$ 
terms where each summand again satisfies 
(\ref{cr4.30}) -- (\ref{cr4.60}) with respect to $\beta + e_m$. 
Using again $|\beta|_1 \geq 1$ we conclude that 
$\partial^{\beta+e_m} G$ can be written in the form
(\ref{cr4.10}) and the corresponding conditions
(\ref{cr4.20}) -- (\ref{cr4.60}) hold.
\end{proof}

\small 

\section{Table of notation used in chapters II and III}
\label{tnot}

\begin{tabular}{ll}
$A$, $A_1$
& (\ref{npc.72}), (\ref{npc.73}) \\

$b$
&theorem \ref{Tsmr3.1} \\

$B_0$, $B_1$, $B_2$ 
&(\ref{npc.80}), (\ref{npc.90}), (\ref{npc.100}) \\

$B_{\rho}(X)$, $B(X, \rho)$
& (\ref{npc.20}) \\

$c$
& (\ref{npc.52}) \\

${\cal C}^{(j)}$
& (\ref{npc.450}), (\ref{npc.460}) \\

$D(\omega)$
& (\ref{ova.85}) \\

$D_N$
& (\ref{npc.54}) \\

$D_P$
& (\ref{npc.480}) \\

$D_V$, $d_V$
& proposition \ref{Plop.1} \\

$D_W$
& lemma \ref{Lenp.1} \\

$D_{\tau}$, $D_K$, $D_E(x)$
& (\ref{npc.110}), (\ref{npc.115}), (\ref{npc.118}) \\

$d_{\tau, c}$, $d_{min}$, $\tilde{d}_{min}$
&(\ref{npc.120}), (\ref{npc.130}), (\ref{npc.140}) \\

$D_{\psi}(k)$
&(\ref{npc.265}) \\

$D_i$, $D_{i, j}$ ($1 \leq i \leq 4$, $1 \leq j \leq 5$)
& (\ref{npc.150}) -- (\ref{npc.260}) \\

$E_{\rho}$, $E_{\delta}$, $E_M$
& (\ref{npc.75}), (\ref{npc.76}), (\ref{npc.77}) \\

$F$
&theorem \ref{Tsmr3.1} \\

$g$ 
& (\ref{smr1.70}) \\

${\cal L}_{\sigma, c}$, ${\cal L}_{\sigma, 1}$
& definition \ref{DFa2} \\

$M_j$, $j \geq 0$
& (\ref{npc.440}) \\

$n_0$
& lemma \ref{hyp} \\

$N_j$, $j \geq 0$
& (\ref{npc.400}), (\ref{npc.410}) \\

$P$ 
&(\ref{ova.110}) \\

${\cal P}$
& (\ref{ova.160}) \\

${\cal POL}$
& (\ref{npc.490}) \\

$q$
&(\ref{npc.78}) \\

$Q$ 
&(\ref{ova.100}) \\

${\cal Q}$ 
& (\ref{ova.150}) \\

$r_{F, b}$
& (\ref{npc.60}) \\

$s$
& (\ref{smr2.50}) \\

${\cal S}$ 
&(\ref{ova.90}) \\

$T^{(j)}(\lambda)$, $T^{(j)}(\theta, \lambda)$
& (\ref{npc.63}), (\ref{npc.64}) \\

$U_{\rho}(Z)$, $U(Z, \rho)$
& (\ref{npc.30}) \\

$V(\omega)$, $V(\theta, \omega)$, $V_l$
& (\ref{smr2.70}), (\ref{ova.520}), (\ref{lop.50}) 
\end{tabular}
\newpage
\begin{tabular}{ll}
$W$ 
& lemma \ref{Lenp.1} \\

$w_{\sigma, c}$, $w_{\sigma, 1}$
& definition \ref{DFa1} \\

$X_{\sigma, c}$, $X_{\sigma, 1}$
& definition \ref{DFa2} \\

$Z_C$, $\tilde{Z}_{C}$
&(\ref{ova.600}), (\ref{ova.610}) \\

$\alpha_k$
& (\ref{npc.70}) \\

$\gamma$
&theorem \ref{Tsmr3.1} \\

$\delta_V$
& proposition \ref{Plop.1} \\

$\delta_j$, $j \geq 0$
& (\ref{npc.430}) \\

$\kappa$
&theorem \ref{Tsmr3.1} \\

$\tilde{\kappa}$
& (\ref{npc.270}) \\

$\lambda^{(0)}$
&(\ref{ova.140}) \\

$\nu$ 
& (\ref{smr2.20}) \\

$\pi^{(j)}_l$
& (\ref{npc.500}) \\

$\rho_j$, $\tilde{\rho}_j$, $j \geq 0$
&(\ref{npc.420}), (\ref{npc.425})   \\

$\sigma_j$
& (\ref{npc.56}) \\

$\tau$
& (\ref{smr2.50}) \\

$\varphi$ 
&(\ref{ova.120}) \\

$\psi$
& section \ref{npc} G) \\

$\omega^{(0)}$ 
&(\ref{smr2.30}), (\ref{smr2.40}) \\

$\Omega$ 
& (\ref{smr2.60}) -- (\ref{smr2.90}) \\

$|\cdot|$, $|\cdot|_1$, $|\cdot|_2$
& section \ref{npc} A) \\

$<\cdot, \cdot>$
& (\ref{npc.40}) \\

$[ \cdot ]$
& (\ref{ova.560}) \\

$\lfloor \cdot \rfloor$ 
& (\ref{npc.50}) 
\end{tabular}

\end{document}